%bgc.tex  The Birth of the Giant Component
% DEAR READER: THIS PAPER IS CURRENTLY BEING REFEREED;
% PLEASE REPORT ANY GLITCHES YOU FIND TO ONE OF THE AUTHORS, IMMEDIATELY!
% e.g. to svante@bellatrix.tdb.uu.se

% YOU NEED TO COPY THE FILE picmac.tex, FOUND IN THIS DIRECTORY; SOME OF YOU
% WILL NEED TO CHANGE "\font\tencirc=lcircle10 \font\tencircw=lcirclew10"
% TO "\font\tencirc=circle10 \font\tencircw=circlew10" IN THAT FILE!

\magnification\magstep1
\baselineskip 14pt
\parskip 0pt plus 2pt

\dimen0=\fontdimen9\tensy \advance\dimen0 .5pt \fontdimen9\tensy=\dimen0
\dimen0=\fontdimen12\tensy \advance\dimen0 .5pt \fontdimen12\tensy=\dimen0
 % that makes a bit more space above and below barline in textstyle fractions

\def\dddots{\mathinner{\mskip1mu\raise1pt\vbox{\kern7pt\hbox{.}}\mskip2mu
    \raise4pt\hbox{.}\mskip2mu\raise7pt\hbox{.}\mskip1mu}}
\def\bib{\par\noindent\hangindent 20pt}
\def\display#1:#2:#3\par{\par\hangindent #1 \noindent
			\hbox to #1{\hfill #2 \hskip .1em}\ignorespaces#3\par}
\def\disleft#1:#2:#3\par{\par\hangindent#1\noindent
			 \hbox to #1{#2 \hfill \hskip .1em}\ignorespaces#3\par}
\def\ra{\rightarrow}
\def\rbar{\relbar\!\!\relbar}
\def\pfbox
 {\hbox{\hskip 3pt\vbox{
 \def\|{\vrule height 4pt depth 2pt}
 \hrule
 \hbox to 5pt{\|\hfill\|}  
 \hrule
 }}\hskip 3pt}
\def\proof{\noindent{\sl Proof}.\enspace}

\def\ovline#1{{\kern.1em\overline{\kern-.1em#1\kern-.1em}\kern.1em}}
\def\MM{\ovline{\ovline M}}
\def\OF{\ovline{\ovline F}}
\def\Pr{\mathop{\rm Pr}\nolimits}
\def\half{{1\over2}}

\mathcode`@="8000
{\catcode`\@=\active \gdef@{\mskip1mu\relax}} % "hair space" in formulas

\input picmac

% Defs for the Bibliography
\def\Bag{1}
\def\BD{2}
\def\BCM{3}
\def\Bi{4}
\def\Bii{5}
\def\Biii{6}
\def\BF{7}
\def\Bor{8}
\def\Bri{9}
\def\Cay{10}
\def\Eis{11}
\def\ERo{12}
\def\ER{13}
\def\FKP{14}
\def\FKG{15}
\def\GJ{16}
\def\CM{17}
\def\Hen{18}
\def\Jan{19}
\def\JL{20}
\def\Kar{21}
\def\AOC{22}
\def\CP{23}
\def\KP{24}
\def\Li{25}
\def\Lii{26}
\def\LW{27}
\def\LPW{28}
\def\Mit{29}
\def\Ram{30}
\def\Ren{31}
\def\Rid{32}
\def\Sei{33}
\def\Sla{34}
\def\Si{35}
\def\Sii{36}
\def\Syl{37}
\def\Vob{38}
\def\Wo{39}
\def\Wsi{40}
\def\Wi{41}
\def\Wsii{42}
\def\Wii{43}
\def\Wiii{44}
\def\Wiv{45}

\centerline{\bf The Birth of the Giant Component}
\bigskip
\centerline{Dedicated to Paul Erd\H{o}s on his 80th birthday}
\bigskip
\centerline{\sl Svante Janson, Donald E. Knuth, Tomasz \L uczak,
 and Boris Pittel}
\smallskip
\centerline{\sl Department of  Mathematics, Uppsala University;}
\centerline{\sl  Computer Science Department, Stanford University;}
\centerline{\sl Department of Discrete Mathematics, Adam Mickiewicz
University; and}
\centerline{\sl Department of Mathematics, The Ohio State University}
\bigskip
{\narrower\smallskip\noindent
{\bf Abstract.}\enspace
Limiting distributions are derived for the sparse connected components that are
present when a random graph on $n$~vertices has approximately $\half  n$
edges. In particular, we show that such a graph consists entirely of trees,
unicyclic components, and bicyclic components with probability approaching
$\sqrt{2\over 3}\,\cosh\sqrt{5\over 18}\approx0.9325$
 as $n\to\infty$. The limiting probability that it consists of trees, unicyclic
components, and at most one other component is approximately
$0.9957$; the limiting probability that it is planar lies between
$0.987$ and $0.9998$. When a random graph evolves and the number of edges 
passes~$\half n$,
its components grow in cyclic complexity according to an interesting
Markov process
whose asymptotic structure is derived. The probability that there never is
more than a single component with more edges than vertices, throughout
the evolution,  approaches $5\pi/18\approx0.8727$.
A~``uniform'' model
of random graphs, which allows self-loops and multiple edges, is shown to lead
to formulas that are substantially simpler than the analogous formulas for
the classical random graphs of Erd{\H o}s and R\'enyi. 
The notions of ``excess''
and ``deficiency,'' which are significant characteristics of the generating
function as well as of the graphs themselves, lead to a mathematically
attractive structural theory for the uniform model. A~general
approach to the study of stopping configurations makes it possible to
sharpen previously obtained estimates in a uniform manner and often to
obtain closed forms for the constants of interest. Empirical results are
presented to complement the analysis, indicating the typical behavior when
$n$ is near 20000.
\footnote{}{\raise 11pt\null
This research 
was supported in part by the National Science Foundation under grant
CCR-86-10181,
and by Office of Naval Research contract
N00014-87-K-0502.}
\smallskip}

\bigskip
\noindent
{\bf 0. Introduction.}\enspace
When edges are added at random to $n$ initially disconnected points, for large~$n$,
a~remarkable transition occurs when the number of edges becomes 
approximately~$\half n$. Erd{\H o}s and R\'enyi [\ER]
studied random graphs with $n$~vertices and ${n\over 2}(1+\mu)$ edges
as $n\to\infty$, and discovered that such graphs almost surely have the
following properties: If $\mu<0$, only small
trees and ``unicyclic'' components
are present, where a unicyclic component is a tree with one additional edge;
moreover, the size of the largest tree 
component is ${(\mu-\ln(1+\mu))^{-1}}\ln n
+O(\log\log n)$. If $\mu =0$, however, the largest component has size 
of order~$n^{2/3}$. And if $\mu >0$, there is a unique ``giant'' component
whose size is of order~$n$; in fact, the size  of this component is
asymptotically~$\alpha n$ when $\mu =-\alpha^{-1}\ln(1-\alpha)-1$. 
Thus, for example, a~random graph with approximately $n\ln 2$ edges will have
a giant component containing $\sim\half n$ vertices.

The research that led to the present paper began in a rather curious
way, as a result of a misunderstanding. In 1988, the students in a class
taught by Richard M. Karp  performed computer experiments in which
graphs with a moderately large number of vertices were generated by adding
one edge at a time. A~rumor spread that these simulations had 
turned up a surprising fact: As each of the
random graphs evolved, the story went, 
never once was there more than a single ``complex''
component; i.e., there never were two or more components present simultaneously
that were neither trees nor unicyclic.
Thus, the first connected component that acquired more edges than vertices was
destined to be the giant component. As more edges were added, this component
gradually swallowed up all of the others, and none of the others ever became
complex before they were swallowed.

Reports of those experiments suggested that a great simplification of
the theory of evolving graphs might be possible.  Could it be that
such behavior occurs almost always, i.e., with probability
approaching~1 as $n\to\infty$? If so, we could hope for the existence
of a much simpler explanation of the fact that a~giant component
emerges during the graph process, and we could devise rather simple
algorithms for online graph updating that would take advantage of the
unique-complex-component phenomenon.  At that time the authors who
began this investigation (DEK and BP) were unaware of Stepanov's
posthumous paper~[\Sii]. We were motivated chiefly by the work of
Bollob\'as [\Bii], who had shown that a component of size $\geq n^{2/3}$
is almost always unique once the number of edges exceeds $\half 
n+2(\ln n)^{1/2}n^{2/3}$; moreover, Bollob\'as proved that such a
component gets approximately 4~vertices larger when each new edge is
added.  His results blended nicely with the unique-complex-component
conjecture.

However, we soon found that the conjecture is false: There is
nonzero probability that a graph 
with $\half n$ edges will contain several pretenders to the giant throne,
and this probability increases when the number of edges is slightly more
than~$\half n$.
We also learned that Stepanov [\Sii] had already obtained similar results.
Thus we could not hope for a theory of random graphs that would be as simple
as the conjecture promised. On the other hand, we learned that
the graph evolution process does satisfy the conjecture with reasonably
high probability; hence algorithms
whose efficiency rests on the assumption of a unique complex component will
not often be inefficient.

Further analysis revealed, in fact, that we must have misunderstood
the initial reports of experimental data. The actual probability that
an evolving graph never has two complex components approaches the limiting
value $5\pi/18\approx 0.8727$; therefore the rumor that got us started could
not have been true. In fact, the computer experiments by Karp's students had
simply reported the state of the graph when exactly $\half n$
edges were present, and at certain other fixed reporting times. A false
impression arose because there is high probability that a random graph with $\half n$
edges has at most one
complex component; indeed, the probability is $0.9957+O(n^{-1/3})$. More
complicated configurations sometimes arise momentarily just after $\half n$
edges are reached. However, the fallacious rumor of 1988 has turned
out to have beneficial effects, because it was a significant catalyst
for the discovery of some remarkably beautiful patterns.

Sections~1--10 of this paper provide a basic introduction to the theory of
evolving graphs and multigraphs, using generating functions as the principal
tool. Two models of graph evolution are presented in section~1, the
``graph process'' and the ``multigraph process.'' Their generating functions
are introduced in section~2, and special aspects of those functions related
to trees and cycles are discussed in section~3. Section~4 explains how to
derive properties of a graph's more complex features by means of
differential equations; the equations are solved for multigraphs in section~5
and for graphs in section~6. The resulting decomposition of multigraphs
turns out to be surprisingly regular. Section~7 explains the regularities and
begins to analyze the algebraic properties of the functions obtained in
section~5. Related results for connected graphs are discussed in section~8.
Section~9 explains the combinatorial significance of the algebraic
structure derived earlier. Finally, section~10 presents a quantitative
lemma about the characteristics of random graphs near the critical
point $\mu=0$, making it possible to derive exact values for many
relevant statistics.

Readers who cannot wait to get to the ``good stuff'' should skim sections 1--10
and move on to section~11, which begins a sequence of applications of the
basic theory. The first step is to analyze the distribution of bicyclic
components; then, in section~12, the same ideas are shown to yield the
joint distribution of all kinds of components. The formulas obtained there
have a simple structure suggesting that the traditional approach of focussing
on connected components is unnecessarily complicated; we obtain a simpler
and more symmetrical theory if we first  consider the {\it excess\/} of edges
over vertices, exclusive of tree components, then look at other properties
like connectedness after conditioning on the excess. Section~13 motivates
this principle, and section~14 derives the probability distribution of a
graph's excess as it passes the critical point. These ideas help to
nail down the probability that a graph with $\half n$ edges is planar,
as shown in section~15.

Section 16 begins the discussion of what may well be the most important
notion in this paper; readers who have time for nothing else are
encouraged to look at Figure~1, which shows the initial stages of the
``big bang.'' The evolution of a graph or multigraph passes through
discrete transitions as the excess increases, and important aspects
of those changes are illustrated in Figure~1; section~17 proves that
this illustration represents a Markov process that characterizes almost
all graph evolutions. The ${5\pi\over18}$ phenomenon alluded to above
is discussed in section~18, which establishes $5\pi\over18$ as an
upper bound for the probability in question. Section~19 shows that, for
small~$n$, the probability of retaining at most one complex component
during the critical stage is in fact greater than~$5\pi\over18$, decreasing
monotonically with~$n$.

The excess of a graph is of principal importance at the critical point, but
a secondary concept called {\it deficiency\/} becomes important shortly
thereafter. A graph with deficiency~0 is called ``clean''; such graphs
are obtained from 3-regular graphs by splitting edges and/or by attaching
trees to vertices of cycles. Section~20 explains how deficiency evolves
jointly with increasing excess. Figure~2, at the end of that section,
illustrates another Markov process that goes on in parallel with
Figure~1. Section~21 shows that most graphs stay
clean until they have acquired approximately $\half n+n^{3/4}$ edges.
Section~22 looks more closely at the moment a graph first becomes unclean.

Section 23 tracks the growth of excess and deficiency as a multigraph
continues to evolve through $\half n+n^{4/5}$, $\half n+n^{5/6}$,
\dots~edges. The excess and deficiency are
 shown to be approximately normally distributed about
certain well-defined values. Specifically, when the number of
edges is ${n\over2}(1+\mu)$, with $\mu=o(1)$, the excess
will be approximately ${2\over3}\mu^3n$ and the deficiency will be
approximately ${2\over3}\mu^4n$. These statistics complement the well-known
fact that the emerging giant component has almost surely grown to
encompass approximately $2\mu n$ vertices.

Sections 24 to 26 develop a theory of ``stopping configurations,'' by
which it is possible to study the first occurrences of various
events during a multigraph's evolution. 
In particular, an explicit formula is derived for
the asymptotic distribution of the time when the excess
first reaches a given value~$r$.
A closed formula is derived for the ``first cycle constant'' of~[\FKP].

Section 27 completes the discussion initiated in sections 17 and~18,
by proving the ${5\pi\over 18}$ phenomenon as a special case of a more
general result about the infinite Markov process in Figure~1.

Finally, section 28 presents empirical data, showing to what extent the
theory relates to practice when $n$ is not too large. Section~29
discusses a number of open questions raised by this work.

\bigbreak
\noindent
{\bf 1. Graph evolution models.}\enspace
We shall consider two ways in which a random graph on $n$~vertices might
evolve, corresponding to sampling with and without replacement. The first
of these, introduced implicitly in [\Bi] and explicitly in [\BF,
proof of Lemma~2.7] and [\FKP],
turns out to be simpler to analyze and simpler to simulate by
computer, therefore more likely to be of importance in applications to
computer science: We generate ordered pairs $\langle x,y\rangle$ repeatedly,
where $1\le x,y\le n$, and add the (undirected) edge
$x\rbar y$ to the graph.
Each ordered pair $\langle x,y\rangle$ occurs with probability $1/n^2$, so
we call this the {\it uniform model\/} of random graph generation. It may
also be called the {\it multigraph process}, because it can  generate
graphs with self-loops $x\rbar x$, and it can also generate 
multiple edges. Notice
that a self-loop $x\rbar x$ is generated with probability $1/n^2$, while an
edge $x\rbar y$ with $x\ne y$ is generated with probability $2/n^2$ because it can
occur either as $\langle x,y\rangle$ or $\langle y,x\rangle$.

The second evolution procedure, introduced by Erd{\H o}s and R\'enyi [\ERo],
is called the {\it permutation model\/} or the {\it graph process}.
In this case we consider all $N={n\choose 2}$ possible edges $x\rbar y$
with $x<y$ and introduce them in random order, with all $N!$~permutations
considered equally likely. In this model there are no self-loops or
multiple edges.

A multigraph $M$ on $n$ labeled vertices can be defined by a symmetric
$n\times n$ matrix of nonnegative integers 
$m_{xy}$, where $m_{xy}=m_{yx}$ is the number of undirected edges
$x\rbar y$ in~$G$. For purposes of analysis, we shall assign
a {\it compensation factor\/}
$$\kappa(M)=
 1\left/\prod_{x=1}^n\left(2^{m_{xx}}\prod_{y=x}^nm_{xy}!\right)\right.
   \eqno(1.1)$$
to $M$; if $m=\sum_{x=1}^n\sum_{y=x}^nm_{xy}$ is the total number of edges,
the number of sequences $\langle x_1,y_1\rangle\langle x_2,y_2\rangle
\,\ldots\,\langle x_m,y_m\rangle$ that lead to~$M$ is then exactly
$$2^m\,m!\,\kappa(M)\,.\eqno(1.2)$$
(The factor $2^m$ accounts for choosing either $\langle x,y\rangle$
or $\langle y,x\rangle$; the $2^{m_{xx}}$ in the denominator of~$\kappa(M)$
compensates for the case $x=y$. The other factor $m!$ accounts for permutations
of the pairs, with $m_{xy}!$ in $\kappa(M)$ to compensate for permutations
between multiple edges.)

Equation (1.2) tells us that $\kappa(M)$ is a natural weighting factor
for a multigraph~$M$, because it corresponds to the relative frequency
with which $M$ tends to occur in applications. For example, consider multigraphs
on three vertices $\{1,2,3\}$ having exactly three edges. The edges will
form the cycle $M_1=\{1\rbar2,\;2\rbar3,\;3\rbar1\}$ much more often than they will
form three identical self-loops $M_2=\{1\rbar1,\;1\rbar1,\;1\rbar1\}$, when the
multigraphs are generated in a uniform way. For if we consider the $3^6$
possible sequences $\langle x_1,y_1\rangle\langle x_2,y_2\rangle\langle x_3,y_3
\rangle$ with $1\le x,y\le 3$, only one of these generates the latter
multigraph, while the cyclic multigraph is obtained in $2^3\,3!=48$ ways.
Therefore it makes sense to assign weights so that
$\kappa(M_2)={1\over48}\kappa(M_1)$, and indeed (1.1) gives
$\kappa(M_1)=1$, $\kappa(M_2)={1\over48}$.

Notice that a given multigraph $M$ is a graph---i.e., it has no loops and no 
multiple edges---if and only if $\kappa(M)=1$.
Notice also that if $M$ consists of several disjoint
components $M_1,\ldots,M_k$, with no edges between vertices of~$M_i$ and~$M_j$
for $i\ne j$, we have
$$\kappa(M)=\kappa(M_1)\,\ldots\,\kappa(M_k)\,.\eqno(1.3)$$

\bigbreak\noindent
{\bf 2. Generating functions.}\enspace
We shall use bivariate generating functions (bgf's) to study labeled 
graphs and multigraphs
and their connected components. If ${\cal F}$ 
is a family of multigraphs with labeled vertices,
 the associated bgf is the formal power series
$$F(w,z)=\sum_{M\in {\cal F}} \kappa(M)\,w^{m(M)}{z^{n(M)}\over n(M)!}\,,\eqno(2.1)$$
where $m(M)$ and $n(M)$ denote the number of edges and the number of
vertices of~$M$. 
We can do many operations on such power series without regard to convergence.
It follows from (1.2) and (2.1) that $m$~steps of the uniform
evolution model on $n$~vertices will produce a multigraph in~${\cal F}$ with
probability
$${2^m\,m!\,n!\over n^{2m}}\,\,[w^mz^n]\,\,F(w,z)\,,\eqno(2.2)$$
where the symbol $[w^mz^n]$ denotes the coefficient of $w^mz^n$ in the
formal power series that follows it. Similarly, if ${\cal F}$ is a family of
graphs with labeled
vertices, the probability that $m$~steps of the permutation model will produce
a graph in~${\cal F}$ is
$${n!\over{N\choose m}}\,[w^mz^n]\,F(w,z)\,,\qquad N={n\choose 2}\,.\eqno(2.3)$$
Formulas (2.2) and (2.3) are asymptotically related by the formula
$${N\choose m}={n^{2m}\over 2^m\,m!}\exp\left(-{m\over n}-{m^2\over n^2}
+O\left({m\over n^2}\right)+O\left({m^3\over n^4}\right)\right)\,,
\qquad 0\le m\le N,\eqno(2.4)$$
which follows from Stirling's approximation.

Incidentally, the exponential factor in (2.4) is the probability that
$m$~steps of the multigraph
 process will produce no self-loops or multiple
edges. When $m=\half n$, this probability is $e^{-3/4}+O(n^{-1})
\approx 0.472$.

When we say that the $n$ vertices of a multigraph are ``labeled,'' it
is often convenient to think of the labeling as an assignment of the
numbers~1 to~$n$. But a strict numeric convention would require us to
recompute the labels whenever vertices are removed or when multigraphs
are combined. The actual value of a label is, in fact, irrelevant;
what really counts is the relative order {\it between\/} labels.
Labeled multigraphs are multigraphs whose vertices have been totally
ordered. In this paper all graphs and multigraphs are assumed to be
labeled, i.e., totally ordered, even when the adjective ``labeled''
is not stated.

The bgf (2.1) is an 
exponential generating function in $z$, and the factor
$\kappa(M)$ is multiplicative 
according to (1.3). Therefore the product of bgf's
$$F_1(w,z)F_2(w,z)\,\ldots\,F_k(w,z)$$
represents ordered $k$-tuples of labeled multigraphs 
$\langle M_1,M_2,\ldots,M_k\rangle$,
each $M_j$ being from family~${\cal F}_j$. Unordered $k$-tuples
$\{M_1,\ldots,M_k\}$ from a common family~${\cal F}$ have the bgf
$F(w,z)^k\!/k!$, if ${\cal F}$ does not include the empty
 multigraph. For example, the bgf for a 3-cycle is $w^3z^3\!/3!$, and
the bgf for two isolated vertices is $z^2\!/2!$; hence the bgf for a
3-cycle and two isolated
 vertices is $(w^3z^3\!/6)(z^2\!/2)=10w^3z^5\!/5!$.
(There are 10 such graphs, one for each choice of the isolated points.)

Let $C(w,z)$ be the bgf for all connected multigraphs, and let $G(w,z)$
be the bgf for the set of all multigraphs. Then we have
$$e^{C(w,z)}=\sum_{k\ge 0}{C(w,z)^k\over k!}=G(w,z)\eqno(2.5)$$
because the term $C(w,z)^k\!/k!$ is the bgf for multigraphs having
exactly $k$~components. 
Similarly, if $\widehat{C}(w,z)$ and $\widehat{G}(w,z)$
are the corresponding bgf's for graphs instead of multigraphs, we have
$$e^{\widehat{C}(w,z)}=\widehat{G}(w,z)\,,\eqno(2.6)$$
a well-known formula due to Riddell [\Rid].
The bgf for all graphs is obviously
$$\widehat{G}(w,z)=\sum_{n\ge 0}(1+w)^{n(n-1)/2}{z^n\over n!}\,.\eqno(2.7)$$
Therefore (2.6) gives us the bgf for connected graphs,
$$\eqalignno{\widehat{C}(w,z)&=\ln\left(1+z+(1+w){z^2\over 2}+
 (1+w)^3{z^3\over 6}+\cdots\,\right)\cr
\noalign{\smallskip}
&=z+w{z^2\over 2}+(3w^2+w^3){z^3\over 6}+\cdots\;.&(2.8)\cr}$$

The bgf $G(w,z)$
for all multigraphs can be found as follows: The coefficient of
$z^n\!/n!$ is $\sum\kappa(M)w^{m(M)}$, summed over multigraphs~$M$ on
$n$~vertices. This is
$$\prod_{x=1}^n \left(\biggl(\sum_{m_{xx}\ge 0}
{w^{m_{xx}}\over 2^{m_{xx}}\,m_{xx}!}\biggr)\prod_{y=x+1}^n\biggl(%
\sum_{m_{xy}\ge 0}{w^{m_{xy}}\over m_{xy}!}\biggr)\right)
=\prod_{x=1}^ne^{w/2}(e^w)^{n-x}=e^{wn^2\!/2}\,.$$
Hence the desired formula is slightly simpler than (2.7):
$$G(w,z)=\sum_{n\ge 0}e^{wn^2\!/2}{z^n\over n!}\,.\eqno(2.9)$$
The corresponding bgf for connected multigraphs is therefore
$$\eqalignno{C(w,z)&=\ln G(w,z)\cr
\noalign{\smallskip}
&=\bigl(1+{\textstyle\half }w+{\textstyle{1\over 8}}w^2+
 {\textstyle{1\over48}}w^3+\cdots\,\bigr)z
+\bigl(w+{\textstyle{3\over 2}}w^2+{\textstyle{7\over 6}}w^3+
\cdots\,\bigr){z^2\over2}\cr
\noalign{\smallskip}
&\qquad\null+\bigl(3w^2+{\textstyle{17\over 2}}
w^3+\cdots\,\bigr){z^3\over6}+\cdots\;.&(2.10)\cr}$$
In this case the coefficient of $w^3z^3$ is ${17\over2}/3!$, because the
connected multigraphs with three edges on three vertices have total
weight~$17\over2$. (The 3-cycle has weight~1; there are 9 multigraphs
obtainable by adding a self-loop to a tree, each of weight~$1\over2$; and
there are six multigraphs obtainable by doubling one edge of a tree,
again weighted by~$1\over2$.)

Notice that expression (2.2) is $[w^mz^n]\,F(w,z)\,/\,[w^mz^n]\,G(w,z)$,
the ratio of the weight of multigraphs in~$\cal F$ to the weight of
all possible multigraphs. Similarly, expression (2.3) is
$[w^mz^n]\,\widehat F(w,z)\,/\,[w^mz^n]\,\widehat G(w,z)$.

It is convenient to group the terms of (2.8) and (2.10)
according to the excess of edges over vertices in connected components. Let
${\cal C}_r$ and $\widehat{\cal C}_r$ denote the families of
connected
multigraphs and
graphs in which there are exactly $r$ more edges than vertices; let $C_r(w,z)$
and $\widehat{C}_r(w,z)$ be the corresponding bgf's. Then we have
$$\eqalignno{C(w,z)&=\sum_rC_r(w,z)=\sum_rw^rC_r(w z)\,,\cr
\noalign{\smallskip}
\widehat{C}(w,z)&=\sum_r\widehat{C}_r(w,z)=\sum_rw^r\widehat{C}_r(w z)\,,&(2.11)\cr}$$
where $C_r(z)$ and $\widehat{C}_r(z)$ are univariate
generating functions for ${\cal C}_r$ and~$\widehat{\cal C}_r$. A~univariate
generating function~$F(z)$ is $\sum\kappa(M)z^n\!/n!$, summed over all graphs
or multigraphs in a given family~$\cal F$. We obtain it from a bgf by
setting $w=1$, thereby ignoring the number of edges. Univariate generating
functions are easier to deal with than bgf's, so we generally try to avoid the
need for two independent variables whenever possible.

\bigbreak\noindent
{\bf 3. Trees, unicycles, and bicycles.}\enspace
Let us say that a connected component has {\it excess\/} $r$ if it
belongs to~${\cal C}_r$, i.e., if it has $r$ more edges than vertices.
A connected graph on $n$~vertices must have at least $n-1$ edges. Hence
$C_r=0$ unless $r\ge -1$. In the extreme case $r=-1$, we have
${\cal C}_{-1}=\widehat{\cal C}_{-1}$, the family of all unrooted trees,
which are {\it acyclic components}.
In the next case $r=0$, the generating functions $C_0$ and $\widehat{C}_0$
represent {\it unicyclic components}, which are trees with an additional
edge. Similarly, $C_1$ and $\widehat{C}_1$ represent {\it bicyclic components}.
In the present paper we shall deal extensively with sparse components
of these three kinds, so it will be convenient to use the special abbreviations
$$\vcenter{\halign{\hfil$#$&$\;#$\hfil\quad&#\hfil\cr
U(z)&=C_{-1}(z)=\widehat{C}_{-1}(z)&for unrooted trees;\cr
\noalign{\smallskip}
V(z)&=C_0(z)\hbox{ and }\widehat{V}(z)
=\widehat{C}_0(z)&for unicyclic components;\cr
\noalign{\smallskip}
W(z)&=C_1(z)\hbox{ and }\widehat{W}(z)
=\widehat{C}_1(z)&for bicyclic components.\cr}}$$
According to a well-known theorem of Sylvester [\Syl] and Borchardt [\Bor],
often attributed erroneously to Cayley [\Cay] although Cayley himself
credited Borchardt,  we have
$U(z)=\sum_{n=1}^\infty n^{n-2}z^n\!/n!$. The other four generating functions
begin as follows:
$$\eqalign{%
V(z)&=\textstyle\half z+{3\over 4}z^2+
{17\over 12}z^3+{71\over 24}z^4+{523\over80}z^5+{899\over60}z^6+\cdots\;;\cr
\noalign{\smallskip}
\widehat{V}(z)&=\textstyle{1\over 6}z^3
+{5\over 8}z^4+{37\over 20}z^5+{61\over 12}z^6
+\cdots\;; \cr
\noalign{\smallskip}
W(z)&=\textstyle{1\over 8}z+{7\over 12}z^2+
{101\over 48}z^3+{83\over 12}z^4+
{12487\over576}z^5+{3961\over60}z^6+\cdots\;; \cr
\noalign{\smallskip}
\widehat{W}(z)&=\textstyle{1\over4}z^4+{41\over 24}z^5+
{95\over 12}z^6+\cdots\;. \cr}$$

All of these generating functions can be expressed succinctly in terms
of the tree function
$$T(z)=\sum_{n\ge 1}n^{n-1}{z^n\over n!}
 =z+z^2+{\textstyle{3\over 2}}z^3+\cdots\,,
\eqno(3.1)$$
which generates {\it rooted\/} labeled trees and satisfies the functional
relation
$$T(z)=ze^{T(z)}\eqno(3.2)$$
due to Eisenstein [\Eis]. Indeed, the relation
$$U(z)=T(z)-{\textstyle\half }T(z)^2\eqno(3.3)$$
is well known, as are the formulas
$$\eqalignno{V(z)&=\half \ln{1\over 1-T(z)}\,,&(3.4)\cr
\noalign{\smallskip}
\widehat{V}(z)&=\half \ln{1\over 1-T(z)}-
\half T(z)-{1\over 4}T(z)^2\,;
&(3.5)\cr}$$
see [\FKP]. 
We can prove (3.4) and (3.5) by noting that the univariate generating
function for connected unicyclic multigraphs whose cycle has length~$k$ is
$${T(z)^k\over 2k}\,;$$
summing over $k\ge 1$  gives (3.4), and summing over $k\ge 3$ gives (3.5).
(If $k=1$, the cycle is a self-loop; hence the multigraph is essentially
a rooted tree and the compensation factor is~$\half $. If $k=2$, the
cycle is a duplicate edge; hence the multigraph is essentially an unordered
pair of rooted trees, and the compensation factor again is~$\half $.
If $k\ge 3$, the unicyclic component is essentially a sequence of $k$~rooted
trees, divided by~$2k$ to account for cyclic order and change of orientation.)

The generating function $\widehat{W}(z)$ was shown by G. N. Bagaev [\Bag]
to be
$$\widehat{W}(z)={T(z)^4\bigl(6-T(z)\bigr)\over 24\bigl(1-T(z)\bigr)^3}\,.
\eqno(3.6)$$
Then E. M. Wright made a careful study of all the generating functions
$\widehat{C}_k(z)$, which he called~$W_k$, in a series of significant papers
[\Wi, \Wii, \Wiii, \Wiv].
We will show below that the bgf for bicyclic connected 
multigraphs is
$$W(z)={T(z)\bigl(3+2T(z)\bigr)\over 24\bigl(1-T(z)\bigr)^3}\,.\eqno(3.7)$$

The coefficients of powers of $1/\bigl(1-T(z)\bigr)$ arise in numerous
applications, so Knuth and Pittel [\KP] began to catalog some of their
interesting properties. For each~$n$ the function $t_n(y)$ defined~by
$${1\over \bigl(1-T(z)\bigr)^y}=\sum_{n\geq 0}t_n(y)\,{z^n\over
n!}\eqno(3.8)$$
is a polynomial of degree $n$ in~$y$, called the {\it tree polynomial}
of order~$n$.
The coefficient of~$y^k$ in $t_n(y)$ is the number of mappings from an
$n$-element set into itself having exactly $k$~cycles. For fixed~$y$
and $n\rightarrow\infty$, we have [\KP, Lemma~2 and (3.16)]
$$t_n(y)={\sqrt{\mskip1mu2\pi}\,n^{n-1/2+y/2}\over
2^{y/2}\Gamma(y/2)}+O(n^{n-1+y/2})\,.\eqno(3.9)$$ 
We can, for example, express the number of connected bicyclic graphs
on $n$~vertices in terms of the tree polynomial~$t_n$, namely
$$\textstyle{{5\over 24}\,t_n(3)-{19\over 24}\,t_n(2)+{13\over
12}\,t_n(1) -{7\over 12}\,t_n(0)+{1\over 24}\,t_n(-1)+{1\over
24}\,t_n(-2)}\,,\eqno(3.10)$$ 
because (3.6) can be rewritten
$$\widehat{W}(z)={5\over 24\bigl(1-T(z)\bigr)^3}
-{19\over 24\bigl(1-T(z)\bigr)^2}+{13\over
12\bigl(1-T(z)\bigr)}-{7\over 12}+{1-T(z)\over 24}+{\bigl(1-T(z)\bigr)^2\over
24}\,.$$
Equation (3.9) tells us that only
the first term ${5\over 24}\,t_n(3)$ of (3.10) is asymptotically
significant. Extensions of (3.9) appear in equations (19.13) and
(19.14) below.

We can also express quantities like (3.10) in terms of Ramanujan's function
[\Ram]
$$\eqalignno{Q(n)&=1+{n-1\over n}+{n-1\over n}\,{n-2\over n}+{n-1\over n}\,{n-2\over
n}\,{n-3\over n}+\cdots\cr
\noalign{\smallskip}
&=\sqrt{{\pi n\over 2}}-{1\over 3}+{1\over 12}\sqrt{{\pi\over 2n}}-
{4\over 135n}+O(n^{-3/2})\,,&(3.11)\cr}$$
which Wright [\Wi] called $1+h(n)/n^n$. For we have
$$t_n(1)=n^n\,;\qquad t_n(2)=n^n(1+Q(n))\,;\qquad
t_n(y+2)=n{t_n(y)\over y}+t_n(y+1)\,,\quad y\ne0.\eqno(3.12)$$
(See [\KP, equations (2.7), (3.14), and (1.9)].) Furthermore, we have
$$[z^n]\,V(z)={\textstyle\half }n^{n-1}Q(n)\,;\eqno(3.13)$$
this follows from a well-known formula of R\'enyi [\Ren].

\bigbreak\noindent
{\bf 4. The cyclic components.}\enspace
For theoretical purposes it proves to be important to partition a multigraph
into its {\it acyclic part}, consisting entirely of isolated vertices or
trees, and its {\it cyclic part}, consisting entirely of components that
each contain at least one cycle. The cyclic part can in turn be partitioned
into the {\it unicyclic part}, consisting entirely of unicyclic components,
and the {\it complex part}, consisting entirely of components that have
more edges than vertices. 
A~multigraph is called cyclic if it equals its cyclic part, complex if
it equals its complex part. In this section and the next, we will
study the generating functions  for cyclic and complex multigraphs.
The formulas turn out to be surprisingly simple, and they will be the
key to much of what follows.

Let $F(w,z)$ be the bgf for all cyclic multigraphs, i.e., for all
multigraphs whose acyclic part is empty.
Formulas (2.5) and (2.11) tell us that
$$F(w,z)=e^{C_0(w,z)+C_1(w,z)+\cdots}=e^{C(w,z)-C_{-1}(w,z)}
 =G(w,z)\,e^{-U(wz)/w};$$
in other words,
$$G(w,z)=e^{U(wz)/w}\,F(w,z)\,.\eqno(4.1)$$
Indeed, this makes sense, because $e^{U(wz)/w}$ is the bgf for all
acyclic multigraphs. We will analyze $F=F(w,z)$ by studying a
linear differential equation satisfied by~$G=G(w,z)$, and seeing that
a similar equation is satisfied by~$F$.

Let $\vartheta_w$ be the differential operator $w{\partial\over\partial w}$,
and let $\vartheta_z$ be~$z{\partial\over\partial z}$. The operator
$\vartheta_w$ corresponds to marking an edge of a multigraph, i.e.,
giving some edge a special label, because $\vartheta_w$ multiplies the
coefficient of $w^mz^n$ by~$m$. Similarly, $\vartheta_z$ corresponds to
marking a vertex, because it multiplies the coefficient of $w^mz^n$ by~$n$.
(For a general discussion of marking, see [\GJ, sections 2.2.24 and
following].) We have
$$\vartheta_w\,G(w,z)\,=\,w\sum_{n\ge0}{n^2\over2}e^{wn^2\!/2}{z^n\over n!}\,
=\,{w\over2}\,\vartheta_z^2\,G(w,z)\,;$$
hence $G$ satisfies the differential equation
$${2\over w}\,\vartheta_w\,G\,=\,\vartheta_z^2\,G\,.\eqno(4.2)$$
Again, this makes sense: The left side represents all multigraphs having
a marked edge and an orientation assigned to that edge, and with the
edge count decreased by~1. The right side represents all multigraphs with an
ordered pair $\langle x,y\rangle$ of marked vertices. Orienting and
discounting an edge is the same as marking two vertices.

We can also write (4.2) in the suggestive form
$$G(w,z)\,=\,e^z+\half \int_0^w\,\vartheta_z^2\,G(w,z)\,dw\,,\eqno(4.3)$$
using the boundary condition $G(0,z)=e^z$. (The generating
function for all multigraphs with no edges is, of course,~$e^z$.)
The operator $\vartheta_z^2$
corresponds to choosing an ordered pair~$\langle x,y\rangle$, and the
operator $\half \int_0^w$ corresponds to disorienting that edge and
blending it into the existing multigraph. (Notice that the English words
``differentiation'' and ``integration'' are remarkably apt synonyms for the
combinatorial operations of marking and blending.)

Most of our work will involve $\vartheta_z$ instead of $\vartheta_w$,
so we shall often write simply $\vartheta$ without a subscript when
we mean $\vartheta_z$. The marking
operator $\vartheta$ has a simple effect on the
generating functions $U(z)$ for unrooted trees and $T(z)$ for rooted trees.
Indeed, we have
$$\vartheta\,U(z)\,=\,T(z)\,,\eqno(4.4)$$
because an unrooted tree with a marked vertex is the same as a rooted tree.
Furthermore
$$\vartheta\,T(z)\,=\,\sum_{k\ge1}T(z)^k\,=\,{T(z)\over1-T(z)}\,,\eqno(4.5)$$
because a rooted tree with a marked vertex is combinatorially
equivalent to an ordered
sequence $\langle T_1,T_2,\ldots,T_k\rangle$ of rooted trees, for some
$k\ge1$. The sequence represents a path of length~$k$ from the marked
vertex to the root, with rooted subtrees sprouting from each point on that
path.

Now let $U=U(wz)/w$ be the function $C_{-1}(w,z)$ that appears in
(4.1), and let $T=T(wz)=C_0(w,z)$. We have
$$\openup1\jot
\eqalign{\vartheta_z U&=z{\partial\over\partial z}{U(wz)\over w}
 = z{wU'(wz)\over w}=z{T(wz)\over wz}={T\over w}\,;\cr
\vartheta_w U&=w{\partial\over\partial w}{U(wz)\over w}
 = w\biggl({zU'(wz)\over w}-{U(wz)\over w^2}\biggr)={T-U\over w}
 ={T^2\over 2w}\,.\cr}$$
Thus
$${2\over w}\,\vartheta_w\,U\,=\,(\vartheta_z\,U)^2\,.\eqno(4.6)$$
In words: ``Orienting and discounting an edge of an unrooted tree is
equivalent to constructing an ordered pair of rooted trees.''

We are now ready to convert (4.2) into a differential equation satisfied
by~$F=F(w,z)$:
$$
\eqalign{\vartheta_w G&=\vartheta_w(e^UF)=(\vartheta_we^U)F+e^U(\vartheta_wF)
=e^U\bigl((\vartheta_wU)F+\vartheta_wF\bigr)\,;\cr
\vartheta_zG&=\vartheta_z(e^UF)=e^U\bigl((\vartheta_zU)F+\vartheta_zF\bigr)\,;\cr
\vartheta_z^2G&=e^U\bigl((\vartheta_z^2U)F+(\vartheta_zU)^2F+
2(\vartheta_zU)(\vartheta_zF)+\vartheta_z^2F\bigr)\,.\cr}$$
Therefore, using (4.6), we have
$${2\over w}\vartheta_wF\;=\;
(\vartheta_z^2U)F\,+\,2(\vartheta_zU)(\vartheta_zF)\,+\,\vartheta_z^2F
\,.\eqno(4.7)$$
And like our other formulas, this one makes combinatorial sense as well
as algebraic sense: The left side tells us that the right side should
yield all ways that the cyclic part of a multigraph can grow, since
${2\over w}\vartheta_wF$ is the number of ways it can go backward one step.
The first term on the right corresponds to marking two vertices of an
unrooted tree (in the acyclic part of the multigraph); joining them will
produce a unicyclic component, thereby increasing the number of components
in~$F$. The middle term corresponds to marking a vertex in some tree of
the acyclic part and another vertex in the cyclic part; joining them will
add new vertices to one of $F$'s existing components. The remaining term
corresponds to marking two vertices in the cyclic part. If such marked
vertices belong to the same component, say a component of excess~$r$,
a new edge between them will change the excess of the component to~$r+1$.
Otherwise, the marked vertices belong to different components, having
respective excesses $r$ and~$s$, possibly with $r=s$; joining them will
merge the components into a new component of excess $r+s+1$.

Similarly, we can proceed to study the bgf $E(w,z)$ for the complex
part of a multigraph, the part whose components all have positive excess.
(The letter~$E$ stands for excess.) We have
$$F(w,z)=e^{V(wz)}\,E(w,z)\,,\eqno(4.8)$$
where $V=V(wz)$ generates unicyclic components. It is easy to verify the
identity
$${2\over w}\vartheta_w V\;=\;\vartheta_z^2U\,+\,2(\vartheta_z U)(\vartheta_zV)
\,,\eqno(4.9)$$
which corresponds to a combinatorially evident fact. Indeed,
$$\vartheta_z^2U={1\over w}{T\over1-T}\,;\qquad\qquad
\vartheta_w V=\vartheta_z V={T\over2(1-T)^2}\,.\eqno(4.10)$$
Therefore we find
$${2\over w}\vartheta_wE\,=\,(\vartheta_z^2V)E+(\vartheta_zV)^2E+
2(\vartheta_zU)(\vartheta_zE)+2(\vartheta_zV)(\vartheta_zE)+\vartheta_z^2E\,.
\eqno(4.11)$$

\bigbreak\noindent
{\bf 5. Enumerating complex multigraphs.}\enspace
To solve the differential equation (4.11), we can first write it in the
form
$${1\over w}(\vartheta_w-T\vartheta_z)E
 \;=\;\half \,e^{-V}\,\vartheta_z^2\,e^V\,E\,.\eqno(5.1)$$
Now we partition $E=E(w,z)$ into terms of equal excess, as we did for
$C(w,z)$ in (2.11):
$$E(w,z)=\sum_rE_r(w,z)=\sum_rw^rE_r(w z)\,.\eqno(5.2)$$
The univariate generating function $E_r(z)$ represents all complex multigraphs
having exactly $r$ more edges than vertices; in particular, $E_0(z)=1$,
since only the empty multigraph is ``complex'' and has excess~0.
Differentiation yields
$$\openup1\jot
\eqalign{\vartheta_wE(w,z)&=\sum_r\bigl(rw^rE_r(wz)+w^r(\vartheta E_r)(wz)
\bigr)\,,\cr
\vartheta_zE(w,z)&=\sum_r w^r(\vartheta E_r)(wz)\,,\cr}$$
where $(\vartheta E_r)(wz)$ here means $\vartheta_z E_r(z)$ with the
argument~$z$ subsequently replaced by $wz$, namely $wzE'_r(wz)$.
Therefore, if we equate the coefficients of $w^{r-1}$ on both sides of
(5.1) and set $w=1$, we obtain a differential recurrence for the
univariate generating functions $E_r=E_r(z)$:
$$\textstyle(r+\vartheta-T\vartheta)E_r\;=\;
\half \,e^{-V}\,\vartheta^2\,e^V\,E_{r-1}\,.\eqno(5.3)$$

It is convenient to introduce a new variable
$$\zeta={T(z)\over 1-T(z)}\eqno(5.4)$$
and to express $E_r$ in terms of $\zeta$ instead of $z$. Note that
$$1+\zeta={1\over1-T(z)}\,;\qquad
T(z)={\zeta\over1+\zeta}\,;\qquad
z={\zeta\over1+\zeta}\exp\left(-\zeta\over1+\zeta\right)\,.\eqno(5.5)$$
Equation (5.3) now takes the form
$$\textstyle\bigl(r+(1+\zeta)^{-1}\vartheta\bigr)E_r
=\half (1+\zeta)^{-1/2}\vartheta^2(1+\zeta)^{1/2}E_{r-1}\,,\eqno(5.6)$$
since $e^V=1/\bigl(1-T(z)\bigr)^{1/2}=(1+\zeta)^{1/2}$ by (3.4).
We will see later that the variable $\zeta$, which represents an
ordered sequence of one or more rooted trees, has important significance
in the study of graphs and multigraphs.

In the $\zeta$ world, with $\vartheta$ still denoting $z{d\over dz}$,
we have the operator equation
$$\vartheta\cdot f(\zeta)=f'(\zeta)\zeta(1+\zeta)^2+f(\zeta)\vartheta\,,
\eqno(5.7)$$
because
$$z{d\zeta\over dz}={T(z)\over1-T(z)}zT'(z)\biggl({1\over T(z)}+{1\over
1-T(z)}\biggr)={T\over(1-T)^3}=\zeta(1+\zeta)^2\,.$$
Equation (5.7) allows us to commute $\vartheta$ with functions of
$\zeta$. For example, we find
$$\eqalign{(1+\zeta)^{-1/2}\vartheta(1+\zeta)^{1/2}
&=\textstyle(1+\zeta)^{-1/2}\bigl(\half (1+\zeta)^{-1/2}\zeta(1+\zeta)^2
+(1+\zeta)^{1/2}\vartheta\bigr)\cr
&=\textstyle\half \zeta(1+\zeta)+\vartheta\,;\cr}$$
hence (5.6) can be rewritten
$$\textstyle\bigl(r+(1+\zeta)^{-1}\vartheta\bigr)E_r
=\half \bigl(\half \zeta(1+\zeta)+\vartheta\bigr)^2E_{r-1}\,.\eqno(5.8)$$
  To simplify the equation even further, we
seek a function $f_r(\zeta)$ such that
$$\textstyle
\vartheta\cdot f_r(\zeta)=
(1+\zeta)f_r(\zeta)\bigl(r+(1+\zeta)^{-1}\vartheta\bigr)
\,;$$
then the differential equation (5.8) will become
$$\textstyle\vartheta\bigl(f_r(\zeta)E_r\bigr)\,=\,\half (1+\zeta)
f_r(\zeta)\bigl(\half \zeta(1+\zeta)+\vartheta\bigr)^2E_{r-1}\,,$$
which can be solved by integration. According to (5.7), the desired
factor $f_r(\zeta)$ is a solution to
$${f'_r(\zeta)\over f_r(\zeta)}={r\over\zeta(1+\zeta)}
={r\over\zeta}-{r\over1+\zeta}\,;$$
so we let $f_r(\zeta)=\zeta^r(1+\zeta)^{-r}$, which incidentally equals
$T(z)^r$. We have derived the equation
$$\vartheta\biggl({\zeta^r\,E_r\over(1+\zeta)^r}\biggr)=\half 
{\zeta^r\over(1+\zeta)^{r-1}}\biggl(\half \zeta(1+\zeta)+\vartheta\biggr)^2
E_{r-1}\,.\eqno(5.9)$$
This differential
equation determines $E_r$ uniquely when $r>0$, given $E_{r-1}$, since
$\zeta^r$ vanishes when $z=0$.

Now all the preliminary groundwork has been laid, and we are ready to
calculate~$E_r$.
We know that $E_0=1$. A bit of experimentation soon reveals a fairly simple
pattern: We can prove by induction on~$r$ that the solution to (5.9) has
the form
$$E_r(z)=\sum_{d=0}^{2r}e_{rd}(1+\zeta)^r\zeta^{2r-d}
=        \sum_{d=0}^{2r}{e_{rd}T(z)^{2r-d}\over\bigl(1-T(z)\bigr)^{
3r-d}}\,,\eqno(5.10)$$
where the coefficients $e_{rd}$ are rational numbers, and where $e_{r(2r)}=0$
for $r>0$. Let $e_{rd}=0$ when $d<0$ or $d>2r$. Assuming that (5.10) holds
for some~$r$, we use (5.7) and (5.8) to compute
$$\openup1\jot\eqalignno{
A_r&=\textstyle\bigl(\half \zeta(1+\zeta)+\vartheta\bigr) E_r\cr
&=\sum_{d=0}^{2r} e_{rd}(1+\zeta)^r\zeta^{2r-d}\zeta(1+\zeta)^2
\biggl({\half \over1+\zeta}+{r\over1+\zeta}+{2r-d\over\zeta}\biggr)\cr
&=\sum_{d=0}^{2r+1} a_{rd}(1+\zeta)^{r+1}\zeta^{2r+1-d}\,,\cr
a_{rd}&=\textstyle(3r+\half -d)e_{rd}+(2r+1-d)e_{r(d-1)}\,;&(5.11)\cr
\noalign{\medskip}
B_r&=\textstyle\bigl(\half \zeta(1+\zeta)+\vartheta\bigr) A_r\cr
&=\sum_{d=0}^{2r+1} a_{rd}(1+\zeta)^{r+1}\zeta^{2r+1-d}\zeta(1+\zeta)^2
\biggl({\half \over1+\zeta}+{r+1\over1+\zeta}+{2r+1-d\over\zeta}\biggr)\cr
&=\sum_{d=0}^{2r+2} b_{rd}(1+\zeta)^{r+2}\zeta^{2r+2-d}\,,\cr
b_{rd}&=\textstyle(3r+{5\over2}-d)a_{rd}+(2r+2-d)a_{r(d-1)}\,.&(5.12)\cr
}$$
Moreover, the left side of Equation (5.9) is a polynomial,
$$\vartheta\bigl(\zeta^r(1+\zeta)^{-r}E_r\bigr)=\vartheta\sum_{d=0}^{2r}e_{rd}
\zeta^{3r-d}=\sum_{d=0}^{2r}(3r-d)e_{rd}(1+\zeta)^2\zeta^{3r-d}\,.$$
The corresponding polynomial on the right-hand side is
$$\textstyle\half \zeta^r(1+\zeta)^{1-r}
\bigl(\half \zeta(1+\zeta)+\vartheta
\bigr)^2E_{r-1}
=\half \sum_db_{(r-1)d}(1+\zeta)^2\zeta^{3r-d}\,;$$
therefore we can complete the induction proof by setting
$$e_{rd}\;=\;{b_{(r-1)d}\over6r-2d}\,,\qquad\qquad 0\le d\le 2r.\eqno(5.13)$$
It is easy to check that $a_{r(2r+1)}=0$ and $b_{r(2r+2)}=0$, hence
$e_{r(2r)}=0$ when $r>0$.

In particular, $a_{00}=\half $, $b_{00}={5\over4}$, $b_{01}=\half $,
and we obtain
$$\textstyle
E_1(z)=(1+\zeta)\bigl({5\over24}\zeta^2+{1\over8}\zeta\bigr)
=\displaystyle\biggl({5\over24\vphantom{\bigl(\bigr)^3}}
{T(z)^2\over\bigl(1-T(z)\bigr)^3}+{1\over8\vphantom{\bigl(\bigr)^2}}{T(z)
   \over\bigl(1-T(z)\bigr)^2}\biggr)\,.
\eqno(5.14)$$
A complex multigraph of excess 1 must consist of a single bicyclic component,
so $E_1(z)$ is the function we called $W(z)$ in~(3.7).
If our only goal had been to compute $W(z)$, we could of course have gotten
this result easily and directly. The more elaborate machinery above has been
developed so that the generating function
 $E_r(z)$ can readily be computed and analyzed for
larger values of~$r$.

\bigbreak\noindent
{\bf 6. Enumerating complex graphs.}\enspace
For graphs instead of multigraphs, the calculations are more intricate,
but it is instructive to look at them and see how they differ.
As in (4.1) and~(4.8), we separate off the cyclic and complex parts
of the bgf by writing
$$\widehat G(w,z)\;=\;e^{U(wz)/w}\,\widehat F(w,z)\,;
\qquad\qquad\widehat F(w,z)\;=\;e^{\widehat V(wz)}\widehat E(w,z).\eqno(6.1)$$
Adding a new edge to a graph means that we want to mark an unordered
pair of {\it distinct\/} vertices, and the operator corresponding to
this is $\half (\vartheta_z^2-\vartheta_z)$. We must also avoid duplicating
an edge that's already present, so we must also subtract~$\vartheta_w$.
Therefore the differential equation satisfied by $\widehat G$ is not
(4.2) but
$${1\over w}\vartheta_w\widehat G=\biggl({\vartheta_z^2-\vartheta_z\over2}-
\vartheta_w\biggr)\,\widehat G\,;\eqno(6.2)$$
and the integral equation corresponding to (4.3) is
$$\widehat G(w,z)=e^z+\int_0^w
\biggl({\vartheta_z^2-\vartheta_z\over2}-
\vartheta_w\biggr)\,\widehat G(w,z)\,dw\,.\eqno(6.3)$$
A computation similar to our derivation of (4.7) now leads to a
differential equation defining~$\widehat F$:
$${1\over w}\vartheta_w\widehat F\;=\;
\biggl(\biggl({\vartheta_z^2-\vartheta_z\over2}-\vartheta_w\biggr)\,U\biggr)
\,\widehat F
\,+\,(\vartheta_zU)(\vartheta_z\widehat F)\,+\,
\biggl({\vartheta_z^2-\vartheta_z\over2}-\vartheta_w\biggr)\,
\widehat F\,.\eqno(6.4)$$

The analog of (5.1) turns out to be
$${1\over w}\bigl(\vartheta_w\widehat E\,-\,T\vartheta_z\widehat E\bigr)
\;=\;e^{-\widehat V}\,\biggl({\vartheta_z^2-\vartheta_z\over2}
-\vartheta_w\biggr)\,e^{\widehat V}\,\widehat E\,;
\eqno(6.5)$$
converting to univariate generating functions $\widehat E_r(w,z)=w^r
\widehat E_r(wz)$ yields
$$\bigl(r+\vartheta-T\vartheta\bigr)\,\widehat E_r\;=\;
e^{-\widehat V}\,\biggl(1-r+{\vartheta^2-3\vartheta\over2}\biggr)\,
e^{\widehat V}\,\widehat E_{r-1}\,.
\eqno(6.6)$$
Again we multiply by the integration factor $\zeta^r\!/(1+\zeta)^r$, but
the differential equation turns out to be rather messy:
$$\vartheta\,\biggl(\biggl({\zeta\over1+\zeta}\biggr)^r\widehat E_r\biggr)=
\biggl({\zeta^r\over(1+\zeta)^{r-1}}\biggr)\biggl(1-r+
{\zeta^4(10+14\zeta+5\zeta^2)\over8(1+\zeta)^2}+
{\zeta^3-3\zeta-3\over2(1+\zeta)}\vartheta+{\vartheta^2\over2}\biggr)\,
\widehat E_{r-1}\,.\eqno(6.7)$$
At least it is linear, and it allows us to compute $\widehat E_r$ for
small~$r$. It turns out that the solution has the form
$$\widehat E_r=
\sum_{d\ge0}\hat e_{rd}{\zeta^{5r-d}\over(1+\zeta)^{2r}}=
\sum_{d\ge0}\hat e_{rd}{T(z)^{5r-d}\over(1-T(z))^{3r-d}}\,,\eqno(6.8)$$
for appropriate coefficients $\hat e_{rd}$. We have, of course,
$\hat e_{00}=1$ and $\hat e_{0d}=0$ for $d\ne0$. When $r>0$, the
values of $\hat e_{rd}$ satisfy the following recurrence, equivalent
to~(6.7):
$$(3r-d)\hat e_{rd}+(6r-d+1)\hat e_{r(d-1)}=
\sum_{j=0}^6 c_j(r-1,d)\hat e_{(r-1)(d-j)}\,,\eqno(6.9)$$
where
$$\eqalign{c_0(r,d)&=(6r-2d+5)(6r-2d+1)/8\,,\cr
c_1(r,d)&=(132r^2+(166-80d)r+45-50d+12d^2)/4\,,\cr
c_2(r,d)&=(398r^2+(584-220d)r+205-160d+30d^2)/4\,,\cr
c_3(r,d)&=(316r^2+(515-160d)r+207-129d+20d^2)/2\,,\cr
c_4(r,d)&=(279r^2+(484-130d)r+208-112d+15d^2)/2\,,\cr
c_5(r,d)&=(13r-3d+10)(5r-d+5)\,,\cr
c_6(r,d)&=(25r^2+(43-10d)r+18-9d+d^2)/2\,.\cr}\eqno(6.10)$$
It is not at all obvious that this recurrence has a solution. We can use
it to compute $\hat e_{rd}$ for $d=0,1,\ldots,3r-1$, but then the
value of $\hat e_{r(3r-1)}$ must satisfy a nontrivial equation when
we set $d=3r$. To get the values of $\hat e_{rd}$ when $d\ge3r$, we
can start by assuming that $\hat e_{rd}=0$ for $d\ge6r$ and work
backward. We will prove later that the recurrence always does have
a solution, and that the last nonzero coefficient for fixed~$r$
can be completely characterized by an almost unbelievable (but true)
formula: If ${s-2\choose2}\le r<{s-1\choose2}$, then
$$\hat e_{r(5r-s)}={{s\choose2}\choose s+r}{1\over s!}\,.\eqno(6.11)$$
Moreover, $\hat e_{rd}=0$ for all $d>5r-s$. Here is a
table of values for small~$r$, in case the reader would like to check
a computer program that is based on the formulas above:
$$\advance\baselineskip5pt
\vbox{\halign{\hfil$#$&&\ \ \hfil$#$\hfil\cr
d=&0&1&2&3&4&5&6&7&8&9&10\cr
\hat e_{0d}=&{1}\cr
\hat e_{1d}=&{5\over {24}}&{1\over 4}\cr
\hat e_{2d}=&{{385}\over {1152}}&{{175}\over {96}}&{{133}\over {32}}
 &{{79}\over {16}}&{{49}\over {16}}&{5\over 6}&{1\over {24}}\cr
\hat e_{3d}=&{{85085}\over {82944}}&{{5005}\over {512}}
 &{{97097}\over {2304}}&{{7777}\over {72}}&{{43621}\over {240}}
 &{{200561}\over {960}}&{{950569}\over {5760}}&{{14001}\over {160}}
 &{{7021}\over {240}}&{{773}\over {144}}&{3\over 8}\cr}}$$

\bigbreak\noindent
{\bf 7. A surprising pattern.}\enspace
The numbers $\hat e_{rd}$ that characterize cyclic graphs of excess~$r$
do not appear
to have any nice mathematical properties. But when we calculate the
corresponding coefficients $e_{rd}$ for multigraphs, as defined in
(5.11)--(5.13), we run into patterns that cry out for explanation.
For example, here is a table showing the values for small~$r$:
$$\advance\baselineskip5pt
\centerline{\vbox{\halign{\hfil$#$&&\ \hfil$#$\hfil\cr
d=&0&1&2&3&4&5&6&7&8&9\cr
e_{0d}=&1\cr
e_{1d}=&{5\over {24}}&{1\over 8}\cr
e_{2d}=&{{385}\over {1152}}&{{35}\over {64}}&{{91}\over {384}}&{1\over {48}}\cr
e_{3d}=&{{85085}\over {82944}}&{{25025}\over {9216}}&{{23023}\over {9216}}
&{{2849}\over {3072}}&{{19}\over {160}}&{1\over {384}}\cr
e_{4d}=&{{37182145}\over {7962624}}&{{11316305}\over {663552}}
&{{3556553}\over {147456}}&{{3658655}\over {221184}}&{{1656083}\over {294912}}
&{{8723}\over {10240}}&{{1969}\over {46080}}&{1\over {3840}}\cr
e_{5d}=&{{5391411025}\over {191102976}}&{{929553625}\over {7077888}}
&{{7994161175}\over {31850496}}&{{8068525465}\over {31850496}}
&{{341105765}\over {2359296}}&{{327803333}\over {7077888}}
&{{1606891}\over {207360}}&{{140569}\over {245760}}
&{{4043}\over {322560}}&{1\over {46080}}\cr}}}$$
Anybody who has played with integers knows that the numerator of
$e_{32}$, 23023, is equal to $7\cdot11\cdot13\cdot23$; moreover, the denominator
is $9216=2^{10}\cdot3^2$. Further experiments show that
the factorization of, say, $e_{55}$, is
$2^{-18}\cdot3^{-3}\cdot11\cdot13\cdot17\cdot19\cdot47\cdot151$.
The occurrence of so many small prime factors cannot be a coincidence!

It is, in fact, easy to see the pattern in the numbers $e_{r0}$,
which satisfy the recurrence
$$e_{r0}={(6r-1)(6r-5)\over24r}e_{(r-1)0}\eqno(7.1)$$
according to rules (5.11)--(5.13). The numbers $\hat e_{r0}$ also
satisfy the same recurrence, according to (6.9) and~(6.10). Therefore we find
$$e_{r0}=\hat e_{r0}={(6r)!\over2^{5r}3^{2r}(3r)!\,(2r)!}\,.\eqno(7.2)$$
But the recurrence defining $e_{rd}$ for $d>0$ is much more complex,
and we have no a priori reason to expect these numbers to have any mathematical
virtues. The following theorem provides an algebraic explanation of 
what is going~on.

\proclaim
Theorem 1. The numbers $e_{rd}$ defined in (5.10) can be expressed as
$$e_{rd}={(6r-2d)!\,P_d(r)\over2^{5r}\,3^{2r-d}\,
 (3r-d)!\,(2r-d)!}\,,\eqno(7.3)$$
where $P_d(r)$ is a polynomial of degree $d$ defined by the formulas
$$\openup2\jot
\eqalignno{P_d(r)&=[z^d]\,F(z)^{2r-d}\,,&(7.4)\cr
F(z)&=3!\sum_{n\ge0}{(4z)^n\over(n+3)!}={6\over(4z)^3}\biggl(e^{4z}
-{(4z)^2\over2}-4z-1\biggr)\,.&(7.5)\cr}$$

\def\F#1#2#3{(r{+}{#1\over#2}{-}{d\over#3})}
\def\FM#1#2#3{(r{-}{#1\over#2}{-}{d\over#3})}
\def\FF#1{(r{+}1{-}{d\over#1})}
\proof
By the duplication and triplication formulas for the Gamma function, expression
(7.3) can also be written
$$e_{rd}=g_{rd}P_d(r),\qquad
g_{rd}={3^r\,\Gamma\F563\,\Gamma\F123\,\Gamma\F163\over
2^{r+d}\,2\pi\,\Gamma\FF2\,\Gamma\F122}\,.\eqno(7.6)$$
Therefore recurrence equation (5.11) becomes
$$\openup1\jot\eqalignno{
a_{rd}&=\textstyle3\F163g_{rd}P_d(r)
 +2\F122g_{r(d-1)}P_{d-1}(r)\cr
&=\textstyle 3\F163g_{rd}A_d(r)\,,\cr
A_d(r)&=\textstyle P_d(r)+{4\over3}P_{d-1}(r)\,.&(7.7)\cr}$$
Similarly, but without as much cancellation, (5.12) becomes
$$\openup1\jot\eqalignno{
b_{rd}&=\textstyle3\F5633\F163g_{rd}A_d(r)
 +2\FF23\F123g_{r(d-1)}A_{d-1}(r)\cr
&=\textstyle {9\over2}g_{r(d-1)}B_d(r)\,,\cr
B_d(r)&=\textstyle \F563\F122A_d(r)
 +{4\over3}\FF2\F123A_{d-1}(r)\,.&(7.8)\cr}$$
Relation (5.13) becomes
$$(3r+3-d)g_{(r+1)d}P_d(r+1)={9\over4}{\FF3\F563\F123\over
\FF2}g_{r(d-1)}P_d(r+1)=\half b_{rd}\,;$$
hence the original recurrence takes the following form:
$$\textstyle\FF3\F563\F123P_d(r+1)=\FF2B_d(r)\,.\eqno(7.9)$$
The boundary conditions are
$$P_d(r)=0\quad\hbox{for $d<0$};\qquad P_0(r)=1;\qquad
P_{2d}(d)=0\quad\hbox{for $d>0$}.\eqno(7.10)$$

It is by no means obvious
 that a polynomial $P_d(r)$ will satisfy (7.7),
(7.8), and (7.9). The key observation that makes everything work is
that a solution to the simpler recurrence
$$\textstyle(r+\half -{d\over3})P_d(r+\half )=
 (r+\half -{d\over2})A_d(r)\eqno(7.11)$$
suffices to solve the more complex one.
 This new recurrence is sort of a
``half step'' between solutions of (7.7), (7.8), and (7.9); it tells
us about multigraphs 
whose excess is an integer plus~$1\over2$, whatever
that may mean.

A solution to (7.11) in the
 extended domain implies a solution to (7.9).
For we will then have
$$\textstyle\FF3\F563\F123P_d(r+1)=\F563\F123\FF2A_d(r+\half )$$
and
$$\eqalign{\textstyle\FF2B_d(r)&=\textstyle\FF2\F563\F122A_d(r)
+{4\over3}\FF2^2\F123A_{d-1}(r)\cr
&=\textstyle\FF2\F563\F123P_d(r+\half )\cr
&\textstyle\qquad
\qquad\null+{4\over3}\FF2\F563\F123P_{d-1}(r+\half )\,.
\cr}$$
Moreover, $P_d({d\over2})=0$ when $d>0$.

We can solve the simultaneous recurrences (7.7) and (7.11) by
constructing solutions to (7.7) that have the desired form (7.4),
namely
$$P_d(r)=[z^d]\,F(z)^{2r-d}\,,\qquad
A_d(r)=[z^d]\,F(z)^{2r-d}\bigl(1+{\textstyle{4\over
3}}\,z\,F(z)\bigr)\,,$$
and noting that the function $F(z)$ of (7.5) satisfies
$$\vartheta F(z)=4z\,F(z)+3-3F(z)\,.\eqno(7.12)$$
Thus we have
$$\openup2\jot
\eqalign{dP_d(r+{\textstyle\half })
&=[z^d]\,\vartheta\bigl(F(z)^{2r+1-d}\bigr)\cr
&=[z^d]\,(2r+1-d)F(z)^{2r-d}\bigl(4zF(z)+3-3F(z)\bigr)\cr
&=(6r+3-3d)\bigl(A_d(r)-P_d(r+{\textstyle\half })\bigr)\,,\cr}$$
and (7.11) holds.\quad\pfbox

\smallskip
Incidentally, the theory of confluent hypergeometric functions
provides us with alternative expressions for the function $F(z)$ in
(7.5). We have, for example,
$$\openup2\jot
\eqalignno{F(z)&=F(1;\,4;\,4z)=3\int_0^1e^{4zt}(1-t)^2\,dt\cr
&={3e^{4z}\over 64z^3}\,\gamma(3,4z)=3e^{4z}\left({1\over 3\cdot
0!}-{4z\over 4\cdot 1!}+{4^2z^2\over 5\cdot 2!}-{4^3z^3\over 6\cdot
3!}+\cdots\right)\,.&(7.13)\cr}$$
The general theory of [\CP] also allows us to write
$$P_d(r)={2r-d\over 2r}\,[z^d]\,G(z)^{2r}\,,\eqno(7.14)$$
where $G(z)=1+z-{1\over 5}z^2+{2\over 15}z^3-{19\over 175}z^4+{2\over
21} z^5-{2018\over 23625}z^6+\cdots$ is defined implicitly by the
relation
$$G\bigl(z\,F(z)\bigr)=F(z)\,.\eqno(7.15)$$

\proclaim
Corollary. For fixed $d\ge0$ we have
$$e_{rd}={3^r\over2^r}{(r+d-1)!\over2\pi\,d!}\bigl(1+O(r^{-1})\bigr)
\eqno(7.16)$$
as $r\to\infty$. Moreover, $e_{rd}$ is a rational number whose
numerator has at most
$$d+O\bigl(d(\log d)^2\!/\log r\bigr)\eqno(7.17)$$
prime factors
 greater than $6r$, and whose denominator has no prime factors
 greater than~$3r$.

\proof
The obvious bounds
$$\openup2\jot
\eqalignno{%
{2r-d\choose d}=[z^d]\,(1+z)^{2r-d}&\leq[z^d]\,F(z)^{2r-d}\cr
&\leq [z^d]\left({1\over 1-z}\right)^{2r-d}={2r-1\choose
d}&(7.18)\cr}$$
tell us that $P_d(r)=(2r)^d\!/d!+O(r^{d-1})$. Formula (7.16) now
follows from (7.3) and Stirling's approximation. (We will derive a
more precise estimate, suitable when $d$ varies with~$r$, in
section~23 below, Lemma~8.)

 All prime factors greater than~$6r$ must
appear as prime factors of $P_d(r)$. We will prove the
 upper bound (7.17) by showing that $m_dP_d(r)$ is an integer, where
$$m_d=5^{\lfloor d/2\rfloor}6^{\lfloor d/3\rfloor}
7^{\lfloor d/4\rfloor}\ldots
=\prod_{k\ge2}(k+3)^{\lfloor d/k\rfloor}\,.\eqno(7.19)$$
It will follow that the denominator of $P_d(r)$ contains no prime
factors greater than $2r+1$, and that if the numerator contains
 $k$~prime factors greater than~$6r$, we have $(6r)^k<m_dP_d(r)\leq
m_d{2r-1\choose d}
<m_d(2r)^d$; i.e., $k\log6r<d\log2r+\log m_d=d\log2r+
O\bigl(d(\log d)^2\bigr)$.

The coefficient of $z^d$ in any power of $F(z)$ is a sum of terms 
$f_1^{k_1}f_2^{k_2}f_3^{k_3}\ldots\,$, where
$f_j=[z^j]\,F(z)={4\over5}{4\over6}\ldots{4\over(j+3)}$ and
$k_1+2k_2+3k_3+\cdots =d$.
Thus, for example, the factor~7 occurs in the denominator 
of $f_1^{k_1}f_2^{k_2}f_3^{k_3}\ldots$ exactly $k_4+k_5+\cdots \leq
d/4$ times.
It follows that the denominator of $P_d$ is a divisor of~$m_d$.\quad\pfbox

\medskip
The estimate (7.17) can be sharpened for small $d$, because $P_d(r)$
always has $(2r-d)$ as a factor when $d>0$. For example,
$$P_1(r)=2r-1,\qquad P_2(r)={(r-1)(10r-7)\over5},\qquad
P_3(r)={(2r-3)(10r^2-21r+10)\over15}.$$
There are no prime factors $>6r$ when $d\le1$, and there is at most
one when $d\le3$.

 Instead of writing
$$E_r(z)=\sum_{d=0}^{2r}e_{rd}{T(z)^{2r-d}\over\bigl(1-T(z)\bigr)^{3r-d}}\,,$$
it is sometimes convenient to use coefficients $e'_{rd}$ such that
$$E_r(z)=\sum_{d=0}^{2r}{e'_{rd}\over\bigl(1-T(z)\bigr)^{3r-d}}\,.\eqno(7.20)$$
The following table shows that the numbers $e'_{rd}$ tend to alternate
in sign:
$$\advance\baselineskip5pt
\centerline{\vbox{\halign{\hfil$#$&&\ \hfil$#$\hfil\cr
d=&0&1&2&3&4&5&6&7&8\cr
e'_{0d}=&1\cr
e'_{1d}=&{5\over {24}}& -{7\over {24}}& {1\over {12}}\cr
e'_{2d}=&{{385}\over {1152}}& -{{455}\over {576}}& {{77}\over {128}}& 
 -{{43}\over {288}}& {1\over {288}}\cr
e'_{3d}=&{{85085}\over {82944}}& -{{95095}\over {27648}}&
 {{119119}\over {27648}}& -{{201355}\over {82944}}&
 {{38623}\over {69120}}& -{{803}\over {34560}}& 
 -{{139}\over {51840}}\cr
e'_{4d}=&{{37182145}\over {7962624}}& -{{40415375}\over {1990656}}& 
 {{141292151}\over {3981312}}& -{{62775713}\over {1990656}}& 
 {{116866321}\over {7962624}}& -{{15867137}\over {4976640}}& 
 {{850003}\over {4976640}}& {{25129}\over {1244160}}& 
 -{{571}\over {2488320}}\cr
% e'_{5d}&{{5391411025}\over {191102976}}& -{{28816162375}\over {191102976}}& 
%  {{5391411025}\over {15925248}}& -{{13125297185}\over {31850496}}& 
%  {{2064417355}\over {7077888}}& -{{7541601353}\over {63700992}}& 
%  {{3902661763}\over {159252480}}& -{{12011857}\over {7962624}}& 
%  -{{11379953}\over {69672960}}& {{2825899}\over {418037760}}& 
%  {{163879}\over {209018880}}\cr
}}}$$
Again, patterns lurk beneath the surface, and there is a prevalence of small
prime factors; for example,
$-e'_{55}={{7541601353}\over {63700992}}=2^{-18}\cdot3^{-5}
\cdot11\cdot13\cdot17\cdot19\cdot23\cdot31\cdot229$. We can in fact prove
the existence of a pattern similar to that of the original 
coefficients~$e_{rd}$:

\proclaim
Corollary. The numbers $e'_{rd}$ defined in (7.20) can be expressed as
$$e'_{rd}={(6r-2d)!\,Q_d(r)\over2^{5r}\,3^{2r-d}\,
 (3r-d)!\,(2r-d)!}\,,\eqno(7.21)$$
where $Q_d(r)$ is a polynomial of degree $d$ for which
$Q_d\bigl({d\over3}-\half \bigr)=0$ when $d>0$.

\proof
By definition, we have
$$e'_{rd}=\sum_{k=0}^d{2r-k\choose d-k}(-1)^{d-k}e_{rk}\,,\eqno(7.22)$$
because the quantity
$T^{2r-k}=\bigl(1-(1-T)\bigr)^{2r-k}$ contributes ${2r-k\choose
d-k}(-1)^{d-k}$ to the coefficient of $(1-T)^{d-3r}$. Now if we plug in
equations (7.3) and (7.21), we find that
$$\eqalignno{\kern-2em
Q_d(r)&=\sum_{k=0}^d{(-1)^{d-k}P_k(r)\over3^{d-k}(d-k)!}{(6r-2k)!\over
(6r-2d)!}{(3r-d)!\over(3r-k)!}\cr
\noalign{\smallskip}
&=\sum_{k=0}^d\left(-{4\over 3}\right)^{d-k}\!\!{3r-k-\half \choose
d-k}P_k(r)=\sum_{k=0}^d\left(-{4\over 3}\right)^{d-k}\!\!{3r-k+{1\over
2}\choose d-k}A_k(r)\,,&(7.23)\cr}$$
clearly a polynomial in $r$ of degree $\leq d$. In fact, the leading term
is $$\sum_k\left(-{4\over3}\right)^{d-k}{(3r)^{d-k}\over(d-k)!}
{(2r)^k\over k!}={(-2)^dr^d\over d!}\,,$$
so the degree is exactly~$d$. If we set $r={d\over
3}-\half $, the sum reduces to $A_d\bigl({d\over 3}-{1\over
2}\bigr)$, which we know is zero for $d>0$ by (7.11).\quad\pfbox

\medskip
It is interesting to try to compute the coefficients $e'_{rd}$ directly,
by proceeding as we did in section~5 but using the variable $\xi=1+\zeta
=\bigl(1-T(z)\bigr)^{-1}$ in place of~$\zeta$. The calculations are
essentially the same, even slightly simpler, until we get to the analog
of equation (5.13); the recurrences that replace (5.11)--(5.13) are
$$\eqalignno{\textstyle a'_{rd}=(3r{+}\half {-}d)e'_{rd}
&\textstyle\null-(3r{+}{3\over2}-d)e'_{r(d{-}1)}\,;&(7.24)\cr
\textstyle (3r{-}d)e'_{rd}-(2r{+}1{-}d)e'_{r(d{-}1)}
&\textstyle=\half \bigl((3r{-}\half {-}d)a'_{(r{-}1)d}
 -(3r{+}\half {-}d)a'_{(r{-}1)(d{-}1)}\bigr)\,.\quad&(7.25)\cr}$$

It appears to be quite difficult to derive (7.21) directly from these
recurrences. The recurrence for $Q_d(r)$, corresponding to equation
(7.9) for $P_d(r)$, turns out to be
$$\eqalignno{
\textstyle(r-{d\over 3})\FM123Q_d(r)
&=\textstyle(r-{d\over 2})\FM122Q_d(r-1)\cr
&\qquad\textstyle\null+{4\over3}\F163\FM123Q_{d-1}(r)\cr
&\qquad\textstyle\null-4(r-{d\over 2})\FM123(r-{d\over 3})
\FM163^{-1}Q_{d-1}(r-1)\cr
&\qquad\textstyle\null+4\F163\FM123Q_{d-2}(r-1)\,,&(7.26)\cr}$$
and we can proceed to solve it for $d=1$, 2, \dots, if we first multiply
both sides by the summation factor $\Gamma(r-{d\over
3})\Gamma(r-\half -{d\over 3})\Gamma(r{+1}-
{d\over2})^{-1}\Gamma\F122^{-1}$. The equation for $d>0$ then takes the form
$$\openup2\jot
\eqalignno{S_d(r)&\textstyle=S_d(r-1)+g_d(r)+g_d(r-\half )\,,\cr
S_d(r)&={\Gamma(r{+}1-{d\over3})\,\Gamma\F123\over
         \Gamma(r{+}1-{d\over2})\,\Gamma\F122}Q_d(r)\,,\cr
g_d(r)&={\Gamma(r{+}1-{d\over 3})\,\Gamma\F123\over
\Gamma\FF2\,\Gamma\F122}f_d(r)\,,\cr}$$
where $f_d(r)=Q_d(r)-{r-d/2\over r-d/3}Q_d(r-\half )$ is a polynomial of
degree $d-1$. For example, $f_1(r)=-{4\over 3}$ and $f_2(r)={8\over
3}r-{4\over 3}$. There is
apparently no analog of the simple relation (7.11) that made everything
work nicely in the theorem above.

A generating function for $Q_d(r)$, analogous to (7.4), can be found
by analyzing (7.23) more carefully. Let $H(z)$ satisfy
$$H(z)=F\bigl(z\,H(z)^{-1/3}\bigr)=1+z+{\textstyle{7\over
15}}\,z^2+{\textstyle{1\over 15}}\,z^3+\cdots\;;\eqno(7.27)$$
then the elementary theory in [\CP] proves that
$$\bigl(x-{\textstyle{d\over
3}}\bigr)\,[z^d]\,H(z)^x=x\,[z^d]\,F(z)^{x-d/3}\,.\eqno(7.28)$$
Hence, by (7.11) and (7.4),
$$A_d(r)={r+\half -{d\over 3}\over r+\half -{d\over
2}}\,[z^d]\,F(z)^{2r+1-d}=[z^d]\,H(z)^{2r+1-2d/3}\,.\eqno(7.29)$$
And (7.23) can therefore be ``summed'':
$$\eqalignno{Q_d(r)
&=\sum_{k=0}^d\,\left(-{4\over 3}\right)^k{3r-d+k+\half \choose
k}\,A_{d-k}(r)\cr
\noalign{\smallskip}
&=\sum_{k=0}^d\,\left({4\over 3}\right)^k{-3r+d-{3\over 2}\choose k}\,
 [z^{d-k}]\,H(z)^{(-2/3)(-3r-3/2+d-k)}\cr
\noalign{\smallskip}
&=[z^d]\,\left({4\over
3}\,z+H(z)^{-2/3}\right)^{-3r-3/2+d}\,.&(7.30)\cr}$$
In particular,
$$\textstyle Q_0(r)=1;\qquad Q_1(r)=-2(r+{1\over6});\qquad
Q_2(r)=2(r-{1\over6})(r-{1\over5}).$$
Although $Q_1(r)=-A_1(r)$ and $Q_2(r)=A_2(r)$,
 we have $Q_3(r)=-A_3(r)+{16\over135}(r-\half )$.

\bigbreak\noindent
{\bf 8. Sparse components.}\enspace
We can readily compute the univariate generating functions $C_1(z)$,
$C_2(z)$, $C_3(z)$, \dots,~$C_r(z)$ for bicyclic, tricyclic, tetracyclic,
\dots,~$(r+1)$-cyclic components, now that we know the simple form of
$E_1(z)$, $E_2(z)$, $E_3(z)$, \dots,~$E_r(z)$, because of the fact that
$$\sum_{r\ge0}w^r\,E_r\;=\;\exp\biggl(\sum_{r\ge1}w^r\,C_r\biggr).\eqno(8.1)$$
Differentiating this formula with respect to $w$ and equating coefficients
of~$w^{r-1}$ leads to the expression
$$r\,E_r\;=\;\sum_{k=1}^r\,k\,C_k\,E_{r-k}\,,\eqno(8.2)$$
from which we may find $C_r$ by calculating
$$C_r\;=\;E_r\,-\,{1\over r}\sum_{k=1}^{r-1}\,k\,C_k\,E_{r-k}\,.\eqno(8.3)$$
Since we know that $E_r=(1+\zeta)^r\sum_{d=0}^{2r-1}e_{rd}\zeta^{2r-d}$ for
$r>0$, it follows by induction that $C_r$ can be written in the same form,
$$C_r=(1+\zeta)^r\sum_{d=0}^{2r-1}c_{rd}\zeta^{2r-d}\,,\eqno(8.4)$$
for appropriate coefficients $c_{rd}$. (The variable $\zeta$ stands
for $T(z)/\bigl(1-T(z)\bigr)$, as in section~5.)
Indeed, relation (8.3) tells us that we can
compute $c_{rd}$ by evaluating a double sum
$$c_{rd}=e_{rd}-{1\over r}\sum_{k=1}^{r-1}k\sum_jc_{kj}e_{(r-k)(d-j)}\,;
\eqno(8.5)$$
the inner sum here is over the range $\max(0,d+1-2r+2k)\le j\le\min(d,2k-1)$,
which is always nonempty for $0<k<r$ except when $d=2r-1$. We always have
$c_{r(2r-1)}=e_{r(2r-1)}=1/\bigl(2^{r+1}(r+1)!\bigr)$. Here is a table of the coefficients
for for small~$r$:
$$\advance\baselineskip5pt
\centerline{\vbox{\halign{\hfil$#$&&\ \ \hfil$#$\hfil\cr
d=&0&1&2&3&4&5&6&7&8&9\cr
c_{1d}=&{5\over {24}}&{1\over 8}\cr
c_{2d}=&{5\over {16}}&{{25}\over {48}}&{{11}\over {48}}&{1\over {48}}\cr
c_{3d}=&{{1105}\over {1152}}&{{985}\over {384}}&{{1373}\over {576}}
&{{515}\over {576}}&{{223}\over {1920}}&{1\over {384}}\cr
c_{4d}=&{{565}\over {128}}&{{12455}\over {768}}&{{26581}\over {1152}}
&{{12227}\over {768}}&{{2089}\over {384}}&{{9583}\over {11520}}
&{{27}\over {640}}&{1\over {3840}}\cr
c_{5d}=&{{82825}\over {3072}}&{{387005}\over {3072}}&{{371195}\over {1536}}
&{{10154003}\over {41472}}&{{121207}\over {864}}&{{519883}\over {11520}}
&{{1573507}\over {207360}}&{{2597}\over {4608}}&{{803}\over {64512}}
&{1\over {46080}}\cr}}}$$

In applications, the leading coefficients $c_{r0}$ of $C_r$ are the most
important, as are the leading coefficients $e_{r0}$ of $E_r$, because
these govern the dominant asymptotic behavior of $[z^n]\,C_r(z)$ and
$[z^n]\,E_r(z)$. Therefore it is convenient to write
$$c_r=c_{r0}, \qquad\qquad e_r=e_{r0}\,.\eqno(8.6)$$
We have seen in (7.2) that there is a simple way to express the numbers $e_r$
in terms of factorials. The values $c_r$ are then easily computed by
using relation (8.3), but with $c_r$ and $e_r$ substituted respectively for
$C_r$ and $E_r$.

Asymptotically speaking, the values of $c_{rd}$ and $e_{rd}$ are
equivalent when $r$ is large.

\proclaim
Theorem 2. For fixed $d\ge0$ we have
$$c_{rd}=e_{rd}\bigl(1+O(r^{-1})\bigr)
={3^r\over2^r}{(r+d-1)!\over2\pi\,d!}\bigl(1+O(r^{-1})\bigr)\eqno(8.7)$$
as $r\to\infty$.

\proof
We know the asymptotic value of $e_{rd}$ from (7.16). To complete the proof,
we need only show that the double sum in (8.5) is $O_d(e_{rd}/r)$, where
$O_d$ implies a bound for fixed~$d$ as $r\to\infty$.

Since $c_{rd}\le e_{rd}$, each term in the double sum is bounded above by
an absolute constant (depending on~$d$) times
$${3^r\over 2^r}{k\over r}{(k+j-1)!\over j!}{(r-k+d-j-1)!\over(d-j)!}=
  {3^r\over 2^r}{k\over r}{(r+d-2)!\over d!}{d\choose j}\bigg/{r+d-2\choose
 k+j-1}\,.$$
We have ${r+d-2\choose k+j-1}\ge r+d-2$ except when $k=1$ and $j=0$
or $k=r-1$ and $j=d$. Therefore all but one term is $O_d(e_{rd}/r^2)$, and
the exceptional term is $O_d(e_{rd}/r)$. There are $O(rd)$ terms altogether,
so the overall double sum is $O_d(e_{rd}/r)$.\quad\pfbox

\medskip
The simple form (8.4) of $C_r(z)$, the generating function for $(r+1)$-cyclic
multigraphs, makes it possible for us to deduce a formula for the
corresponding graph-based function $\widehat C_r(z)$, which turns out
to be only about 50\% more complicated. In fact, we will prove
a result that applies to the generating functions for infinitely many
models of random graphs, including both $G(w,z)$ and $\widehat G(w,z)$
as special cases.

Our starting point for this calculation is the formal power series relation
$$ \widehat G(w,z)=G(\,\ln(1+w),z/\sqrt{\,1+w}\,)\,.\eqno(8.8)$$
which is an immediate consequence of (2.7) and (2.9). It follows that
$$ \widehat C(w,z)=C(\,\ln(1+w),z/\sqrt{\,1+w}\,)\,.\eqno(8.9)$$
We can therefore obtain a near-polynomial formula for $\widehat C_r(z)$
as a special case of the following result.

\proclaim
Theorem 3. If $f(w)=1+f_1w+f_2w^2+\cdots$
 and $g(w)=1+g_1w+g_2w^2+\cdots$ are
arbitrary formal power series with $f(0)=g(0)=1$, and if
$$\widetilde C(w,z)=C\left(wf(w),z{g(w)\over f(w)}\right)=\sum_rw^r
\widetilde C_r(wz)\,,\eqno(8.10)$$
where $C$ is the bgf (2.10) for connected multigraphs,
then there exist coefficients $\tilde c_{rd}$ such that
$$\widetilde{C}_r(z)=\sum_{d=0}^{3r+2}\tilde c_{rd}\zeta^{3r+2-d}
(1+\zeta)^{-2}=\sum_{d=0}^{3r+2}\tilde c_{rd}{T(z)^{3r+2-d}\over
\bigl(1-T(z)\bigr)^{3r-d}}\eqno(8.11)$$
for all $r>0$.

\proof
Consider Ramanujan's function $Q(n)$ of (3.11), which has the asymptotic
value $\sqrt{\pi n\over2}+O(1)$ as $n\to\infty$. Following Knuth [\AOC], 
we shall say that a function $s(n)$ of the form $p(n)+q(n)Q(n)$
is a {\it semipolynomial\/} when $p$ and~$q$ are polynomials. The
{\it degree\/} of a semipolynomial is computed by assuming that $Q(n)$
is of degree~$1\over2$. For example, $3+2n+(1+n)Q(n)$
is a semipolynomial of degree~${3\over 2}$. More formally,
if $d$ is any nonnegative integer, the semipolynomial $p(n)+q(n)Q(n)$
has degree $\le \half d$ if and only if $p$ has degree $\le \half d$
and $q$ has degree $<\half d$.

The formulas (3.12) of section~3, taken from [\KP], show that
generating functions of the form $F(z)=\sum_{k=1}^d
a_k \big/\bigl(1-T(z)
\bigr)^k$ are precisely those whose coefficients satisfy
$$[z^n]\,F(z)={n^ns(n)\over n!}$$
where $s(n)$ is a semipolynomial of degree $\le\half (d-1)$.

Consider now the expansion
$$\sum_rw^rf(w)^rC_r\bigl(zwg(w)\bigr)=\sum_rw^r\,\widetilde{C}_r(wz)$$
which follows 
from (8.10) and (2.11). We will study how each term on the left contributes
to terms on the right. First, when $r=-1$ we have
$$\eqalign{{U\bigl(zwg(w)\bigr)\over w\,f(w)}
&=\sum_{n\ge 1}{n^{n-2}z^nw^{n-1}(1+g_1w+\cdots\,)^n\over n!\,(1+f_1w+\cdots\,)}\cr
\noalign{\smallskip}
&=\sum_{n\ge 1}{n^{n-2}z^nw^{n-1}\bigl(1+np_0(n)w+np_1(n)w^2
+\cdots\,\bigr)\over n!\,(1+f_1w+\cdots\,)}\cr}$$
where each $p_l(n)$ is a polynomial of degree $\leq l$. The effect is to 
make $\widetilde C_{-1}(z)=U(z)$, and to
contribute a linear combination of $U(z)$, $T(z)$, and 
$\bigl(1-T(z)\bigr)^{-1}$,
$\ldots\,$, $\bigl(1-T(z)\bigr)^{-2l+1}$ to $\widetilde{C}_l(z)$ for
each $l\ge0$. Next, when $r=0$ we have
$$V\bigl(zwg(w)\bigr)=\half \,
\sum_{n\ge 1}n^{n-1}Q(n)z^nw^n\bigl(1+np_0(n)w
+np_1(n)w^2+\cdots\,\bigr)\,;$$
this contributes $V(z)$ to $\widetilde C_0(z)$ and
a linear combination of $\bigl(1-T(z)\bigr)^{-1},\ldots,
\bigl(1-T(z)\bigr)^{-2l}$ to $\widetilde C_l(z)$ for each $l>0$.
Finally, when $r>0$ we have, by (5.11),
$$w^r f(w)^r
C_r\bigl(zwg(w)\bigr)=\sum_{n\geq 0}{n^ns(n)\over n!}z^nw^{n+r}
\bigl(1+np_0(n)w+np_1(n)w^2+\cdots\,\bigr)f(w)^r\,,$$
where $s(n)$ is a semipolynomial of degree $\le {3\over 2}r-\half $. This
contributes a linear combination of $\bigl(1-T(z)\bigr)^{-1},\ldots,
\bigl(1-T(z)\bigr)^{-2l-r}$ to $\widetilde{C}_l(z)$ for each $l\ge r$. The
proof of (8.11) is complete, because $U(z)=
\half \zeta(2+\zeta)/(1+\zeta)^2$ and $T(z)=\zeta/(1+\zeta)$.\quad\pfbox

\medskip
Incidentally, our proof shows that the only contribution to the 
coefficient of the ``leading term''
$T(z)^{3l+2}/\bigl(1-T(z)\bigr)^{3l}$ of $\widetilde{C}_l(z)$ comes
from $C_l(z)$ itself.
Therefore $\widetilde{C}_r(z)$ and $C_r(z)$ have identical leading
coefficients. In particular, $\hat c_{r0}=c_{r0}=c_r$.
We will see below that this gives the same asymptotic characteristics
to the limiting distribution of component types in the uniform and permutation
models when $m\approx \half n$.

Theorem 3 justifies our earlier assertion that the recurrence (6.9)--(6.10)
for $\hat e_{rd}$ has a solution. The coefficients $\hat c_{rd}$ can be
computed from those coefficients $\hat e_{rd}$ using the relation
$\widehat C_r=\widehat E_r-{1\over r}\sum_{k=1}^{r-1}k\widehat C_k
\widehat E_{r-k}$; but that makes $\widehat C_r$ a polynomial of
degree~$5r$ with denominator $(1+\zeta)^{2r}$, so the numerator and
denominator must be divided by $(1+\zeta)^{2r-2}$. A simpler recurrence
for $\widehat C_r$ was found by Wright~[\Wi], who proved Theorem~3 in
the special case $\widetilde C_r=\widehat C_r$ by a different method.
Translated into the notation of the present paper, Wright's recurrence is
$$\vartheta\left(\zeta\over1+\zeta\right)^r\widehat C_r=
{\zeta^r\over2(1+\zeta)^{r-1}}\biggl(\,\sum_{j=0}^{r-1}(\vartheta\widehat C_j)
(\vartheta\widehat C_{r-1-j})\;+\;\bigl(\vartheta^2-3\vartheta-2(r-1)\bigr)
\widehat C_{r-1}\biggr)\,,\quad r>0,\eqno(8.12)$$
with $\vartheta\widehat C_0=\half \zeta^3(1+\zeta)^{-1}$. As we saw
for the related sequence $\widehat E_r$ in section~6, it isn't obvious
that this recurrence has a solution of the desired form
$$\widehat{C}_r(z)=\sum_{d=0}^{3r+2}\hat c_{rd}\zeta^{3r+2-d}
(1+\zeta)^{-2}=\sum_{d=0}^{3r+2}\hat c_{rd}{T(z)^{3r+2-d}\over
\bigl(1-T(z)\bigr)^{3r-d}}\,,\eqno(8.13)$$
when $r>0$. Theorem 3 provides an algebraic proof, while Wright proved
the existence by a combination of algebraic and combinatorial methods that
we will consider in the next section. Here is a table of the first
few values of the coefficients:
$$\advance\baselineskip5pt
\centerline{\vbox{\halign{\hfil$#$&&\enspace\hfil$#$\hfil\cr
d=&0&1&2&3&4&5&6&7&8&9&10&11&12\cr
\hat c_{1d}=&{5\over {24}}&{1\over 4}\cr
\hat c_{2d}=&{5\over {16}}&{{55}\over {48}}&{{73}\over {48}}
&{3\over 4}&{1\over {24}}\cr
\hat c_{3d}=&{{1105}\over {1152}}&{{395}\over {72}}
&{{15131}\over {1152}}&{{2399}\over {144}}&{{8303}\over {720}}
&{{557}\over {144}}&{3\over 8}\cr
\hat c_{4d}=&{{565}\over {128}}&{{26165}\over {768}}&{{133651}\over {1152}}
&{{523789}\over {2304}}&{{80573}\over {288}}&{{317611}\over {1440}}
&{{77773}\over {720}}&{{89}\over 3}&{{839}\over {240}}&{1\over {12}}\cr
\hat c_{5d}=&{{82825}\over {3072}}&{{67005}\over {256}}&{{1770535}\over {1536}}
&{{31448897}\over {10368}}&{{438258631}\over {82944}}&{{1146749}\over {180}}
&{{86265}\over {16}}&{{304411}\over {96}}&{{25180997}\over {20160}}
&{{109627}\over {360}}&{{781}\over {20}}&{{439}\over {240}}&{1\over {120}}\cr
}}}$$

Notice that $\hat c_{rd}=0$ for sufficiently large values of~$d$; we do
not have to go all the way up to $d=3r+2$. In fact, we will see in the
next section that the final nonzero coefficient is $\hat c_{r(3r+2-s)}$
when ${s-2\choose 2}\le r<{s-1\choose2}$, and it has the value
exhibited in~(6.11).

The asymptotic value of the leading coefficients $\hat c_{r0}=c_{r0}=c_r$ has
an interesting history. Wright~[\Wiii] gave a complicated argument establishing
that $\hat c_{r0}$ is asymptotically $\bigl({3\over2}\bigr)^r(r-1)!$ times
a certain constant, for which he obtained the numerical value 0.159155.
Stepanov~[\Si] independently computed the numerical value `0,46\dots' for
three times the constant; the approximation 0.48 would have been more
accurate, but Stepanov was perhaps conjecturing that the true value would
be ${1\over3}+{1\over\pi}\bigl(\sqrt3+\ln(2-\sqrt3\,)\bigr)\approx 0.46546$,
which he announced at the same time in connection with another problem
concerning the size of the largest component when the centroid is
removed from a random tree.
Wright's constant was identified as $1/2\pi$ by G.~N. Bagaev and E.~F.
Dmitriev~[\BD], who presented without proof a list of asymptotic expressions
for the solution of several
related enumeration problems. Lambert Meertens independently found a proof
in 1986, but did not publish it; his approach was reported later
in~[\BCM]. A detailed analysis was
also carried out by V.~A. Vobly{\u\i}~[\Vob], who obtained a number of
interesting auxiliary formulas. In particular, if we write $c(z)=c_1z+c_2z^2
+c_3z^3+\cdots\,$, Vobly{\u\i} proved the formal power series relation
$$\vartheta c(z)=-{1\over6}+{1\over3z}\biggl(1-{I_{-2/3}(1/3z)\over
I_{1/3}(1/3z)}\biggr)\,.\eqno(8.14)$$
In other words, he proved that the coefficients $c_r$ show up in the
asymptotic series
$${I_{-2/3}(1/3z)\over I_{1/3}(1/3z)}\sim1-{z\over2}-3c_1z^2-6c_2z^3
-9c_3z^3-\cdots\,,\eqno(8.15)$$
as $z\to0$. This is interesting because the left-hand side can also be
expressed as a continued fraction
$$2z+{1\over\displaystyle 8z+
          {\strut 1\over\displaystyle 14z+
            {\strut 1\over\displaystyle 20z+
              {\strut 1\over\vphantom{{1^1\over1}}26z+\cdots}}}}\,,\eqno(8.16)$$
using the standard recurrence $zI_{\nu+1}(z)=zI_{\nu-1}(z)-2\nu I_\nu(z)$
for the modified Bessel functions~$I_\nu(z)$. In the course of his
investigation, Vobly{\u\i} noticed that the coefficients of $e^{c(z)}$ have
a simple form, although he did not mention their combinatorial significance;
these are the numbers we have called~$e_r$. He gave the formulas
$${2^r\over3^r}e_r=(-1)^r(1/3,r)={\Gamma(r+5/6)\,\Gamma(r+1/6)
\over2\pi r!}\,,\eqno(8.17)$$
which are equivalent to (7.2). Here $(\nu,r)$ denotes Hankel's symbol,
$$(\nu,r)={1\over r!}\prod_{k=1}^r\textstyle(\nu+k-\half )(\nu-k+\half )
\,.$$

\bigbreak\noindent
{\bf 9. Structure of complex multigraphs.}\enspace
The generating functions $E_r$, $C_r$, $(1+\zeta)^{2r}\widehat E_r$, and
$(1+\zeta)^2\widehat C_r$ are polynomials in~$\zeta$, and these polynomials
have a combinatorial interpretation that provides considerable insight into
what is happening as a graph or
multigraph evolves. The inner structure in the case of $\widehat C_r$
was studied by Wright in his original paper [\Wi]; we will see that his
results for graphs become
simpler when we consider the analogous results for multigraphs.

Let $M$ be a cyclic multigraph of excess~$r$, i.e., any multigraph with
no acyclic components, having $r$ more edges than vertices. We
can ``prune''~$M$ by repeatedly cutting off any vertex of
degree~1 and the edge leading to that vertex;
this eliminates as many edges as vertices, so the pruned
multigraph~$\ovline{M}$ still has excess~$r$. Each vertex of~$\ovline{M}$
has degree at least~2. Such multigraphs are called {\it smooth}.

Conversely, given any smooth multigraph~$\ovline{M}$, we obtain all
multigraphs~$M$ that prune down to~it by simply sprouting a
tree from each vertex of~$\ovline{M}$ (i.e., identifying that vertex with
the root of a rooted tree). Since $T(z)$ is the generating function
for rooted trees, it follows that
$$F_r(z)=\ovline{F}_r\bigl(T(z)\bigr)\,,\eqno(9.1)$$
where $F_r(z)$ is the generating function for all cyclic multigraphs of excess~$r$
and $\ovline{F}_r$ is the generating function for all
smooth multigraphs of excess~$r$. Thus, for example, we must have
$$\ovline{F}_1(z)={1\over 24}\,z(3+2z)/(1-z)^{7/2}\,,\eqno(9.2)$$
because we know from (3.4), (4.8), (5.2), and (5.14) that
$$F_1(z)=
e^{V(z)}E_1(z)={1\over
24}\,T(z)\bigl(3+2T(z)\bigr)/\bigl(1-T(z)\bigr)^{7/2}.$$
The coefficient of~$z^n$ in $\ovline{F}_1(z)$ is the sum of
$\kappa(\ovline{M})$ over all multigraphs~$\ovline{M}$ on
$n$~labeled vertices having $n+1$ edges and all vertices of degree~2
or more, divided by~$n!$.  For example, the coefficient of~$z$ is 1/8;
this is the compensation factor of the multigraph with a single vertex~$x$
and two loops from~$x$ to itself. The coefficient of~$z^2$ is
${25\over 48}={25\over 24}/2!$; the smooth labeled multigraphs
$$\unitlength=8pt
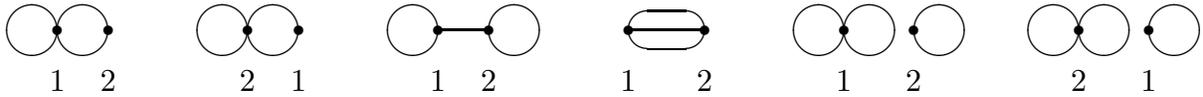
\beginpicture(5,3)(0,0)
\put(2,0){\makebox(0,0){1}}
\put(4,0){\makebox(0,0){2}}
\put(2,2){\disk{.4}}
\put(4,2){\disk{.4}}
\put(1,2){\circle2}
\put(3,2){\circle2}
\endpicture
\qquad
\beginpicture(5,3)(0,0)
\put(2,0){\makebox(0,0){2}}
\put(4,0){\makebox(0,0){1}}
\put(2,2){\disk{.4}}
\put(4,2){\disk{.4}}
\put(1,2){\circle2}
\put(3,2){\circle2}
\endpicture
\qquad
\beginpicture(6,3)(0,0)
\put(2,0){\makebox(0,0){1}}
\put(4,0){\makebox(0,0){2}}
\put(2,2){\disk{.4}}
\put(4,2){\disk{.4}}
\put(1,2){\circle2}
\put(5,2){\circle2}
\put(2,2){\line(1,0)2}
\endpicture
\qquad
\beginpicture(5,3)(0,0)
\put(1,0){\makebox(0,0){1}}
\put(4,0){\makebox(0,0){2}}
\put(1,2){\disk{.4}}
\put(4,2){\disk{.4}}
\put(2.5,2){\oval(3,1.5)}
\put(1,2){\line(1,0)3}
\endpicture
\qquad
\beginpicture(6.75,3)(0,0)
\put(2,0){\makebox(0,0){1}}
\put(4.75,0){\makebox(0,0){2}}
\put(2,2){\disk{.4}}
\put(4.75,2){\disk{.4}}
\put(1,2){\circle2}
\put(3,2){\circle2}
\put(5.75,2){\circle2}
\endpicture
\qquad
\beginpicture(6.75,3)(0,0)
\put(2,0){\makebox(0,0){2}}
\put(4.75,0){\makebox(0,0){1}}
\put(2,2){\disk{.4}}
\put(4.75,2){\disk{.4}}
\put(1,2){\circle2}
\put(3,2){\circle2}
\put(5.75,2){\circle2}
\endpicture
$$
have compensation factors ${1\over 4}$, ${1\over 4}$, ${1\over 4}$,
${1\over 6}$, $1\over16$, and $1\over16$, respectively, summing
to~$25\over24$.

The smooth multigraph $\ovline{M}$ obtained by repeatedly pruning $M$ is
called the {\it core\/} of~$M$ (see [\Lii]). Let $\ovline{F}$ be any
family of smooth multigraphs, and let $F$ be the set of all cyclic
multigraphs whose core is a member of~$\ovline{F}$. The argument
that proves (9.1) also proves that the univariate and bivariate
generating functions for~$F$ and~$\ovline{F}$ are related by the
equations
$$F(z)=\ovline{F}\bigl(T(z)\bigr)\,;\qquad
F(w,z)=\ovline{F}\bigl(w,\,T(wz)/w\bigr)\,.\eqno(9.3)$$
In particular we have
$\widehat E_r(z)=
\ovline{\widehat E}_r\bigl(T(z)\bigr)$, where $\ovline{\widehat E}_r$
counts all smooth graphs of excess~$r$ having no unicyclic components.
This relationship accounts for the curious formula (6.11) about
the last nonvanishing coefficient $\hat e_{rd}$; we can reason as follows:
The minimum number of vertices among all graphs of
excess~$r$, when ${s-2\choose2}\le r<{s-1\choose2}$, is~$s$, because a
graph on $s-1$ vertices has at most $s-1\choose2$ edges and ${s-1\choose2}
<s-1+r$. The coefficient of the minimum power of~$\zeta$ in $\widehat E_r=
\ovline{\widehat E}_r\bigl(\zeta/(1+\zeta)\bigr)$
therefore comes entirely from the ${s(s-1)/2\choose s+r}$ graphs on
$s$ labeled vertices having exactly $s+r$ edges. All such graphs
are smooth.

When $M$ has no unicyclic components
we can go beyond pruning to another kind of vertexectomy
that we will call {\it cancelling\/}: If any vertex has degree~2, we
can remove it and splice together the two edges that it formerly
touched. Repeated application of this process on any smooth
multigraph~$\ovline{M}$ of excess~$r$ will lead to a multigraph
$\MM$ of excess~$r$ in which every vertex has degree~3 or
more. (A~self-loop $\langle x,x\rangle$ is assumed to  contribute~2 to
the degree of~$x$. A~vertex with a self-loop will be connected to at
least one other vertex, because there are no unicycles,
so we will never cancel it.) The
multigraph~$\MM$ can be called {\it reduced}.
Only the middle two multigraphs of the six pictured above are reduced.

There are only finitely many reduced multigraphs of excess~$r$. For
if such a multigraph has $n$~vertices of degrees $d_1,d_2,\ldots,d_n$,
it has $n+r=\half \,(d_1+d_2+\cdots +d_n)\geq {3\over 2}\,n$
edges, hence $n\leq 2r$. The extreme case $n=2r$ occurs if and only if
the multigraph is 3-regular, i.e., every vertex has degree exactly~3.
We will see later that such regularity is, in fact, normal: The complex
components of a random graph or multigraph with $\half @n+o(n^{3/4})$ edges
almost always reduce to components that are 3-regular.

The reduced multigraph $\MM$ obtained by pruning and cancellation from
a given complex multigraph~$M$ is called the {\it kernel\/} of~$M$
(see [\Lii]).
Our immediate goal is to find the generating function for all smooth
multigraphs~$\ovline{M}$ without unicyclic components
that have  a given reduced
multigraph~$\MM$ as their kernel. For this it will be convenient to
introduce another representation of a multigraph~$M$: We label
both the vertices and the edges, and we assign an arbitrary
orientation to each edge,
 thereby obtaining a {\it directed edge-labeled\/}
multigraph. Let $V=V(M)$ be the set of vertex labels and $E=E(M)$ the
set of edge labels. Each edge $e\in E$ has a {\it dual
edge\/}~$\overline{e}$, and
 $\ovline{E}$ is the set of all dual edges. The
multigraph~$M$ is then represented as a mapping~$M$ from 
$E\cup\ovline{E}$
to~$V$, with the interpretation that each directed edge~$e$ runs from
$M(e)$ to $M(\overline{e})$. The dual of~$\overline{e}$, namely 
$\overline{\overline{e}}$,
is~$e$; thus $\overline{e}$ runs from $M(\overline{e})$ to~$M(e)$.

If the vertex labels are $1,\ldots,n$ and if the edge labels are
$1,\ldots,m$, the multigraph mapping~$M$ takes the set
$\{1,\ldots,m,\overline{1},\ldots,\overline{m}\}$ 
into the set $\{1,\ldots,n\}$.
Any such mapping is equivalent to a sequence $\langle
x_1,y_1\rangle\,\ldots,\,\langle x_m,y_m\rangle$ of ordered pairs
generated by the multigraph process of section~1, where $x_k=M(k)$ and
$y_k=M(\overline{k})$. 

The number of different mappings $M$ corresponding to a given 
multigraph~$M$ is $2^mm!\,\kappa(M)$, where $\kappa$ is the
compensation factor defined in
 (1.1).  This holds because $2^mm!$ is the
number of ways to orient the edges and to assign edge labels, 
and $\kappa$ accounts
for duplicate assignments that leave us with the same mapping~$M$. 

Duplicate assignments can be treated
 more formally as follows. A {\it
signed permutation\/}~$\sigma$ of 
a set~$E$ and its dual~$\ovline{E}$ is
a permutation of $E\cup\ovline{E}$ with the property that
$\sigma\overline{e}=\overline{\sigma e}$ for all~$e$. (The group of all
signed permutations on a set of $m$~elements is conventionally called
the hyperoctahedral group~{\bf B}$_m$; it is the group of all $2^mm!$
symmetries of an $m$-cube.) Given a multigraph represented as a
mapping~$M$ from $E\cup\ovline{E}$ to~$V$, an {\it edge
 automorphism\/} is a
signed permutation~$\sigma$ of $E\cup\ovline{E}$
 with the property that $M(\sigma e)=M(e)$.

 It is easy to see that the number of  edge automorphisms of~$M$ is
$1/\kappa(M)$. 
Such a mapping~$\sigma$ must be the product of one of the
$2^{m_{xx}}(m_{xx})!$ signed permutations of the $m_{xx}$~self-loops
from~$x$ to~$x$, for each~$x$, times one of the $(m_{xy})!$
signed permutations of the
 $m_{xy}$~edges from~$x$ to~$y$, for each $x<y$. 
Edge automorphisms are the automorphisms of multigraphs with
labeled vertices and unlabeled edges; this explains why $\kappa(M)$
is used as a weighting function for each~$M$ in the generating
functions we have been discussing.

We are now ready 
to prove a basic lemma about multigraphs, motivated by but
noticeably simpler than the corresponding result for graphs obtained
by Wright [\Wi]:

\proclaim
Lemma 1. If $\MM$ is a reduced multigraph having
$\nu$~vertices, $\mu$~edges, and compensation factor~$\kappa$, the
generating function for all smooth, complex 
 multigraphs~$\ovline{M}$ that reduce to~$\MM$ under cancellation is
$${\kappa\,z^{\nu}\over (1-z)^{\mu}\,\nu!}\,.\eqno(9.4)$$

\proof
This result is ``intuitively obvious,'' but it requires a formal proof
to ensure that everything is counted properly in the presence of
compensation factors. We assume that $\MM$ is represented
by a fixed mapping from edges and dual edges to vertices, where the
set of edge labels is $\{[1],\ldots,[\mu]\}$ and the set of vertex
labels is $\{(1),\ldots,(\nu)\}$. The dual of edge~$[j]$ will be
denoted by $\ovline{[j]}=[-j]$. The given multigraph mapping can be
represented as a function~$M$ from $\{-\mu,\ldots,-1,1,\ldots,\mu\}$
to $\{1,\ldots,\nu\}$, such that edge~$[j]$ runs from
$\bigl(M(j)\bigr)$ to $\bigl(M(-j)\bigr)$ and edge~$[-j]$ runs from
$\bigl(M(-j)\bigr)$ to $\bigl(M(j)\bigr)$. Square brackets and round
parentheses are used notationally here in order to distinguish edge labels
from vertex labels, although $M$ is a function from integers to integers.

Let $s_n$ be the coefficient of $z^n$ in $z^{\nu}/(1-z)^{\mu}$. This
quantity~$s_n$  is the number of solutions $\langle
n_1,\ldots,n_{\mu}\rangle$ to the equation
$$n_1+\cdots +n_{\mu}=n-\nu\eqno(9.5)$$
in nonnegative integers. Let $m-\mu=n-\nu$; then $m$ is the number of
edges in an $n$-vertex multigraph that cancels to~$\MM$.

We will construct $2^mm!\,n!\,s_n/\nu!$ sequences of ordered pairs $\langle
x_1,y_1\rangle\,\ldots\,\langle x_m,y_m\rangle$ of
integers $1\leq x_j,y_j\leq n$ such that (a)~every constructed
sequence defines a smooth multigraph that cancels
to~$\MM$; (b)~every sequence that defines such a smooth
multigraph is constructed exactly $1/\kappa$~times. This will prove the
lemma, because of (2.2). As noted earlier, constructing a sequence $\langle
x_1,y_1\rangle\,\ldots\,\langle x_m,y_m\rangle$ is equivalent to
constructing a map~$\ovline{M}$ from $\{-m,\ldots,-1,1,\ldots,m\}$ into
$\{1,\ldots,n\}$, if we let $x_j=\ovline{M}(j)$ and $y_j=\ovline{M}(-j)$.

The construction is as follows. For each ordered solution $\langle
n_1,\ldots,n_{\mu}\rangle$ to (9.5), we effectively insert $n_j$~new
vertices into edge~$[j]$, thereby undoing the effect of cancellation.
Formally, we construct a set of $m$~edge labels
$$E=\{\,[j,k]\mid 1\leq j\leq \mu,\; 0\leq k\leq n_j\,\}\eqno(9.6)$$
and a set of $n$ vertex labels
$$V=\{\,(i)\mid 1\leq i\leq \nu\,\}\;\;\cup\;\;\{\,(j,k)\mid 1\leq j\leq
\mu,\;
1\leq k\leq n_j\,\}\,.\eqno(9.7)$$
Edge $[j,k]$ runs from vertex $(j,k)$ to vertex $(j,k+1)$, where we
define for convenience
$$(j,0)=\bigl(M(j)\bigr)\,,\qquad
(j,n_j+1)=\bigl(M(-j)\bigr)\,.\eqno(9.8)$$ Thus the original
edge~$[j]$ from $\bigl(M(j)\bigr)$ to $\bigl(M(-j)\bigr)$ has become a
sequence of $n_j+1$ edges $[j,0]\,\ldots\,[j,n_j]$ between the same
two vertices, with intermediate vertices $(j,1),\ldots,(j,n_j)$.

The dual of edge $[j,k]$ will be denoted by $-[j,k]$. We also define
$$[-j,k]=-[j,n_{\vert j\vert}-k]\,,\qquad (-j,k)=(j,n_{\vert
j\vert}+1-k)\,;\eqno(9.9)$$ 
this means that the original edge $[-j]$ has become the edge sequence
$[-j,0]\,\ldots\,[-j,n_j]$, which is the reverse of
$[j,0]\,\ldots\,[j,n_j]$. Edge $[-j,k]$ runs from $(-j,k)$ to
$(-j,k+1)$.

To complete the construction, let $f$ be any one-to-one mapping
from~$V$ to $\{1,\ldots,n\}$ that preserves the order of the original
labels $(1),\ldots,(\nu)$; and let $g$ be any {\it signed
bijection\/} from $\ovline{E}\cup E$ 
to $\{-m,\ldots,-1,1,\ldots,m\}$.
$\bigl($A~signed bijection is a one-to-one correspondence such that
$g(\overline{e})=-g(e)$.$\bigr)$ Then we define
$$\ovline{M}\bigl(g([j,k])\bigr)=f\bigl((j,k)\bigr)\,,\eqno(9.10)$$
for all $[j,k]$ in $\ovline{E}\cup E$. This mapping~$\ovline{M}$
corresponds to a sequence $\langle x_1,y_1\rangle \,\ldots\,\langle
x_m,y_m\rangle$ that defines a multigraph~$\ovline{M}$ on
$\{1,\ldots,n\}$, as stated above. We have constructed
$2^mm!\,n!\,s_n/\nu!$ such sequences, since there are $2^mm!$ choices
for~$g$ and $n!/\nu!$ for~$f$, given any solution $\langle
n_1,\ldots,n_{\mu}\rangle$ to (9.5).

It is clear that 
$\ovline{M}$ is a smooth multigraph on~$n$ vertices that
cancels to the given reduced multigraph~$\MM$,
 and that every such 
$\ovline{M}$ is
constructed at least once. We need to verify that every
mapping~$\ovline{M}$ is obtained exactly $1/\kappa$~times among the
$2^mm!\,n!\,s_n/\nu!$ constructed mappings. 

Suppose $\ovline{M}$ has been constructed from $(\langle
n_1,\ldots,n_{\mu}\rangle,f,g)$, and suppose $\sigma$ is one of the
$1/\kappa$ edge
automorphisms of~$\MM$. We will define a new construction
$(\langle n'_1,\ldots,n'_{\mu}\rangle,f',g')$ that produces the same
mapping~$\ovline{M}$. Our notational conventions allow us to regard
$\sigma$ as a permutation of $\{-\mu,\ldots,-1,1,\ldots,\mu\}$,
where 
$$\sigma(-j)=-\sigma j\qquad{\rm and}\qquad M(\sigma j)=M(j)\,.\eqno(9.11)$$ 
The new construction is defined by
$$\vcenter{\halign{$\hfil#\;$&$#\hfil$\qquad&$#\hfil$\quad&$#\hfil$\cr
n'_j&=n_{\vert\sigma(j)\vert}\,,&1\leq j\leq\mu\,;\cr
\noalign{\smallskip}
f'\bigl((i)\bigr)&=f\bigl((i)\bigr)\,,&1\leq i\leq\nu\,;\cr
\noalign{\smallskip}
f'\bigl((j,k)'\bigr)&=f\bigl((\sigma j,k)\bigr)\,,&1\leq j\leq
\mu\,,&1\leq k\leq n'_j\,;\cr
\noalign{\smallskip}
g'([j,k]')&=g([\sigma j,k])\,,&1\leq j\leq \mu\,,&0\leq k\leq n'_j\,;\cr
\noalign{\smallskip}
\ovline{M}{}'\bigl(g'([j,k]')\bigr)&=f'\bigl((j,k)'\bigr)\,,
&1\leq\vert j\vert\leq\mu\,,&0\leq k\leq n'_{\vert
j\vert}\,.\cr}}\eqno(9.12)$$ 
Here $(j,k)'$ and $[j,k]'$ are the new vertex and edge labels
corresponding to $\langle n'_1,\ldots,n'_{\mu}\rangle$; they are
defined in (9.6)--(9.9).

It is easy to verify that the definitions in (9.12) imply validity of
the same formulas for the whole range of~$j$ and~$k$ values:
$$\vcenter{\halign{$\hfil#\;$&$#\hfil$\qquad&$#\hfil$\quad&$#\hfil$\cr
f'\bigl((j,k)'\bigr)&=f\bigl((\sigma j,k)\bigr)\,,
&1\leq \vert j\vert\leq \mu\,,&0\leq k\leq n'_{\vert
j\vert}+1\,;\cr
\noalign{\smallskip}
g'([j,k]')&=g([\sigma j,k])\,,&1\leq \vert j\vert\leq\mu\,,
&0\leq k\leq n'_{\vert j\vert}\,.\cr}}\eqno(9.13)$$
For example, if $j>0$ we have
$$f'\bigl((j,0)'\bigr)=f'\bigl(\bigl(M(j)\bigr)\bigr)
=f\bigl(\bigl(M(j)\bigr)\bigr) 
=f\bigl(\bigl(M(\sigma j)\bigr)\bigr)=f\bigl((\sigma
j,0)\bigr)\,;$$
\vskip-18pt
$$\eqalign{f'\bigl((j,n'_j+1)'\bigr)=f'\bigl(\bigl(M(-j)\bigr)\bigr)
&=f\bigl(\bigl(M(-j)\bigr)\bigr)\cr
\noalign{\smallskip}
&=f\bigl(\bigl(M\bigl(\sigma(-j)\bigr)\bigr)\bigr)
=f\bigl(\bigl(M(-\sigma
j)\bigr)\bigr)=f\bigl((\sigma_j,n'_j+1)\bigr)\,;\cr
\noalign{\smallskip\medskip}
f'\bigl((-j,k)'\bigr)=
f'\bigl((j,n'_j{+}1{-}k)'\bigr)&=f\bigl((\sigma j,n'_j{+}1{-}k)\bigr)\cr
\noalign{\smallskip}
&=f\bigl((\sigma j,n_{\vert\sigma j\vert}{+}1{-}k)\bigr)
=f\bigl((-\sigma
j,k)\bigr)=f\bigl(\bigl(\sigma(-j),k\bigr)\bigr)\,.\cr}$$
Therefore if $l$ is any value in $\{-m,\ldots,-1,1,\ldots,m\}$, we can
verify that $\ovline{M}{}'(l)=\ovline{M}(l)$, as follows:
There are unique~$j$
and~$k$ such that $l=g([\sigma j,k])$. Hence $l=g'([j,k]')$, and
$$\ovline{M}{}'(l)=f'\bigl((j,k)'\bigr)=f\bigl((\sigma
j,k)\bigr)=\ovline{M}(l)\,.$$ 

Conversely, if $(\langle n'_1,\ldots,n'_{\mu}\rangle,f',g')$ is
another construction that makes $\ovline{M}{}'(l)=\ovline{M}(l)$ for
all~$l$, we can reverse this process and find a unique edge
automorphism~$\sigma$ satisfying all the conditions of
(9.12). Exactly $\nu$ of the
vertices of $\ovline{M}=\ovline{M}{}'$ have degree $\geq 3$, since
$\MM$ is reduced; these are the images
under~$f$ and~$f'$ of $(1),\ldots,(\nu)$, and they have the same order
in~$\ovline{M}$. Therefore $f'\bigl((i)\bigr)=f\bigl((i)\bigr)$ for
$1\leq i\leq\nu$. 

Let $l=g'([j,0])$. Since $\ovline{M}{}'(l)=f'\bigl((j,0)\bigr)=f'
\bigl(\bigl(M(j)\bigr)\bigr)=f\bigl(\bigl(M(j)\bigr)\bigr)$, we
know that $\ovline{M}(l)$ must be a vertex of degree $\geq 3$, so
there must be a value~$j'$ (either positive or negative) such that
$l=g([j',0])$. This rule defines $\sigma j=j'$. We have
$\ovline{M}(l)=f\bigl((j',0)\bigr)
=f\bigl(\bigl(M(j')\bigr)\bigr)$, hence $M(\sigma j)=M(j)$.

Let us say that the edge $[j,k]'$ of $\ovline{M}{}'$ corresponds to the
edge $[j',k']$ of~$\ovline{M}$ if $g'([j,k]')=g[j',k']$. We have
defined $\sigma j$ for $1\leq j\leq \mu$ in such a way that $[j,0]'$
corresponds to $[\sigma j,0]$. Suppose we know that $[j,k]'$
corresponds to $[\sigma j,k]$ for some $k<n'_j$; then $-[j,k]'$ also
corresponds to $-[\sigma j,k]$. Also
$\ovline{M}{}'\bigl(g'(-[j,k]')\bigr) =
\ovline{M}{}'\bigl(g'([-j,n'_j-k]')\bigr) =f'\bigl((-j,{n'_j-k})'\bigr)
=f'\bigl((j,k+1)'\bigr)=\ovline{M}{}'\bigl(g'([j,k+1]')\bigr)$ is a
vertex~$v$ of degree~2 in~$\ovline{M}$, which therefore equals
$\ovline{M}\bigl(g(-[\sigma j,k])\bigr)=f\bigl((-\sigma
j,n_{\vert\sigma j\vert}-k)\bigr)$. Consequently we have $k<n_{\vert\sigma
j\vert}$, $f'\bigl((j,k+1)'\bigr)=f\bigl((\sigma j,k+1)\bigr)$, and
$v=\ovline{M}{}'\bigl(g'([j,k+1]')\bigr)=\ovline{M}\bigl(g([\sigma
j,k+1])\bigr)$. Now $[j,k+1]'$ must correspond to $[\sigma j,k+1]$,
since there is only one value $l\neq -g'([j,k]')$ such that
$\ovline{M}(l)=v$. In this way we prove inductively that $[j,k]'$
corresponds to $[\sigma j,k]$ for $0\leq k\leq n'_j$, and that
$n'_j=n'_{\vert\sigma j\vert}$. Hence (9.12) holds.\quad\pfbox

\medskip
Let $\OF$ be a family of reduced multigraphs, and let $\ovline{F}$
be the family of all smooth complex multigraphs 
that reduce under cancellation to a member of~$\OF$. The bivariate
generating functions of~$\ovline{F}$ and~$\OF$ are then related by
the equation
$$\ovline{F}(w,z)=\OF\bigl(w/(1-wz),z\bigr)\,,\eqno(9.14)$$
because Lemma 1 establishes this relation in the case that $\OF$ has only
one member. Equation (9.14) says simply that every edge in~$\OF$,
represented by~$w$, is to be replaced by a sequence of one or more
edges, represented by $w/(1-wz)=w+w^2z+w^3z^2+\cdots\,$; perhaps this
means that Lemma~1 is indeed obvious and that the lengthy proof was
unnecessary. It is, however, comforting to know that a formal
verification is possible, when one is beginning to learn the power of
generating function techniques. And somehow, examples of multigraphs
with numerous self-loops and repeated edges do seem to mandate a
formal proof, because compensation factors change when edges are
manipulated. 

As an example of Lemma 1, let us derive explicitly the generating
function $\ovline E_1(z)=\ovline{C}_1(z)$ for all smooth
bicyclic multigraphs. All
such multigraphs cancel to a reduced multigraph of excess~1,
which can have at most 2~vertices and 3~edges. There are only three
possibilities,
$$\unitlength=10pt
\beginpicture(4,2)(0,0)
\put(2,1){\disk{.4}}
\put(1,1){\circle2}
\put(3,1){\circle2}
\endpicture
\,,\qquad
\beginpicture(6,2)(0,0)
\put(2,1){\disk{.4}}
\put(4,1){\disk{.4}}
\put(2,1){\line(1,0)2}
\put(1,1){\circle2}
\put(5,1){\circle2}
\endpicture
\,,\qquad
\beginpicture(4,2)(0,0)
\put(.5,1){\disk{.4}}
\put(3.5,1){\disk{.4}}
\put(2,1){\oval(3,1.5)}
\put(.5,1){\line(1,0)3}
\endpicture
\,,\eqno(9.15)$$
having $\kappa={1\over 8}$, ${1\over 4}$, and ${1\over 6}$,
respectively. Therefore
$$\ovline{E}_1(z)=\ovline{C}_1(z)
={z\over 8(1-z)^2}+{z^2\over 8(1-z)^3}+{z^2\over
12(1-z)^3}={z(3+2z)\over 24(1-z)^3}\,,\eqno(9.16)$$
in agreement with (9.2). Wright [\Wi] states that there are 15
connected, unlabeled, 
reduced multigraphs of excess~2, and 107 of excess~3. 

If a reduced multigraph of excess $r$ has exactly $2r-d$ vertices,
we will say that it has {\it deficiency\/}~$d$. A~reduced multigraph
of deficiency~0 is 3-regular; we will call such multigraphs {\it clean}.

\proclaim
Corollary. The coefficient $e_{rd}$ in (5.10) and (7.3) is $(2r-d)!^{-1}
\sum \kappa(\MM)$, summed over all  reduced, labeled  multigraphs
$\MM$ of excess~$r$ and deficiency~$d$.
The coefficient $c_{rd}$ 
in (8.4) can be obtained in the same way, but restricting the sum
to connected multigraphs.\quad\pfbox

\medskip
This corollary leads to a completely different proof of Theorem~1,
because it allows us to obtain formula (7.3) for~$e_{rd}$ by a
combinatorial counting argument. Consider a reduced multigraph that
has exactly $d_k$~vertices of degree~$k$, for each $k\geq 3$; then
$d_3+d_4+\cdots =n$ and $3d_3+4d_4+\cdots =2m$. We can calculate
$\sum\kappa(\MM)$ over all such~$\MM$ by counting the number of relevant
sequences $\langle x_1,y_1\rangle\,\ldots\,\langle x_m,y_m\rangle$ and
dividing by~$2^mm!$; and the number of ways to choose $\langle
x_1,y_1\rangle\,\ldots\,\langle x_m,y_m\rangle$ is clearly a product
of multinomial coefficients,
$${(2m)!\over 3!^{d_3}\,4!^{d_4}\,\ldots\;}\quad{n!\over
d_3!~\,d_4!\,\ldots\;}\,,$$ 
since the first factor is the number of ways to partition $2m$ slots
into $d_k$~labeled classes of size~$k$ for each~$k$, and the second factor
counts the assignments of vertex labels to those classes. To obtain
all reduced multigraphs of excess~$r$ and deficiency~$d$, we sum over
all sequences of nonnegative integers $\langle
d_3,d_4,\ldots\,\rangle$ such that $\sum_{k\geq 3}d_k=2r-d$ and
$\sum_{k\geq 3}kd_k=6r-2d$, or equivalently
$$\sum_{k\geq 3}(k-3)d_k=d\qquad{\rm and}\qquad\sum_{k\geq
3}(k-2)d_k=2r\,.$$ 
Let
$$f_{cd}=\left.\sum\left\{\,\prod_{k\geq 3}\,{1\over
k!^{d_k}\,d_k!}\;\right\vert\;\sum_{k\geq 3}(k-3)d_k=d\quad{\rm
and}\quad
\sum_{k\geq 3}(k-2)d_k=c\right\}\,.\eqno(9.17)$$
We have just proved that
$$e_{rd}={(6r-2d)!\over 2^{3r-d}(3r-d)!}\,f_{(2r)d}\,.\eqno(9.18)$$
And we can readily calculate a bivariate generating function for the
coefficients~$f_{rd}$: 
$$\openup2\jot
\eqalign{\sum_{r,d\geq 0}f_{rd}w^dz^r
&=\sum_{d_3,d_4,\ldots\geq 0}\;\prod_{k\geq
3}\,{w^{(k-3)d_k}z^{(k-2)d_k}\over k!^{d_k}\,d_k!}\cr
&=\prod_{k\geq 3}\;\sum_{d_k\geq 0}\left({w^{k-3}z^{k-2}\over
k!}\right)^{d_k}\,{1\over d_k!}\cr
&=\prod_{k\geq 3}\exp(w^{k-3}z^{k-2}\!/k!)\cr
&=\exp\left(w^{-3}z^{-2}\sum_{k\geq 3}\,{(wz)^k\over k!}\right)
=\exp\biggl({z\over 6}\,F\left({wz\over 4}\right)\biggr)\,,\cr}$$
where $F$ is the function defined in (7.5). Comparing (9.18) to (7.3)
now yields the promised proof of (7.4):
$$\openup2\jot
\eqalign{P_d(r)
&=2^{2r+d}3^{2r-d}(2r-d)!\,f_{(2r)d}\cr
&=2^{2r+d}3^{2r-d}(2r-d)!\,[w^dz^{2r}]\,\exp\bigl(z\,F(wz/4)/6\bigr)\cr
&=2^{2r+d}3^{2r-d}(2r-d)!\,[w^dz^{2r-d}]\exp\bigl(z\,F(w/4)/6\bigr)\cr
&=2^{2d}(2r-d)!\,[w^dz^{2r-d}]\exp\bigl(z\,F(w/4)\bigr)=[w^d]\,F(w)^{2r-d}\,.\cr}$$

These observations also allow us to express $e_{rd}$ in the suggestive
form
$$e_{rd}={1\over 2^{3r-d}(3r-d)!}\,{6r-2d\brace 2r-d}_{\geq
3}\,,\eqno(9.19)$$ 
where ${m\brace n}_{\geq 3}$ denotes the number of ways to partition
an $m$-element set into $n$~subsets, each containing at least
3~elements. The asymptotic behavior of the integers
$2^{3r-d}(3r-d)!\,e_{rd}$ will therefore be analogous to the
 asymptotic behavior of Stirling numbers. 

Lemma 1 captures the combinatorial essence of the generating functions
for all complex multigraphs. We can
 obtain a similar generating function for  graphs instead of multigraphs,
but we must work a bit harder, and the formulas are not as attractive.
The following improvement over Wright's
original treatment~[\Wi] is based on an approach suggested by V.~E.
Stepanov~[\Sii].

\proclaim
Lemma 2. Let $\MM$ be a reduced multigraph having
$\nu$~vertices, $\mu$~edges, compensation factor~$\kappa$, and $\mu_{xy}$
edges between $x$ and~$y$ for $1\le x\le y\le\nu$. The
generating function for all smooth, complex  graphs~$\ovline{G}$ that lead
to~$\MM$ under cancellation is
$${\kappa\,z^{\nu}\over (1-z)^{\mu}\,\nu!}\,P(\MM,z)\,,
\eqno(9.20)$$
where
$$P(\MM,z)=\prod_{x=1}^\nu\left(z^{2\mu_{xx}}
\prod_{y=x+1}^\nu z^{\mu_{xy}-1}\bigl(\mu_{xy}-(\mu_{xy}-1)z\bigr)\right)
\eqno(9.21)$$
is a polynomial in $z$ such that $P(\MM,1)=1$.

\proof
We argue as in Lemma 1, but we must restrict the solutions $\langle n_1,\ldots,
n_\mu\rangle$ of (9.5) to cases that produce a graph instead of a multigraph.
Thus, each $n_j$ that corresponds to a self-loop must be $\ge2$, so we
use $z^2\!/(1-z)$ instead of $1/(1-z)$ in the contribution that $n_j$ makes
to the overall generating function. A subsequence $\langle n_j,\ldots,
n_{j+k-1}\rangle$ that corresponds to $k=\mu_{xy}$ edges between
distinct vertices $x<y$ must have the property that at most one of
$\langle n_j,\ldots,n_{j+k-1}\rangle$ is zero; hence we use
$${z^k\over(1-z)^k}+{kz^{k-1}\over(1-z)^{k-1}}
\;=\;{z^{k-1}\bigl(k-(k-1)z\bigr)\over(1-z)^k}$$
instead of $1/(1-z)^k$ in its contribution. The net effect is to
multiply the previous generating function by $P(\MM,z)$.\quad\pfbox

\medskip
Replacing $z$ by $T(z)$ gives the generating function for all graphs
that prune and cancel to $\MM$. For example,
the generating function $\widehat E_1(z)=\widehat C_1(z)=\widehat W(z)$
of (3.6) can be read off from~(9.15): It is
$${T(z)^5\over8\bigl(1-T(z)\bigr)^2}+
  {T(z)^6\over8\bigl(1-T(z)\bigr)^3}+
  {T(z)^4\bigl(3-2T(z)\bigr)\over 12\bigl(1-T(z)\bigr)^3}\,.\eqno(9.22)$$

The degree of the polynomial $P(\MM,z)$ is the total number of ``penalty
points'' of $\MM$, where each self-loop costs two penalty points, and where
each cluster of $\mu_{xy}>1$ multiple edges between distinct vertices
costs~$\mu_{xy}-1$. If $\MM$ is a graph, the degree is zero and
$P(\MM,z)=1$. At the other extreme, if all edges of~$\MM$ are self-loops,
the degree is~$2\mu$.

The quantity $T(z)^\nu/\bigl(1-T(z)\bigr)^\mu$ becomes $\zeta^\nu
(1+\zeta)^{\mu-\nu}$, when we express it in terms of the variable
$\zeta=T(z)/\bigl(1-T(z)\bigr)$ introduced in section~5;
the quantity 
$P\bigl(\MM,T(z)\bigr)$ becomes $P\bigl(\MM,\zeta/(1+\zeta)\bigr)$.
If we restrict consideration to connected multigraphs of excess~$r$, we
get rational functions of $\zeta$ with denominator $(1+\zeta)^{r+2}$;
this denominator occurs when there are $(r+1)$ self-loops in~$\MM$.
However, we have seen in Theorem~3 that the
denominator of $\widehat C_r$ is always a divisor of $(1+\zeta)^2$.
There seems to be no easy combinatorial explanation for the cancellation
that occurs when the contributions of different~$\MM$ are added together.
Some of the properties of connected graphs are easier to derive by
combinatorics, others are easier to derive by algebra.

The actual coefficients of $P\bigl(\MM,\zeta/(1+\zeta)\bigr)$
do not make any significant difference asymptotically, when graphs are sparse;
we will see later that the asymptotic behavior
as $\zeta\to\infty$ is what counts,
hence we only need to know that $P(\MM,1)=1$. We observed earlier that
the leading coefficients $\hat e_{r0}$ and $e_{r0}$ of $\widehat E$ and~$E$
are equal, as are the leading coefficients $\hat c_{r0}$ and $c_{r0}$. Now
Lemma~2 shows in fact that each reduced multigraph~$\MM$ makes the
same contribution to the leading coefficient for graphs as it
does for multigraphs.

\bigbreak\noindent
{\bf 10. A lemma from contour integration.}\enspace
Studies of random graphs that have $m\approx\half n$ edges are traditionally
broken into two cases, the ``subcritical'' case where $m<\half n$
and the ``supercritical'' case where $m>\half n$. It is desirable, however,
to have estimates of probabilities that hold uniformly for all~$m$
in the vicinity of $\half n$, passing smoothly from one side to the
other. The following lemma, based on techniques introduced in~[\FKP],
will be our key tool for the computation of probabilities.

\proclaim
Lemma 3. If $m=\half n(1+\mu n^{-1/3})$ and if $y$ is any real constant,
we have
$${2^m\,m!\,n!\over (n-m)!\,n^{2m}}\,[z^n]\,{U(z)^{n-m}\over 
 \bigl(1-T(z)\bigr)^y}
=\sqrt{\mskip1mu2\pi}\,A(y,\mu)\,n^{y/3-1/6}+O\bigl((1+
\vert\mu\vert^B)n^{y/3-1/2}\bigr)
\eqno(10.1)$$
uniformly for $|\mu|\le n^{1/12}$, where $B=\max(4,\,{9\over 2}-y)$ and
$$A(y,\mu)={e^{-\mu^3\!/6}\over3^{(y+1)/3}}\sum_{k\ge0}
{\bigl(\half 3^{2/3}\mu\bigr)^k\over k!\,\Gamma\bigl((y+1-2k)/3\bigr)}
\,.\eqno(10.2)$$
As $\mu\to-\infty$, we have
$$A(y,\mu)={1\over\sqrt{\mskip1mu2\pi}\,|\mu|^{y-1/2}}
\biggl(1-{3y^2+3y-1\over 6|\mu|^3}
+O(\mu^{-6})\biggr)\,;\eqno(10.3)$$
as $\mu\to+\infty$, we have
$$A(y,\mu)={e^{-\mu^3\!/6}\over2^{y/2}\mu^{1-y/2}}
\biggl({1\over\Gamma(y/2)}+{4\mu^{-3/2}\over3\sqrt2\,\Gamma(y/2-3/2)}
+O(\mu^{-2})\biggr)\,.\eqno(10.4)$$
Moreover, (10.1) can be improved to
$${2^m\,m!\,n!\over (n-m)!\,n^{2m}}\,[z^n]\,{U(z)^{n-m}\over 
 \bigl(1-T(z)\bigr)^y}
=\sqrt{\mskip1mu2\pi}\,A(y,\mu)\,n^{y/3-1/6}\bigl(1+O(\mu^4n^{-1/3})\bigr)
\eqno(10.5)$$
if $|\mu|$ goes to infinity with $n$ while remaining $\leq n^{1/12}$.

\proof
First we need to derive some auxiliary results about the function~$A$. If
$\alpha$ is any positive number, we define a path $\Pi(\alpha)$
in the complex plane that consists of the following three straight line
segments:
$$\openup1\jot
s(t)=\cases{-e^{-\pi i/3}\,t,&for $-\infty<t\le-2\alpha$;\cr
\alpha+it\sin\pi/3,&for $-2\alpha\le t\le+2\alpha$;\cr
e^{+\pi i/3}\,t,&for $+2\alpha\le t<+\infty$.\cr}\eqno(10.6)$$
Now we define
$$A(y,\mu)={1\over2\pi i}\int_{\Pi(1)} s^{1-y}e^{K(\mu,s)}\,ds\,,
\eqno(10.7)$$
where $K(\mu,s)$ is the polynomial
$$K(\mu,s)={(s+\mu)^2(2s-\mu)\over6}=
{s^3\over3}+{\mu s^2\over2}-{\mu^3\over6}\,.\eqno(10.8)$$
Our first goal is to show that $A(y,\mu)$ satisfies (10.2), (10.3), and~(10.4).

To get (10.2), we make the substitution $u=s^3\!/3$. As $s$
traverses~$\Pi(1)$, the variable~$u$ traverses an interesting
contour~$\Gamma$ that begins at $-\infty$ and hugs the lower edge of the
negative axis, then circles the origin counterclockwise and returns
to~$-\infty$ along the upper edge of the axis. On this contour~$\Gamma$ we
have Hankel's well-known formula for the reciprocal Gamma function,
$${1\over \Gamma(z)}={1\over 2\pi i}\int_{\Gamma}{e^u\,du\over u^z}\,.$$
(See, for example,
[\Hen, Theorem 8.4b].) So we can expand (10.7) into an absolutely convergent
series, after substituting $3^{1/3}u^{1/3}$ for~$s$:
$$\openup3\jot
\eqalign{\int_{\Pi(1)} s^{1-y}e^{K(\mu,s)}\,ds
&={e^{-\mu^3\!/6}\over3^{(y+1)/3}}\int_\Gamma
{e^u\exp\bigl(\half 3^{2/3}\mu u^{2/3}\bigr)\,du\over u^{(y+1)/3}}\cr
&={e^{-\mu^3\!/6}\over3^{(y+1)/3}}\int_\Gamma
\sum_{k\ge0}{\bigl(\half 3^{2/3}\mu\bigr)^ke^u\,du\over k!\, u^{(y+1-2k)/3}}\,.\cr
}$$
Interchanging summation and integration, and applying Hankel's
formula,  gives~(10.2).

To get (10.3) and (10.4), we note first that the integral (10.7) can be
taken over any path $\Pi(\alpha)$, not just $\Pi(1)$, because
$e^{K(\mu,s)}$ has no singularities. Moreover, we can ``straighten out''
the path $\Pi(\alpha)$, changing it to a single straight line from
$\alpha-i\infty$ to $\alpha+i\infty$, if $\alpha$ is sufficiently
large. For we can readily verify that the integrand is exponentially small on
any large circular arc $s=Re^{i\theta}$, as $|\theta|$ increases from
$\pi/3$ to the angle where $R\cos\theta=\alpha$: The real part of
$s^3$ is $R^3\cos3\theta$, which increases from $-R^3$ to $4\alpha^3-
3R^2\alpha$; and the real part of~$s^2$ lies between $-R^2$ and
$-R^2\!/2$. Hence the real part of $K(\mu,s)$ will be at most $-cR^2$
for some positive $c=c(\alpha)$ on the entire arc, whenever $\alpha>0$ and
$\alpha>-\half \mu$; this will make $s^{1-y}e^{K(\mu,s)}$ exponentially
small.

If $\mu$ is negative, let $\alpha=-\mu$; then
$$\openup5\jot
\eqalignno{A(y,-\alpha)
&={1\over2\pi}\int_{-\infty}^\infty(\alpha+it)^{1-y}e^{K(-\alpha,
\alpha+it)}\,dt\cr
&={1\over2\pi\sqrt\alpha}\int_{-\infty}^\infty(\alpha+it/\sqrt\alpha\,)^{1-y}
e^{K(-\alpha,\alpha+it/\sqrt\alpha\,)}\,dt\cr
&={1\over2\pi\alpha^{y-1/2}}\int_{-\infty}^\infty\biggl(1+{it\over
 \alpha^{3/2}}\biggr)^{1-y}
e^{-t^2\!/2-it^3\!/(3\alpha^{3/2})}\,dt\,,&(10.9)\cr
}$$
and we can find the asymptotic value of the remaining integral by
using Laplace's standard technique of 
``tail-exchange'' (see [\CM, section~9.4]):
$$\openup3\jot
\eqalign{&
\int_{-\infty}^\infty\biggl(1+{it\over\alpha^{3/2}}\biggr)^{1-y}
 e^{-t^2\!/2-it^3\!/(3\alpha^{3/2})}\,dt\cr
&\hskip5em=\int_{-\alpha^\epsilon}^{\alpha^\epsilon}
\biggl(1+{it\over\alpha^{3/2}}\biggr)^{1-y}
 e^{-t^2\!/2-it^3\!/(3\alpha^{3/2})}\,dt\;
+\;O\bigl(e^{-\alpha^{2\epsilon}\!/3}
\bigr)\cr
&\hskip5em=\int_{-\alpha^\epsilon}^{\alpha^\epsilon}
e^{-t^2\!/2}\biggl(1+{(1-y)\mkern1mu it\over
\alpha^{3/2}}-{it^3\over 3\alpha^{3/2}}
+O(\alpha^{6\epsilon-3})\biggr)dt\;+\;O\bigl(e^{-\alpha^{2\epsilon}\!/3}
\bigr)\cr
&\hskip5em=\sqrt{\mskip1mu2\pi}\;
+\;O\bigl(\alpha^{6\epsilon-3}\bigr)\,.\cr}$$
If we expand the integrand further,
 to terms that are $O(\alpha^{12\epsilon
-6})$, we obtain
$$A(y,-\alpha)={1\over\sqrt{\mskip1mu2\pi}\alpha^{y-1/2}}
\biggl(1\,-\,{3y^2+3y-1\over6\alpha^3}\,+\,O\bigl(\alpha^{12\epsilon
-6}\bigr)\biggr)\,.$$
The method can clearly be extended, in principle,
to give a complete
 asymptotic series in powers of~$\alpha^{-3}$, beginning
as shown in (10.3).

We also want to know the
 asymptotic value of $A(y,\mu)$ as $\mu\to+\infty$,
and for this we need to work a bit harder. A combination of the methods
we have used to prove (10.2) and (10.3) will establish (10.4). The idea
now is to integrate on the path $\mu^{-1}+it/\sqrt\mu\,$:
$$\openup3\jot
\eqalign{A(y,\mu)
&={e^{K(\mu,\mu^{-1})}\over 2\pi\sqrt\mu}
\int_{-\infty}^\infty\biggl(\mu^{-1}+{it\over\sqrt\mu}\biggr)^{\!1-y}
\exp\bigl(\textstyle it(\mu^{-1/2}+\mu^{-5/2})\cr
&\textstyle\hskip22em{}\null-t^2\bigl(\half +\mu^{-2}\bigr)
-{1\over3}it^3\mu^{-3/2}\bigr)\,dt\cr
&={e^{K(\mu,\mu^{-1})}\over2\pi\mu^{1-y/2}}
\int_{-\infty}^\infty(\mu^{-1/2}+it)^{1-y}
e^{-t^2\!/2}g(it,\mu)\,dt\cr
&={e^{K(\mu,\mu^{-1})}\over2\pi i\mu^{1-y/2}}
\int_{-\infty i}^{\infty i}(v+\mu^{-1/2})^{1-y}
e^{v^2\!/2}g(v,\mu)\,dv\,,\cr
}$$
where the last step replaces $it$ by $v$. We can distort the path of~$v$
so that it crosses the positive real axis, and then replace $v^2\!/2$ by~$u$
to get Hankel's contour~$\Gamma$ again:
$$\openup3\jot
\eqalign{A(y,\mu)
&={e^{K(\mu,\mu^{-1})}\over2\pi i\mu^{1-y/2}}
\int_\Gamma(\sqrt{\mskip1mu2u}+\mu^{-1/2})^{1-y}
e^ug(\sqrt{\mskip1mu2u},\mu)\,{du\over\sqrt{\mskip1mu2u}}\cr
&={e^{K(\mu,\mu^{-1})}\over2^{1+y/2}\pi i\mu^{1-y/2}}
\int_\Gamma(1+(2\mu u)^{-1/2})^{1-y}
u^{-y/2}e^ug(\sqrt{\mskip1mu2u},\mu)\,du\,.\cr
}$$

For definiteness we can stipulate 
that the contour $\Gamma$ lies entirely on the
negative axis, except for a circular loop about~0 with a radius of~1.
When $u$ is on the negative axis, 
say $u=-t$, the quantity $\sqrt{\mskip1mu2u}$ will
be $-i\sqrt{\mskip1mu2t}$ 
on the first part of~$\Gamma$ and $+i\sqrt{\mskip1mu2t}$ on the~last,
so we will have
$$g(\sqrt{\mskip1mu2u},\mu)=\exp\bigl(\textstyle
\mp i\sqrt{\mskip1mu2t}\mkern1mu(\mu^{-1/2}+\mu^{-5/2})
-2t\mu^{-2}\pm{1\over3}i(2t)^{3/2}
 \mu^{-3/2}\bigr)\,.$$
On the portions of $\Gamma$ for which $|u|\ge\mu^\epsilon$, the integrand
is superpolynomially small;\footnote{*}{``Superpolynomially small'' means that
it approaches zero faster than any negative power of the argument.} hence
$$\int_\Gamma\;=\;\int_{\Gamma[\mu^\epsilon]}\,+\,O\bigl(e^{-\mu^\epsilon\!/2}
\bigr)\,,$$
where $\Gamma[\mu^\epsilon]$ is the subcontour that runs along the
lower edge of the negative axis from $-\mu^\epsilon$ to the circle
$u=e^{i\theta}$ and back to $-\mu^\epsilon$ on the top edge
of the axis. On~$\Gamma[\mu^\epsilon]$ we have
$$\eqalign{(2\mu u)^{-1/2}&=O(\mu^{-1/2})\,,\cr
g(\sqrt{\mskip1mu2u},\mu)&=
1+\sqrt{\mskip1mu2u}\,\mu^{-1/2}+O(\mu^{\epsilon-1})\,;\cr}$$
and $\int_{\Gamma}\vert u^{-y/2}e^u\vert\,du$ exists. Hence
$$\openup2\jot
\eqalign{&\int_{\Gamma[\mu^\epsilon]}\bigl(1
+(2\mu u)^{-1/2}\bigr)^{1-y}u^{-y/2}
e^u g\bigl(\sqrt{\mskip1mu2u},\mu\bigr)\,du\cr
&\qquad=\int_{\Gamma[\mu^\epsilon]}\bigl(u^{-y/2}
+{1\over\sqrt{\mathstrut\smash{2\mu}}}\bigl(
(1-y)u^{-(y+1)/2}
+2u^{-(y-1)/2}\bigr)\bigr)e^u\,du+O(\mu^{\epsilon-1})\cr
&\qquad=2\pi i\biggl({1\over\Gamma(y/2)}+
{1\over\sqrt{\mathstrut\smash{2\mu}}}\biggl({1-y\over
\Gamma\bigl((y+1)/2\bigr)}
+{2\over\Gamma\bigl((y-1)/2\bigr)}\biggr)\biggr)+
O(\mu^{\epsilon-1})\,.\cr}$$
The coefficient of
 $\mu^{-1/2}$ vanishes, because $\Gamma\bigl((y+1)/2\bigr)
=\half (y-1)\Gamma\bigl((y-1)/2\bigr)$. 
We can use the same method to
expand the integrand further, obtaining~(10.4).

Notice that
 $1/\Gamma(y/2)$ or $1/\Gamma(y/2-3/2)$ may be zero, but not both.
Therefore (10.4) gives the asymptotically leading term of $A(y,\mu)$
in all cases.

Whew---we have worked
 pretty hard to establish (10.2)--(10.4), and we still
haven't begun to 
tackle the main assertion of the lemma. Fortunately, the
work we have done
 so far will help streamline the rest of the proof. The
next step is to analyze the factor at the left of~(10.1); a routine
application of Stirling's approximation shows that
$${2^m\,m!\,n!\over (n-m)!\,n^{2m}}\;=\;
\sqrt{\mskip1mu2\pi n}\,
2^{n-m}\,e^{-\mu^3\!/6-n}\biggl(1+O\biggl({1+\mu^4\over
 n^{1/3}}\biggr)\biggr)\,,\eqno(10.10)$$
uniformly for $|\mu|\le n^{1/12}$ as $n\to\infty$, when
$m={n\over2}(1+\mu n^{-1/3})$.

Now we turn to the other parts of (10.1). Equation (3.2) implies that 
$T$ has an analytic continuation in which 
$T(ze^{-z})=z$ for $|z|<1$. Hence, by (3.3) and
Cauchy's formula for $[z^n]\,f(z)$, we can substitute $\tau=ze^{-z}$
and get
$$\eqalignno{[z^n]\,{U(z)^{n-m}\over \bigl(1-T(z)\bigr)^y}
&={1\over 2\pi i}
\oint {U(\tau)^{n-m}\,d\tau\over\bigl(1-T(\tau)\bigr)^y
\,\tau^{n+1}}\cr
\noalign{\smallskip}
&={e^n\,2^{m-n}\over 
2\pi i}\oint (1-z)^{1-y}e^{nh(z)}\,{dz\over z}\,,&(10.11)\cr}$$
where
$$\eqalignno{h(z)
&=z-1-{m\over n}\ln z+\left(1-{m\over n}\right)\ln\,(2-z)\cr
\noalign{\smallskip}
&=z-1-\ln z-\left(1-{m\over n}\right)
\ln\,{1\over 1-(z-1)^2}\,.&(10.12)\cr}$$
The contour in (10.11) should keep 
$|z|<1$. Notice that $h(1)=h'(1)=0$; if
$m=\half n$ we 
also have $h''(1)=0$. This triple zero accounts for the
procedure we shall
 use to investigate the value of (10.11) for large~$n$.

Let $\nu=n^{-1/3}$, and let $\alpha$ be the positive solution to
$$\mu\;=\;\alpha^{-1}\,-\,\alpha\,.\eqno(10.13)$$
We will evaluate (10.11) on the path $z=e^{-(\alpha+it)\nu}$,
where $t$ runs from $-\pi n^{1/3}$ to $\pi n^{1/3}$:
$$\oint f(z)\,{dz\over z}=i\nu\int_{-\pi n^{1/3}}^{\pi n^{1/3}}
f(e^{-(\alpha+it)\nu})\,dt\,.\eqno(10.14)$$
It will turn out 
that the main contribution to the value of this integral
comes from the vicinity of $t=0$.

The magnitude of 
$e^{h(z)}$ depends on $\Re h(z)$.\footnote{\dag}{$\Re(x+iy)=x$
denotes the real part of the complex number $x+iy$.}
If $z=\rho e^{i\theta}$, we have
$$\Re h(\rho e^{i\theta})
=\rho \cos\theta-1-{m\over n}\,\ln \rho+\half 
\left(1-{m\over n}\right)\ln\,(4-4\rho\cos\theta +\rho^2)\,.\eqno(10.15)$$
The derivative with respect 
to $\theta$ is $-\rho g(\theta)\sin\theta$, where
$$g(\theta)=1-{2\left(1-{m\over n}\right)\over 
4-4\rho\cos\theta+\rho^2}
\ge {(2-\rho)^2-2\left(1-{m\over 
n}\right)\over4-4\rho\cos\theta+\rho^2}
\,;\eqno(10.16)$$
and $g(\theta)$ is positive when $\rho=e^{-\alpha\nu}$, because
$2(1-{m\over n})=1-\mu\nu<1+\alpha\nu<(2-e^{-\alpha\nu})^2$. (We
always have $0<\alpha\nu<2$ when $\vert\mu\vert\leq n^{1/12}$, and it
is not difficult to verify that $(2-e^{-x})^2>1+x$ when $0<x<2$.)
Hence $\Re h(e^{-(\alpha+it)\nu})$ decreases as $|t|$
increases, and $|e^{nh(z)}|$ has its
 maximum on the circle $z=e^{-(\alpha+it)\nu}$
when $t=0$.

Looking further at $nh(e^{-s\nu})$, we have the asymptotic estimate
$$n\,h(e^{-s\nu})={\textstyle{1\over 3}}\,s^3+
\textstyle\half \mu s^2
+O\bigl((\mu^2s^2+s^4)\nu
\bigr)\,,\eqno(10.17)$$
uniformly in any
 region such that $|s\nu|\le c$ where $c<\ln 2$. This follows
from (10.12), using the expansion
$$\ln\,{1\over 1-(e^u-1)^2}=u^2+u^3+O(u^4)\,,\qquad |u|\le c\,.$$
We also have
$$(1-e^{-s\nu})^{1-y}=
s^{1-y}\nu^{1-y}\bigl(1+O(s\nu)\bigr)\,.\eqno(10.18)$$
Therefore if $f(z)=(1-z)^{1-y}e^{nh(z)}$ is the integrand of (10.11)
and (10.14), we have
$$e^{-\mu^3\!/6}f(e^{-s\nu})=\nu^{1-y}s^{1-y}e^{K(\mu,s)}
\bigl(1+O(s\nu)+O(\mu^2s^2\nu)+O(s^4\nu) \bigr)\,,\eqno(10.19)$$
when $s=O(n^{1/12})$. 
(This restriction on~$s$ ensures that $\mu^2s^2\nu$
and $s^4\nu$ are
bounded, 
hence the $O$~terms of (10.17) can be moved out of the exponent.)

The exponent $K(\mu,s)$ in (10.19), when $s=\alpha+it$, is
$$\textstyle K(\alpha^{-1}-\alpha,\alpha+it)=
\bigl(\half \alpha^{-1}-{1\over6}\alpha^{-3}\bigr)+it-
\half (\alpha+\alpha^{-1})t^2-{1\over3}it^3\,.$$
The real part is bounded above by ${1\over3}-t^2$, for all $\alpha>0$,
since $3\alpha^{-1}-\alpha^{-3}\le2\le\alpha+\alpha^{-1}$, with
equality iff~$\alpha=1$. Hence the integrand $f(e^{-s\nu})$ becomes
superpolynomially small when $\vert t\vert$ grows, and we have
$$\openup3\jot
\eqalign{{e^{-\mu^3\!/6}\over 2\pi i}\,\oint f(z)\,{dz\over z}
&={\nu e^{-\mu^3\!/6}\over 2\pi}\,
\int_{-n^{1/12}}^{n^{1/12}} f\bigl(e^{-(\alpha+it)\nu})\,dt
 +O\bigl(e^{-(\alpha+\alpha^{-1})n^{1/6}\!/3}\bigr)\cr
&={\nu^{2-y}\over 2\pi i}\,\int_{\alpha-n^{1/12}i}^{\alpha+n^{1/12}i}
s^{1-y}e^{K(\mu,s)}\,ds+O(\nu^{3-y}R)
 +O\bigl(e^{-(\alpha+\alpha^{-1})n^{1/6}\!/3}\bigr)\cr
&=\nu^{2-y}A(y,\mu)+O(\nu^{3-y}R) 
+O\bigl(e^{-\max(2,\vert\mu\vert)n^{1/6}\!/3}\bigr)\,,\cr}$$
where $s=\alpha+it$ and
$$R=\int_{-\infty}^\infty\bigl(\,|s^{2-y}|
+\mu^2|s^{3-y}|+|s^{5-y}|\,\bigr)
\bigl|e^{K(\mu,s)}\bigr|\,dt=R_1+R_2+R_3.$$
The lemma will be proved if we can show that $R=O\bigl(1+\mu^B)$
and that $R/A(y,\mu)=O(\mu^4)$ as $|\mu|\to\infty$.

To show that each remainder integral $R_1$, $R_2$, $R_3$ is small, we
will let  $s=\alpha+iu/\beta$, where $u=\beta t$ and
$$\beta=\sqrt{\alpha+\alpha^{-1}}\,.\eqno(10.20)$$
Notice that when
 $\mu\le0$ we have $\alpha\ge1$ and $\alpha=|\mu|+O\bigl(
|\mu|^{-1}\bigr)$;
 when $\mu\ge0$ we have $0<\alpha\le1$ and $\alpha^{-1}
=\mu+O(\mu^{-1})$. Therefore in both cases
$$\beta=|\mu|^{1/2}+
O\bigl(|\mu|^{-1/2}\bigr)\qquad\hbox{as $|\mu|\to\infty$.}
\eqno(10.21)$$

The first remainder, $R_1$, is
$$\int_{-\infty}^\infty |\alpha+it|^{2-y}\,|e^{K(\mu,\alpha+it)}|\,dt
={e^{\alpha^{-1}\!/2-\alpha^{-3}\!/6}\over\beta}
\int_{-\infty}^{\infty}\,\biggl|\mkern1mu\alpha
+{iu\over\beta}\biggr|^{2-y}
e^{-u^2\!/2} \,du\,.$$
If $\mu<0$, we have $\alpha\beta\ge\sqrt2$, hence
$$R_1\le {O(1)\,\alpha^{2-y}\over\beta}
\int_{-\infty}^\infty\max\biggl(1,\,
\biggl|1+{iu\over\sqrt2}\biggr|^{2-y}\biggr)e^{-u^2\!/2}\,du\,;$$
and the integral exists,
 so this is $O\bigl(|\mu|^{3/2-y}\bigr)$ by (10.21).
Similarly, $R_2=O\bigl(|\mu|^{2+5/2-y}\bigr)$ when $\mu<0$, and
$R_3=O\bigl(|\mu|^{9/2-y}\bigr)$.

On the other hand, when $\mu>0$ we have $\alpha\beta\le\sqrt2$, and
we need to be more cautious. Instead of letting $t$ run from $-\infty$
to~$+\infty$ through real values in the derivation above, we will
distort the path slightly near the origin, so that $t$ passes through
the point $-i/\beta$ and so that $\beta s=\alpha\beta +iu$ 
never has magnitude less than~1. $\bigl($We used
 essentially the same sort of contour when deriving
(10.4).$\bigr)$ Then $u$ passes
 through the point~$-i$, and we have
$$R_1\le {O(1)\,e^{-\mu^3\!/6}\over\beta^{\,3-y}}\int_{-\infty}^\infty
\max\bigl(1,\,\bigl|\sqrt2+iu\bigr|^{2-y}\bigr)e^{-u^2\!/2}\,du\,.$$
We therefore have $R_1=O(e^{-\mu^3\!/6}\mu^{y/2-3/2})$; similarly,
        $R_2=O(e^{-\mu^3\!/6}\mu^{y/2-4/2+2})$ and
        $R_3=O(e^{-\mu^3\!/6}\mu^{y/2-6/2})$.
From (10.4) we know that $A(y,\mu)$ grows at least as
fast as $e^{-\mu^3\!/6}\mu^{y/2-5/2}$.
 So in this case the remainders behave even
better than we have claimed in (10.5), although the error term
$O(\mu^4\!/n^{1/3})$ is still necessary because of~(10.10).\quad\pfbox

\medskip
If we differentiate the integral (10.7) with respect to $s$ and with
respect to~$\mu$, we obtain a recurrence relation for $A(y,\mu)$ and
a formula for the derivative:
$$\eqalignno{(y-2)A(y,\mu)&=\mu A(y-2,\mu)+A(y-3,\mu)\,;&(10.22)\cr
A'(y,\mu)&=
\textstyle\half A(y-2,\mu)-\half \mu^2A(y,\mu)\,.&(10.23)\cr
}$$
(The prime here denotes differentiation with respect to the second
argument,~$\mu$. The derivative with respect to~$y$ could also be worked
out; but it depends on the derivative of the Gamma function in a rather
complicated way, and it is not expressible directly in terms of~$A$
itself.)

The derivative is more easily investigated if we define
$$B(y,\mu)=e^{\mu^3\!/6}\,A(y,\mu)\,.\eqno(10.24)$$
Then
$$\eqalignno{(y-2)B(y,\mu)&=\mu B(y-2,\mu)+B(y-3,\mu)\,;&(10.25)\cr
B'(y,\mu)&=\textstyle\half B(y-2,\mu)\,.&(10.26)\cr
}$$
It is easy to verify that the infinite series of (10.2) satisfies these
relations. Repeated application of (10.25) and (10.26) leads to a
third-order differential equation for $B=B(y,\mu)$:
$$8B'''-4\mu^2B''+2\mu(2y-9)B'-(y-2)(y-5)B=0\,.\eqno(10.27)$$

We can see from (10.22) that, for any fixed $\mu\ge0$, there are
infinitely many negative values of~$y$ such that $A(y,\mu)=0$.
For if $y<0$ and there is no root between $y-1$ and~$y$, then
$A(y-1,\mu)$ and $A(y,\mu)$ have the same sign; hence $A(y+2,\mu)$
has the opposite sign, and there's a root between $y$ and~$y+2$.
Therefore we cannot use equation (10.5) until $|\mu|$ is
sufficiently large, at least not when $y<0$ and $\mu\ge0$.

Lemma 3 implies the nonobvious inequality
 $A(y,\mu)\ge0$ for all $y\ge0$, since $A(y,\mu)$ is proportional to
the limiting value of 
 the coefficients of $U(z)^{n-m}/\bigl(1-T(z)\bigr)^y$, and these
coefficients  are nonnegative. Moreover,
$A(y,\mu)$ is strictly positive for $y\geq 2$ and all~$\mu$. For if
$y\geq 2$ and $A(y,\mu_0)=0$, we have $B(y,\mu_0)=0$; but
$B'(y,\mu)\geq 0$ by (10.26), hence we must have $B(y,\mu)=0$ for all
$\mu\leq\mu_0$, which is impossible because $B(y,\mu)$ is a
nonconstant analytic function of~$\mu$ by (10.2).

When $y=1$ there is a ``closed form'' in terms of the Airy function:
$$A(1,\mu)=e^{-\mu^3\!/12}{\rm Ai}(\mu^2\!/4)\,;\eqno(10.28)$$
this is proved in [\FKP, (A.12) and (A.19)]. If we differentiate
(10.28) with respect to~$\mu$, taking note of the fact that (10.22) gives
$$A(-1,\mu)=-\mu A(0,\mu)\,,\eqno(10.29)$$
we find
$$\textstyle A(0,\mu)=-\half \mu e^{-\mu^3\!/12}{\rm Ai}(\mu^2\!/4)
- e^{-\mu^3\!/12}{\rm Ai}'(\mu^2\!/4)\,.\eqno(10.30)$$
Therefore in particular,
$$\textstyle e^{\mu^3\!/12}A(1,\mu)\qquad\qquad\hbox{and}\qquad\qquad
  e^{\mu^3\!/12}\bigl(A(0,\mu)+{\textstyle\half }\,\mu\,A(1,\mu)\bigr)$$
are even functions of $\mu$. The well-known relations between Airy
functions and Bessel functions,
$${\rm
Ai}(z)={1\over\pi}\sqrt{z\over3}\,K_{1/3}\biggl({2\over3}z^{3/2}\biggr),
\qquad
{\rm Ai}' 
 (z)={-z\over \pi\,\sqrt{3}}\,K_{2/3}\biggl({2\over3}z^{3/2}\biggr),$$
yield the additional formulas
$$\eqalignno{
A(1,\mu)&={e^{-\mu^3\!/12}\mu\over2\pi\sqrt3}\,
K_{1/3}\!\left(\mu^3\over12\right)\,,
&(10.31)\cr
A(0,\mu)+{\mu\over2}A(1,\mu)&=A(3,\mu)-{\mu\over2}A(1,\mu)=
{e^{-\mu^3\!/12}\mu^2\over4\pi\sqrt3}\,
K_{2/3}\!\left(\mu^3\over12\right)\,.&(10.32)\cr
}$$

Since we know $A(y,\mu)$ for $y=-1$, 0, and 1, we can use (10.22) to
determine $A(y,\mu)$ for all negative integers~$y$, and for $y=3$
as indicated in (10.32). But a new idea is needed if we hope to have a
closed form when $y=2$. It is possible to express $A(2,\mu)$ as an
infinite sum of Bessel functions,
$$A(2,\mu)={1\over3}\biggl(e^{-\mu^3\!/6}+e^{-\mu^3\!/12}\biggl(\sum_{k\ge0}
(-1)^k\Bigl(I_{k+1/3}\Bigl({\mu^3\over12}\Bigr)-I_{k+2/3}\Bigl({\mu^3\over12}
\Bigr)\Bigr)\biggr)\biggr)\,,\eqno(10.33)$$
but this may be as close to a closed form as possible unless we use
general hypergeometric functions.
Equation (10.33) follows from (10.2) and the hypergeometric identity
$$\eqalignno{&F({\textstyle{1\over
2}+a,\,1+2a-b-c;\,1+2a-b,\,1+2a-c;\,2z})\cr
\noalign{\smallskip}
&\qquad\qquad={e^z\,\Gamma(a)\over (z/2)^a}\,\sum_{k\geq 0}\,
{(-1)^k(2a)^{\overline{k}}b^{\overline{k}}
c^{\overline{k}}(k+a)I_{k+a}(z)\over
(1+2a-b)^{\overline{k}}\,(1+2a-c)^{\overline{k}} k!}
&(10.34)\cr}$$
[\Sla, equation (2.8)]; here $x^{\overline{k}}$ denotes
$\Gamma(x+k)/\Gamma(x)$, and we obtain (10.33) by setting
$z=\mu^3\!/12$, $(a,b,c)=\bigl({1\over 3},{1\over 3},1\bigr)$ and
$\bigl({2\over 3},{2\over 3},1\bigr)$.

The facts that $K_{1/3}(z)=3^{-1/2}\pi\bigl(I_{-1/3}(z)-I_{1/3}(z)\bigr)$,
               $K_{2/3}(z)=3^{-1/2}\pi\bigl(I_{-2/3}(z)-I_{2/3}(z)\bigr)$,
and $e^{-z}=I_0(z)+2\sum_{k\ge1}(-1)^kI_k(z)$ suggest that we look for
an identity of the form
$$\eqalignno{\kern-2em A(y,\mu)
&=\left({\mu\over 2}\right)^{2-y}e^{-\mu^3\!/12}\sum_{k\geq 0}a_k(y)
I_{(k+y-2)/3}\left({\mu^3\over 12}\right)\cr
\noalign{\smallskip}
&={e^{-\mu^3\!/6}\over 3^{(y-2)/3}}\,\sum_{k\geq
0}{a_k(y)\over\Gamma\bigl({k+y+1\over3}\bigr)}\left({\mu\over
2\cdot 3^{1/3}}\right)^k\!F\!\left({2k{+}2y{-}1\over 6}@;
 {2k{+}2y{-}1\over 3}@;{\mu^3\over 6}\right).&(10.35)\cr}$$
Any formal power series in $\mu$ has such an expansion, for all
$y>-1$. But the coefficients $a_k(y)$ do not appear to have a simple
form except in the cases already mentioned. We have
$$\displaylines{a_0(y)={1\over 3}\,,\quad a_1(y)={y-1\over 3}\,,\quad
a_2(y)={y(y-3)\over 6}\,,\quad a_3(y)={(y^2-1)(y-6)\over 18}\cr
\noalign{\smallskip}
a_4(y)={(y-1)(y+2)(y^2-11y+12)\over 72}\,,\quad a_5(y)=
{(y+3)y(y-1)(y-3)(y-14)\over 360}\,.\cr}$$

Splitting (10.2) into three sums according to the value of $k\bmod3$ yields
a closed form for $A(y,\mu)$ in terms of general hypergeometric series:
$$\eqalignno{A(y,\mu)&=e^{-\mu^3\!/6}\left({1\over
3^{(y+1)/3}\,\Gamma\bigl((y{+}1)/3\bigr)}\; F\left({2{-}y\over
6}\,,\,{5{-}y\over
6}\,;\, {1\over 3}\,,\,{2\over 3}\,;\,{\mu^3\over 6}\right)\right.\cr
\noalign{\smallskip}
&\qquad\qquad\null+{1\over 3^{(y-1)/3}\,\Gamma\bigl((y{-}1)/3\bigr)}
{\mu\over 2}\,
F\left({4{-}y\over
6}\,,\,{7{-}y\over 6}\,;\,{2\over 3}\,,\,{4\over 3}\,;\,{\mu^3\over
6}\right)\cr
\noalign{\smallskip}
&\qquad\qquad\null+\left.{1\over3^{(y-3)/3}\Gamma\bigl((y{-}3)/3\bigr)}
{\mu^2\over8}\,F\left({6{-}y\over 6}\,,\,{9{-}y\over
6}\,;\,{4\over 3}\,,\,{5\over 3}\,;\,{\mu^3\over
6}\right)\right).&(10.36)\cr}$$ 

\bigbreak\noindent
{\bf 11. Application to bicyclic components.}\enspace
Now we are ready to begin using the basic theoretical results of the
preceding sections. We will start by considering the case when the
parameter~$\mu$ of Lemma~3 is very small, say $\mu=O(n^{-1/3})$. Then
there are $m=\half n+O(n^{1/3})$ edges.

\proclaim
Theorem 4. The probability that a random graph or multigraph with $n$~vertices
and $\half n+O(n^{1/3})$ edges has exactly $r$~bicyclic components,
and no components of higher cyclic order, is
$$\left({ 5\over 18}\right)^r\sqrt{{2\over 3}}\;{1\over
(2r)!}+O(n^{-1/3})\,.\eqno(11.1)$$

\proof
(The special case $r=0$ and $m=\half n$ of this theorem was Corollary 9
of~[\FKP].)
Consider first the case of random multigraphs,
since this case is simpler. If there are $n$~vertices,
$m$~edges, $r$~bicyclic components, and no components with higher cyclic
order, there must be exactly $n-m+r$ acyclic components. The probability
of such a configuration, according to (2.2), is therefore
$${2^m\,m!\,n!\over n^{2m}}\;[z^n]\;{U(z)^{n-m+r}\over (n-m+r)!}\;
e^{V(z)}\;{W(z)^r\over r!}\,,\eqno(11.2)$$
where $U(z)$, $V(z)$, $W(z)$ are the generating functions (3.3), (3.4), and
(3.7). Now
$$W(z)={5\over 24}\;{1\over\bigl(1-T(z)\bigr)^3}-{7\over 24}\;
{1\over\bigl(1-T(z)\bigr)^2}+
{1\over 12}\;{1\over\bigl(1-T(z)\bigr)}\,,\eqno(11.3)$$
using the coefficients $e'_{1d}$ of~(7.20);
so we see that $W(z)^r$ is a polynomial of degree $3r$ in $\bigl(1-T(z)\bigr)^{-1}$,
with leading coefficient $\bigl({5\over 24}\bigr)^r$. Lemma~3 tells us that the
leading term of $W(z)^r$ is the only significant one, asymptotically
speaking, because the other terms
contribute at most $n^{-1/3}$ times as much as the leading term.
We can also write
$$U(z)^r=2^{-r}\bigl(1-\bigl(1-T(z)\bigr)^2\bigr)^r\,;\eqno(11.4)$$
this allows us to replace $U(z)^r$ by $2^{-r}$ in (11.2).
Since $e^{V(z)}=\bigl(1-T(z)\bigr)^{-1/2}$, the value of (11.2) is
$${(n-m)!\over (n-m+r)!\,r!\,2^r}\;{\sqrt{\mskip1mu2\pi}\,
n^r\over 3^{r+1/2}\Gamma(r+1/2)}\;\left({5\over 24}\right)^r
\bigl(1+O(n^{-1/3})\bigr)\,.$$
This simplifies to (11.1) using the fact that
$${(n-m)!\over (n-m+r)!}={2^r\over
n^r}\,\bigl(1+O(rn^{-2/3}+r^2n^{-1})\bigr)\,,$$ 
and using a special case of the duplication formula
for the Gamma function,
$$\Gamma(r+1/2)={(2r)!\,\sqrt{\pi}\over 4^rr!}\,.
\eqno(11.5)$$

On the other hand if we are dealing with random graphs we must replace (11.2) by
$${n!\over {n(n-1)/2\choose m}}\;[z^n]\;{U(z)^{n-m+r}\over
(n-m+r)!}\;e^{\widehat{V}(z)}\;{\widehat{W}(z)^r\over r!}\,,\eqno(11.6)$$
where $\widehat{V}(z)$ and $\widehat{W}(z)$ appear in (3.5) and (3.6).
Again we have $\widehat{W}(z)={5\over 24}\bigl(1-T(z)\bigr)^{-3}$ plus less
significant terms, so $\widehat{W}(z)$ produces an effect similar to $W(z)$.
But $\widehat{V}(z)=V(z)-\half T(z)-{1\over 4}T(z)^2$; so we now want
the coefficient of $[z^n]$ in an expression proportional to
$${U(z)^{n-m}\over \bigl(1-T(z)\bigr)^{3r+1/2}}\;e^{-T(z)/2-T(z)^2\!/4}\,,$$
which has an exponential factor not covered by Lemma~3. The proof of
Lemma~3 shows, however, that this exponential factor simply changes the
result by a factor of $e^{-3/4}+O(n^{-1/3})$:
We multiply (10.18) by $\exp(-e^{-s\nu}/2-e^{-2s\nu}/4)=e^{-3/4}+O(s\nu)$.

Furthermore, (11.6) contains a factor $e^{+3/4}$ to cancel the $e^{-3/4}$, because
of (2.4). Therefore the leading term of the asymptotic probability for graphs
is the same as it was for multigraphs.\quad\pfbox

\proclaim
Corollary. The probability that a random graph or multigraph with $n$~vertices
and $\half n$ edges has only acyclic, unicyclic, and bicyclic components is
$$\sqrt{2\over 3}\cosh\sqrt{5\over 18}+O(n^{-1/3})\approx 0.9325\,.\eqno(11.7)$$

\proof
The sum over $r$ of the estimate made in Theorem 4 clearly gives a lower
bound, so we must prove that it is also an upper bound. That sum can
be written 
$${2^m\,m!\,n!\over (n-m)!\,n^{2m}}\;[z^n]\,U(z)^{n-m}f_{n-m}(z)\,,$$
where
$$f_l(z)=\sum_{r\ge 0}{l\mskip1mu!\over (l+r)!}\;{\bigl(U(z)W(z)\bigr)^r\over r!}\;e^{V(z)}
\,.\eqno(11.8)$$
If we look at the proof of Theorem 4, and the proof of Lemma 3 on which it is
based, we see that the calculations all depend on $f_l(ze^{-z})$, where
$|z|\le e^{-\nu}$ and $\nu=n^{-1/3}$. In this region,
$$|T(ze^{-z})|\le e^{-\nu}\,,\qquad
 |1-T(ze^{-z})|\ge \nu+O(\nu^2)\,.\eqno(11.9)$$
Thus the sum $f_{n-m}(z e^{-z})$ converges {\it uniformly\/} for all~$n$
and all $|z|\le e^{-\nu}$. Uniform convergence allows us to interchange
summation and integration.
(Notice that the function $h(z)$ in the proof of Lemma 3, which influences
the behavior of the integrand most strongly as $n\to\infty$, is independent of~$r$.)
\quad\pfbox

\medskip
Another proof of (11.7) will be given below.

\bigbreak\noindent
{\bf 12. Components of higher cyclic order.}\enspace
Now let's consider components that are tricyclic, tetracyclic, etc. (Notice
that tricyclic components correspond to $C_2(z)$, not
$C_3(z)$, in the notation of section~2;  our notation has mathematical
advantages, but it is slightly out of phase
with the traditional terminology.)

\proclaim
Theorem 5. The probability that a random graph or multigraph with $n$ vertices
and $\half n+O(n^{1/3})$ edges has exactly $r_1$~bicyclic components,
$r_2$~tricyclic components, $\ldots\,$, $r_q$~$(q+1)$-cyclic components, and
no components of higher cyclic order, is
$$\biggl({4\over3}\biggr)^{\!r}
\sqrt{2\over 3}\;{c_1^{r_1}\over r_1!}\;
{c_2^{r_2}\over r_2!}\;\cdots\;
{c_q^{r_q}\over r_q!}\;
{r!\over 
(2r)!}+O(n^{-1/3})\,,\eqno(12.1)$$
where $r=r_1+2r_2+\cdots +qr_q$ and the constants $c_j$ are defined in (8.6).

\proof
If there are $n$ vertices and $m$ edges, there
must be exactly $n-m+r$ acyclic 
components. So we can argue as in Theorem~4 to find
$$\displaylines{\qquad{2^m\,m!\,n!\over n^{2m}}\;[z^n]\;
{U(z)^{n-m+r}\over (n-m+r)!}\;e^{V(z)}\;{C_1(z)^{r_1}\over r_1!}\;
{C_2(z)^{r_2}\over r_2!}\;\cdots\;{C_q(z)^{r_q}\over r_q!}\hfil\cr
\noalign{\smallskip}
\hfil ={c_1^{r_1}\over r_1!}\;
{c_2^{r_2}\over r_2!}\;\cdots\;
{c_q^{r_q}\over r_q!}\;{ \sqrt{\mskip1mu2\pi}\over 3^{r+1/2}\Gamma(r+1/2)}
+O(n^{-1/3})\,.\qquad\cr}$$
Formula (12.1) now follows from (11.5) as before.\quad\pfbox

\medskip
Let's illustrate the consequences of Theorem~5 by computing the
limiting probabilities for small values of the parameters
$(r_1,r_2,\ldots,r_q)$. Here is a list of all configurations with
$r_1+r_2+\cdots +r_q>1$ that occur with limiting probability $.000005$
or more, showing the probabilities rounded to five decimal places:
$$\vcenter{\halign{$\hfil#\;$&$#$\hfil\qquad
&$\hfil#\;$&$#$\hfil\qquad
&$\hfil#\;$&$#$\hfil\cr
[2]&=.00263&[0,2]&=.00008&[1,0,0,0,0,1]&=.00002\cr
[1,1]&=.00105&[1,0,0,0,1]&=.00004&[2,1]&=.00001\cr
[1,0,1]&=.00031&[0,1,1]&=.00003&[0,1,0,1]&=.00001\cr
[1,0,0,1]&=.00010&[3]&=.00002&[1,0,0,0,0,0,1]&=.00001\cr}}$$
(The notation $[2]$ stands for the case $r_1=2$, $r_2=r_3=\cdots =0$;
similarly $[r_1,\ldots,r_q]$ implies that there are no complex
components of cyclic order greater than $q+1$.)

The sum of these probabilities, .00431, is nicely balanced by $\sqrt{2/3}$
plus the sum of probabilities when there is only one complex component,
i.e., when $r_q=1$ and all other $r$'s are zero:
$$\eqalign{.81650&+.11340+.03780+.01547+.00678+.00307+.00141\cr
&\quad\null+.00066+.00031+.00015+.00007+.00003+.00002+.00001\,;\cr}$$
this comes to $.99568=.99999-.00431$.

Suppose ${\cal R}$ is any countably infinite set of configurations
$[r_1,r_2,\ldots,r_q]$, where $q$ might be unbounded. We would like to
prove that a random graph or multigraph with approximately $\half n$ edges lies
in~${\cal R}$ with limiting probability
$$\sum\{\,P[r_1,r_2,\ldots,r_q]\bigm| [r_1,r_2,\ldots,r_q]\in {\cal
R}\,\}\,, \eqno(12.2)$$
where $P[r_1,r_2,\ldots,r_q]$ is the limiting value stated in
Theorem~5. The technique we used to prove (11.7)
does not apply, because the infinite sums over which integration takes
place might not converge uniformly when $q$ is unbounded.

However, we are obviously justified in claiming that (12.2) is a {\it
lower\/} bound for the stated probability, because the sum over any
finite subset of~${\cal R}$ yields a lower bound. 

We will prove below that the sum of $P[r_1,r_2,\ldots,r_q]$ over all
possible configurations $[r_1,r_2,\ldots,r_q]$ is~1. Consequently, the
sum (12.2) must in fact be the limiting probability of a random graph
or multigraph
being in~${\cal R}$, not just a lower bound. If (12.2) were too low, we
would not obtain~1 by adding the complementary probabilities
$P[r_1,r_2,\ldots,r_q]$ for $[r_1,r_2,\ldots,r_q]\not\in{\cal R}$.
This observation will lead to the promised ``second proof'' of~(11.7),
if we also sum less significant terms to obtain the error bound $O(n^{-1/3})$.

\bigbreak\noindent
{\bf 13. Excess Edges.}\enspace
The notion of ``excess'' was used somewhat informally in the introductory
sections of this paper. Let us now define it formally, saying that
the {\it excess\/} of a graph or multigraph is the number of edges
plus the number of acyclic components, minus the number of vertices.
Thus a $(q+1)$-cyclic component has excess~$q$, when $q\geq 0$. If a
graph or multigraph
has $r_1$~bicyclic components, $r_2$~tricyclic components, etc.,
then it has excess $r=r_1+2r_2+3r_3+\cdots\;$.

If $G$ and $G'$ are graphs on the same vertices, and if $G\cup G'$ and
$G\cap G'$ denote the graphs obtained by taking the union and
intersection of their edges, the excesses satisfy
$$r(G)+r(G')\leq r(G\cup G')+r(G\cap G')\,.$$
For we can start with empty graphs and insert the edges of $G\cap G'$,
preserving equality. Then if we insert an edge of $G\setminus G'$ or
of $G'\setminus G$, each side of the inequality increases by either~0
or~1; and the left side cannot increase by~1 unless the right side
does also. For example, if the left side increases by~1 when we add an
edge of $G\setminus G'$, the endpoints of that edge are in non-trees
of~$G$, so they surely are in non-trees of $G\cup G'$.

We have seen in Theorem 5 that the limiting joint probability distribution
of the random variables $(r_1,r_2,\ldots\,)$ in a large random graph or
multigraph with
approximately $\half n$ edges has the form
$${c_1^{r_1}\over r_1!}\;{c_2^{r_2}\over r_2!}\;\cdots\;{c_q^{r_q}\over r_q!}
\;f(r)\,,\eqno(13.1)$$
where $r=r_1+2r_2+\cdots+qr_q$ is the excess of the graph and $r_l=0$ for $l>q$.
Indeed, this is not surprising, if we look at the problem in another
way.

Let $\cal S$ be the set of all multigraphs of configuration $[r_1,r_2,
\ldots,r_q]$, and let $S(w,z)$ be its bgf. The probability that a given
multigraph with $m$~edges and $n$~vertices lies in~$\cal S$ is then
$$\Pr_{mn}({\cal S})={[w^mz^n]\,S(w,z)\over[w^mz^n]\,G(w,z)}\,.\eqno(13.2)$$
We can also express this as
$$\Pr_{mn}({\cal S})=\Pr_{mn}({\cal S} \mid r)\Pr_{mn}({\cal E}_r)\,,
\eqno(13.3)$$
where $\Pr_{mn}({\cal S}\mid r)$ means the probability of obtaining an
element of~$\cal S$ given that the excess is~$r$, and $\Pr_{mn}({\cal E}_r)$
is the probability that a random multigraph has excess~$r$:
$$\Pr_{mn}({\cal S}\mid r)={[w^mz^n]\,S(w,z)\over[w^mz^n]\,e^{U(w,z)+V(w,z)}
E_r(w,z)}\,,\quad
  \Pr_{mn}({\cal E}_r)={[w^mz^n]\,e^{U(w,z)+V(w,z)}
E_r(w,z)\over[w^mz^n]\,G(w,z)}\,.\kern-14pt\eqno(13.4)$$
Since all elements of $\cal S$ have excess $r$, we can compute $[w^mz^n]\,
S(w,z)$ with univariate generating functions:
$$\openup2\jot
\eqalign{S(w,z)
&=e^{U(w,z)+V(w,z)}{C_1(w,z)^{r_1}\over r_1!}
{C_2(w,z)^{r_2}\over r_2!}\ldots{C_q(w,z)^{r_q}\over r_q!}\cr
&=e^{U(wz)/w+V(wz)}{\bigl(wC_1(wz)\bigr)^{r_1}\over r_1!}
{\bigl(w^2C_2(wz)\bigr)^{r_2}\over r_2!}\ldots
{\bigl(w^qC_q(wz)\bigr)^{r_q}\over r_q!}\cr
&=e^{U(wz)/w+V(wz)}w^r{C_1(wz)^{r_1}\over r_1!}
{C_2(wz)^{r_2}\over r_2!}\ldots{C_q(wz)^{r_q}\over r_q!}\,;\cr}$$
hence
$$[w^mz^n]\,S(w,z)=[z^n]\,{U(z)^{n+r-m}\over(n+r-m)!}e^{V(z)}S(z)\,,
\eqno(13.5)$$
if we let
$$S(z)={C_1(z)^{r_1}\over r_1!}
{C_2(z)^{r_2}\over r_2!}\ldots{C_q(z)^{r_q}\over r_q!}\,.$$
Similarly
$$[w^mz^n]\,e^{U(w,z)+V(w,z)}E_r(w,z)
=[z^n]\,{U(z)^{n+r-m}\over(n+r-m)!}e^{V(z)}E_r(z)\,.$$

A multigraph with $m$ edges, $n$ vertices, and excess $r>0$ has $t=n+r-m$
components that are trees (including isolated vertices). Suppose it has
$n_1$~vertices in complex components and $n_0$ vertices in trees and
unicyclic components. Then
$$\eqalign{\Pr_{mn}({\cal S}\mid r)
={\displaystyle[z^n]\,{U(z)^t\over t!_{\mathstrut}}e^{V(t)}S(z)\over
  \displaystyle[z^n]\,{U(z)^{t^{\mathstrut}}\over t!}e^{V(t)}E_r(z)}
&={\displaystyle\sum_{n_0+n_1=n}\bigl([z^{n_0}]\,U(z)^t e^{V(z)}\bigr)
\bigl([z^{n_1}]\,S(z)\bigr)\over
\displaystyle\sum_{n_0+n_1=n}\bigl([z^{n_0}]\,U(z)^t e^{V(z)}\bigr)
\bigl([z^{n_1}]\,E_r(z)\bigr)}\cr
\noalign{\medskip}
&=\sum_{n_1}\Pr({\cal S}\mid r,n_1)\Pr_{mn}(n_1\mid {\cal E}_r)\,,\cr}$$
where
$$\openup2\jot
\eqalignno{
\Pr({\cal S}\mid r,n_1)&={[z^{n_1}]\,S(z)\over[z^{n_1}]\,E_r(z)}\,;&(13.6)\cr
\Pr_{mn}(n_1\mid{\cal E}_r)&={\bigl([z^{n-n_1}]\,U(z)^t e^{V(z)}\bigr)
\bigl([z^{n_1}]\,E_r(z)\bigr)\over[z^n]\,U(z)^te^{V(z)}E_r(z)}\,.
&(13.7)\cr}$$

Thus, $\Pr(\cal S)$ has been expressed in terms of a simple ratio (13.6),
the number of multigraphs consisting of precisely $r_j$ components
of excess~$j$ for $1\le j\le q$, divided by the number of complex
multigraphs of excess~$r$. We know from section~9 that
there are coefficients~$s_d$ such that
$$S(z)={s_0 T(z)^{2r}\over\bigl(1-T(z)\bigr)^{3r}}+
       {s_1 T(z)^{2r-1}\over\bigl(1-T(z)\bigr)^{3r-1}}+\cdots+
       {s_{2r-1}T(z)\over\bigl(1-T(z)\bigr)^{r+1}}\,.$$
Indeed, section 9 
tells us that $s_d$ is $\sum \kappa(\MM)/(2r-d)!$,
summed over all reduced multigraphs of configuration $[r_1,r_2,\ldots,
r_q]$ having exactly $2r-d$ vertices. We can also write
$$S(z)={s'_0\over\bigl(1-T(z)\bigr)^{3r}}+
       {s'_1\over\bigl(1-T(z)\bigr)^{3r-1}}+\cdots+
       {s'_{2r}\over\bigl(1-T(z)\bigr)^r}\,,\eqno(13.8)$$
letting $s'_d=\sum_k{2r-k\choose d-k}(-1)^{d-k}s_k$ as in (7.22). Therefore,
$$n!\,[z^n]\,S(z)=s'_0t_n(3r)+s'_1t_n(3r-1)+\cdots+s'_{2r}t_n(r)\,,$$
expressing the relevant number of multigraphs in terms of the tree
polynomials (3.8); and (3.9) tells us that
$$n!\,[z^n]\,S(z)=
s'_0{\sqrt{\mskip1mu2\pi}\,n^{n-1/2+3r/2}\over2^{3r/2}\Gamma(3r/2)}
\bigl(1+O(n^{-1/2})\bigr)\,.$$
Similarly, we have
$$n!\,[z^n]\,E_r(z)=
e'_{r0}{\sqrt{\mskip1mu2\pi}\,n^{n-1/2+3r/2}\over2^{3r/2}\Gamma(3r/2)}
\bigl(1+O(n^{-1/2})\bigr)\,.$$
Therefore the ratio (13.6) is
$$\Pr({\cal S}\mid r,n_1)={s'_0\over e'_{r0}}\bigl(1+O(n_1^{-1/2})\bigr)\,;$$
and we can sum over $n_1$ to get
$$\Pr_{mn}({\cal S})=\biggl({s'_0\over e'_{r0}}+O(\epsilon)\biggr)\Pr_{mn}
({\cal E}_r)\,,\eqno(13.9)$$
where $\epsilon$ is the expected value of $n_1^{-1/2}$ in the probability
distribution~(13.7).

Moreover, the leading coefficient is
$$s'_0=s_0={c_1^{r_1}\over r_1!}{c_2^{r_2}\over r_2!}\ldots
           {c_q^{r_q}\over r_q!}\,;\eqno(13.10)$$
and $e'_{r0}$ is just $e_r$, the sum 
of (13.10) over all configurations $[r_1,r_2,
\ldots,r_q]$ with $r_1+2r_2+\cdots+qr_q=r$.
This derivation explains why we obtained a formula of the form (13.1)
in Theorem~5.

With graphs instead of multigraphs, the same considerations apply, but
we must add more terms to the formulas. For example, (13.8) becomes
$$\widehat S(z)={\hat s'_0\over\bigl(1-T(z)\bigr)^{3r}}+
       {\hat s'_1\over\bigl(1-T(z)\bigr)^{3r-1}}+\cdots+
       \hat s'_{3r}+
       \hat s'_{3r+1}\bigl(1-T(z)\bigr)+
       \hat s'_{3r+2}\bigl(1-T(z)\bigr)^2\,.\eqno(13.11)$$
The leading coefficient $\hat s'_0$ is the same as $s_0$, so the
asymptotic behavior is the same as before, if we assume that $m$ is
large enough to make the expected value of~$n_1^{-1/2}$ approach zero.

We can estimate the expected value of $n_1^{-1/2}$ by finding the
expected value of
$${[z^{n_1}]\,S(z)-(s_0/e_r)E_r(z)\over[z^{n_1}]\,E_r(z)}\,;\eqno(13.12)$$
indeed, this expected value is the true error in the approximation (13.9),
so it is even more relevant than the expected value of~$n_1^{-1/2}$.
Since $S(z)-(s_0/e_r)E_r(z)$ can be expressed as
 $(s'_1-e'_{r1}s_0/e_r)/\bigl(1-T(z)\bigr)
^{3r-1}$ plus less significant terms, the desired expected value times
$\Pr_{mn}({\cal E}_r)$ is obtained by applying Lemma~3 as we did in
the proof of Theorem~4, but with $3r$ replaced by $3r-1$. The result,
when $m$ is near~$\half n$, is proportional to~$n^{-1/3}$.

The expected value of $n_1^k$ can be computed if we replace $S(z)$ by
$\vartheta^kE_r(z)$ in these formulas, because $[z^n]\,\vartheta^kE_r(z)
=n^k[z^n]\,E_r(z)$. This has the effect of changing the leading term
from $e_r\big/\bigl(1-T(z)\bigr)^{3r}$ to $(3r)(3r+2)\ldots(3r+2k-2)e_r
\big/\bigl(1-T(z)\bigr)^{3r+2k}$, so the result when $m$ is near~$\half n$
is proportional to~$n^{2k/3}$. We have proved

\proclaim
Corollary. If $m=\half n(1+\mu n^{-1/3})$ and $|\mu|\le n^{1/12}$,
the $k$th  moment ${\rm E}_{mn}(n_1^k\mid r)$ of the number of
vertices in complex components, given that the total excess is~$r$, is
$$\openup2\jot
\tabskip\centering
\halign to\displaywidth{
$\displaystyle{#}$\hfil\tabskip0pt
&\quad#\hfil\quad\tabskip\centering
&\llap{$#$}\tabskip0pt\cr
\alpha_{kr}{\Gamma(r+\half )\over3^{2k/3}\Gamma(r+\half +{2\over3}k)}
n^{2k/3}\bigl(1+O(\mu)+O(n^{-1/3})\bigr)\,,
&if $\mu=O(1)$;&(13.13)\cr
\alpha_{kr}{n^{2k/3}\over\mu^{2k}}\bigl(1+O(|\mu|^{-3})+O(\mu^4n^{-1/3})
\bigr)\,,
&if $\mu\to-\infty$;&(13.14)\cr
\alpha_{kr}{n^{2k/3}\mu^k\over2^k}{\Gamma({3\over2}r+{1\over4})\over
\Gamma({3\over2}r+{1\over4}+k)}\bigl(1+O(\mu^{-1})+O(\mu^4n^{-1/3})\bigr)\,,
&if $\mu\to+\infty$;&(13.15)\cr}
$$
here $\alpha_{kr}=(3r)(3r+2)\ldots(3r+2k-2)$.

\proof
These expressions are $\alpha_{kr}$ times the ratios of formulas
(10.2), (10.3), and (10.4) when $y=3r+\half +2k$ to their values
when $y=3r+\half $.\quad\pfbox

\medskip
Notice that when $m$ is approximately
$\half n-n^{3/4}$, the probable value of~$n_1$
is proportional to $n^{2/3-2/12}=n^{1/2}$; when $m\approx \half n
+n^{3/4}$, it is proportional to $n^{2/3+1/12}=n^{3/4}$. These are the
extreme cases $|\mu|=n^{1/12}$ at the limits of Lemma~3's range.

We can use formula (13.3) whenever $\cal S$ is a collection of multigraphs
whose complex components have total excess~$r$. We can use formula
(13.6) whenever $\cal S$ also places no restriction on its
non-complex (acyclic and unicyclic) components. For example, we
can determine the conditional probability that a random graph with
$\half n$ edges has a bicyclic component of each of the three
types in~(9.15), given that it has excess~1. The generating functions
$S(z)$ for the three cases are respectively
${1\over8}T/(1-T)^2$, ${1\over8}T^2\!/(1-T)^3$, ${1\over12}T^2\!/(1-T)^3$;
so the respective conditional probabilities are
$$O(n^{-1/3}),\qquad {3\over5}+O(n^{-1/3}),\qquad
{2\over5}+O(n^{-1/3})\,.\eqno(13.16)$$

All probabilities that are conditional on excess~$r$ must, of course,
be multiplied by $\Pr_{mn}({\cal E}_r)$, the probability that a
random multigraph has excess~$r$.  Lemma~3 and the method of Theorem~5
make this easy to compute:

\proclaim
Corollary. A graph or multigraph with $m=\half n(1+\mu n^{-1/3})$
edges and $n$ vertices has excess~$r$ with probability
$$\Pr_{mn}({\cal E}_r)=\sqrt{\mskip1mu2\pi}\,e_r\,A(3r+{\textstyle
\half },\mu)+O\biggl({1+\mu^4\over n^{1/3}}\biggr)\,,\eqno(13.17)$$
uniformly for $|\mu|\leq n^{1/12}$ as $n\to\infty$, where $e_r=e_{r0}$
is given by~(7.2) and $A(y,\mu)$ is given by~(10.2). When $\mu\to
-\infty$, the probability is $O(|\mu|^{-3r})$; when $\mu\to+\infty$
it is $O(\mu^{3r/2}e^{-\mu^3\!/6})$.\quad\pfbox

(The special case $r=0$ in (13.17), without the error bound, was
found by Britikov~[\Bri], who proved that a random graph has excess~0
with probability approaching $\sqrt{2\pi}@A(\half,\mu)$,
for fixed $\mu$ as $n\to\infty$.)

Here is a table that shows how the probabilities of having excess~$r$
change as the graph or multigraph
evolves past the critical point $m=\half n$:
$$\vcenter{\halign{$#$\hfil\quad%
&\hfil$#$\hfil\enspace%
&\hfil$#$\hfil\enspace%
&\hfil$#$\hfil\enspace%
&\hfil$#$\hfil\enspace%
&\hfil$#$\hfil\enspace%
&\hfil$#$\hfil\enspace%
&\hfil$#$\hfil\enspace%
&\hfil$#$\hfil\enspace%
&\hfil$#$\hfil\enspace%
&\hfil$#$\hfil\enspace%
&\hfil$#$\hfil\cr
&r=0&r=1&r=2&r=3&r=4&r=5&r=6&r=7&r=8&r=9&r=10\cr
\noalign{\smallskip}
\mu=-3&.994&.006&.000&.000&.000&.000&.000&.000&.000&.000&.000\cr
\mu=-2&.983&.015&.001&.000&.000&.000&.000&.000&.000&.000&.000\cr
\mu=-1&.947&.043&.008&.002&.000&.000&.000&.000&.000&.000&.000\cr
\mu=0&.816&.113&.040&.017&.007&.003&.001&.001&.000&.000&.000\cr
\mu=1&.475&.179&.115&.077&.052&.035&.023&.015&.010&.007&.004\cr
\mu=2&.100&.082&.085&.086&.084&.079&.073&.066&.058&.051&.043\cr
\mu=3&.003&.004&.007&.010&.013&.017&.020&.024&.028&.031&.034\cr}}$$
The mean excess is approximately .308, 1.544, 6.364, 19.009, for
$\mu=0,1,2,3$. 

In this paper we are interested mainly in graphs or multigraphs with
approximately $\half n$ edges, but it is instructive to consider
also the formulas that arise when $m$ is somewhat smaller. The excess
is then almost surely zero. In fact, we can obtain a formula that has
a much better error bound than (13.17), in the case $r=0$ and
$\mu<-n^{-\epsilon}$: If we set $\lambda=2m/n$, and if $m<{1\over
2}n-n^{2/3+\epsilon}$, the probability of excess~0 can be shown to be
exactly
$$\eqalignno{&\kern-4em{2^mm!\,n!\over (n-m)!\,n^{2m}}\;[z^n]\;{U(z)^{n-m}\over
\bigl(1-T(z)\bigr)^{1/2}}\cr
\noalign{\smallskip}
&\kern-4em\quad ={S(m)S(n)\over\sqrt{2\pi}\,S(n-m)}\oint\left(
1-{it\beta\over1-\lambda}\right)^{\!\!1/2}\!\!
\left(1+{it\beta\over\lambda}\right)^{\!\!-1}
\!\!\!e^{h(n,\lambda,t)-t^2\!/2}\,dt\,,\quad&(13.18)\cr}$$
where
$$\eqalignno{S(n)&={n!\,e^n\over n^n\,\sqrt{2\pi n}}=1+O\left({1\over
n}\right)\,,&(13.19)\cr
\noalign{\smallskip}
\beta&=\sqrt{\lambda(2-\lambda)\over(1-\lambda)n}\,\,,&(13.20)\cr
\noalign{\smallskip}
h(n,\lambda,t)&=n@h(\lambda+it\beta)-n@h(\lambda)+t^2\!/2&(13.21)\cr
&={n\over 2}\,\sum_{k\geq 3}\,{(it\beta)^k\over k}\,\left(
{(-1)^k\over\lambda^{k-1}}-{1\over (2-\lambda)^{k-1}}\right)
\,,&(13.22)\cr}$$
and the contour of integration makes $z=\lambda+it\beta$ traverse the circle
$\vert z\vert=\lambda$ as $t$ varies. The function $h(z)$ in (13.21)
is the function defined in (10.12). We are
essentially simplifying the proof of Lemma~3 by choosing a path of
integration through the saddle point $z=\lambda$, as in the proof of
Theorem~4 in [\FKP]. The proof of that theorem justifies restricting $t$
to a neighborhood of zero, so that the tail-exchange method can be
applied as in the derivation following (10.9).
It follows that the probability of excess~0 is
$1-O\bigl(n^2\!/(n-2m)^3\bigr)$ for $m$ in the stated range. We have
in fact the estimate
$${\rm Pr}_{mn}({\cal E}_0)=1-{5\over
24}\,\alpha^{-3}\bigl(1+O(\alpha^{-3})+O(\alpha
n^{-1/3})\bigr)\eqno(13.23)$$
when $m=\half n(1-\alpha\,n^{-1/3})$, uniformly for $(\ln n)^2\leq
\alpha\leq \half n^{1/3}$.

It is interesting to note that the tail-exchange method can be used to extend
(13.23) to an asymptotic series in $\alpha^{-1}$ and $\alpha n^{-1/3}$,
although the integral (13.18) actually diverges if we let $t$ run through
all real values from $-\infty$ to $+\infty$ instead of describing the
stated contour. Indeed, the magnitude of the integrand in (13.18) for
large real values of~$|t|$ is approximately $|t|^{n-2m-1/2}$.

\bigbreak\noindent
{\bf 14. Probability distribution of the excess.}\enspace
One way to check our calculations is to verify that the
probabilities in (13.17) sum to~1.
Thus we want to prove that
$$\sum_{k\geq 0}\,\sqrt{{2\pi\over 3}}\,{(\half \,3^{2/3}\mu)^k\over
k!}\,\sum_{r\geq 0}\,{(6r)!\over
2^{5r}3^{3r}(2r)!\,(3r)!\,\Gamma(r+1/2-2k/3)}=e^{\mu^3\!/6}\,.\eqno(14.1)$$
The inner sum is a hypergeometric series whose sum is known;
$${1\over\Gamma(1/2-2k/3)}\,F\left({1\over 6},{5\over
6}\,;\,\half -{2k\over 3}\,;\,{1\over
2}\right)={2^{1/2+2k/3}\,\sqrt{\pi}\over
\Gamma\bigl((1-k)/3\bigr)\,\Gamma\bigl((2-k)/3\bigr)}\,.\eqno(14.2)$$ 
Indeed, the special hypergeometric
$$f(a,b,z)=F(a,1-a;b;z)\,,$$
which is related to a Legendre function,
satisfies
$${f(a,b,{1\over 2_{\mathstrut}})\over \Gamma(b)}={2^{1-b}\,\sqrt{\pi}\over
\Gamma\bigl(\half (a+b)\bigr)\,\Gamma\bigl({1\over
2}(1-a+b)\bigr)}\,.\eqno(14.3)$$
This well-known relation can be obtained by applying
Euler's identity $$F(a,b;c;z)=(1-z)^{c-a-b}F(c-a,c-b;c;z)$$ and
Gauss's identities
$$\eqalign{
F(a,b;c;1)&=\bigl(\Gamma(c-a-b)\Gamma(c)\bigr)\big/
\bigl(\Gamma(c-a)\Gamma(c-b)\bigr)\,,\cr
F(2a,2b;a+b+\half ;z)&=F\bigl(a,b;a+b+\half ;
4z(1-z)\bigr)\,,\cr
}$$
which can be found, for example, in [\CM, (5.92), (5.111), exercise 5.28]):
$$\textstyle(1-z)^{1-b}F(a,1-a;b\mskip1mu;z)=F(b-a,b+a-1;b\mskip1mu;z)=
F\bigl(\half b-\half a,\half b+\half a-\half ;b\mskip1mu;
4z(1-z)\bigr)\,;
$$ we obtain (14.3) by letting $z\to\half $.

The sum (14.2) vanishes except when $k=3m$, and in this case the $k$\/th
term on the left of (14.1) reduces to simply $(\mu^3\!/6)^m\!/m!$
because of the formula
$$\Gamma({\textstyle{1\over 3}}-m)\,\Gamma({\textstyle{2\over
3}}-m)=3^{3m-1/2}2\pi\;{m!\over (3m)!}\,.\eqno(14.4)$$ 
Hence (14.1) is true. It is remarkable that so much of
nineteenth century mathematics has turned out to be relevant to
the study of random graphs.

When $\mu=0$, the generating function for the limiting
probabilities of excess~$r$ turns out to have a closed form: It is
$$\sum_{r\geq 0}\,\left(4\over3\right)^{\!r}\sqrt{2\over 3}\,{e_rr!\,z^r\over
(2r)!}=\sqrt{2\over3}\,F\left({1\over 6}\,,
\,{5\over 6}\,;\,\half \,;\,{z\over
2}\right) 
={2\cos\bigl({2\over 3_{\mathstrut}}\,{\rm arcsin}\,\sqrt{z/2}\,\bigr)\over
\sqrt{\mathstrut\mskip1mu6-3z}^{\mathstrut}}\,. \eqno(14.5)$$
From this expression it is easy to calculate the limiting value of the
mean excess when $m=\half n$, namely $\half -3^{-3/2}\approx
0.308$. The variance, similarly, is ${23\over 27}-3^{-3/2}$.

The limiting mean excess when the number of edges is ${1\over
2}n(1+\mu n^{-1/3})$ does not seem to have a simple closed form,
although we can express it as a hypergeometric series and find the
asymptotic value. Suppose we insert the factor~$z^r$ into the
left-hand side of (14.1). Then the left-hand side of (14.2) becomes
$${1\over \Gamma(1/2-2k/3)}\,F\left({1\over 6}\,,\,{5\over 6}\,;\,
\half -{2k\over 3}\,;\,{z\over 2}\right)
={1\over \Gamma(1/2-2k/3)}\,
f\left({1\over6}\,,\half -{2k\over3}\,,{z\over2}\right)\,.\eqno(14.6)$$
To evaluate the derivative of such a function at~$\half $, we
can use the identity
$$z(1-z)f'(a,b,z)=\left(az+{1-a-b\over 2}\right)f(a,b,z)-{1-a-b\over
2}\, f(-a,b,z)\,,\eqno(14.7)$$
which is readily verified by checking that the coefficients of~$z^n$
agree on both sides. To get the mean value of~$r$, we want to
differentiate (14.6) with respect to~$z$ and set $z=1$; and according to
(14.7), this is equivalent to replacing (14.6) by
$${1\over \Gamma(1/2-2k/3)}\,\bigl(({\textstyle\half +{2k\over
3}})f({\textstyle{1\over 6}\,,\,\half -{2k\over 3}\,,\,{1\over
2}})-({\textstyle{1\over 3}+{2k\over 3}})f(-{\textstyle{1\over 6}\,,\,
\half -{2k\over 3}\,,\,\half })\bigr)\,.\eqno(14.8)$$

Again, $f({1\over 6},\half -{2k\over 3},\half )$ vanishes
unless $k=3m$.  The contribution to the mean from this half of (14.8) is
just what we had when we were summing the probabilities, but with an
additional factor of $(\half +{2k\over 3})$; so it is
$$e^{-\mu^3\!/6}\sum_{m\geq 0}({\textstyle{1\over
2}}+2m)\,{(\mu^3\!/6)^m\over m!}=\half +{\mu^3\over
3}\,.\eqno(14.9)$$
The other half of (14.8) is, however, more complicated, since all values
of~$k$ make a contribution. According to (14.3), we want to evaluate
$$\sum_{k\geq 0}\,\sqrt{{2\pi\over 3}}\,{(\half 3^{2/3}\mu)^k\over
k!}\,
{2^{1/2+2k/3}\,\sqrt{\pi}\,(1/3+2k/3)\over\Gamma(1/6-k/3)\,\Gamma(5/6-k/3)}
=\Sigma_0+\Sigma_1+\Sigma_2\,,$$
where $\Sigma_j$ is a hypergeometric series corresponding to $k=3m+j$:
$$\eqalignno{%
\Sigma_0&={1\over 3\,\sqrt{3}\,}\,F\left({5\over 6}\,,\,{7\over
6}\,;\,{1\over 3}\,,\,{2\over 3}\,;\,{\mu^3\over 6}\right)\,;\cr
\noalign{\smallskip}
\Sigma_1&=-{1\over \sqrt{3}\,}\,
{\mu\,\sqrt{\pi}\over 6^{1/3}\Gamma(5/6)}
\,F\left({7\over 6}\,,\,{3\over
2}\,;\,{2\over 3}\,,\,{4\over 3}\,;\,{\mu^3\over 6}\right)\,;\cr
\noalign{\smallskip}
\Sigma_2&=-{5\,\sqrt{3}\over 2}\,{\mu^2\over 6^{2/3}}\,
{\sqrt{\pi}\over\Gamma(1/6)}
\,F\left({3\over 2}\,,\,{11\over
6}\,;\,{4\over 3}\,,\,{5\over 3}\,;\,{\mu^3\over
6}\right)\,.&(14.10)\cr}$$ 
As $z\rightarrow +\infty$, such hypergeometric series satisfy the
asymptotic formula
$$F(a,b\mskip1mu;c,d\mskip1mu;z)
={\Gamma(c)\Gamma(d)\over\Gamma(a)\Gamma(b)}\,z^{\delta}e^z
\left(1+{\delta(a+b-1)-ab+cd\over z}+O(z^{-2})\right)\,,\eqno(14.11)$$
where $\delta=a+b-c-d$; this follows by plugging the right-hand side
into the differential equation
$$\vartheta(\vartheta +c-1)(\vartheta +d-1)F=z(\vartheta +a)(\vartheta
+b)F$$ 
satisfied by the left. We obtain
$$\eqalignno{%
e^{-\mu^3\!/6}\Sigma_0&={1\over 3\,\sqrt{3}\,}\,{\Gamma({1\over
3})\Gamma({2\over 3})\over\Gamma({5\over 6})\Gamma({7\over 6})}\,
\left({\mu^3\over 6}+{1\over 4}+O(\mu^{-3})\right)\,;\cr
\noalign{\smallskip}
e^{-\mu^3\!/6}\Sigma_1&=-{1\over \sqrt{3}\,}\,
{\sqrt{\pi}\over\Gamma({5\over 6})}\,
{\Gamma({2\over
3})\Gamma({4\over 3})\over\Gamma({7\over 6})\Gamma({3\over 2})}\,
\left({\mu^3\over 6}+{1\over 4}+O(\mu^{-3})\right)\,;\cr
\noalign{\smallskip}
e^{-\mu^3\!/6}\Sigma_2&=-{5\sqrt{3}\over 2}\,
{\sqrt{\pi}\over \Gamma({1\over 6})}\,
{\Gamma({4\over
3})\Gamma({5\over 3})\over\Gamma({3\over 2})\Gamma({11\over 6})}\,
\left({\mu^3\over 6}+{1\over 4}+O(\mu^{-3})\right)\,;&(14.12)\cr}$$
therefore $e^{-\mu^3\!/6}(\Sigma_0+\Sigma_1+\Sigma_2)=({2\over
3}-{4\over 3}-{4\over 3})\bigl({\mu^3\over 6}+{1\over
4}+O(\mu^{-3})\bigr)$.
Subtracting this from (14.9), and using computer algebra to refine the
estimate further, gives us the answer we seek:

\proclaim
Theorem 6. The expected value of the excess, when there are ${1\over
2}n(1+\mu n^{-1/3})$ edges, approaches
$${2\over 3}\mu^3+1+{5\over 24}\mu^{-3}+{15\over
16}\mu^{-6}+O(\mu^{-9})\,, \eqno(14.13)$$
for fixed $\mu\geq\delta >0$ as $n\rightarrow \infty$.\quad\pfbox

\medskip
This method of calculation shows also that the variance will be
$O(\mu^6)$ and the $k$\/th moment will be $O(\mu^{3k})$; each
derivative of (14.6) can do no worse than multiply by~$\mu^3$, because
of (14.7).

Incidentally, 
 the $O(\mu^{-3})$ terms in all three equations of (14.12)
turn out to equal ${5\over 48}\mu^{-3}+{15\over
32}\mu^{-6}+O(\mu^{-9})$, and this is no coincidence. We have, in fact,

$$\eqalignno{
{\Gamma({5\over 6})\Gamma({7\over 6})\over\Gamma({1\over
3})\Gamma({2\over 3})}\,
F\left({5\over 6}\,,\,{7\over 6}\,;\,{1\over 3}\,,\, {2\over
3}\,;z^3\right)
&\sim{\Gamma({7\over 6})\Gamma({9\over 6})\over\Gamma({2\over
3})\Gamma({4\over 3})}\,z\,F\left({7\over 6}\,,\,{3\over
2}\,;\,{2\over 3}\,,\,{4\over 3}\,;z^3\right)\cr
\noalign{\smallskip}
&\sim {\Gamma({3\over 2})\Gamma({11\over 6})\over\Gamma({4\over
3})\Gamma({5\over 3})}\,z^2F\left({3\over 2}\,,\,{11\over
6}\,;\,{4\over 3}\,,\,{5\over 3}\,;z^3\right)\,,&(14.14)\cr}$$
in the sense that all three functions have the same asymptotic series
$\sum s_kz^{-3k}e^{z^3}$ as $z\rightarrow \infty$. 
This follows because all three functions satisfy the same differential
equation, and because their asymptotic behavior depends only on the
differential equation except for a constant of proportionality. It is
well known that the general hypergeometric functions
$F(a_1,\ldots,a_m;\,b_1,\ldots,b_n;\,z)/\Gamma(b_1)\,
\ldots\,\Gamma(b_n)$
and
$z^{1-b_1}F(a_1+1-b_1,\ldots,a_m+1-b_1;\,b_2+1-b_1,
\ldots,b_n+1-b_1,2-b_1;\,z)/\Gamma(b_2+1-b_1)\,
\ldots\,\Gamma(b_n+1-b_1)$
both satisfy the differential equation
$\vartheta(\vartheta+b_1-1)\,\ldots\,(\vartheta+b_n-1)F=z(\vartheta
+a_1)\,\ldots\,(\vartheta +a_m)F$. In the case of (14.14), even more
is true: The three asymptotically equivalent functions shown there can
be written respectively as ${1\over 3}\bigl(G(z)+G(\omega z)
+G(\omega^2z)\bigr)$,
${1\over 3}\bigl(G(z)+\omega^2G(\omega z)+\omega G(\omega^2z)\bigr)$, 
${1\over 3}\bigl(G(z)
+\omega G(\omega z)+\omega^2G(\omega^2z)\bigr)$, where
$$G(z)={1\over\sqrt{3}}\,\sum_{k\geq 0}\,{\Gamma(3/2+k)z^k\over
\Gamma(1/2+k/3)\,k!}\eqno(14.15)$$
and $\omega=e^{2\pi i/3}$.

\bigbreak\noindent
{\bf 15. Deficiency, planarity, and complexity.}\enspace
The calculations in the preceding section can be combined with the
structure theory of section~9 to yield the following general result.

\proclaim
Theorem 7. Let $\MM$ be a reduced multigraph
of excess~$r$ and deficiency~$d$, i.e., a reduced multigraph having $2r-d$
vertices and $3r-d$ edges. The probability that the complex part of a random
graph or multigraph reduces to $\MM$ is asymptotically
$${\sqrt{2\pi}\,\kappa(\MM)\over (2r-d)!}
A(3r+{\textstyle\half }-d,\mu)\,n^{-d/3}\eqno(15.1)$$
when there are $\half n(1+\mu n^{-1/3})$ edges and $n$ vertices,
$|\mu|=o(n^{1/12})$. Here $\kappa(\MM)$ denotes the compensation
factor (1.1), and $A(y,\mu)$ is defined in~(10.2).
The sum of (15.1) over all $\MM$ of deficiency~0 is~1.
For each $d\ge0$, the probability that a random multigraph
has deficiency $\ge d$ is $O\bigl((1+\mu^4)^dn^{-d/3}\bigr)$,
uniformly in $n$ and~$\mu$.

\proof
When $d=0$, this theorem is a consequence of the
corollary following~(9.16), together with (13.17) and (14.1).

When $d>0$, 
(15.1) is clear, but we need two auxiliary results of independent interest
before we can prove the desired uniform estimate.

\proclaim
Lemma 4. Let $E_{r(\geq d)}$ denote the generating function for all
complex multigraphs of excess~$r$ whose deficiency is at least~$d$.
Then
$${e_{rd}T\over (1-T)^{3r-d}}-{(2r-d-1)e_{rd}T\over
(1-T)^{3r-d-1}}\leq E_{r(\geq d)}\leq{e_{rd}T\over
(1-T)^{3r-d}}\,,\eqno(15.2)$$ 
where inequality between generating functions means that the
coefficients of every power of~$z$ obey the stated relation.

\proof
The claim is trivial when $d=2r-1$; and it is true when $r=1$, because
$E_1={5\over 24}\,T(1-T)^{-3}-{1\over 12}\,T(1-T)^{-2}$. The lower
bound is easily seen to be a lower bound on
$e_{rd}T^{2r-d}/(1-T)^{3r-d}$ itself.

The proof of the upper bound now proceeds by induction on $r$. Let
$$E'_r=\sum_{k=0}^{d-1}e_{rk}(1+\zeta)^r\zeta^{2r-k}+e_{rd}\zeta(1+\zeta)
^{3r-d-1}\,,\eqno(15.3)$$
in the notation of section 5. We want to prove that $E_r\leq E'_r$; it
suffices to show, by (5.8), that
$$\bigl(r+(1-T(z))\vartheta\bigr)E'_r=
\bigl(r+(1+\zeta)^{-1}\vartheta\bigr)E'_r\geq{\textstyle\half }
\,\bigl({\textstyle{1\over
2}}\,\zeta(1+\zeta)+\vartheta\bigr)^2E'_{r-1}\,,\eqno(15.4)$$
considering both sides as generating functions in powers of~$z$.
Proceeding as in (5.11) and (5.12) to form
$$\textstyle{A'_r=\bigl({1\over
2}\,\zeta(1+\zeta)+\vartheta\bigr)E'_r\,,\qquad
B'_r=\bigl({1\over
2}\,\zeta(1+\zeta)+\vartheta\bigr)A'_r\,,}$$
a bit of algebra shows that when $0\leq d\leq 2r-3$ we have
$$\eqalignno{&\bigl(r+(1+\zeta)^{-1}\vartheta\bigr)E'_r-\textstyle\half 
B'_{r-1}\cr
\noalign{\smallskip}
&\qquad =\half \zeta(1+\zeta)^{r+1}\biggl(\sum_{k\geq 0}\zeta^{2r-d-2-k}
\bigl((\alpha_k+\beta_k)e_{rd}-
(\gamma_k+\delta_k+\epsilon_k)e_{(r-1)d}\bigr)\cr
\noalign{\smallskip}
&\qquad\qquad\qquad\null-(2r-1-d)^2e_{(r-1)(d-1)}\zeta^{2r-d-2}\biggr)\,,
&(15.5)\cr}$$
where
$$\displaylines{\alpha_k=2(3r-d){2r-d-2\choose k+1}\,,\qquad
\beta_k=2(r+1){2r-d-2\choose k}\,;\cr
\noalign{\smallskip}
\gamma_k=\bigl(3r-{\textstyle\half }-d\bigr)
         \bigl(3r-{\textstyle{5\over 2}}-d\bigr){2r-d-3\choose k+1}\,,\qquad
\delta_k=\bigl(9r-{\textstyle{13\over 2}}-3d\bigr){2r-d-3\choose
k}\,,\cr
\noalign{\smallskip}
\epsilon_k={2r-d-3\choose k-1}\,.\cr}$$
Obviously $\beta_ke_{rd}\geq \epsilon_ke_{(r-1)d}$,
since $e_{rd}\geq e_{(r-1)d}$. And the inequality $9r-{13\over
2}-3d\leq (3r-\half -d)(3r-{5\over 2}-d)$ for $0\leq d\leq 2r-3$ yields
$$(\gamma_k+\delta_k)e_{(r-1)d}\leq \bigl(3r-{\textstyle{1\over
2}}-d\bigr)\bigl(3r-{\textstyle{5\over 2}}-d\bigr){2r-d-2\choose k+1}
e_{(r-1)d}\le \alpha_ke_{rd}\,.$$
In fact, (5.11)--(5.13) imply that
$$2(3r-d)e_{rd}\ge \bigl(3r-{\textstyle{1\over
2}}-d\bigr)\bigl(3r-{\textstyle{5\over 2}}-d\bigr)e_{(r-1)d}
+(2r-d)(2r-1-d)e_{(r-1)(d-1)}\,;$$
so (15.5) is a polynomial in $\zeta$ with nonnegative coefficients,
and thus a power series in $z$ with nonnegative coefficients, proving
(15.4). The case $d=2r-2$
needs to be handled separately, but it offers no
difficulty.\quad\pfbox

\proclaim Lemma 5.
There exists a constant $\epsilon>0$ such that, for every fixed $d\ge0$,
a random multigraph with $n$ vertices and $m={n\over2}(1+\mu n^{-1/3})$
edges has excess~$r$ and deficiency $\ge d$ with probability
$$\cases{O(\mu^{4d-3/2}n^{-d/3}e^{-\epsilon(r-{2\over3}\mu^3)^2\!/\mu^3}),
&if $r\le\mu^3$,\cr
O(n^{-d/3}e^{-\epsilon r}),& if $r\ge\mu^3$,\cr}$$
uniformly in $n$, $r$, and $\mu$ when $\mu\le n^{1/12}$.

\proof
Let $p_{rd}=p_{rd}(n,\mu)$ be the stated probability. It suffices to
prove the lemma when $\mu\ge1$ and $r\ge1$. For if $r=d=0$, the result
follows from Lemma~3;
and $p_{0d}=0$ when $d>0$. On the other hand, if $\mu<1$ we
have $p_{rd}(n,\mu)\le \sum_{j=r}^\infty p_{jd}(n,1)$.

Using Lemma 4 and arguing as in the proof of Lemma 3, equation (10.11),
we obtain
$$\openup2\jot
\eqalign{p_{rd}&={2^mm!\,n!\over n^{2m}}\;[z^n]\;{U^{n-m+r}\over
(n-m+r)!}\,e^V\,E_{r(\ge d)}\cr
&\le{2^mm!\,n!\over n^{2m}}\;[z^n]\;{U^{n-m+r}\over(n-m+r)!}
{e_{rd}\,T\over(1-T)^{3r-d+1/2}}\cr
&={2^mm!\,n!\,e_{rd}\,e^n@2^{m-n-r}\over(n-m+r)!\,n^{2m}\,2\pi i}\oint
\left(z(2-z)\over1-z\right)^{\!r}(1-z)^{d-2r+1/2}e^{nh(z)}dz\,,\cr}$$
with $h(z)$ as in (10.12), and where the integral is taken around a
circle $z=\rho e^{i\theta}$ with $0<\rho<1$. On this circle, both
$\bigl|(2-z)/(1-z)\bigr|$ and $|1-z|^{-1}$ attain their maxima at
$z=\rho$. Moreover, by (10.16) we have ${d\over d\theta}\Re h=-\rho@
g(\theta)\sin\theta$, where $g(\theta)>\bigl((2-\rho)^2-1\bigr)/9>
{2\over9}(1-\rho)$; therefore
$$\Re h(\rho e^{i\theta})\le h(\rho)+{2\over9}(1-\rho)\rho@(\cos\theta-1)
\le h(\rho)-{4\over9\pi^2}\rho@(1-\rho)\theta^2,\qquad\hbox{for
$|\theta|\le\pi$.}$$
Now $p_{rd}=0$ if $d\ge2r$, because we are assuming that $r\ge1$.
Hence $d-2r+1/2<0$, and the contour integral including the factor
$1/(2\pi i)$ is less than
$$\eqalign{&{\rho\over2\pi}\left(\rho@(2-\rho)\over1-\rho\right)^{\!r}(1-\rho)^{
d-2r+1/2}e^{nh(\rho)}\int_{-\pi}^{\pi}\exp\left(-{4n\rho@(1-\rho)\over
9\pi^2}\theta^2\right)d\theta\cr
&\qquad<{3\over4}\sqrt{\pi\over n}\rho^{r+1/2}(2-\rho)^r(1-\rho)^{d-3r}
e^{nh(\rho)}\,.\cr}$$

In the following argument, unspecified positive constants will be denoted
by $\epsilon_1$, $\epsilon_2$, \dots, while positive numbers that may
depend on~$d$ will be denoted by $C_1$, $C_2$, \dots~. Let $\nu=n^{-1/3}$.
If we apply (10.10) to the coefficient in front of the contour integral,
and if we use the estimate
$${(n-m+r)!\over(n-m)!}>(n-m)^r=\left(n\over2\right)^{\!r}(1-\mu\nu)^r
>\left(n\over2\right)^{\!r}e^{-2\mu\nu r},$$
which is valid when $\mu\nu\le\half$, we obtain the upper bound
$$p_{rd}\le C_1@e_{rd}@n^{-r}\rho^r@(2-\rho)^r(1-\rho)^{d-3r}
e^{nh(\rho)-\mu^3\!/6+2\mu\nu r},\eqno(15.6)$$
where $\rho$ is any number between 0 and 1.

Suppose now that $r\le12\mu^3$, and set $\rho=1-\xi\mu\nu$,
$r={2\over3}x\mu^3$. If $\xi=O(1)$, we have
$$n@h(1-\xi\mu\nu)={1\over3}\xi^3\mu^3+\half\xi^2\mu^3+O(1)$$
as in (10.17). Therefore, since $\rho@(2-\rho)<1$,
$$\eqalign{p_{rd}
&\le C_2@e_{rd}@n^{-r}(\xi\mu\nu)^{d-3r}e^{\xi^3\mu^3\!/3+\xi^2\mu^3\!/2-
\mu^3\!/6}\cr
\noalign{\smallskip}
&\le C_3@3^r@2^{-r}r^{r+d-1/2}e^{-r}n^{-d/3}@(\xi\mu)^{d-3r}
e^{\xi^3\mu^3\!/3+\xi^2\mu^3\!/2-\mu^3\!/6}\cr
\noalign{\smallskip}
&=C_3@r^{d-1/2}n^{-d/3}@(\xi\mu)^d e^{k(x,\xi)\mu^3\!/6},
\quad k(x,\xi)=2\xi^3+3\xi^2-1+4x\ln\left(x\over e\xi^3\right)\,,\cr}$$
by (7.16) and Stirling's formula. Given $x$ between 0 and 18, we minimize
$k(x,\xi)$ by letting $\xi$ be the positive root of $\xi^3+\xi^2=2x$;
notice that this makes $\xi\le3$, justifying our assumption that
$\xi=O(1)$. The minimum $k(x,\xi)$ satisfies
$$\openup1\jot
\eqalign{k(x,\xi)&=\xi^2-1+2(\xi^3+\xi^2)\ln\left(1+\xi^{-1}\over2\right)\cr
&\le\xi^2-1+2(\xi^3+\xi^2)\,{\xi^{-1}-1\over2}=-(\xi+1)(\xi-1)^2.\cr}$$
We also have $|\xi-1|\ge\epsilon_1|x-1|$, hence $k(x,\xi)\le
-\epsilon_2(x-1)^2$. Our estimates have shown that
$$p_{rd}\le C_5 @r^{d-1/2}n^{-d/3}\mu^d e^{-\epsilon_2(x-1)^2\mu^3\!/6},$$
when $r={2\over3}x\mu^3\le12\mu^3$, so the first half of the lemma
has been proved.

When $r\ge12\mu^3$, let us set $\rho=1-\eta$ and $r=yn$. In this case we will
in fact prove the lemma for a much larger range of~$\mu$, assuming only
that $\mu\nu\le\delta$, when $\delta$ is a suitably small constant. If
$\delta\le{1\over5}$ we can assume that $0<y<{3\over5}$, since
$m$ is at most ${1+\delta\over2}n$ and since $p_{rd}=0$ when $r\ge m$.
Using (7.16) and (15.6) again, we find
$$p_{rd}\le C_6@ r^{d-1/2}\eta^d e^{nl(y,\eta)-\mu^3\!/6+2r\mu\nu},$$
where
$$\openup1\jot\eqalign{l(y,\eta)
&=y\ln\left(3y(1-\eta^2)\over2e\eta^3\right)+h(1-\eta)\cr
&=y\ln\left(3y(1-\eta^2)\over2e\eta^3\right)-\eta-\ln(1-\eta)+
{1-\mu\nu\over2}\ln(1-\eta^2)\,.\cr}$$

Given $y$, the minimum value of $l(y,\eta)$ occurs when $y=\eta^2(\eta+
\mu\nu)/(3-\eta^2)$. However, we do not need to find the exact minimum,
in order
to achieve the upper bound in the lemma; it will suffice to be close to the
minimum when $y$ is small. Therefore we choose $\eta$ in such a way that
the calculations will be relatively simple:
$$y\;=\;{2\eta^3\over3(1-\eta^2)}\,.\eqno(15.7)$$
With this choice, we always have $\eta<{3\over4}$; and
$$l(y,\eta)=f(\eta)=-{2\eta^3\over3(1-\eta^2)}-\eta-\ln(1-\eta)
+{1-\mu\nu\over2}\ln(1-\eta^2)\,.$$
If we set $\eta=\mu\nu$, this function $f(\eta)$ reduces to
$$\sum_{k=1}^\infty\eta^{2k+1}\left({1\over2k}+{1\over2k+1}-{2\over3}\right)
<{\eta^3\over6}={\mu^3\over 6n}\,.$$
On the other hand, the actual value of $\eta$ must be larger than $2\mu\nu$,
because $2\mu\nu$ is too small to satisfy (15.7):
$${2(2\mu\nu)^3\over3\bigl(1-(2\mu\nu)^2\bigr)}\le
{16(\mu\nu)^3\over3(1-{4\over25})}={400\mu^3\over63n}<{r\over n}=y\,.$$
When $\eta>\mu\nu$ we have
$$f'(\eta)=-{\eta\bigl(\eta^3+3\eta^2\mu\nu+3(\eta-\mu\nu)\bigr)\over
3(1-\eta^2)^2}<-\eta(\eta-\mu\nu)<-(\eta-\mu\nu)^2\,;$$
hence when $\eta$ satisfies (15.7) we have
$$l(y,\eta)<{\mu^3\over6n}-{(\eta-\mu\nu)^3\over3}
\le{\mu^3\over6n}-\epsilon_3@\eta^3
\le{\mu^3\over6n}-\epsilon_4@y\,.$$
We have proved that
$$p_{rd}\le C_7@r^d y^{d/3}e^{-\epsilon_4r+2r\mu\nu},$$
and this is at most $C_8n^{-d/3}e^{-\epsilon_5r}$ if $\delta$ is
less than $\half\epsilon_4$.\quad\pfbox

\medskip
Returning to the proof of Theorem 7,
its final claim now follows for $\mu\le n^{1/12}$ by summing the
upper bounds of Lemma~5 over all values of~$r$.
The claim is trivial when $\mu>n^{1/12}$.\quad\pfbox

\medskip
As remarked earlier, the fact that (15.1) sums to 1 allows us to compute
asymptotic probabilities of any collection of graphs or multigraphs obtained
as a union over an infinite set of reduced multigraphs, as long as
at least one multigraph in the set is clean (has deficiency zero).
We simply sum the individual probabilities, neglecting unclean
cases.

One corollary of Theorem 7 is the fact that a random graph with
$\half (n+\mu n^{2/3})$ edges is clean with probability $1-o(1)$
whenever $\mu=o(n^{1/12})$. Stepanov proved this for $\mu\leq 0$ [\Sii,
Theorem~3] and conjectured that it would also hold for positive~$\mu$.
His conjecture was proved for all fixed~$\mu$ by \L uczak, Pittel, and
Wierman [\LPW].

Erd\H os and R\'enyi remarked in their pioneering paper [\ER, \S8] that,
if $x$ is any real number, the probability that a graph with
$\half n+xn^{1/2}$ edges is nonplanar ``has a positive lower limit,
but we cannot calculate its value. It may even be~1, though this
seems unlikely.''  They gave no proof that the limiting probability
is positive, and their remark was embedded in a section of~[\ER] that
contains a technical error (see [\LW]); but a proof of their assertion
was found later by Stepanov~[\Sii, Corollary~2 following~(10)]. In the other
direction, the fact that nonplanarity occurs with probability strictly
less than~1 follows from the fact that a graph with $\half n+o(n^{2/3})$
edges has excess~0 with probability~$\sqrt{2\over3}$, as observed
in~[\FKP, Corollary~8].

We are now in a position to make a more precise estimate of the
probability in question.

\proclaim
Theorem 8. The probability that a graph with $\half n+o(n^{2/3})$
edges is nonplanar approaches a limit~$\rho$ as $n\to\infty$, where
$$0.000229\le\rho\le0.012926.\eqno(15.8)$$

\proof
The condition $m=\half n+o(n^{2/3})$ is equivalent to saying that
$\mu=o(1)$ when $m=\half n(1+\mu n^{-1/3})$, so we can let $\mu=0$
in the asymptotic formulas above. By Theorem 7, the constant~$\rho$
is the sum $\sum\sqrt{\mskip1mu2\pi}A(3r+\half ,0)\kappa(\MM)/(2r)!=
\sum\sqrt{2\over3}\bigl({4\over3}\bigr)^{r}r!\,\kappa(\MM)/(2r)!^2$
over all nonplanar, reduced, labeled, clean multigraphs~$\MM$, where
$r=r(\MM)$ is the excess of~$\MM$.

A clean multigraph cannot contain a subgraph that is homeomorphic to the
complete graph~$K_5$, i.e., a subgraph that cancels to~$K_5$, because
$K_5$ has deficiency~5. Adding an edge to any multigraph increases the
excess by 0 or~1 and increases the deficiency by 0, 1, or~2 (see
section~20 below for further discussion); hence all subgraphs of a clean
multigraph are clean. Indeed, this argument implies that a
random graph with $\half n+o(n^{2/3})$ edges has probability $O(n^{-5/3})$
of containing a~$K_5$.

Therefore, if a sparse graph or multigraph is nonplanar, its nonplanarity
comes almost surely from a subgraph that cancels to the complete
bipartite graph~$K_{3,3}$, which is clean and has excess~3.

One way to
obtain bounds on $\rho$ is to restrict consideration to reduced
multigraphs whose components all have excess $\le3$. If such a multigraph
contains a $K_{3,3}$, it corresponds only to nonplanar graphs; if it
does not, it corresponds only to planar graphs. The difference between the
upper and lower bounds so obtained is the probability that a random
graph of $\half n+o(n^{2/3})$ edges has at least one component
of excess~$\ge4$, i.e., that at least one component is more than tetracyclic.

The multigraph $K_{3,3}$ has compensation factor~1, because it is a
graph, and its vertices can be labeled in $\half {6\choose 3}=10$
different ways. Thus it contributes only
${10\over6!}={1\over 72}$ to the constant 
$c_3={1105\over1152}$ that accounts for all clean
connected multigraphs of excess~3.

Let $f_r=[z^r]\,
\exp(c_1z+c_2z^2+c_3z^3)$ and $g_r=[z^r]\,\exp\bigl(c_1z+
c_2z^2+(c_3-{1\over72})z^3\bigr)$. Then the quantities
$$p=\sum_{r\ge0}\sqrt{2\over3}
\left(4\over3\right)^{\!r}f_r{r!\over(2r)!}
\qquad\hbox{and}\qquad
  q=\sum_{r\ge0}\sqrt{2\over3}
\left(4\over3\right)^{\!r}g_r{r!\over(2r)!}
$$
are respectively the probability that a sparse graph has all
components of excess~$\le3$ and the probability that, moreover,
no component cancels to~$K_{3,3}$. These series converge rapidly and
lead to the numerical bounds $p-q$ and $1-q$ in~(15.8).\quad\pfbox

\medskip
It is interesting to study the expected number E$n_1$ of vertices in complex
components, as a function of~$\mu$, because it will be the expected
number of vertices in the giant component when $\mu$ increases. We have
${\rm E}n_1=\sum_r{\rm E}(n_1@|@ r)\Pr({\cal E}_r)$.
By (13.17) and the
remarks preceding (13.13), each term in this sum can be approximated, to
within relative error $O\bigl((1+\mu^4)n^{-1/3}\bigr)$, by
$3r\sqrt{2\pi}e_r@A(3r+{5\over2},\mu)@n^{2/3}$. Let us, for simplicity,
assume that $\mu$ is bounded. Then the proof of Lemma~5 is easily
modified to show that the $r$th term of the sum is $O\bigl(n^{2/3}(r+1)
e^{-\epsilon r}\bigr)$, uniformly in $n$ and~$r$. Thus, by dominated
convergence, E$n_1=\bigl(f(\mu)+o(1)\bigr)n^{2/3}$, where
$$f(\mu)=\sum_{r\geq 0}3r\sqrt{2\pi}\,e_r\,A(3r+{\textstyle{5\over
2}},\mu)\,.\eqno(15.9)$$
Equation (10.23) tells us that
$$\textstyle{\half \,\mu^2f(\mu)+f'(\mu)=\half \,\sum_{r\geq
0}3r\sqrt{2\pi}\,e_r\,A(3r+\half ,\mu)={3\over
2}\,g(\mu)\,,}\eqno(15.10)$$
where $g(\mu)$ is the expected value of $r$; we calculated $g(\mu)$ in
the discussion leading up to Theorem~6. Thus, we obtain the estimate
$$f(\mu)=2\mu-\mu^{-2}-{\textstyle{27\over 8}}\,\mu^{-5}-
{\textstyle{495\over 16}}\,\mu^{-8}+O(\mu^{-11})\,,\eqno(15.11)$$
for $\mu\geq\delta>0$, by combining (15.9)
with the asymptotic formula for $g(\mu)$ in (14.13).

We can express $f(\mu)$ in ``closed hypergeometric form'' by
proceeding as in (14.9) and (14.10). The result is
$$\eqalignno{f(\mu)&=-{2^{-2/3}\pi\over 3^{7/6}\,\Gamma\bigl({2\over
3}\bigr)}\;e^{-\mu^3\!/6}+\mu-{\mu\over 4}e^{-\mu^3\!/6}\,F\left({1\over3}\,;
\,{4\over3};\,{\mu^3\over 6}\right)\cr
\noalign{\smallskip}
&\qquad\null+e^{-\mu^3\!/6}\left({2^{1/3}\,\sqrt{\pi}\over
3^{7/6}\,\Gamma\bigl({7\over 6}\bigr)}\; F\left({1\over
2}\,,\,{5\over
6}\,;\, {1\over 3}\,,\,{2\over 3}\,;\,{\mu^3\over 6}\right)\right.\cr
\noalign{\smallskip}
&\qquad\qquad\qquad\qquad\null-{\mu\over 2\,\sqrt{3}}\,
F\left({5\over
6}\,,\,{7\over 6}\,;\,{2\over 3}\,,\,{4\over 3}\,;\,{\mu^3\over
6}\right)\cr
\noalign{\smallskip}
&\qquad\qquad\qquad\qquad\null+\left.{3^{1/6}\,\sqrt{\pi}\,\mu^2\over
2^{7/3}\,\Gamma\left({5\over 6}\right)}\,
F\left({7\over 6}\,,\,{3\over
2}\,;\,{4\over 3}\,,\,{5\over 3}\,;\,{\mu^3\over
6}\right)\right)\,.&(15.12)\cr}$$ 
It is instructive to compare this expression with alternative formulas for the
same quantity obtained in [\LPW] by a different method:
$$\eqalignno{f(\mu)&={1\over\sqrt{2\pi}}\int_0^\infty\biggl(\sum_{r\ge1}f_rx^{3r/2}\biggr)e^{G(x,\mu)}dx\cr
&=\mu+{1\over
\sqrt{2\pi}}\,\int_0^{\infty}\,{1-e^{G(x,\mu)}\over x^{3/2}}\,dx -
{1\over 4}\,\int_0^{\infty}e^{G(x,\mu)}\,dx\,.&(15.13)\cr}$$
Here $G(x,\mu)=\bigl((\mu-x)^3-\mu^3\bigr)/6$, and $f_rn^{n+(3r-1)/2}$
is Wright's asymptotic estimate~[\Wiii] for the number of connected
graphs with excess~$r$.

\bigbreak\noindent
{\bf 16. Evolutionary paths.}\enspace
Consider any graph or multigraph that evolves by starting out with
isolated vertices and then by acquiring edges one at a time. Initially
its excess is~0; then each new edge either preserves the current
excess or increases it by~1. We observed in section~4, following~(4.7),
that a new edge
 augments the excess if and only if both of its endpoints
are currently in the cyclic part. We observed in section~13 that
many interesting statistics about random graphs can be usefully
represented in terms of probabilities that are conditional on
the graph having a given excess. Therefore it is natural to look more
closely at the way a graph changes character as its excess grows.

Every evolution of a graph or multigraph traces a  path from
left to right in the following diagram, which shows the beginning of
an infinite partial ordering of all possible configurations
$[r_1,r_2,\ldots,r_q]$:

% illustration for "Birth of the Giant Component"

\font\bmit=cmmib10

\def\bmath{\textfont0=\tenbf\textfont1=\bmit}
\unitlength=1cm
\centerline{\beginpicture(13,10)(0,-5)
\put(0,0){\bmath\makebox(0,0){$[0]$}}
\put(3,0){\bmath\makebox(0,0){$[1]$}}
\put(6,1){\bmath\makebox(0,0){$[0,1]$}}
\put(6,-1){\bmath\makebox(0,0){$[2]$}}
\put(9,2){\bmath\makebox(0,0){$[0,0,1]$}}
\put(9,0){\bmath\makebox(0,0){$[1,1]$}}
\put(9,-2){\bmath\makebox(0,0){$[3]$}}
\put(12,4){\bmath\makebox(0,0){$[0,0,0,1]$}}
\put(12,2){\bmath\makebox(0,0){$[1,0,1]$}}
\put(12,0){\bmath\makebox(0,0){$[0,2]$}}
\put(12,-2){\bmath\makebox(0,0){$[2,1]$}}
\put(12,-4){\bmath\makebox(0,0){$[4]$}}
\def\\#1/#2//{\makebox(0,0){\smash{\raise12pt\hbox{$({#1\over#2})$}}}}%
\put(0,0){\\1/1//}
\put(3,0){\\1/1//}
\put(6,1){\\72/77//}
\put(6,-1){\\5/77//}
\put(9,2){\\15912/17017//}
\put(9,0){\\1080/17017//}
\put(9,-2){\\25/17017//}
\put(12,4){\\7029504/7436429//}
\put(12,2){\\318240/7436429//}
\put(12,0){\\77760/7436429//}
\put(12,-2){\\10800/7436429//}
\put(12,-4){\\125/7436429//}
\def\\#1/#2//{\makebox(0,0){$#1\over#2$}}%
\put(0.4,0){\line(1,0){.6}}
\put(1.2,0){\\1/1//}
\put(1.4,0){\vector(1,0){1.3}}
\put(3.3,.1){\line(3,1){.9}}
\put(4.5,.5){\\72/77//}
\put(4.8,.6){\vector(3,1){.7}}
\put(3.3,-.1){\line(3,-1){.9}}
\put(4.5,-.6){\\5/77//}
\put(4.8,-.6){\vector(3,-1){.8}}
\put(6.6,1.2){\line(3,1){.5}}
\put(7.4,1.6){\\216/221//}
\put(7.7,1.567){\vector(3,1){.6}}
\put(6.2,.8){\line(3,-1){.4}}
\put(6.9,0.6){\\5/221//}
\put(7.2,.467){\vector(3,-1){1.1}}
\put(6.3,-.95){\line(3,4){.4}}
\put(7.0,-.15){\\72/221//}
\put(7.15,.18){\vector(3,4){1.0}}
\put(6.4,-1.016){\line(3,1){.5}}
\put(7.25,-.85){\\144/221//}
\put(7.6,-.616){\vector(3,1){.9}}
\put(6.4,-1.233){\line(3,-1){.4}}
\put(7.15,-1.5){\\5/221//}
\put(7.5,-1.6){\vector(3,-1){1.1}}
\put(9.6,2.4){\line(3,2){.4}}
\put(10.3,2.8){\\432/437//}
\put(10.6,3.033){\vector(3,2){.8}}
\put(9.65,2){\line(1,0){.2}}
\put(10.2,2.1){\\5/437//}
\put(10.5,2){\vector(1,0){.8}}
\put(9.6,.6){\line(3,4){.4}}
\put(10.3,1.4){\\144/437//}
\put(10.5,1.8){\vector(3,4){1.3}}
\put(9.6,.3){\line(3,2){.5}}
\put(10.475,0.85){\\216/437//}
\put(10.8,1.1){\vector(3,2){.7}}
\put(9.475,0){\line(1,0){.2}}
\put(10.0,0.07){\\72/437//}
\put(10.3,0){\vector(1,0){1.1}}
\put(9.06,-.24){\line(3,-2){.34}}
\put(9.7,-.5){\\5/437//}
\put(9.883,-.783){\vector(3,-2){1.5}}
\put(9.6,-1.4){\line(3,4){.7}}
\put(10.55,-0.4){\\216/437//}
\put(10.575,-.1){\vector(3,4){1.3}}
\put(9.3,-2){\line(1,0){.4}}
\put(10.05,-2){\\216/437//}
\put(10.4,-2){\vector(1,0){1.0}}
\put(9.15,-2.3){\line(3,-2){.55}}
\put(9.9,-2.8){\\5/437//}
\put(10.2,-3.0){\vector(3,-2){1.4}}
\put(12.9,4.3){\vector(3,2){.5}} % corrected
\put(12.9,3.8){\vector(3,-1){.5}}
\put(12.7,2.3){\vector(1,1){.6}} % corrected
\put(12.7,2.0){\vector(3,1){.6}} % corrected
\put(12.7,1.9){\vector(3,-2){.6}}
\put(12.65,1.8){\vector(3,-4){.65}}
\put(12.6,.15){\vector(3,1){.7}}
\put(12.6,0){\vector(1,0){.7}}
\put(12.6,-.15){\vector(3,-1){.7}}
\put(12.5,-1.84){\vector(1,1){.8}}
\put(12.5,-1.93){\vector(2,1){.8}}
\put(12.5,-2.0){\vector(1,0){.8}}
\put(12.5,-2.1){\vector(2,-1){.8}}
\put(12.5,-2.2){\vector(1,-1){.8}}
\put(12.3,-3.9){\vector(2,1){1.0}}
\put(12.3,-4.0){\vector(1,0){1.0}}
\put(12.3,-4.2){\vector(3,-2){1.0}} % corrected
\endpicture
}
%\bigskip
{\narrower\narrower\smallskip\noindent
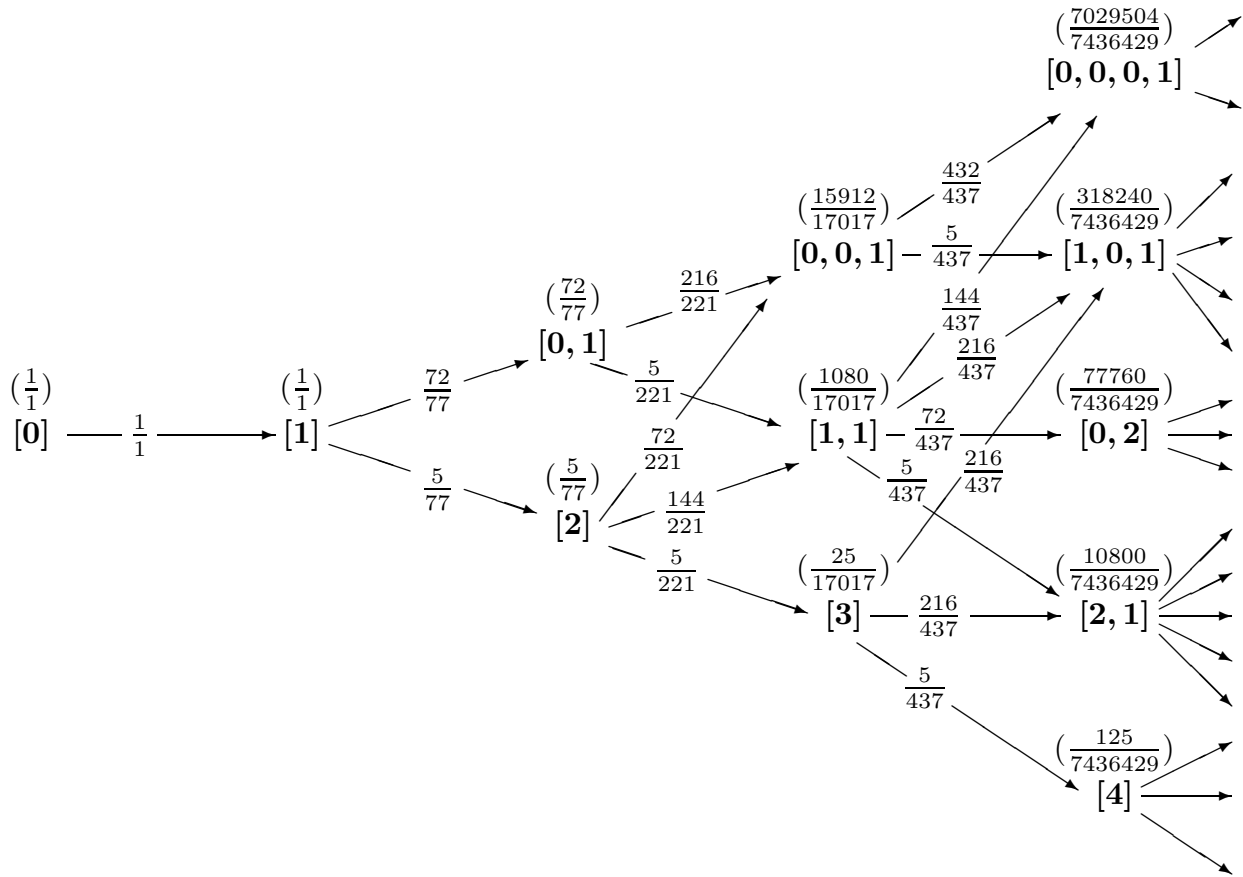
{\bf Figure 1.}\enspace The evolution of complex components. \ Each
configuration $[r_1,r_2,\ldots,r_q]$ stands for a graph or multigraph with
$r_1$~bicyclic components, $r_2$~tricyclic components, \dots,
$r_q$~($q+1$)-cyclic components. As a graph evolves, its excess
$r_1+2r_2+3r_3+\cdots\,$ increases in unit steps,
and the configurations follow a path from left to right in this
partial ordering.
\smallskip}

\bigskip
\noindent
When complex components begin to form, they follow a path in
this diagram, with the indicated transition probabilities. The upper
path is followed most frequently; on this path
there is a unique complex component that will become
the ``giant.'' Parenthesized ratios
are the probabilities of reaching a given configuration.
At the moment the excess first reaches~2, the configuration must
either be $[0,1]$ (one tricyclic component) or~$[2]$ (two bicyclic
components). When the excess goes from~2 to~3, we go either from
$[0,1]$ to $[0,0,1]$ or $[1,1]$, or from~$[2]$ to $[0,0,1]$, [1,1],
or~$[3]$; and so on. Each configuration $[r_1,r_2,\ldots,r_q]$
corresponds to a partition of the excess $r=r_1+2r_2+\cdots +qr_q$.
The fraction in parenthesis shown above each configuration in Figure~1
is the limiting probability $c_1^{r_1}c_2^{r_2}\ldots
c_q^{r_q}\!/(r_1!\,r_2!\,\ldots\, r_q!\,e_r)$ that a random graph of
excess~$r$ has configuration $[r_1,r_2,\ldots,r_q]$. This is the
limiting probability that the infinite path traced out in the
infinite extension of Figure~1 will pass through
$[r_1,r_2,\ldots,r_q]$ during the evolution of a random graph or
multigraph on a large number of vertices.

A random graph almost always acquires nearly $\half n$~edges
before taking the first step from~[0] to~[1] in Figure~1. Indeed, the
uniform estimate (13.17), with $\mu=-n^{1/21}$, implies that the
probability of excess~$r$ when $m=\half n\exp(-n^{-2/7})$ is of
order~$n^{-r/7}$.

The fractions shown on arcs leading between configurations are
transition probabilities, namely the limiting probabilities that a
random graph of configuration $[r_1,r_2,\ldots,r_q]$ will go to another
specified configuration when its excess next changes. For example, a~random
graph in configuration~[2],
having two bicyclic components and no other complex components,
will proceed next to configuration~[1,1] with
probability~${144\over 221}$. These transition probabilities have a
fairly simple characterization:

\proclaim
Theorem 9. Let $r_1+2r_2+\cdots +qr_q=r$ and
$\delta_1+2\delta_2+3\delta_3+\cdots =1$. The asymptotic probability
that a random graph or multigraph of configuration
$[r_1,r_2,\ldots,r_q]$, having no acyclic components,  will change to
configuration
$[r_1+\delta_1,r_2+\delta_2,\ldots,r_q+\delta_q,\delta_{q+1},\ldots\,]$
when a random edge is added, can be computed as follows:
$$\vcenter{\halign{$\hfil#\hfil$\qquad&$\hfil#\hfil$\cr
\hbox{\rm Nonzero $\delta$'s}&{\rm Probability}\cr
\noalign{\medskip}
\delta_1=1&{5\over 4}/(3r+\half )(3r+{5\over 2})\cr
\noalign{\smallskip}
\delta_j=-1,\;\delta_{j+1}=1&9j(j+1)r_j/(3r+\half )(3r+{5\over 2})\cr
\noalign{\smallskip}
\delta_j=-2,\;\delta_{2j+1}=1&9j^2r_j(r_j-1)/(3r+\half )(3r+{5\over 2})\cr
\noalign{\smallskip}
\delta_j=-1,\;\delta_k=-1,\;\delta_{j+k+1}=1,\;j<k%
&18j\,k\,r_jr_k/(3r+\half )(3r+{5\over 2})\cr}}$$
In all other cases, the probability is~0. The estimates are correct to
within $O(n^{-1/2})$ when there are $n$~vertices.

\proof
As usual, it is easiest to consider first the uniform multigraph
process. We know that the generating function for
the cyclic multigraphs under consideration is
$$S(z)=e^{V(z)}\;{C_1(z)^{r_1}\over r_1!}\;{C_2(z)^{r_2}\over
r_2!}\;\ldots\; {C_q(z)^{r_q}\over r_q!}\;;\eqno(16.1)$$ 
the number of such multigraphs, weighted as usual by the compensation
factor~(1.1), is~$[z^n]\,S(z)$. We also know from (3.4) that
$V(z)=-\half \ln\bigl(1-T(z)\bigr)$, hence
$$e^{V(z)}={1\over \bigl(1-T(z)\bigr)^{1/2}}\,.$$

We observed in section~4 that the operator $\vartheta=z{d\over dz}$
corresponds to ``marking'' or singling out a
particular vertex. The function $\vartheta^2S(z)$ can therefore be
regarded as the generating function for multigraphs of configuration
$[r_1,r_2,\ldots,r_q]$ together with an ordered pair of marked
vertices $\langle x,y\rangle$. When $S(z)$ is a product $A(z)B(z)$,
the familiar relation
$$\vartheta^2\bigl(A(z)B(z)\bigr)
=\bigl(\vartheta^2A(z)\bigr)B(z)+2\bigr(\vartheta A(z)\bigr)
\bigl(\vartheta
B(z)\bigr)+A(z)\bigl(\vartheta^2B(z)\bigr)\eqno(16.2)$$
has a natural combinatorial interpretation: The product $A(z)B(z)$
stands for ordered pairs of graphs, generated respectively by $A(z)$
and $B(z)$, with no edges between them; the first term
$\bigl(\vartheta^2A(z)\bigr)B(z)$ of (16.2) corresponds to cases when
both of the marked vertices $\langle x,y\rangle$ are in the graph
generated by~$A(z)$; the last term corresponds to cases when both~$x$
and~$y$ belong to the $B(z)$ graph. The middle term $2\bigl(\vartheta
A(z)\bigr)\bigl(\vartheta B(z)\bigr)$ corresponds to the cases where
$x$ is in~$A$ and~$y$ is in~$B$ or vice versa.

We can use this idea in connection with (16.1) to understand what
happens when the graph gains a new edge.
The coefficient of~$z^n$ in~$\vartheta^2 S(z)$ represents all possibilities
$\langle x,y\rangle$; we can divide this
into cases by writing
$$\vartheta^2S(z)=S(z)\biggl(\sum_{0\leq j\leq
q}\,{\vartheta^2f_j(z)\over f_j(z)}+2\sum_{0\leq j<k<q}\,
{\vartheta f_j(z)\over f_j(k)}\,{\vartheta f_k(z)\over
f_k(z)}\biggr)\eqno(16.3)$$
where $f_0(z)=e^{V(z)}$ and $f_j(z)=C_j(z)^{r_j}\!/r_j!$ for $j\geq 1$.
A~term like $S(z)\bigl(\vartheta^2 f_j(z)\bigr)/f_j(z)$, say, then
corresponds to cases where $x$ and~$y$ both belong to $(j+1)$-cyclic
components. 

Each of the factors $f_j(z)$ is a linear combination of powers of the
quantity $\xi=1+\zeta=1/(1-T(z))$. For example, $f_0(z)=\xi^{1/2}$ and
$f_1(z)={5\over 24}\xi^3-{7\over 24}\xi^2+{1\over 12}\xi$,
according to (3.4) and (11.3). Hence it is easy to
compute $\vartheta f_j$ and~$\vartheta^2 f_j$, using rule (4.5):
$$\eqalignno{%
\vartheta(\xi^\alpha)&=\alpha\xi^{\alpha+2}-\alpha\xi^{\alpha+1}\,;\cr
\noalign{\smallskip}
\vartheta^2(\xi^{\alpha})
&={\alpha(\alpha+2)\xi^{\alpha+4}}-{\alpha(2\alpha+3)
\xi^{\alpha+3}}+{\alpha(\alpha+1)\xi^{\alpha+2}}\,.&(16.4)\cr}$$
The overall function $S(z)$ has the form $\xi^{3r+1/2}P(\xi^{-1})$ for
some polynomial~$P$, with $P(0)\neq 0$; hence the coefficient
$[z^n]\,S(z)$ is $t_n(3r+\half )P(0)\bigl(1+O(n^{-1/2})\bigr)/n!$ by
(3.8) and~(3.9).
It follows from (16.4) that $\vartheta^2S(z)=\xi^{3r+9/2}Q(\xi^{-1})$
for some polynomial~$Q$,
where $Q(0)=(3r+\half )(3r+{5\over 2})P(0)$. Hence
$$n^2={[z^n]\,\vartheta^2S(z)\over [z^n]\,S(z)}={\textstyle
(3r+\half )(3r+{5\over 2})}\,{t_n(3r+{9\over 2})\over t_n(3r+{1\over
2})}\,
\bigl(1+O(n^{-1/2})\bigr)\,.\eqno(16.5)$$
The transition probabilities we wish to compute are the fractions of
$(3r+\half )(3r+{5\over 2})$ that occur when
$\vartheta^2$~operates on individual factors of~$S(z)$.

For example, consider first the term
$S(z)\bigl(\vartheta^2f_0(z)\bigr)/f_0(z)$ of (16.3). This corresponds
to the case where both~$x$ and~$y$ belong to a cyclic component
(possibly the same one), thereby creating a new bicyclic component;
thus it corresponds to having $\delta_1=1$ and all other $\delta_j=0$.
 In this case
$[z^n]\,S(z)\bigl(\vartheta^2f_0(z)\bigr)/f_0(z)\sim {1\over
2}\cdot{5\over 2}\,t_n(3r+{9\over 2})P(0)/n!$, and the latter is asymptotically
${5\over 4}/(3r+\half )(3r+{5\over 2})$ of the total
$[z^n]\,\vartheta^2S(z)$. 

The term $2S(z)\bigl(\vartheta f_0(z)\bigr)\bigl(\vartheta
f_j(z)\bigr)/f_0(z)f_j(z)$, similarly, gives the probability that a
vertex from a cyclic component joins with a $(j+1)$-cyclic component;
this occurs with probability $2(\half )(3j\,r_j)/(3r+{1\over
2})(3r+{5\over 2})$. The net effect on components corresponds to
$\delta_j=-1$, $\delta_{j+1}=+1$.

There is also another way to get $\delta_j=-1$ and $\delta_{j+1}=+1$,
namely if both~$x$ and~$y$ belong to the same $(j+1)$-cyclic
component. The probability of this case works out to be
$(3j)(3j+2)r_j/(3r+{1\over
2})(3r+{5\over 2})$; hence the total transition probability for
$\delta_j=-1$ and $\delta_{j+1}=+1$ is $9j(j+1)r_j/(3r+{1\over
2})(3r+{5\over 2})$ as stated in the theorem. 

Notice that
$$\vartheta^2C_j^{r_j}=r_jC_j^{r_j-1}(\vartheta^2C_j)+r_j(r_j-1)C_j^{r_j-2}
(\vartheta C_j)^2\,.\eqno(16.6)$$
We have just taken care of the first term; the second term corresponds
to vertices~$x$ and~$y$ in distinct~$C_j$'s, when the new edge makes
$\delta_j=-2$ and $\delta_{2j+1}=+1$. The probability is
$9j^2r_j(r_j-1)/(3r+\half )(3r+{5\over 2})$. 

Finally, the term $2S(z)\bigl(\vartheta f_j(z)\bigr)(\vartheta
f_k(z)\bigr)/f_j(z)f_k(z)$ of (16.3) represents a case that
occurs with probability
$2(3j\,r_j)(3k\,r_k)/(3r+\half )(3r+{5\over 2})$ and corresponds
to $\delta_j=\delta_k=-1$, $\delta_{j+k+1}=+1$.

If we are working with the graph process instead of the multigraph
process, we must use $\widehat{C}_j(z)$ instead of $C_j(z)$; but
$f_0(z)$ is still essentially of degree $-1/2$ in~$\xi^{-1}$, and
$f_j(z)$ is still of degree~$-3j$, so the asymptotic calculations work
out as before.

However, in a
 random graph we must use the operator $\half (\vartheta^2_z
-\vartheta_z)-\vartheta_w$ instead of $\vartheta^2_z$, and we must
work with bivariate generating functions, as discussed in section~6.
The bgf corresponding to (16.1) is
almost univariate, however:
$$\hat{F}(w,z)=w^r\,e^{\hat{V}(wz)}\,{C_1(wz)^{r_1}\over
r_1!}\,{C_2(wz)^{r_2}\over r_2!}\,\cdots\,{C_q(wz)^{r_q}\over
r_q!}\,.$$
It is not difficult to see that the effect of $\vartheta_z^2$ swamps
the effects of~$\vartheta_z$ and~$\vartheta_w$, asymptotically, so
the multigraph analysis carries through.\quad\pfbox

\medskip
One amusing consequence of Theorem 9 is that we can use it to discover
and prove formula~(7.2) for
 the numbers~$e_r$ in a completely different way. The probability of
reaching the configuration~$[r]$, consisting of $r$~bicyclic
components and none of higher cyclic order, is $c_1^r/(r!\,e_r)$. The
only way to reach this configuration, when $r>0$, is from
$[r-1]$, and the transition probability is
$${{5\over 4}\over
 (3r-{5\over 2})(3r-\half )}={c_1^r\!/(r!\,e_r)\over
c_1^{r-1}/\bigl((r-1)!\,e_{r-1}\bigr)}\,.$$
Since $c_1=5/24$, we have $e_r=(6r-5)(6r-1)e_{r-1}/24r$, and (7.2)
follows by induction. This indirect method is probably the simplest
way to deduce the fact that Wright's constant is~$1/(2\pi)$.

\bigbreak\noindent
{\bf 17. A near-Markov process.}\enspace
We proved in Theorem~9 that the transition probabilities shown in
Figure~1 are the limiting probabilities, averaged over all multigraphs,
that a multigraph reaching a particular state will take a particular
step as its excess increases. But we did not prove that those
transition probabilities are independent of past history. For all we
know, the path taken to a particular configuration during the
evolution of a random graph might strongly influence the probability
distribution of its next leap forward. The next theorem addresses this
question. 

\proclaim Theorem 10. For any fixed $R$, an evolving random graph or
multigraph almost surely carries out a random walk in the first 
$R$~levels of the partial ordering shown in Figure~1, with transition
probabilities that approach the limiting values derived in Theorem~9.

\proof
As in previous proofs, it suffices to consider random multigraphs.
We will show that the transition probabilities have the asymptotic
behavior of Theorem~9 for all random multigraphs that remain
clean---i.e., for all multigraphs that reduce, under the pruning and
cancelling algorithms of Section~9, to 3-regular multigraphs~$\MM$
having $2r$~vertices and $3r$~edges, when the excess is $r\leq R$. We
know from Theorem~7 that the multigraph will be clean with probability
$1-O\bigl((1+\mu^4)n^{-1/3}\bigr)$;
 and we know from (13.17) that the probability of
excess~$r$ becomes superpolynomially small as the number of edges
passes~${n\over 2}$. So the excess almost surely increases past any
given value before a large multigraph becomes unclean. For example, if
$\mu\rightarrow\infty$ with $\mu=o(n^{1/12})$, the probability of
excess $\leq R$ approaches zero while the probability of remaining
clean is $1-o(1)$.

The proof for clean multigraphs 
is not as trivial as might be expected: Multigraphs that
follow a given path to $[r_1,r_2,\ldots,r_q]$ in the partial ordering
are {\it not\/} uniformly distributed, among all multigraphs whose
complex parts are enumerated by the generating function
$$e^V(C_1^{r_1}/r_1!)(C_2^{r_2}/r_2!) \ldots(C_q^{r_q}/r_q!)$$ 
assumed
in the proof of Theorem~9. Past history does affect the frequency of
certain types of components. For example, a~tricyclic component that
prunes and cancels to~$K_{3,3}$ cannot evolve along the path
$[1]\ra[2]\ra[0,0,1]$; removing any edge of~$K_{3,3}$ leaves a
connected graph.

Let's try to clarify the situation by working an example. Consider
the reduced multigraph
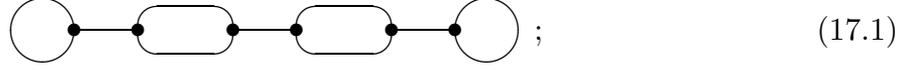
$$\advance\abovedisplayshortskip-\baselineskip \unitlength=10pt
\beginpicture(16,2)(0,.7)
\put(2,1){\disk{.4}}
\put(1,1){\circle2}
\put(2,1){\line(1,0)2}
\put(4,1){\disk{.4}}
\put(5.5,1){\oval(3,1.5)}
%\put(4,1){\line(1,0)3}
\put(7,1){\disk{.4}}
\put(7,1){\line(1,0)2}
\put(9,1){\disk{.4}}
\put(10.5,1){\oval(3,1.5)}
%\put(9,1){\line(1,0)3}
\put(12,1){\disk{.4}}
\put(12,1){\line(1,0)2}
\put(14,1){\disk{.4}}
\put(15,1){\circle2}
\endpicture
\ \,;
\eqno(17.1)$$
suppose we wish to compute the transition probabilities for
multigraphs of excess~3 that prune and cancel to (17.1) after
following the path $[1]\ra[0,1]\ra[0,0,1]$. The generating function
for all such multigraphs, assuming that there are no acyclic components,
would be ${1\over 32}\,e^VT^6\!/(1-T)^9$, if we did not specify the past
history $[1]\ra [0,1]\ra[0,0,1]$; but it turns out to be only ${8\over
9}$ as much when we prescribe the history. The reason is, intuitively,
that (17.1) has 9~edges, and a multigraph with history
$[1]\ra[0,1]\ra[0,0,1]$ can reduce to it only if the ``middle''
edge is not the last to be completed. The latter event happens with
probability~${8\over 9}$. 

A formal proof of the ${8\over 9}$ phenomenon can be given as follows.
The generating function $e^VT^6(1-T)^9$ expands to
$e^VT^6\sum_{n_1,n_2,\ldots,n_9\geq 0}T^{n_1}T^{n_2}\ldots T^{n_9}$;
 the individual terms represent the insertion of $\langle
n_1,\ldots,n_9\rangle$ vertices into the nine edges of (17.1), after
which a tree is sprouted at each vertex. The resulting multigraph will
have $n$~vertices and $m=n+3$ edges; there will be $6+n_1+\cdots +n_9$
root vertices and $9+n_1+\cdots +n_9$ ``critical'' edges on paths
between root vertices. Suppose we color each critical edge with one of
9~colors, corresponding to the original edge of (17.1) from which it
was subdivided. Then among the $m!$~permutations of edges that could
generate any such multigraph, exactly ${8\over 9}$ have the property
that the last critical edge has some color besides the ``middle''
color. (This follows by symmetry between $n_1,n_2,\ldots,n_9$.) Such
permutations are precisely those for which the history will be $[1]\ra
 [0,1]\ra [0,0,1]$; hence we obtain (17.1) with exactly ${8\over 9}$
times its overall probability, given that history.

It turns out that there are 17 unlabeled clean, connected, reduced multigraphs
of excess~3; and exactly 6 of them occur with weight~${8\over 9}$ when
the past history is $[1]\ra[0,1]\ra[0,0,1]$. Those 6 occur with
weight~${1\over 9}$ when the past history is
$[1]\ra[2]\ra[0,0,1]$, and the other 11 do not occur at all in that
case.

In general, given any $\MM$ that can arise for a given past history,
there will be a fraction $\beta >0$ such that each multigraph reducing
to~$\MM$ arises $\beta$~times as often with the given history as it
does overall. The reason is a slight generalization of the method by
which we proved the ${8\over 9}$~phenomenon: Each permutation of
colors of critical edges is equally likely to be the sequence of last
appearances in a random permutation of $n_1+n_2+\cdots+n_{3r}$
critical edges, and such permutations determine the past history. The
generating function for~$\MM$ will then be a constant multiple of
$e^V\bigl(T^{2r_1}/(1-T)^{3r_1}\bigr)\bigl(T^{2r_2}/(1-T)^{3r_2}\bigr)\ldots
\allowbreak
{\bigl(T^{2r_q}/(1-T)^{3r_q}\bigr)}$. Hence the asymptotic transition
probabilities will be the same for every feasible~$\MM$, exactly as
calculated in Theorem~9.\quad\pfbox

\bigbreak\noindent
{\bf 18. An emerging giant.}\enspace
The classic papers of Erd{\H o}s and R\'enyi [\ERo, \ER] tell us that
an evolving graph almost surely develops a single giant component,
which eventually is surrounded by only a few trees and later by only
isolated vertices, until the entire graph becomes connected. Thus
there will be a time when the graph reaches some configuration
$[0,0,\ldots,0,1]$ on the top line of Figure~1 and stays on that top
line ever afterward.

Indeed, the most 
probable path in Figure~1 is the one that goes directly from [1] to
$[0,1]$ to $[0,0,1]$ and so~on, never leaving the top line. The first
transition probability is~${72\over 77}$, the next is~${216\over
221}$, and subsequent steps are ever more likely to stay in line. In
such cases we can see the ``seed'' around which the giant component is
forming, before that component has become in any way gigantic. (The
complex components of any given finite excess almost always have only
$O(n^{2/3})$ vertices, a~vanishingly small percentage of the total;
each step at the beginning of
 Figure~1 occurs after adding about $n^{2/3}$ more edges.)

If we assume that the transition probabilities in Figure~1 are exact,
the overall probability that an evolving graph adheres strictly to the
top line---never having more than one complex component throughout its
entire evolution---is
$$\prod_{r=1}^{\infty}\,{r(r+1)\over (r+{1\over 6})(r+{5\over
6})}={\Gamma({7\over 6})\,\Gamma({11\over 6})\over\Gamma(1)\,\Gamma(2)}
={5\over 36}\,\Gamma\!\left({1\over 6}\right)\Gamma\!\left({5\over
6}\right)={5\pi\over 18}\,.\eqno(18.1)$$
Numerically, this comes to 0.8726646, roughly 7 times out of every~8.

Is ${5\pi\over 18}$ the true limiting probability that an evolving
graph or multigraph
never acquires two simultaneous components of positive excess,
throughout its evolution? We can at least prove that ${5\pi\over 18}$
is an upper bound. For we know from Theorem~10 that an evolving graph
will hug the top line of Figure~1 for at least $R$~steps with
probability
$$\prod_{r=1}^R\,{r(r+1)\over (r+{1\over 6})(r+{5\over 6})}={5\pi\over
18} +O(R^{-1})+O(n^{-1/3})\eqno(18.2)$$
for any fixed $R$, as $n\ra\infty$.

It is natural to conjecture that ${5\pi\over 18}$ is also a lower
bound, because a large component tends to propagate itself as soon as
it becomes large enough. Still, it is conceivable that a random graph
might have a tendency to leave the top line briefly when it first
becomes unclean. The transition probability for remaining on the top
line becomes strictly less than $r(r+1)/(r+{1\over 6})(r+{5\over 6})$
when the graph has a positive deficiency.
For example, suppose the initial bicyclic component is already
unclean; it will then correspond to the double self-loop of (9.15). We
know from (13.16) that this case arises with probability
$O(n^{-1/3})$. But if it does occur, the generating function for the
complex part will be a constant multiple of $T/(1-T)^2$ instead of
$T^2\!/(1-T)^3$, so the proof technique of Theorem~9 will yield a
transition probability from~[1] to $[0,1]$ of only~${8\over 9}$
instead of~${72\over 77}$. In general, when the deficiency is~$d$, the
asymptotic transition probability drops to 
$$\textstyle{(r-{d\over 3})(r-{d\over
3}+1)/(r-{d\over 3}+{1\over 6})(r-{d\over 3}+{5\over 6})\,.}$$
 This
probability estimate is, moreover, valid only when the excess is
reasonably small as a function of~$n$; 
otherwise the trees that sprout from the pruned
multigraph~$\ovline{M}$ will not be large enough to assert their
asymptotic behavior. 

\bigbreak\noindent
{\bf 19. A monotonicity property.}\enspace
During the time when an evolving graph or multigraph
 stays clean, we can show that the asymptotic
top-line transition probabilities $r(r+1)/\allowbreak
{(r+{1\over 6})}{(r+{5\over
6})}$ are in fact {\it lower\/} bounds for the correct
(non-asymptotic) probabilities. More precisely, the proof of
Theorem~9 shows that the true transition probability is a ratio of
expressions involving the tree polynomials $t_n(y)$, when there are
$n$~vertices in the cyclic part of the multigraph. We will prove that
this ratio decreases monotonically to $r(r+1)/\allowbreak
(r+{1\over 6})(r+{5\over 6})$ as $n\ra\infty$.

First we need to prove an auxiliary result about tree polynomials that
is interesting in its own right. Let us generalize the definition
of $t_n(y)$ in (3.8) by introducing a new parameter $m\geq 0$:
$${T(z)^m\over
\bigl(1-T(z)\bigr)^y}=\sum_{n=0}^{\infty}t_{m,n}(y)\,{z^n\over
n!}\,.\eqno(19.1)$$ 
Thus
$$t_{m,n}(y)=\sum_{j=0}^m{m\choose j}(-1)^jt_n(y-j)\eqno(19.2)$$
is the $m$\/th backward difference of $t_n(y)$.

\proclaim Lemma 6. Let $m$ be a nonnegative integer. For any fixed
integer $n>m$ and arbitrary real $y>0$, the ratio
$t_{m,n+1}(y)/t_{m,n}(y)$ is an increasing function of~$y$.
Equivalently, for fixed $y>0$ and any integer $n>m$, the ratio
$t'_{m,n}(y)/t_{m,n}(y)$ is an increasing function of~$n$.

\proof
The two statements of the lemma are clearly equivalent, because
$t_{m,n}(y)$ is positive when $y>0$ and $n>m$. 

Equation (2.12) of [\KP] states that
$$t_n(y)=n^{n-1}\sum_{k\geq 0}\,{y^{\overline{k+1}}\over k!}\;
{(n-1)^{\underline{k}}\over n^k}\,,\eqno(19.3)$$
where $x^{\overline{k}}$ means $x(x+1)\ldots(x+k-1)$ and
$x^{\underline{k}}$ means $x(x-1)\ldots(x-k+1)$. Therefore, by (19.2),
$$\eqalignno{t_{m,n}(y)&=n^{n-1}\sum_{k\geq
m-1}\,(k+1)\,{y^{\overline{k+1-m}}\over
(k+1-m)!}\;{(n-1)^{\underline{k}}\over n^k}\cr
\noalign{\smallskip}
&=n^{n-m}\sum_{k=0}^{n-m}\,(k+m)\,{y^{\overline{k}}\over
k!}\;{(n-1)^{\underline{k+m-1}}\over n^k}\,.&(19.4)\cr}$$
It follows that the inequality
$t'_{m,n}(y)/t_{m,n}(y)<t'_{m,n+1}(y)/t_{m,n+1}(y)$ is equivalent to
$${\sum_{k=0}^N\,a_k\alpha_k\over\sum_{k=0}^N\,b_k\alpha_k} >
{\sum_{k=0}^N\,a_k\beta_k\over\sum_{k=0}^N\,b_k\beta_k}\,,\eqno(19.5)$$
where $N=n+1-m$ and
$$\eqalignno{a_k=(k+m)\,{y^{\overline{k}}\over k!}\,,\qquad&
b_k=(k+m)\,{d\over dy}\;{y^{\overline{k}}\over k!}\,,\cr
\noalign{\smallskip}
\alpha_k=(n-1)^{\underline{k+m-1}}\,n^{n-m-k}\,,\qquad&\beta_k=
n^{\underline{k+m-1}}\,(n+1)^{n+1-m-k}\,.&(19.6)\cr}$$

The following condition is sufficient to prove (19.5), assuming
positive denominators:
$${a_0\over b_0}>{a_1\over b_1}>\cdots >{a_N\over
b_N}\qquad\hbox{and}\qquad
{\alpha_0\over\beta_0}>{\alpha_1\over\beta_1}
>\cdots>{\alpha_N\over\beta_N}\,.\eqno(19.7)$$
 For we have
$$\eqalignno{\sum_{k=0}^Nb_k\beta_k
&\sum_{j=0}^Na_j\alpha_j-\sum_{j=0}^Nb_j\alpha_j\sum_{k=0}^Na_k\beta_k\cr
\noalign{\smallskip}
&=\sum_{0\leq j<k\leq
n}(b_ka_j-b_ja_k)(\beta_k\alpha_j-\beta_j\alpha_k) >0\,.&(19.8)\cr}$$
(Historical note: Inequality (19.5) under condition (19.7) goes back at
least to Seitz in 1936 [\Sei]; see [\Mit, Section~2.5, Theorem~4],
where a supplementary condition is needed: The product of the
denominators must be positive. In linearly ordered discrete probability
space, the inequality is equivalent to saying that
$E\bigl(f(X)g(X)\bigr)\geq E\bigl(f(X)\bigr)E\bigl(g(X)\bigr)$
whenever $f$ and~$g$ are increasing functions of the random
variable~$X$. This inequality is, in turn, a~special case of the
celebrated FKG inequality [\FKG], which applies to certain partially
ordered probability spaces. The equality in (19.8), which reduces to
Lagrange's identity when we set $a_k=\alpha_k$ and $b_k=\beta_k$, is
the Binet-Cauchy identity for det~$AB$ when $A$ is a matrix of size
$2\times n$ and $B$ is $n\times 2$.)

And (19.7) is not difficult to verify, under the substitutions (19.6).
We have
$$\eqalign{{b_{k+1}\over a_{k+1}}&={1\over y}+{1\over y+1}+\cdots
+{1\over y+k}={b_k\over a_k}+{1\over y+k}\,;\cr
\noalign{\medskip}
{\alpha_{k+1}\over\alpha_k}&={n-k-m\over n}<{n-k-m+1\over
n+1}={\beta_{k+1}\over\beta_k}\,.\cr}$$
(When $m=0$ we omit the terms for $k=0$.)\quad\pfbox

\medskip
Assume now that the cyclic part of a random multigraph contains
$n$~vertices. The ``top line'' transition probability from a single
clean component of excess~$r$ to a single component of excess $r+1$ is
$1-p_{nr}$, where $p_{nr}$ is the probability that a new bicyclic
component will be formed. By the argument of Theorem~9,
$$p_{nr}={[z^n]\,\bigl(\vartheta^2 V(z)\bigr)S(z)\over [z^n]\,
\vartheta^2\bigl(V(z)S(z)\bigr)}\,,\eqno(19.9)$$
where $V(z)=1/\bigl(1-T(z)\bigr)^{1/2}$ is the generating function for
unicyclic components and $S(z)=T(z)^{2r}/\bigl(1-T(z)\bigr)^{3r}$ is a
prototypical generating function for clean components of excess~$r$.
We want to show that $p_{nr}$ is an increasing function of~$n$, since
we want $1-p_{nr}$ to be decreasing.

Let's work on a simpler problem first, showing that
$$q_{nr}={[z^n]\,\bigl(\vartheta\,A(z)\bigr)S(z)\over
 [z^n]\,\vartheta\bigl(A(z)S(z)\bigr)}\eqno(19.10)$$ 
is an increasing function of $n$ whenever
$$A(z)={T(z)^a\over\bigl(1-T(z)\bigr)^b}\,,\qquad b>{3\over
2}\,a\,.\eqno(19.11)$$ 
Here $a$ is a nonnegative integer; we will assume that $n\geq 2r+a$,
so that the denominator of (19.10) is nonzero. We have
$$\eqalign{\vartheta A(z)&={b\,T(z)^a\over
\bigl(1-T(z)\bigr)^{b+2}}-{(b-a)\,T(z)^a\over\bigl(1-T(z)\bigr)^{b+1}}\,,\cr
\noalign{\smallskip}
\vartheta\bigl(A(z)S(z)\bigr)
&={(3r+b)T(z)^{2r+a}\over\bigl(1-T(z)\bigr)^{3r+b+2}}-{(r+b-a)T(z)^{2r+a}\over
\bigl(1-T(z)\bigr)^{3r+b+1}}\,;\cr}$$
hence
$$\eqalign{q_{nr}&={b\,t_{2r+a,n}(3r+b+2)-(b-a)\,t_{2r+a,n}(3r+b+1)\over
(3r+b)\,t_{2r+a,n}(3r+b+2)-(r+b-a)\,t_{2r+a,n}(3r+b+1)}\cr
\noalign{\smallskip}
&={b\over 3r+b}\,\biggl(1-{r(2b-3a)/(3rb+b^2)\over
\left(\displaystyle{{t_{2r+a,n}(3r+b+2)\over t_{2r+a,n}(3r+b+1)}-
{r+b-a\over 3r+b}}\right)}\biggr)\,.\cr}$$
Since the coefficients of $t_{2r+a,n}(y)$ are nonnegative, we have
$$t_{2r+a,n}(3r+b+2)/t_{2r+a,n}(3r+b+1)\geq 1>(r+b-a)/(3r+b)\,.$$
 It
follows that $q_{nr}$ is increasing iff
$${t_{2r+a,n}(3r+b+2)\over
t_{2r+a,n}(3r+b+1)}<{t_{2r+a,n+1}(3r+b+2)\over
t_{2r+a,n+1}(3r+b+1)}\,.\eqno(19.12)$$
And (19.12) does hold, because $t_{2r+a,n+1}(y)/t_{2r+a,n}(y)$ is an
increasing function of~$y$ by Lemma~6.

Incidentally, this argument also shows that $q_{nr}$ is constant when
$b={3\over 2}a$ and decreasing when $0<b<{3\over 2}a$. 

Now to prove that $p_{nr}$ is increasing, we can write
$$\eqalign{p_{nr}&={[z^n]\,\bigl(\vartheta^2V(z)\bigr)S(z)\over 
 [z^n]\,\vartheta\bigl(\bigl(\vartheta
V(z)\bigr)S(z)\bigr)}\;{[z^n]\,\vartheta\bigl(\bigl(\vartheta
V(z)\bigr)S(z)\bigr)\over [z^n]\,\vartheta^2\bigl(V(z)S(z)\bigr)}\cr
\noalign{\medskip}
&={[z^n]\,\bigl(\vartheta^2V(z)\bigr)S(z)\over 
 [z^n]\,\vartheta\bigl(\bigl(\vartheta
V(z)\bigr)S(z)\bigr)}\;{[z^n]\,\bigl(\vartheta V(z)\bigr)S(z)\over
 [z^n]\,\vartheta\bigl(V(z)S(z)\bigr)}\,.\cr}$$
The first factor is of type $q_{nr}$ if we put $A(z)=\vartheta
V(z)=\half \,T(z)/\bigl(1-T(z)\bigr)^{5/2}$; here $a=1$,
$b={5\over 2}$, so $q_{nr}$ is increasing. The second factor is of
type $q_{nr}$ if we put $A(z)=V(z)$; here $a=0$, $b=\half $, and
again $q_{nr}$ is increasing. We have proved

\proclaim Theorem 11. The probability that a clean random multigraph
of excess $r>0$ will not acquire a new bicyclic component when its
excess next changes is strictly greater than the limiting value
$r(r+1)/(r+{1\over 6})(r+{5\over 6})$. \quad\pfbox

\bigskip
Theorem 11 gives further support to the ${5\pi\over 18}$ conjecture of
Section~18, because ${5\pi\over 18}$ was shown there to be an upper
bound. If the top-line transition probability were always strictly
greater than $r(r+1)/(r+{1\over 6})(r+{5\over 6})$, we could establish
${5\pi\over 18}$ as a lower bound. However, Theorem~11 does not prove
the conjecture, because the probability becomes smaller than
$r(r+1)/(r+{1\over 6})(r+{5\over 6})$ when a graph becomes unclean.

Incidentally, when the number of edges gets large, we may need
asymptotic formulas for $t_n(y)$ that are valid when $y$ goes to
infinity with~$n$. Formula (3.9) can be extended to
$$t_n(y)={\sqrt{2\pi}\,n^{n-1/2+y/2}\over
2^{y/2}\,\Gamma(y/2)}\,\bigl(1+O(y^{3/2}n^{-1/2})\bigr)\,,\eqno(19.13)$$
uniformly for $1\leq y\leq n^{1/3}$, using the proof technique of
Lemma~3. Still larger values of $y$ can be handled by using the saddle
point method to derive the following general estimate:
$$t_{a\mskip-1mu\lambda n,n}(\lambda n+b)={n!\,
e^{n\rho}\rho^{(a\mskip-1mu\lambda
-1)n}\lambda^{(1-b)/2}\over 2\,\sqrt{\pi n}\,(1-\rho)^{\lambda n}}\,
\bigl(1+O(\sqrt{\lambda}\,)+O(1/\sqrt{\lambda
n}\,)\bigr)\,,\eqno(19.14)$$
for fixed $a$ and $b$ as $\lambda\ra 0$ and $\lambda n/(\log
n)^2\ra\infty$, where
$$\rho=1+c\lambda-\sqrt{\lambda(1+c^2\lambda)}=1-\sqrt{\lambda}+c\lambda
-{c^2\over 2}\,\lambda^{3/2}+O(\lambda^{5/2})\,,\quad c={1-a\over
2}\,.\eqno(19.15)$$
For example, to estimate $t_{2r,n}(3r)$ when $r=n^{1/2}$, we can use
(19.14) with $a={2\over 3}$, $b=0$, and $\lambda=3n^{-1/2}$. The
complicated dependence on~$\rho$ can also be expressed as
$${e^{n\rho}\rho^{(a\lambda -1)n}\over (1-\rho)^{\lambda n}}
=\exp\bigl(n\bigl(\textstyle{1-\half \lambda \ln\lambda +{1\over
2}\lambda +({1\over 3}-a)\lambda^{3/2}-{1\over
4}a^2\lambda^2+O(\lambda^{5/2})}\bigr)\bigr)\,,\eqno(19.16)$$ 
which is sufficiently accurate if $\lambda\leq n^{-1/4}$.

\bigbreak\noindent
{\bf 20. The evolution of uncleanness.}\enspace
We get further insight into the behavior of an evolving multigraph by
studying how its reduced multigraph~$\MM$ changes as the excess
increases. Let's review the theory of Section~9 in light of what we
have learned since then. The generating function for the cyclic part
of all multigraphs having excess~$r$ and deficiency~$d$ is
$$E_{rd}(z)=e_{rd}\,{T(z)^{2r-d}\over
\bigl(1-T(z)\bigr)^{3r-d+1/2}}\,.\eqno(20.1)$$
We can interpret it as follows, ignoring the constant factor $e_{rd}$
for a moment: There is a reduced multigraph~$\MM$ having $\nu=2r-d$
vertices and $\mu=3r-d$ edges; each vertex has degree $\geq 3$, where
a self-loop is considered to add~2 to the degree. We can obtain all
cyclic multigraphs~$M$ that reduce to~$\MM$ by a two-step process. First we
insert~0 or more vertices of degree~2 on each edge; and we also
construct any desired number of cycles, as separate components. 
All of the newly
constructed vertices, including the vertices in the cycles, have
degree~2. This first step creates a set of multigraphs with the
univariate generating function $z^{\nu}(1-z)^{-\mu}(1-z)^{-1/2}$,
because each edge subdivision corresponds to $(1-z)^{-1}$, and because
the cycles are generated by $\exp(\half z+{1\over 4}z^2+{1\over
6}z^3+\cdots\,)=(1-z)^{-1/2}$. Now we proceed to step two, which
sprouts a rooted tree from every vertex; this changes~$z$ to $T(z)$ in
the generating function.

The excess increases by 1 when we add a new edge $\langle x,y\rangle$
to~$M$. How does the new edge change~$\MM$? A~moment's thought shows
that~$\MM$ will gain 2, 1, or~0 vertices; this means the deficiency
will either stay the same or it will increase by~1 or~2.

In fact there is a nice algebraic and quantitative way to understand
what happens, in terms of the generating function. Again we consider a
two-step process: First we choose a vertex~$x$ of~$M$; this means we
apply the marking operator~$\vartheta$ to the generating function.
There are three cases: The marked vertex either belongs to a tree
attached to one of the $\nu$~special vertices of~$\MM$, or it belongs
to a tree attached to a vertex within one of the $\mu$~edges, or it
belongs to a tree attached to a vertex in some cycle. We represent
Case~1 by attaching a ``half-edge'' to the existing vertex; we
represent Case~2 by introducing a new vertex into the split edge and
attaching a half-edge to~it; we represent Case~3 by introducing a new
vertex with a self-loop and attaching a half-edge to~it.

A half-edge is like an edge but it touches only one vertex.
For example, if $\MM$ is the multigraph~$K_4$, the symbolic
representations of the three possible outcomes of step~1 are
$$\unitlength=10pt
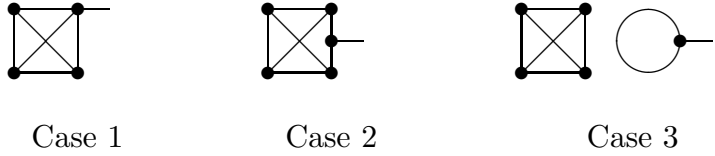
\beginpicture(4,4)(0,0)
\put(0,4){\disk{.4}}
\put(0,2){\disk{.4}}
\put(2,2){\disk{.4}}
\put(2,4){\disk{.4}}
\put(0,2){\line(1,0)2}
\put(0,4){\line(1,0)2}
\put(0,2){\line(0,1)2}
\put(2,2){\line(0,1)2}
\put(0,2){\line(1,1)2}
\put(2,2){\line(-1,1)2}
\put(2,4){\line(1,0)1}
\put(2,0){\makebox(0,0){Case 1}}
\endpicture
\qquad\qquad
\beginpicture(4,4)(0,0)
\put(0,4){\disk{.4}}
\put(0,2){\disk{.4}}
\put(2,2){\disk{.4}}
\put(2,4){\disk{.4}}
\put(0,2){\line(1,0)2}
\put(0,4){\line(1,0)2}
\put(0,2){\line(0,1)2}
\put(2,2){\line(0,1)2}
\put(0,2){\line(1,1)2}
\put(2,2){\line(-1,1)2}
\put(2,3){\disk{.4}}
\put(2,3){\line(1,0)1}
\put(2,0){\makebox(0,0){Case 2}}
\endpicture
\qquad\qquad
\beginpicture(7,4)(0,0)
\put(0,4){\disk{.4}}
\put(0,2){\disk{.4}}
\put(2,2){\disk{.4}}
\put(2,4){\disk{.4}}
\put(0,2){\line(1,0)2}
\put(0,4){\line(1,0)2}
\put(0,2){\line(0,1)2}
\put(2,2){\line(0,1)2}
\put(0,2){\line(1,1)2}
\put(2,2){\line(-1,1)2}
\put(5,3){\disk{.4}}
\put(4,3){\circle2}
\put(5,3){\line(1,0)1}
\put(3.5,0){\makebox(0,0){Case 3}}
\endpicture
$$
Let's call this augmented multigraph~$\MM{}'$. 

A cyclic multigraph $M'$ with a marked vertex can be reduced by
attaching a half-edge to the marked vertex, then pruning all vertices
of degree~1 and cancelling all vertices of degree~2. Conversely, the
marked cyclic multigraphs that reduce to a given~$\MM{}'$ are obtained
by adding zero or more vertices to each edge ({\it including\/} the
half edge), also adding cycles, then sprouting trees from each vertex.
Thus the generating function for~$M'$ in Case~1 is
$$e_{rd}\,{\nu\,T(z)^{\nu}\over
\bigl(1-T(z)\bigr)^{\mu+3/2}}\,;\eqno(20.2)$$ 
the $\nu$ in the numerator accounts for the number of vertices that
can be chosen, and the extra $\bigl(1-T(z)\bigr)$ in the denominator
accounts for the new half-edge. The generating function for~$M'$ in
Case~2 is
$$e_{rd}\,{\mu\,T(z)^{\nu+1}\over\bigl(1-T(z)\bigr)^{\mu+5/2}}\,;\eqno(20.3)$$
now we have $\mu$ edges that can be split, and we include an
additional $T(z)$ in the numerator for the new vertex and an
additional $\bigl(1-T(z)\bigr)^2$ in the denominator for the new
half-edge and the additional split edge. Finally, the generating
function for~$M'$ in Case~3 is
$$e_{rd}\,{{1\over
2}\,T(z)^{\nu+1}\over\bigl(1-T(z)\bigr)^{\mu+5/2}}\,;
\eqno(20.4)$$
as in Case 2, the diagram has gained one vertex and two edges. The
factor~$\half $ is due to the compensation factor~$\kappa$ of a
self-loop. 

If our calculations are correct, the sum $(20.2)+(20.3)+(20.4)$ should
be the result of applying~$\vartheta$ to the overall generating
function (20.1). And sure enough, 
$$\vartheta\,{T(z)^{\nu}\over\bigl(1-T(z)\bigr)^{\mu+1/2}}=
{\nu\,T(z)^{\nu}\over\bigl(1-T(z)\bigr)^{\mu+3/2}}+
{(\mu+\half )\,T(z)^{\nu+1}\over\bigl(1-T(z)\bigr)^{\mu+5/2}}\,;
\eqno(20.5)$$
 everything checks out fine.

The next step, choosing $y$, is the same, except that now we mark a
vertex of~$M'$ and obtain~$M''$. The transition from~$\MM{}'$ to~$\MM{}''$
again leads to three cases; we attach another half-edge and possibly
split an existing edge or add a new self-loop. In particular, we might
split the half-edge of~$\MM{}'$. The change in the generating function
is once again represented by (20.5), but this time $\nu$ and $\mu$
have to be adjusted to equal the number of vertices and edges
of~$\MM{}'$. The left term of (20.5) therefore becomes
$${\nu^2T(z)^{\nu}\over\bigl(1-T(z)\bigr)^{\mu+5/2}}+{\nu(\mu+{3\over
2})\,T(z)^{\nu+1}\over \bigl(1-T(z)\bigr)^{\mu+7/2}}\,,\eqno(20.6)$$
and the right term becomes
$${(\mu+\half )(\nu+1)\,T(z)^{\nu+1}\over\bigl(1-T(z)\bigr)^{\mu+7/2}}+
{(\mu+\half )(\mu+{5\over
2})\,T(z)^{\nu+2}\over\bigl(1-T(z)\bigr)^{\mu+9/2}}\,.\eqno(20.7)$$

Notice that the first term of (20.5) corresponds to the case that the
deficiency increases by~1 when $x$ is chosen, while the second term
corresponds to the case where the deficiency stays the same.
Similarly, the first terms of (20.6) and (20.7) correspond to an
increase in deficiency when $y$ is chosen, after $x$ has already been
marked.

By looking at the coefficients  of these generating functions we can
see why the deficiency rarely increases unless the total number of
vertices in the cyclic part is not much larger than~$\nu$. Suppose we
change the generating function to 
$$F(z,s)={T(z)^{\nu}\over\bigl(1-s\,T(z)\bigr)^{\mu+1/2}}\,;$$
then
$${[z^n]\;{\partial\over\partial s}\;F(z,s)\big\vert_{s=1}\over
 [z^n]\,F(z,s)\vert_{s=1}}$$
will be the average number of tree-root vertices that appear within
the edges of~$\MM$. For fixed~$\nu$ and~$\mu$ as $n\ra\infty$ this
number~is 
$${[z^n]\,(\mu+{1\over
2})\,T(z)^{\nu+1}\bigl(1-T(z)\bigr)^{-\mu-3/2}\over [z^n]\,T(z)^{\nu}
\bigl(1-T(z)\bigr)^{-\mu-1/2}}=
{(\mu+\half )t_n(\mu+{3\over 2})\over t_n(\mu+{1\over
2})}\,\bigl(1+O(n^{-1/2})\bigr)\,,\eqno(20.8)$$
which is approximately $\sqrt{\mathstrut\mu n}$ by (3.9), when $\mu$
is large.  Thus, there are about
$\sqrt{\mathstrut\mu n}$ tree roots, only $\nu$ of which will increase the
deficiency when chosen; almost all choices of~$x$ and~$y$ will fall in
trees that add new vertices to~$\MM$ and~$\MM{}'$.

If we replace one of the factors $T(z)$ in the numerator of the
generating function by $\vartheta T(z)=T(z)/\bigl(1-T(z)\bigr)$, we
multiply the coefficient of~$z^n$ by the average size of a rooted
tree; we find that each rooted tree contains about $\sqrt{n/\mu}$
vertices.

The number $n$ in these calculations has been the number of vertices
in the cyclic part of a multigraph, and the number~$\mu$ is~$3r$. Let's
return to our other notational convention, where $n$ is the total
number of vertices in the evolving multigraph and $m={n\over 2}(1+\mu
n^{-1/3})$ is the total number of edges. Recall that the average
excess~$r$ grows as~${2\over 3} \mu^3$, for $\mu\leq n^{1/12}$; the size of the
cyclic part, similarly, has order $\mu n^{2/3}$. The probability that
a random new edge falls in the cyclic part (and therefore increases
the excess) is therefore of order $(\mu n^{2/3}\!/n)^2=\mu^2n^{-2/3}$;
we must add about $n^{2/3}\!/\mu^2$ more edges before the excess
increases. And when it does, the probability of choosing a ``bad''~$x$
or~$y$, making the new multigraph unclean, is the ratio of~$2r$ to the {\it
total\/} number of tree roots, which is of order
$${2r\over \sqrt{3r(\mu n^{2/3})}}\approx{\mu^3\over\sqrt{\mu^4n^{2/3}}}
=\mu n^{-1/3}\,.$$
We will probably have to do $n^{1/3}\!/\mu$ augmentations of excess,
 adding $(n^{2/3}\!/\mu^2)(n^{1/3}\!/\mu)=n/\mu^3$ more edges, 
before we reach an unclean
multigraph. That is why the multigraph tends to stay clean until $\mu=n^{1/12}$,
as asserted in Theorem~7.

After $x$ and $y$ are chosen to form the endpoints of a new edge,
a~third step takes place: This new  edge is merged or integrated with
the other edges. 
Symbolically, the two half-edges for~$x$ and~$y$ are now spliced together.
We can complete our study of how the generating function changes at
the time of excess augmentation by considering this third and final
step.

It is easiest to consider the {\it inverse\/} of the final step,
namely the operation of marking an edge whose removal would {\it
decrease\/} the excess. Such an edge must be in the complex part, not
the acyclic or unicyclic
part. The operator that corresponds to marking an arbitrary
edge in a complex multigraph of excess~$r$ is $r+\vartheta$, because this
multiplies the coefficient of~$z^n$ by $r+n$, the total number of
edges. However, we also need to figure out the generating function for
``insignificant'' edges, edges whose removal would leave the excess
unchanged. Such edges can be described by an ordered pair consisting
of a rooted tree and a multigraph of excess~$r$ with a marked vertex; one
end of the edge is attached to the marked vertex and the other end is
attached to the root of the tree. Thus the appropriate operator for
insignificant edges is $T(z)\vartheta$. Altogether we find that the
generating function that corresponds to marking a significant edge, given
a family of complex multigraphs of excess~$r$, is $r+\vartheta
-T\vartheta$. We also should multiply this by two, because we assign
an orientation to the edge with the ordered pair $\langle x,y\rangle$.

When the operator $2\bigl(r+\vartheta-T(z)\vartheta\bigr)$ is applied
to a generating function of the form
$T(z)^{\nu}/\bigl(1-T(z)\bigr)^{\mu}$, 
with
$\mu=\nu+r$, we get
$$\eqalignno{2\bigl(r+\bigl(1-T(z)\bigr)\vartheta\bigr)\,{T(z)^{\nu}\over
\bigl(1-T(z)\bigr)^{\mu}}
&={2\bigl((r+\nu)T(z)^{\nu}\bigr)\over \bigl(1-T(z)\bigr)^{\mu}}+
{2\mu\,T(z)^{\nu+1}\over\bigl(1-T(z)\bigr)^{\mu+1}}\cr
\noalign{\smallskip}
&={2\mu\,T(z)^{\nu}\over \bigl(1-T(z)\bigr)^{\mu+1}}\,.&(20.9)\cr}$$
Therefore the inverse operation we seek, which merges an ordered
$\langle x,y\rangle$ into the set of existing edges, takes
$${T(z)^{\nu}\over \bigl(1-T(z)\bigr)^{\mu+1}}\;\longmapsto\;
{1\over 2\mu}\;{T(z)^{\nu}\over
\bigl(1-T(z)\bigr)^{\mu}}\,.\eqno(20.10)$$
For example, the first term of (20.6) will go into
$${\nu^2\,T(z)^{\nu}\over 2(\mu+1)\bigl(1-T(z)\bigr)^{\mu+3/2}}\,.$$
(First we multiply by $\bigl(1-T(z)\bigr)^{1/2}$ to get rid of the
unicyclic components, then we apply the inverse operation (20.10), then
we put the unicyclic components back.) 

Altogether we find that the generating function
$T(z)^{2r-d}/\bigl(1-T(z)\bigr)^{3r-d+1/2}$ for cyclic multigraphs of
excess~$r$ and deficiency~$d$ makes the following contributions to the
generating functions for cyclic multigraphs of excess $r+1$ and
deficiencies~$d$, $d+1$, and $d+2$, according to (20.6), (20.7), and
(20.10):
$$\eqalignno{&{(6r-2d+5)(6r-2d+1)\over 8(3r-d+3)}\;{T(z)^{2r+2-d}\over
\bigl(1-T(z)\bigr)^{3r+3-d+1/2}}\cr
\noalign{\smallskip}
&\qquad\null+{\bigl((2r-d)(6r-2d+3)+(2r-d+1)(6r-2d+1)\bigr)\over
4(3r-d+2)}\; {T(z)^{2r+1-d}\over\bigl(1-T(z)\bigr)^{3r+2-d+1/2}}\cr
\noalign{\smallskip}
&\qquad\null+{(2r-d)^2\over
2(3r-d+1)}\;{T(z)^{2r-d}\over\bigl(1-T(z)\bigr)^{3r+1-d+1/2}}\,.&(20.11)\cr}$$
This is essentially the same as the recurrence relation for $e_{rd}$
in (5.11)--(5.13).

We can illustrate the observations of this section by introducing
another partial ordering analogous to Figure~1. Every evolving graph
or multigraph traces a path in Figure~2, just as it does in Figure~1; but in
Figure~2 the state $(r,d)$ represents excess~$r$ and deficiency~$d$.
Fractions in brackets above each state are the coefficients~$e_{rd}$
of the generating function (5.10). Fractions on the arrows are not
transition probabilities but rather the amounts by which each
generating function coefficient affects the coefficients at the next
level; these fractions are the coefficients in (20.11). 

\vfill\eject
\unitlength=1cm
\centerline{\beginpicture(10,9)(3,-5)
\put(3,0){\bmath\makebox(0,0){$(0,0)$}}
\put(6,1){\bmath\makebox(0,0){$(1,0)$}}
\put(6,-1){\bmath\makebox(0,0){$(1,1)$}}
\put(9,3){\bmath\makebox(0,0){$(2,0)$}}
\put(9,1){\bmath\makebox(0,0){$(2,1)$}}
\put(9,-1){\bmath\makebox(0,0){$(2,2)$}}
\put(9,-3){\bmath\makebox(0,0){$(2,3)$}}
\put(12,4.8){\bmath\makebox(0,0){$(3,0)$}}
\put(12,3){\bmath\makebox(0,0){$(3,1)$}}
\put(12,1){\bmath\makebox(0,0){$(3,2)$}}
\put(12,-1){\bmath\makebox(0,0){$(3,3)$}}
\put(12,-3){\bmath\makebox(0,0){$(3,4)$}}
\put(12,-4.8){\bmath\makebox(0,0){$(3,5)$}}
\def\\#1/#2//{\makebox(0,0){\smash{\raise12pt\hbox{$[{#1\over#2}]$}}}}%
\put(3,0){\\1/1//}
\put(6,1){\\5/24//}
\put(6,-1){\\1/8//}
\put(9,3){\\385/1152//}
\put(9,1){\\35/64//}
\put(9,-1){\\91/384//}
\put(9,-3){\\1/48//}
\def\\#1/#2//{\makebox(0,0){$#1\over#2$}}%
\put(3.5,.233){\line(3,1){.7}}
\put(4.5,.5){\\5/24//}
\put(4.7,.633){\vector(3,1){.7}}
\put(3.5,-.233){\line(3,-1){.74}}
\put(4.5,-.6){\\1/8//}
\put(4.7,-.633){\vector(3,-1){.8}}
\put(6.6,1.5){\line(3,2){.5}}
\put(7.4,2.0){\\77/48//}
\put(7.7,2.234){\vector(3,2){.8}}
\put(6.5,1){\line(1,0){.3}}
\put(7.1,1.1){\\39/20//}
\put(7.4,1){\vector(1,0){.95}}
\put(6.35,0.7){\line(3,-2){.5}}
\put(7.05,0.3){\\1/2//}
\put(7.21,.06){\vector(3,-2){1.19}}
\put(6.4,-.75){\line(3,2){.4}}
\put(7,-.4){\\9/8//}
\put(7.21,-.21){\vector(3,2){1.29}}
\put(6.6,-1){\line(1,0){.5}}
\put(7.45,-1){\\17/16//}
\put(7.75,-1){\vector(1,0){.7}}
\put(6.4,-1.3){\line(3,-2){.55}}
\put(7.15,-1.75){\\1/6//}
\put(7.3,-1.9){\vector(3,-2){1.2}}
\put(9.6,3.5){\vector(3,2){1.8}}
\put(9.65,3){\vector(1,0){1.7}}
\put(9.4,2.6){\vector(3,-2){1.9}}
\put(9.6,1.5){\vector(3,2){1.9}}
\put(9.65,1){\vector(1,0){1.7}}
\put(9.5,0.6){\vector(3,-2){1.9}}
\put(9.6,-.5){\vector(3,2){1.9}}
\put(9.65,-1){\vector(1,0){1.7}}
\put(9.5,-1.4){\vector(3,-2){1.9}}
\put(9.6,-2.5){\vector(3,2){1.9}}
\put(9.65,-3){\vector(1,0){1.8}}
\put(9.5,-3.4){\vector(3,-2){1.9}}
\endpicture
}
\bigskip
{\narrower\narrower\smallskip\noindent
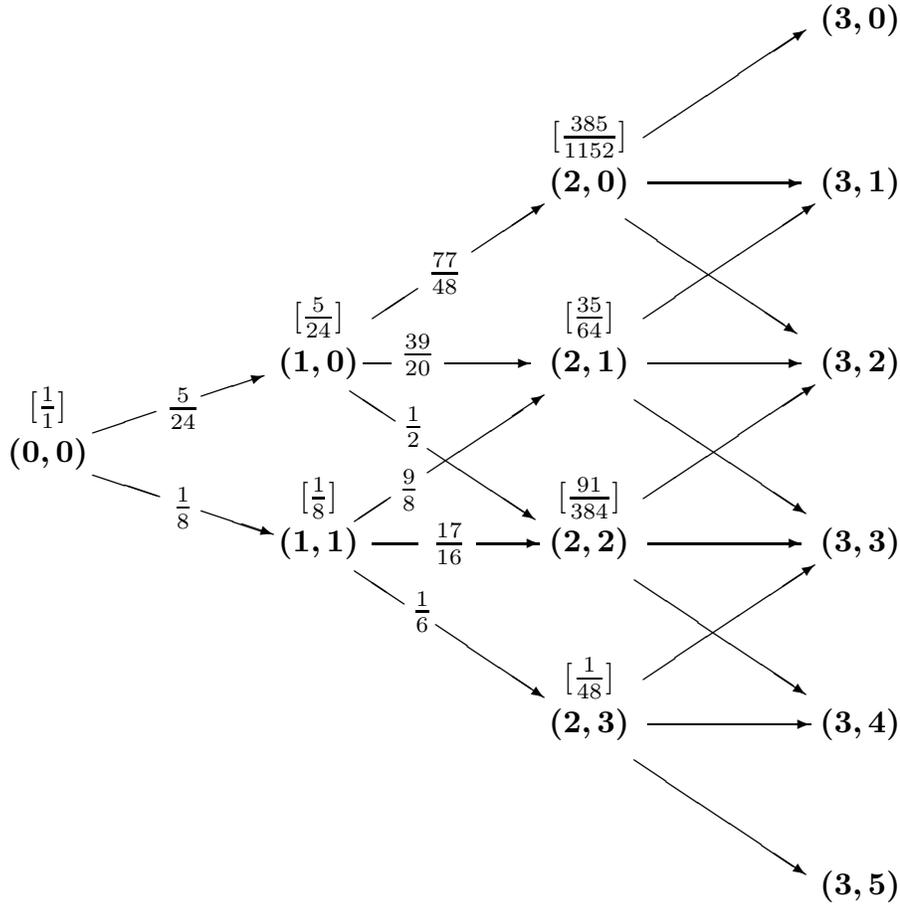
{\bf Figure 2.}\enspace The evolution of deficiency.
Each configuration $(r,d)$ stands for a graph or multigraph
whose complex part reduces to a multigraph with $2r-d$ vertices and
$3r-d$ edges, when vertices of degrees~1 and~2 are eliminated. A~graph
or multigraph
with deficiency~0 is called ``clean''; the reduced multigraphs in such
cases are 3-regular. When $r$ is small, each unit increase in
deficiency occurs with probability of order $n^{-1/3}$; therefore most
random graphs stay clean until $r$ is quite large.
\smallskip}

\bigbreak\noindent
{\bf 21. Waiting for uncleanness.}\enspace
We have seen that a graph almost surely stays clean while it has
$\half (n+\mu n^{2/3})$ edges, as long as $\mu$ is $o(n^{1/12})$.
What happens when $\mu$ gets a bit larger? Another contour integral
provides the answer; in this one, we rescale~$\mu$ in preparation for
the appearance of the giant component, but we allow~$\mu$ to be small
enough that there is  a substantial overlap with the estimate
(10.1) of Lemma~3.

\proclaim Lemma 7. If $m=\half (n+\mu n)$ and $r={2\over
3}\mu^3 n+\rho\sqrt{\mu^3n}$, we have
$$\eqalignno{&{2^mm!\,n!\,e_r\over n^{2m}(n-m+r)!}\;[z^n]\;
{U(z)^{n-m+r}T(z)^{2r}\over \bigl(1-T(z)\bigr)^{3r+y}}\cr
\noalign{\nobreak\smallskip}
&\qquad
=B(y,\mu,\rho,n)\exp\bigl(O\bigl((1+ \vert\rho\vert^3)\mu^{-3/2}n^{-1/2}+
(1+\vert\rho\vert)\mu^{5/2}n^{1/2}\bigr)\bigr)\,,&(21.1)\cr}$$
where
$$B(y,\mu,\rho,n)=\sqrt{{3\over 20\pi n}}\,\mu^{-1-y}\exp\left(-{2\over
3}\,\mu^4n- {3\over 20}\,\rho^2\right)\,,\eqno(21.2)$$
uniformly for $n^{-1/3}\log n\leq\mu\leq n^{-1/5}$,
$\vert\rho\vert\leq{2\over 3}
\mu^{3/2}n^{1/2}$, and fixed~$y$ as $n\ra\infty$.

\proof This is the sort of lemma for which computer algebra really pays
off. We can begin by using Stirling's approximation to show that
$$\eqalignno{&\log\left({2^mm!\,n!\, e_r\over
n^{2m}(n-m+r)!\,2^{n-m+r}}\right) 
=\textstyle{-n+3r\ln\mu-{5\over 6}\mu^3n}\cr
\noalign{\smallskip}
&\qquad\qquad\null-\textstyle{{3\over 2}\ln\mu+\half \ln{3\over 2}+{2\over
3}\mu^4n+{3\over 4}\rho^2}\cr
\noalign{\smallskip}
&\qquad\qquad\null+O\bigl((1+\vert\rho\vert^3)\mu^{-3/2}n^{-1/2}
+(1+\vert\rho\vert)\mu^{5/2}n^{1/2}\bigr)\,.\quad&(21.3)\cr}$$

Now we express the remaining factor by using the trick of (10.11):
$$[z^n]\;{(2U(z))^{n-m+r}T(z)^{2r}\over\bigl(1-T(z)\bigr)^{3r+y}}={1\over
2\pi i}\oint (1-z)^{1-y}e^{g(z)}\,{dz\over z}\,,\eqno(21.4)$$
where
$$g(z)=nz+(3r-m)\ln z-3r\ln(1-z)+(n-m+r)\ln(2-z)\,.\eqno(21.5)$$
As before we can show that the asymptotic value of the
integral depends only on the behavior
of the integrand near $z=1$. This time we need not worry about a
three-legged saddle point, because we are sufficiently far from the
critical region near $\mu=0$. A~good path of integration turns out to
be $z=1-\alpha+it\mu^{-1/2}n^{-1/2}$, where 
$\alpha=\mu-{2\over 3}\mu^2+{3\over 5}\rho\mu^{-1/2}n^{-1/2}$. Indeed,
some beautiful cancellation occurs in the most significant terms:
$$\eqalignno{g(1-\alpha+it\mu^{-1/2}n^{-1/2})
&=g(1-\alpha)-{\textstyle{5\over
2}}t^2+O\bigl((\mu^{5/2}n^{1/2}+\mu^{-3/2}n^{-1/2}\rho^2)t\bigr)\cr
\noalign{\smallskip}
&\qquad\null+
O\bigl(\bigl((1+\vert\rho\vert)\mu^{-3/2}n^{-1/2}+\mu\bigr)t^2\bigr)\,,
&(21.6)\cr}$$
when $\vert t\vert\leq\log n$. The $O$~bounds follow from the fact
that the power series for $\log z$, $\log(1-z)$, and $\log(2-z)$
converge in the stated ranges.

The other factors of the integrand, besides $e^{g(z)}$, are 
$$(1-z)^{1-y}\,{dz\over
z}=\mu^{1-y}i\mu^{-1/2}n^{-1/2}\,dt\sum_{k=0}^{\infty}
\mu^k\beta^{k+1-y}\,,$$
where $\beta=(\alpha-it\mu^{-1/2}n^{-1/2})/\mu=1-{2\over
3}\mu+({3\over 5}\rho-it)\mu^{-3/2}n^{-1/2}$. We can now write the
integral as a factor independent of~$t$ times
$$\int
e^{-5t^2\!/2}(1+\gamma_1t+\gamma_2t^2+\cdots\,)\,dt\,,\eqno(21.7)$$
where the $\gamma$'s are functions of $\mu$ and~$\rho$, and the series is
convergent for $\vert t\vert\leq \log n$. The integrand is
superpolynomially small when $t=\pm\log n$; hence we can bound the error
terms for $\vert t\vert\leq\log n$, then integrate from~$-\infty$
to~$\infty$, showing that (21.7) is
$$\sqrt{{2\pi\over
5}}\,\bigl(1+O(\mu^{5/2}n^{1/2}+(1+\rho^2)\mu^{-3/2}n^{-1/2})\bigr)\,.
\eqno(21.8)$$ 

Finally we observe that the other factors nicely cancel the leading
terms of (21.3); only (21.1) and (21.2) are left. 
The overall formula (21.1) has a weaker estimate than (21.8) because
Stirling's approximation (21.3) is more sensitive to the value
of~$\rho$ and because of the term $g(1-\alpha)$.\quad\pfbox

Notice that Lemma 7 matches the first estimate of Lemma~5, which
says that the asymptotic probability of excess~$r$ is like
that for a normal distribution with mean ${2\over3}\mu^3$ and variance
of order~$\mu^3$, as long as $r=O(\mu^3)$. On the other hand,
the extreme tails for larger values of~$r$ are not as small as they
would be in a normal distribution; they decrease only as shown in the
second estimate of Lemma~5. For example, with probability $100^{-m}$
all edges will join vertices in the first $n/10$ vertices; so
there will be at least $0.9n$ isolated vertices, and the excess will
be at least $m-n+0.9n>0.4n$.

\proclaim Theorem 12.
The probability that a random multigraph with $n$~vertices and
$m=\half (n+\mu n)$ edges is clean, when $0\leq\mu\leq
n^{-1/5}$, is 
$$\exp\bigl(-{\textstyle{2\over
3}}\mu^4n+O\bigl(\bigl(\mu^{5/2}n^{1/2}\log
n+\mu^{-3/2}n^{-1/2}(\log n)^3\bigr)\bigr)\,.\eqno(21.9)$$

\proof
The probability decreases as $\mu$ increases. Therefore we need to
verify the result only for $\mu$ greater than $n^{-3/11}$ or so, when
the error estimate $\mu^{-3/2}n^{-1/2}(\log n)^3$ does not swamp the
main term $\exp(-{2\over 3}\mu^4n)=1-{2\over 3}\mu^4n+O(\mu^8n^2)$. 

Formula (21.1) is the probability that a random graph or multigraph
with $m$~edges is clean
and has excess~$r$, if we set $y=\half $. That probability is
superpolynomially small unless $\vert\rho\vert\leq\log n$,
because of the term $-\rho^2$ in the exponent. 
Extremely large values of~$\rho$, not covered by the hypotheses of
Lemma~7, are also negligible. 
Therefore we can sum
over~$r$ by integrating over~$\rho$ from $-\log n$ to
$+\log n$; and we can then extend the integral
from~$-\infty$ to~$\infty$ without changing its asymptotic value. Hence the
probability of cleanliness is
$$n^{-1/2}\sqrt{{3\over
20\pi}}\,\mu^{-3/2}\exp\bigl(-{\textstyle{2\over 3}}
\mu^4n\bigr)\int_{-\infty}^{\infty}e^{-3\rho^2\!/20}\sqrt{\mu^3n}\,
d\rho=\exp\bigl(-{\textstyle{2\over 3}}\mu^4n\bigr)\,,$$
plus the error term. Another nice bit of cancellation.\quad\pfbox

\proclaim
Corollary. The average number of edges added to an evolving multigraph
until it first becomes unclean is
$$\postdisplaypenalty10000
{\textstyle\half }n+{3^{1/4}\Gamma\bigl({1\over 4}\bigr)\over
2^{13/4}}\,n^{3/4}+O(n^{8/11+\epsilon})\,,\eqno(21.10)$$ 
and the standard deviation is of order $n^{3/4}$.

\proof
The stated average number is $\sum_{m\geq 0} p_m$, where $p_m$ is the
probability in the theorem. When $\mu\leq 0$, the probability of
uncleanliness is $O(n^{-1/3})$ by Theorem~7, so the sum for $0\leq
m<\half n$ is $\half n-O(n^{2/3})$. When $0\leq\mu\leq
n^{-3/11}(\log n)^{6/11}$, the probability of uncleanliness is
$O(n^{-1/11} (\log n)^{24/11})$ by (21.9); 
after that  the error is negligible in
comparison with the integral
$${\textstyle\half }n\int_0^{\infty}e^{-(2/3)\mu^4n}\,d\mu
={\textstyle\half }n^{3/4}\,{\textstyle{1\over 4}\,({3\over
2})}^{1/4} \int_0^{\infty}e^{-u}u^{-3/4}\,du=cn^{3/4}\,,$$
where $c$ is the coefficient of $n^{3/4}$ in (21.10). This proves
(21.10).

The expected value of $m^2$ at the stopping time is $\sum_{m\geq
0}(2m+1)p_m$, 
and we need to be especially careful when evaluating this sum; the
simple estimate $p_m=1-O(n^{-1/3})$ for $m\leq \half n$ will not
do, because it will obliterate significant terms by adding
$O(n^{5/3})$. Appropriate accuracy is maintained by computing
the expected value of $(m-\half n)^2$, which is
$${n^2\over 4}+\sum_{m\geq
0}(2m+1-n)p_m=\sum_{m=0}^{n/2}(n-2m)(1-p_m)
+\sum_{m=n/2}^{\infty}(2m-n)p_m+O(n)\,.$$
We can show that the terms for $m\leq \half n$ are now negligible,
because the cleanliness probability~$p_m$ is bounded below by the
probability that a multigraph with $m$~edges has excess~0. Therefore
$1-p_m=O\bigl(n^2\!/(n-2m)^3\bigr)$ when $m\leq m_0={1\over
2}n-n^{2/3+\epsilon}$, by the remarks preceding (13.23); and
$$\sum_{m=0}^{n/2}(n-2m)(1-p_m)=\sum_{m=0}^{m_0}O\left({n^2\over
(n-2m)^2}\right)
+\sum_{m=m_0}^{n/2}O(n-2m)=O(n^{4/3-\epsilon})+O(n^{4/3+2\epsilon})\,.$$
The other terms can be approximated by
$$\sum_{m=n/2}^{\infty}(2m-n)p_m=
\int_0^{\infty}{n^2\mu\over2}e^{-(2/3)\mu^4n}\,d\mu+O(n^{16/11+\epsilon})\,,$$
with an error estimate coming from the range $0\leq \mu\leq
n^{-3/11+\epsilon}$ as before. It follows that the variance is
asymptotic to this integral minus the square of $\bigl((21.10)-{1\over
2}n\bigr)$, namely $\bigl(3^{1/2}\Gamma({1\over
2})2^{-7/2}-c^2\bigr)n^{3/2}$.

Incidentally, the value of $c$ is approximately 0.50155, and the
standard deviation is approximately $0.1407n^{3/4}$.\quad\pfbox

\medskip
Once a graph begins to get dirty, its deficiency rises rapidly. For
fixed~$d$ we can estimate the probability of excess~$r$ and
deficiency~$d$ by taking $y=\half -d$ and multiplying (21.1)
by~$r^d\!/d!$, because of (7.16). The fact that (21.1) has $T(z)^{2r}$ in
the numerator instead of $T(z)^{2r-d}$ is unimportant, since
$T(z)^{2r}=T(z)^{2r-d}\sum{d\choose k}\bigl(T(z)-1\bigr)^k$. We obtain
a probability about $({2\over 3}\mu^4n)^d\!/d!$ times as large as before,
but this is damped rapidly by the factor $\exp(-{2\over 3}\mu^4n)$ when
$\mu$ becomes greater than~$n^{-1/4}$. We will look further at the
growth of deficiency in section~23.

\bigbreak\noindent
{\bf 22. A closer look.}\enspace
The structure theory of section~20 gives us more detailed information
about what happens when an evolving multigraph first changes from
clean to unclean. We learned in that section that the process of
adding a new edge $\langle x,y\rangle$ can be broken into three parts,
namely the introduction of half-edges at~$x$ and~$y$ followed by the
joining of those two edges. The deficiency can increase by~1 during
each of the first two stages. 

The probability that a clean graph becomes potentially 
deficient when a half-edge
is attached to~$x$ is the probability that the image $\MM{}'$ of the
half-edge after pruning and cancellation does not create a new vertex
not in~$\MM$. According to the analysis of section~20, the expected
number of times this happens is 
$$p_1(n)=\sum_m\,{2^mm!\,n!\over
n^{2m+1}}\,[w^mz^n]\,G_1(w,z)\,,\eqno(22.1)$$
\vskip-10pt
$$\eqalignno{%
G_1(w,z)&=e^{U(wz)/w}\sum_{r\geq 0}e_rw^r\,{2r\,T(wz)^{2r}\over
\bigl(1-T(wz)\bigr)^{3r+3/2}}\cr
\noalign{\smallskip}
&=\textstyle{{5\over 12}w^3z^2+{5\over
24}\,(2w^3+13w^4)\,z^3+\cdots\;.}
&(22.2)\cr}$$
The factor $2r$ covers the deficient choices of $x$, as in the first
term of (20.5).

Actually (22.2) is an overestimate, because some apparently bad
choices of~$x$ are ``false alarms.'' If the half-edge of~$x$ does not
add a vertex to~$\MM$, there's still a possibility that $y$ will be
chosen in the acyclic part; then the new edge $\langle x,y\rangle$
will not increase the excess and the multigraph will still be clean.
The expected number of false alarms is
$$p'_1(n)=\sum_m\,{2^mm!\,n!\over n^{2m+2}}\,[w^mz^n]\;{T(wz)\over
w}\,G_1(w,z)\,.\eqno(22.3)$$

The multigraph becomes unclean when $y$ is chosen if the half-edge
for~$y$ prunes and cancels to a reduced multigraph~$\MM{}''$ having
the same $2r+1$ vertices as~$\MM{}'$. This occurs with probability
$$p_2(n)=\sum_m\,{2^mm!\,n!\over
n^{2m+2}}\,[w^mz^n]\,G_2(w,z)\,,\eqno(22.4)$$
\vskip-10pt
$$\eqalignno{%
G_2(w,z)&=e^{U(wz)/w}\sum_{r\geq 0}e_rw^r\,
{(3r+\half )(2r+1)\,T(wz)^{2r+1}\over
\bigl(1-T(wz)\bigr)^{3r+7/2}}\cr
\noalign{\smallskip}
&=\textstyle{\half wz+{1\over
4}\,(2w+9w^2)z^2+\cdots\;.}
&(22.5)\cr}$$

Consequently we must have
$$p_1(n)-p'_1(n)+p_2(n)=1\eqno(22.6)$$
for all $n$; this identity is a nontrivial property of the bivariate
generating functions $G_1(w,z)$ and $G_2(w,z)$. When $n=6$, for
example, computer calculations show that
$$\displaylines{p_1(6)={10288260775\over 22039921152}\approx 0.4668\,;\qquad
p'_1(6)={38865625\over 612220032}\approx 0.0635\,;\cr
\noalign{\smallskip}
p_2(6)={13150822877\over 22039921152}\approx 0.5967\,.\cr}$$

We can use Lemma 7 to calculate the approximate values of these
quantities when $n$ is large, ignoring extreme terms not covered by
that lemma:
$$\eqalignno{%
p_1(n)&\sim{1\over
n}\int_0^{\infty}\int_{-\infty}^{\infty}2rB({\textstyle{3\over
2}},\mu,\rho,n)(\mu^{3/2}n^{1/2}\,d\rho)({\textstyle{1\over
2}}n\,d\mu)\cr
\noalign{\smallskip}
&=2^{-7/4}3^{-1/4}\Gamma({\textstyle{3\over 4}})\,n^{1/4}\,;&(22.7)\cr
\noalign{\medskip}
p'_1(n)&=\sum_m\,{2^{m-1}(m-1)!\,n!\over
n^{2m}}\,[w^mz^n]\,T(wz)G_1(w,z)\cr
\noalign{\smallskip}
&\sim{1\over
2}\int_0^{\infty}m^{-1}\int_{-\infty}^{\infty}
2rB({\textstyle{3\over
2}},\mu,\rho,n)(\mu^{3/2}n^{1/2}\,d\rho)({\textstyle{1\over
2}}n\,d\mu)\cr 
\noalign{\smallskip}
&\sim 2^{-7/4}3^{-1/4}\Gamma({\textstyle{3\over 4}})n^{1/4}\,;&(22.8)\cr
\noalign{\medskip}
p_2(n)&\sim {1\over n^2}\int_{n^{-1/3}}^{\infty}\,\int_{-\infty}^{\infty}
\textstyle{(3r+\half )(2r+1)B({7\over
2},\mu,\rho,n)(\mu^{3/2}n^{1/2}\,d\rho)(\half n\,d\mu)}\cr
\noalign{\smallskip}
&\sim\half\,.&(22.9)\cr}$$
Notice that $p_1(n)$ and $p'_1(n)$ are unbounded, so they must be
regarded as expected values (not probabilities). But $p_1(n)-p'_1(n)$
is the probability of a ``true alarm.''
As we might have guessed, the transition from clean to unclean occurs
about half the time when $x$ is chosen, half the time when $y$ is
chosen.

\bigbreak\noindent
{\bf 23. Giant growth.}\enspace
We know from the classical theory [\ER] that a giant component will
emerge when the number of edges is ${n\over 2}(1+\mu)$ for a positive
constant~$\mu$. 
The classical theory deals with graphs, but the same phenomenon will
occur with multigraphs, because random graphs are generated by the
multigraph process if we discard self-loops and duplicate edges;
discarded edges do not affect the size of components, and
comparatively few edges are discarded until the graph has gotten
rather dense (see [\Bi]).

Instead of relying on the classical theory, we can also deduce the
existence of a giant component by studying the generating function
$G(w,z)$. The proof
 is indirect: First we count the vertices that lie in trees
and unicyclic components, showing that there probably aren't too many
of those. Then we show that it is improbable to have two distinct
complex components.

The first part is easy, because there is a simple closed form for the
expected number of vertices in trees. If we mark just the vertices in
trees of size~$k$, by differentiating the generating function
$$G(w,z)\,\exp\bigl(-k^{k-2}w^{k-1}z^k\!/k!+k^{k-2}w^{k-1}z^ks^k\!/k!\bigr)$$
with respect to $s$ and setting $s=1$, we see that the expected number
of such vertices is just
$$\eqalign{{2^mm!\,n!\over n^{2m}}\;[w^mz^n]\;{k^{k-1}\over
k!}\,w^{k-1}z^kG(w,z)
&={2^mm!\,n!\,k^{k-1}\over n^{2m}\,k!}\;[w^{m-k+1}z^{n-k}]\,G(w,z)\cr
\noalign{\smallskip}
&={2^mm!\,n!\,k^{k-1}\over n^{2m}\,k!}\;{(n-k)^{2(m-k+1)}\over
2^{m-k+1}(m-k+1)!\,(n-k)!}\;;\cr}$$
this can be written
$${k^{k-1}\over k!}\;{2^{k-1}m^{\underline{\smash{k-1}}}
{\mkern1mu}n^{\underline{k}}\over
(n-k)^{2k-2}}\,\left({n-k\over n}\right)^{2m}\eqno(23.1)$$
in terms of falling factorial powers
$x^{\underline{k}}=x(x-1)\,\ldots\,(x-k+1)$. 

Asymptotically, we have $n^{\underline{k}}=n^k\bigl(1+O(k^2\!/n)\bigr)$
and $(n-k)^k=n^k\bigl(1+O(k^2\!/n)\bigr)$ for all~$k$; also
$(1-k/n)^n=e^{-k}\bigl(1+O(k^2\!/n)\bigr)$ for $k\leq\sqrt{n}$ and
$(1-k/n)^n\leq e^{-k}$ for $k\leq n$. If $\mu$ is a nonzero constant,
$\mu >-1$,  and if
$m={n\over 2}(1+\mu)$, expression (23.1) is
$${n\over 1+\mu}\;{k^{k-1}\over
k!}\,(1+\mu)^ke^{-k(1+\mu)}\bigl(1+O({\textstyle{{k^2\over
n}}})\bigr)\eqno(23.2)$$
for $k\leq\sqrt{n}$; and it is superpolynomially small when
$k=\sqrt{n}$, because it is $O\bigl(\bigl((1+\mu)
e^{-\mu}\bigr)^kk^{1/2}\bigr)$ and $(1+\mu)e^{-\mu}<1$. 
It is also superpolynomially small when $k>\sqrt{n}$, because we will
prove in section~27 below that a continuous approximation of the quantity
$$e^k\,\sqrt{{m-k\over m}}\,\sqrt{{n-k\over
n}}\;{2^k\,m^{\underline{k}}\,n^{\underline{k}}\over (n-k)^{2k}}\;
\left({n-k\over n}\right)^{2m}\eqno(23.3)$$
decreases when $k$ increases.

Let $\sigma$ be defined by the formula
$$(1+\mu){\mkern1mu}e^{-\mu}=(1-\sigma){\mkern1mu}e^{\sigma}\,,
\qquad\sigma=\mu+O(\mu^2)\,.\eqno(23.4)$$
Then $\sigma$ is the quantity called $1-x\bigl({1\over
2}(1+\mu)\bigr)$ in [\ER], and we have
$$\eqalign{\sum_{k\geq 1} {k^{k-1}\over k!}\;\bigl((1+\mu){\mkern1mu}
e^{-(1+\mu)}\bigr)^k
&=\sum_{k\geq 1} {k^{k-1}\over k!}\,\bigl((1-\sigma){\mkern1mu}
e^{-(1-\sigma)}\bigr)^k\cr
\noalign{\smallskip}
&=T\bigl((1-\sigma)e^{-(1-\sigma)}\bigr)=1-\sigma\,,\cr}$$
when $\mu$ is positive. 
By summing (23.2) over all~$k$, we conclude that the expected total
number of vertices in trees is
$${1-\sigma\over 1+\mu}\,n+O(\sigma^{-3})\,;
\eqno(23.5)$$ 
the error term $O(\sigma^{-3})$ here comes from summing
$\vartheta^2T\bigl((1-\sigma)e^{-(1-\sigma)}\bigr)$, which brings a
factor of~$k^2$ into each term.

For example, if $1+\mu=\ln 4$, we have $1-\sigma=\ln2$, because
${1\over 4}\ln 4=\half \ln 2$. When the number of edges reaches
$n\ln 2$ the expected number of vertices in trees will be ${1\over
2}n$.
And in general when the number of edges reaches ${n\over 2x}\ln{1\over
1-x}$, the expected number of vertices in trees will be $(1-x)n$, for
$0<x<1$.

The expected number of vertices in unicyclic components can be found
in a similar way, by differentiating 
$$G(w,z)\,e^{-V(wz)+V(wzs)}$$
with respect to $s$ and setting $s=1$. The generating function is
$${1\over
2}\;{T(wz)\over\bigl(1-T(wz)\bigr)^2}\,G(w,z)
=\bigl(\vartheta V(wz)\bigr)\,G(w,z)\,,\eqno(23.6)$$
and we have
$${T(z)\over\bigl(1-T(z)\bigr)^2}=\sum_{k\geq 1}\,{k^kQ(k)\over
k!}\,z^k\eqno(23.7)$$ 
by (3.12). The expected number of vertices belonging to 
unicyclic components of size~$k$  therefore can be expressed in closed form,
analogous to (23.1):
$$\eqalignno{\half \;&{k^kQ(k)\over k!}\;{2^mm!\,n!\over
n^{2m}}\;[w^{m-k}z^{n-k}]\,G(w,z)\cr
\noalign{\smallskip}
&\qquad =\half \;{k^kQ(k)\over
k!}\;{2^km^{\underline{k}}{\mkern1mu}n^{\underline{k}}\over
(n-k)^{2k}}\,\left({n-k\over n}\right)^{2m}\,.&(23.8)\cr}$$
Summing over $k$, and breaking the sum into two
parts $k\leq\sqrt{n}$ and $k>\sqrt{n}$ as above, now yields
$$\half \,\sum_{k\geq 1}\,{k^kQ(k)\over
k!}\,\bigl((1+\mu){\mkern1mu}
e^{-(1+\mu)}\bigr)^k \bigl(1+O({k^2\over n})\bigr)={1-\sigma\over
2\sigma^2}+O(\sigma^{-6}n^{-1})\,.\eqno(23.9)$$
(We will obtain sharper bounds in section~27.)

We have assumed in this discussion that $\mu$ is a constant. But our
relatively coarse asymptotic arguments are in fact valid if $\mu$
varies with~$n$, provided that it is not too small. Relation (23.4)
defines~$\sigma$ as an analytic function of~$\mu$,
$$\eqalignno{\sigma&=\mu-{2\over 3}\mu^2+{4\over 9}\mu^3-{44\over
135}\mu^4+{104\over 405}\mu^5-{40\over 189}\mu^6\cr
\noalign{\smallskip}
&\qquad\null+{7648\over 42525}\mu^7-{2848\over 18225}\mu^8+{31712\over
229635}\mu^9-{23429344\over 189448875}\mu^{10}+\cdots\;,&(23.10)\cr}$$
where the power series converges for $\vert\mu\vert <1$. The
quantity $\bigl((1+\mu){\mkern1mu}e^{-\mu}\bigr)^k$ is superpolynomially small
for $k=\sqrt{n}$ if $\mu$ is at least, say, $n^{-1/4}\log n$. We are
therefore justified in using (23.5)+(23.9) as the expected number  of
vertices in non-complex components whenever $\mu\geq n^{-1/4}\log n$.

Suppose $\mu=n^{-1/4}\log n$. Then the expected number of vertices in
unicyclic components is approximately ${1\over
2}\sigma^{-2}\sim\half \mu^{-2}=\half n^{1/2}(\log n)^{-2}$,
and a similar argument proves that the expected value of the square of
this number is approximately ${5\over 4}\sigma^{-4}\sim{5\over
4}n(\log n)^{-4}$. So 
 the probability of choosing two vertices in unicyclic components is
approximately ${5\over
4}n^{-1}(\log n)^{-4}$. This probability decreases steadily as $m$
increases, but even if it stayed fixed we would have to add
about ${4\over 5}n(\log n)^4$ more edges before hitting two unicyclic vertices,
i.e., before creating a new bicyclic component. By that time the
expected number of vertices in trees and unicyclic components will be
nearly zero, so the multigraph will almost surely contain no such
vertices. Therefore, if there is only one complex component present when
$\mu=n^{-1/4}\log n$, there will almost surely be only one complex
component from that time on; it will become gigantic. (We will obtain
sharper results in section 27; see Lemma~9 and its corollary.)

Let's look more closely at what happens as the giant component
develops. According to (23.5), it will have approximately
$$\left(1-{1-\sigma\over 1+\mu}\,\right)\,n={\mu+\sigma\over
1+\mu}\,n=2\mu n+O(\mu^2n)\eqno(23.11)$$
vertices when $m={n\over 2}(1+\mu)$; this is substantially larger than
the number $\half \mu^{-2}+O(\mu^{-1})$ of unicyclic vertices.
When $m$ increases by~1, the value of $\mu n$ increases by~2, so
(23.11) increases by~4. Notice that (23.11) agrees with the leading
term of (15.11).

We saw in section~21 that the expected excess $r$ is approximately
${2\over 3}\mu^3n$ when $m={n\over 2}(1+\mu)$, at least for
$0\leq\mu\leq n^{-1/4}$. We will prove momentarily that this
relationship continues to hold as long as $\mu$ remains $o(1)$; 
 but before giving the proof, let's look at the situation
heuristically. The probability that a new edge increases the excess is
the probability that 
 both of its endpoints lie in the cyclic part, namely
$(2\mu)^2$. The change in~$r$ with respect to~$m$ is
$(dr/d\mu)(d\mu/dm)=(2\mu^2n)(2/n)$, and this too is $(2\mu)^2$. So the
relation $r={2\over 3}\mu^3n$ is consistent with (23.11) when $\mu$ is
not too large.

The expected value of the deficiency $d$ turns out to be approximately
${2\over 3}\mu^4n$, about $\mu$ times~$r$. Heuristic 
justification comes from the considerations of section~20: When a new
edge $\langle x,y\rangle$ falls in the cyclic part, the probability
that $x$ is ``bad'' (in the sense that it increases the deficiency)
will be the number of reduced vertices $2r-d$ divided by the square
root of $3r-d$ times the size of the complex part $\bigl($see the remarks
following (20.8)$\bigr)$. So it will be
approximately ${4\over 3}\mu^3n$ divided by
$\bigl((2\mu^3n)(2\mu n)\bigr)^{1/2}$, namely ${2\over 3}\mu$. The same
holds for~$y$. Hence the expected increase in~$d$, given that $r$
increases, is ${4\over 3}\mu$. And the derivative of ${2\over
3}\mu^4n$ with respect to~$\mu$ is indeed ${4\over 3}\mu$ times the
derivative of ${2\over 3}\mu^3n$.

In order to carry out a rigorous proof as $\mu$ increases from
$n^{-1/4}$ to $n^{-1/5}$ to $n^{-1/6}$ and so on, we need to track the
full asymptotic spectrum of the behavior of~$r$ and~$d$, not
using just the leading terms. It turns out that $r$ and~$d$ are approximately
given by the following
joint functions of~$\mu$ and~$\sigma$, whose asymptotic series can be
computed from (23.10):
$$\eqalignno{r_{\mu}&={\mu^2-\sigma^2\over 2(1+\mu)}\,n\,;&(23.12)\cr
\noalign{\smallskip}
d_{\mu}&={3\mu^2-3\sigma^2-\sigma(\mu+\sigma)^2\over
2(1+\mu)}\,n\,.&(23.13)\cr}$$ 
Notice that the numerators of both $r_{\mu}$ and $d_{\mu}$ are
divisible by $(\mu+\sigma)n$, so $r_{\mu}$ and~$d_{\mu}$ are multiples
of the formula $(\mu+\sigma)n/(1+\mu)$ for giant component size
(23.11). The quantity $\mu+\sigma$ can also, incidentally, be expressed
as $\ln(1+\mu)-\ln(1-\sigma)$.

These values $r_{\mu}$ and $d_{\mu}$ also have a surprising relation
to the confluent hypergeometric series $F(z)=F(1;\,4;\,4z)$ of (7.5).
It is not difficult to check that
$$\openup2\jot
\eqalignno{%
F\bigl((\mu+\sigma)/4\bigr)={6e^{\mu+\sigma}\over (\mu+\sigma)^3}
\left({2r_{\mu}-d_{\mu}\over n}\right)
&={6(1+\mu)\over (1-\sigma)(\mu+\sigma)^3}\left({2r_{\mu}-d_{\mu}\over
n}\right)\,;&(23.14)\cr
{\vartheta F\bigl((\mu+\sigma)/4\bigr)\over
F\bigl((\mu+\sigma)/4\bigr)}
&={d_{\mu}\over 2r_{\mu}-d_{\mu}}\,.&(23.15)\cr}$$

The quantities $r_\mu$ and $d_\mu$ are not the exact expected values
of $r$ and~$d$. Indeed, the exact expected values are rational
numbers, when $m$ and $n$ are integers, while $\sigma$ is always
irrational when $\mu$ is rational. But we will prove that
the distributions of $r$ and~$d$ are approximately normal with
expectations $r_{\mu}$ and~$d_{\mu}$.

Before we can prove such a claim,
we need to improve the estimate of~$e_{rd}$ in (7.16),
because that estimate was derived only for fixed~$d$.

\proclaim
Lemma 8.
Let $F(z)$ be the function defined in (7.5). If $r\ra\infty$ and if 
$d$ varies in such a way that $d/r\ra 0$, the polynomial
$P_d(r)=[z^d]\,F(z)^{2r-d}$ satisfies
$$P_d(r)={F(s)^{2r-d}\over s^d}\;{(d/e)^d\over
d!}\,\left(1+O\left({d\over r}\right)\right)\,,\eqno(23.16)$$
where $s$ is the solution to $\vartheta F(s)/F(s)=d/(2r-d)$.

\proof
We have
$$P_d(r)={1\over 2\pi i}\oint {F(z)^{2r-d}\over z^d}\;{dz\over z}
={1\over 2\pi i}\oint e^{(2r-d)f(z)}\,{dz\over z}\,,$$
where $f(z)=\ln F(z)-\bigl(d/(2r-d)\bigr)\ln z$, integrated on the
circle $\vert z\vert=s$. By hypothesis, $\vartheta f(s)=0$. Using the
expansion formula
$$f(se^t)=\sum_{k=0}^n\,{t^k\over
k!}\,\vartheta^kf(s)+\int_0^t\,{x^n\over n!}\,
\vartheta^{n+1}f(se^{t-x})\,dx\eqno(23.17)$$
with $n=2$ and $t=i\theta$, we obtain
$$f(se^{i\theta})=f(s)-{\textstyle\half }\theta^2\vartheta^2
f(s)+O(\theta^3s)$$
because $\vert\vartheta^3 f(se^{i\theta})\vert=O(s)$. If $d\ra\infty$, the
contour integral is
$$\openup2\jot
\eqalignno{&{1\over
2\pi}\,\int_{-\pi}^{\pi}\exp\bigl((2r-d)\bigl(f(s)-{\textstyle{1\over
2}}\theta^2\vartheta^2f(s)+O(\theta^3s)\bigr)\bigr)\,d\theta\cr
&\qquad ={1\over 2\pi\sqrt{d}}\,\int_{-\pi\sqrt{d}}^{\pi\sqrt{d}}
\exp\bigl((2r-d)f(s)-t^2\!/2+O(t^2d/r)+O(t^3d^{-1/2})\bigr)\,dt\cr
&\qquad ={F(s)^{2r-d}\over s^d\sqrt{2\pi
d}\,}\bigl(1+O(d/r)+O(d^{-1/2})\bigr)\,,&(23.18)\cr}$$
because $\vartheta^2f(s)=s+O(s^2)=d/(2r-d)+O(d^2\!/r^2)$.
The terms $O(t^2d/r)$ and $O(t^3d^{-1/2})$ can safely be moved out of
the exponent because they are bounded when $\vert t\vert\leq d^{1/6}$
and $\vert t\vert\leq\sqrt{r/d}$. Larger values of $\vert t\vert$ are
unimportant in the integral because of the factor $e^{-t^2\!/2}$, and
because the relation
$$F(z)=3\int_0^1(1-u)^2e^{4zu}\,du$$
implies that $\vert F(z)\vert\leq F(\Re z)$; once the real part is
sufficiently small, we can neglect the remaining part of the path.

Equation (23.18) does not match (23.16) perfectly, although it would
be sufficient for the applications considered below. To derive the
sharper estimate claimed in (23.16) when $d$ is small, we can apply
(23.17) to $f(z)-z$ instead of to $f(z)$, obtaining
$$\openup2\jot
\eqalign{f(se^{i\theta})-se^{i\theta}&=f(s)-s-i\theta
s+O(\theta^2s^2)\,;\cr
f(se^{i\theta})&=f(s)+s(e^{i\theta}-i\theta-1)+O(\theta^2s^2)\cr
&=f(s)+{d\over 2r-d}(e^{i\theta}-i\theta-1)+O\left({\theta^2d^2\over
r^2}\right)\,.\cr}$$
The contour integral without the $O$ term can be evaluated exactly,
$$\openup2\jot
\eqalign{&{1\over 2\pi}\int_{-\pi}^{\pi}
\exp\bigl((2r-d)\bigl(f(s)+(e^{i\theta}-i\theta
-1)d/(2r-d)\bigr)\bigr)\,d\theta\cr
&\qquad ={F(s)^{2r-d}\over s^d}\cdot{1\over 2\pi}\int_{-\pi}^{\pi}
e^{(e^{i\theta}-1)d}d\theta/e^{i\theta d}\cr
&\qquad ={F(s)^{2r-d}\over s^d}\,[z^d]\,e^{(z-1)d}={F(s)^{2r-d}\over
s^d}\;{(d/e)^d\over d!}\,.\cr}$$
The $O$ term contributes a relative error of $O(d/r)$, because we have
$$\openup2\jot
\eqalign{\int_{-\pi}^{\pi}\bigl\vert\exp
\bigl((e^{i{\theta}}-1-i\theta)d\bigr)\bigr\vert
\theta^2\,d\theta
&=\int_{-\pi}^{\pi}e^{(\cos\theta-1)d}\theta^2\,d\theta\cr
&\leq\int_{-\pi}^{\pi}
e^{-c\theta^2d}\theta^2\,d\theta=O(d^{-3/2})\,,\cr}$$
where $c=2/\pi^2$.\quad\pfbox

\proclaim
Theorem 13.
The joint distribution of the
excess $r$ and deficiency $d$ of a random multigraph with
$m={n\over 2}(1+\mu)$ edges is approximately normal about the
expected values~$r_{\mu}$ and~$d_{\mu}$ in (23.12) and (23.13),
with zero covariance.
More precisely, there exists $\epsilon>0$ such that if
$$r=r_{\mu}+\rho\sqrt{\mu^3n}\,,\qquad
d=d_{\mu}+\delta\sqrt{\mu^4n}\,,\eqno(23.19)$$
the probability that a random multigraph has excess $r$ and
deficiency~$d$ is
$${3\over 4\pi\sqrt{5}\,\mu^{7/2}n}\,\exp\left(-{3\over
20}\rho^2-{3\over 4}\delta^2
+O\left((1+\vert\rho\vert+\vert\delta\vert)^2\,\mu^{1/2}
+{1+\vert\rho\vert^3\over(\mu^3n)^{1/2}}
+{1+\vert\delta\vert^3\over (\mu^4n)^{1/2}}\right)\right),\eqno(23.20)$$
when $n^{-1/4}\leq\mu\leq\epsilon$ and $n\ra\infty$, uniformly for
$\vert\rho\vert\leq\half \,\sqrt{n\mu^3}$ and
$\vert\delta\vert\leq \half \,\sqrt{n\mu^4}$.

\proof
Before proving formula (23.20), we can verify that its leading factor
 yields total
probability~1 when integrated over all values of~$r$ and~$d$
near~$r_{\mu}$ and~$d_{\mu}$: The integral over~$d$ gives a factor of
$\sqrt{4\pi\mu^4n/3}$, and the integral over~$r$ gives a factor of
$\sqrt{\mskip1mu20\pi\mu^3n/3}$.

Let $r$ and~$d$ be given by (23.19); the
probability of excess~$r$ and deficiency~$d$ is then
$${2^mm!\,n!\,e_{rd}\over n^{2m}(n-m+r)!\,2^{n-m+r}}\;[z^n]\;
{\bigl(2-T(z)\bigr)^{n-m+r}T(z)^{n-m+3r-d}\over
\bigl(1-T(z)\bigr)^{3r-d+1/2}}\,.\eqno(23.21)$$ 
We find the coefficient of $z^n$ by evaluating a contour integral as
in (10.11) and (21.4); it is
$${1\over 2\pi i}\oint e^{g(z)}(1-z)^{1/2}\,{dz\over
z}\,,\eqno(23.22)$$
\vskip-15pt
$$g(z)=nz+(3r-d)\bigl(\ln z-\ln(1-z)\bigr)+r\ln(2-z)-m\ln
z+(n-m)\ln(2-z)\,.\eqno(23.23)$$

The key to this theorem is the fact that, when $\rho=\delta=0$, there
is a saddle point at $z=1-\sigma$:
$$\eqalignno{{g'(1-\sigma)\over n}&=1+{\mu+\sigma\over
2(1+\mu)}\,\left({\sigma(\mu+\sigma)\over
1-\sigma}+{\sigma(\mu+\sigma)\over \sigma}-{\mu-\sigma\over
1+\sigma}\right)\cr
\noalign{\smallskip}
&\hskip10em\null- {1+\mu\over 2(1-\sigma)}-{1-\mu\over
2(1+\sigma)}=0\,.&(23.24)\cr}$$
Moreover, $g''(1-\sigma)=5\mu n+O(\mu^2n)$ in that case. 
If we integrate on the path $z=1-\sigma+it/\sqrt{\mu n}$, as we did in
Lemma~7 (section~21), the logarithm of the result will be
$$g(1-\sigma)+\ln\,{2^mm!\,n!\,e_{rd}\over n^{2m}(n-m+r)!\,
2^{n-m+r}}\,\sqrt{1\over 10\pi n}
+O(\mu)+O(\mu^{-3/2}n^{-1/2})\,,$$
where $r=r_{\mu}$ and $d=d_{\mu}$. The relevant quantity~$s$ needed in
Lemma~8 is 
$$s={\mu+\sigma\over 4}\eqno(23.25)$$ 
because of (23.15). The evaluation of the stated logarithm
is tedious, but it can be done in a reasonable amount of time with
computer assistance, using some simplifications such as
$$3r-d={\sigma(\mu+\sigma)^2\over 2(1+\mu)}\,n\,,\qquad
n-m+r={1-\sigma^2\over 2(1+\mu)}\,n\,.$$
The term $\ln\bigl((6r-2d)!/(3r-d)!\bigr)$ from (7.3) can be evaluated as
$(3r-d)\ln(3r-d)-3r+d+(6r-2d)\ln 2+\half \ln 2+O(\mu)$. It is not
difficult to verify that the terms involving $n\ln n$ cancel. There
are three terms involving $n\ln\mu$, namely $(3r-d)\ln\mu^2$, $-d\ln
\mu$, and $-(2r-d)\ln\mu^3$, coming respectively from within
expansions of $(3r-d)\ln(3r-d)$, $-d\ln s$, and
 $-(2r-d)\ln(2r-d)$; there are two other terms,
$-(3r-d)\ln\sigma$ from within $g(1-\sigma)$ and $+(3r-d)\ln\sigma$
from within $(3r-d)\ln(3r-d)$, which also cancel. The most difficult
part of the computation is the sum of about 16~terms that are rational
functions in~$\mu$ and~$\sigma$, times~$n$; these too sum to zero,
using relations (23.14). The net result is that the complicated
logarithm sums to $\ln 3-2\ln 2-\ln\pi-\half \ln 5-{7\over
2}\ln\mu -\ln n+O(\mu)+O(\mu^{-4}n^{-1})$; this proves the theorem
when $\rho=\delta=0$.

For the case of general $\rho$ and $\delta$ the calculations are
similar but even worse. We now choose the integration path
$$z=1-\sigma-{\textstyle{3\over 5}}\,\rho/\sqrt{\mu n}+it/\sqrt{\mu
n}\,;\eqno(23.26)$$
the first-order effects of $\rho$ and $\delta$ then cancel out, and
the second-order effects contribute $-{3\over 20}\rho^2-{3\over
4}\delta^2$ to the logarithm of the result.\quad\pfbox

\noindent\bigbreak
{\bf 24. A waiting game.}\enspace
Now let's consider a little game. Start with an empty multigraph and
add edges repeatedly at random until either (1)~two different complex
components are present; or (2)~the multigraph is unclean. Case~1
represents the event ``we have left the top line of Figure~1 before
leaving the top line of Figure~2.''

Let $G_0(w,z)$ be the bgf for all multigraphs such that the game has
not yet stopped. Then
$$\sum_m\,{2^mm!\,n!\over n^{2m}}\;[w^mz^n]\,G_0(w,z)\eqno(24.1)$$
is the expected running time of the game. We have
$$G_0(w,z)=e^{U(wz)/w}\sum_r w^rK_r(wz)\,,\eqno(24.2)$$
where $K_r(z)$ generates all clean cyclic multigraphs, weighted by the
probability that they will arise as the cyclic part of a multigraph
occurring during the game. 

We learned in section 17 how to compute weighting factors that account
for the history of transitions in Figure~1 among clean multigraphs;
and we learned more specifically in section~20 how these coefficients
arise as a multigraph gains random edges. In consequence, we can
conclude that $K_r(z)=k_rT(z)^{2r}/\bigl(1-T(z)\bigr)^{3r+1/2}$, where
$k_1=e_1={5\over 24}$ and the later coefficients obey the rule
$$k_{r+1}={3\over 2}rk_r\,.\eqno(24.3)$$
Here's why: Given $k_rT^{2r}/(1-T)^{3r+1/2}$, the generating function
for a clean vertex~$x$ is
$${\textstyle{1\over
2}}k_rT^{2r+1}/(1-T)^{3r+5/2}+3rk_rT^{2r+1}/(1-T)^{3r+5/2}\,,\eqno(24.4)$$
where the first term corresponds to cases where $x$ is in the
unicyclic part. Similarly, given the generating function ${1\over
2}k_rT^{2r+1}/(1-T)^{3r+5/2}$ after $x$ is chosen to be unicyclic, the
generating function for a clean unicyclic~$y$ is
$${\textstyle{5\over 2}\,{1\over
2}}k_rT^{2r+2}/(1-T)^{3r+9/2}\,;\eqno(24.5)$$ 
here ${5\over 2}=1+1+\half $, for choosing $y$ on the half-edge
to~$x$, or on the self-loop attached to that half-edge, or in a
different unicyclic component. We obtain a new bicyclic component if
and only if both~$x$ and~$y$ were unicyclic. Therefore the generating
function for cases where the game continues is
$$\bigl((3r+\textstyle{5\over 2})(3r+\textstyle{1\over
2})-\textstyle{5\over 4}\bigr)\,k_rT^{2r+2}/(1-T)^{3r+9/2}\,.$$
As in (20.9) and (20.10), we multiply by $(1-T)/(6r+6)$ to account for
merging $\langle x,y\rangle$ with the existing edges. This proves
(24.3).

Equation (24.3) implies, of course, that
$$k_r={5\over 36}\,\left({3\over 2}\right)^r(r-1)!\,.\eqno(24.6)$$
Comparing this to the case $d=0$ of (7.16), we have
$$k_r={5\pi\over 18}\,e_r\bigl(1+O(r^{-1})\bigr)\,.\eqno(24.7)$$
Therefore the similar 
calculations of section~22, where we found that
$p_1(n)-p'_1(n)=1-p_2(n)\sim\half $,  tell us that {\sl the game will
stop in Case~(2) with probability\/}~${5\pi\over 18}$. This provides
further evidence in support of the top-line conjecture that was made
in section~18.

We can now try to compute the expected time for the game to be
completed, but it appears to be quite complicated. The contribution to
(24.1) from a given $m$ and~$r$ can be obtained by changing $e_r$
to~$k_r$ in (13.17) when $m$ and~$r$ are not too large; this means we
want to evaluate
$$\sum_{k\geq 0}\sqrt{2\pi\over 3}\;{(\half \,3^{2/3}\mu)^k\over
k!}\,\left({1\over\Gamma(1/2-2k/3)}
+\sum_{r\geq 1}\,{5\over 36}\,\left(\half \right)^r\;
{(r-1)!\over \Gamma(r+1/2-2k/3)}\right)\eqno(24.8)$$
in place of (14.1), representing $e^{\mu^3\!/6}$ times the probability
that the game is still alive after $m$~edges. 
 The inner sum is known to be ${5\over 72}$ times
$$\eqalignno{\sum_{r\geq 0}\left(\half \right)^r\!{r!\over
\Gamma(r+3/2-2k/3)}
&={1\over\Gamma(3/2-2k/3)}\,F\left(1,1;{3\over 2}-{2k\over 3}\,;{1\over
2}\,\right)\cr
\noalign{\smallskip}
&={1\over\Gamma(1/2-2k/3)}\,\left(\psi
\left({3\over 4}-{k\over 3}\right)-\psi\left({1\over 4}-{k\over
3}\right)\right)\,,\quad&(24.9)\cr}$$
so it has the value $\sqrt\pi$ when $k=0$. (Here, as usual, $\psi(z)=
\Gamma'(z)/\Gamma(z)$.) Further study of (24.8)
should prove to be interesting.

\bigbreak\noindent
{\bf 25. Waiting time in general.}\enspace
Bivariate generating functions provide a useful tool for studying the
``first occurrences'' of particular graphs or multigraphs, as shown in
[\FKP]. The special problems considered in that paper can be put into
the following general framework.

Let ${\cal S}$ be any collection of multigraphs, with bgf $S(w,z)$.
Suppose we wish to study the first time that an evolving multigraph on
$n$~vertices does not lie in~${\cal S}$. If $[z^n]\, S(0,z)=0$, the
empty graph on $n$~vertices is not in~${\cal S}$, so the process never
gets started. Otherwise, the probability that an evolving multigraph lies
in~${\cal S}$ when it has $m-1$ edges but not when it has $m$ is
$${2^{m-1}(m-1)!\,n!\over
n^{2m}}\;[w^mz^n]\,(w\vartheta_z^2-2\vartheta_w)S(w,z)\,.\eqno(25.1)$$
The proof is simple, by definition of the operators $\vartheta_z$
and~$\vartheta_w$, because the probability in question is
$${2^{m-1}(m-1)!\,n!\over
n^{2(m-1)}}\;[w^{m-1}z^n]\,S(w,z)-{2^mm!\,n!\over
n^{2m}}\;[w^mz^n]\,S(w,z)\,.$$

For convenience we shall write
$$\nabla\,S(w,z)=(w\vartheta_z^2-2\vartheta_w)\,S(w,z)\,;\eqno(25.2)$$
we call $\nabla S$ the bgf for ``stopping configurations,'' while $S$ itself
is the bgf for ``going configurations.''

The operator $\Phi_n$, introduced in [\FKP], is 
$$\Phi_nF(w,z)=\sum_{m=1}^{\infty}\,{2^{m-1}(m-1)!\,n!\over
n^{2m}}\;[w^mz^n]\,F(w,z)\,.\eqno(25.3)$$
Equations (25.1)--(25.3) imply that $\Phi_n\nabla S(w,z)$ is the
probability that a stopping configuration will be encountered when
some edge is added to an initially empty multigraph. A~similar
operator 
$$\widehat{\Phi}_n\,\widehat{F}(w,z)=\sum_{m=1}^{\infty}\,{n!\over
2m{n(n-1)/2\choose m}}\;[w^mz^n]\,\widehat{F}(w,z)\eqno(25.4)$$
for graphs instead of multigraphs is considered in [\FKP],
but we will restrict consideration to multigraphs for simplicity. (As
one might expect from section~6, we should use the operator
$$\widehat{\nabla}=w(\vartheta_z^2-\vartheta_z-2\vartheta_w)-2\vartheta_w\eqno(25.5)$$
in place of $\nabla$ when defining stopping configurations for the
graph process.)

Several examples will help clarify these definitions and demonstrate
their usefulness. Since the bgf $G(w,z)$ for all multigraphs satisfies
$\vartheta_z^2G=2w^{-1}\vartheta_wG$, equation (4.2), we have $\nabla
G=0$; this, of course, is obvious, because there are no stopping
configurations when all multigraphs are permitted.

\medskip\noindent
{\bf Example 1.}\enspace
 Let $S(w,z)$ be the bgf for all multigraphs having nothing
but self-loops. Clearly $S(w,z)=e^{ze^{w/2}}$, because $ze^{w/2}$ is
the bgf for a single vertex with nothing but self-loops.
Formula (25.2) now tells us that
$$\nabla S(w,z)=wz^2e^we^{ze^{w/2}}\,,\eqno(25.6)$$
because $\vartheta_z^2S=z^2e^wS+ze^{w/2}S$ and
$\vartheta_wS={w\over 2}ze^{w/2}S$.
 Thus, by (25.1), the probability that an
evolving multigraph first fails to lie in~$S$ when it acquires the
$m$\/th edge is
$$\eqalign{%
&{2^{m-1}(m-1)!\,n!\over n^{2m}}\;[w^mz^n]\,wz^2e^we^{ze^{w/2}}\cr
\noalign{\smallskip}
&\qquad={2^{m-1}(m-1)!\,n!\over n^{2m}}\;[w^m]\;{we^we^{(n-2)w/2}\over
(n-2)!}\cr
\noalign{\smallskip}
&\qquad={2^{m-1}(m-1)!\,n(n-1)\over
n^{2m}}\;[w^{m-1}]\,e^{nw/2}={n(n-1)\over n^{m+1}}\,.\cr}$$
And sure enough, $n^{1-m}-n^{-m}$ is obviously the probability that
a sequence of edges $\langle x_1,y_1\rangle\,\ldots\,\langle
x_m,y_m\rangle$ will have $x_1=y_1$, \dots, $x_{m-1}=y_{m-1}$, $x_m\neq
y_m$.

\medskip\noindent
{\bf Example 2.}\enspace
Let $S(w,z)$ be the bgf for all acyclic multigraphs, namely
$e^{U(w,z)}=e^{U(wz)/w}$. The formulas 
$$\vartheta_z U=w^{-1}T\,,\qquad
\vartheta_z^2U=w^{-1}T/(1-T)\,,\qquad
\vartheta_wU={\textstyle\half }w^{-1}T^2\eqno(25.7)$$
were derived in section 4, and we have
$$\vartheta_z^2e^F=
(\vartheta_z^2\,F)e^F+(\vartheta_z\,F)^2e^F\eqno(25.8)$$
for any $F=F(w,z)$; hence
$$\nabla e^U={T\over 1-T}\,e^U\,.\eqno(25.9)$$
These are the stopping configurations that define the appearance of
the first cycle in an evolving multigraph. The term $T^ke^U$
corresponds to a first cycle of length~$k$; therefore if we replace
$T^k$ by~$kT^k$ and sum over all stopping times, we get an expression
for the expected length of the first cycle,
$$\Phi_n\,{T\over (1-T)^2}\,e^U\,.\eqno(25.10)$$
This was one of the main problems studied in [\FKP], where it was shown
that the expected length is proportional to~$n^{1/6}$ although the
standard deviation is proportional to~$n^{1/4}$.

\medskip\noindent
{\bf Example 3.}\enspace
Let $S(w,z)=U(w,z)$ be the bgf for unrooted trees. This is a perverse
example, thrown in primarily because
(25.7) gives us the information we need to calculate
$$\eqalignno{\nabla U&={T\over 1-T}-{T^2\over w}\cr
\noalign{\smallskip}
&=wz+(-w+2w^2)z^2+(-2w^2+{\textstyle{9\over 2}}\,w^3)z^3
+\cdots\;.&(25.11)\cr}$$
What is the meaning of these negative coefficients? 

The example does make sense, if we rephrase our interpretation of
(25.1). The exact meaning of
$${2^{m-1}(m-1)!\,n!\over n^{2m}}\;[w^mz^n]\,\nabla S(w,z)$$
is, ``the probability that an 
evolving multigraph leaves ${\cal S}$ when the
$m$\/th~edge is added, minus the
 probability that it enters $\cal S$ when the
$m$th edge is added.'' In our example, $U(0,z)=z$; when there are
two or more vertices, the empty multigraph is not a tree, but it can
become one later. The bgf for becoming a tree is $w^{-1}T^2$,
corresponding to an ordered pair of rooted trees with $m-1$ edges. The
bgf for adding a new edge $\langle x,y\rangle$ to a tree is
$\sum_{k\geq 1}T^k$, where the term~$T^k$ corresponds to cases where
$x$ and~$y$ are at distance~$k$. (Each appearance of $T=T(wz)$
includes an implicit edge touching the tree root, because $w$ and~$z$
appear with equal powers in every term.)

Example 3 cautions us to interpret the operators~$\nabla$ and~$\Phi_n$
a bit more carefully. In general, we have the identity
$$\Phi_n\nabla S(w,z)=n!\;[z^n]\;S(0,z)-\lim_{m\rightarrow\infty}\,
{2^mm!\,n!\over n^{2m}}\;[w^mz^n]\,S(w,z)\,,\eqno(25.12)$$
for any bgf $S(w,z)$ such that the limit exists, because
$$\Phi_n\nabla
S(w,z)=\sum_{m=1}^{\infty}\,\left({2^{m-1}(m-1)!\,n!\over
n^{2(m-1)}}\;[w^{m-1}z^n]\,S(w,z)-{2^mm!\,n!\over
n^{2m}}\;[w^mz^n]\,S(w,z)\right)\,.$$ 
A sufficient condition for the limit to exist is that the coefficients
of $\nabla S(w,z)$ are nonnegative. A~sufficient condition for the
coefficients to be nonnegative is that $S(w,z)$ should represent a
family of multigraphs~${\cal S}$ with the property that the deletion
of any edge preserves membership in~${\cal S}$.

\medskip\noindent
{\bf Example 4.}\enspace
Let $S(w,z)=G(w,z)-C(w,z)$ be the bgf for all disconnected
multigraphs. The stopping configurations now represent the first time
an evolving multigraph becomes connected. Since $G(w,z)=e^{C(w,z)}$,
we have
$$\eqalign{\vartheta_wC&=\vartheta_w\ln G=(\vartheta_wG)/G\,;\cr
\vartheta_zC&=\vartheta_z\ln G=(\vartheta_zG)/G\,;\cr
\vartheta_z^2C&=(\vartheta_z^2G)/G-(\vartheta_zG)^2\!/G^2\,;\cr}$$
hence
$$\nabla S=\nabla G-\nabla C=w(\vartheta_zC)^2\,.\eqno(25.13)$$
Of course! This is an edge that joins an ordered pair of vertices
marked in distinct components.

\medskip\noindent
{\bf Example 5.}\enspace
Let $S(w,z)$ be any bgf of the form
$$S(w,z)=e^{U(w,z)+V(w,z)}H(w,z)\,.\eqno(25.14)$$
Then we can use (25.7) and (4.9) to compute
$$\nabla
S=e^{U+V}\bigl((2T\vartheta_z-2\vartheta_w+we^{-V}\vartheta_z^2
e^V)\,H\bigr)\,.\eqno(25.15)$$ 
For example, when $S(w,z)=G(w,z)$, the left side of (25.15) is zero,
and $H(w,z)$ is the bgf we have called $E(w,z)$. Equating the right
side of (25.15) to zero gives the differential equation (5.1) that we
originally used to compute $E(w,z)$.

In the special case $H(w,z)=1$, the stopping configurations correspond
to the first time an evolving multigraph acquires a bicyclic
component, i.e., the time when its excess changes from~0 to~1. This is
another problem that was considered in [\FKP], where it was shown that
the expected number of unicyclic components present at the time is
${1\over 6}\ln n+O(1)$. 
If we express $H$ in terms of univariate generating functions,
$$H(w,z)=\sum_{r\geq 0}w^rH_r(wz)\,,\eqno(25.16)$$
then (25.15) can be written
$$\nabla S=e^{U+V}\sum_{r\geq 1}w^r\nabla H_r(wz)\,,\eqno(25.17)$$
where the univariate function $H_r(z)$ is related to (5.3):
$$\nabla
H_r(z)=e^{-V}\vartheta^2e^VH_{r-1}(z)-2\bigl(r+(1-T)\vartheta\bigr)
H_r(z) \,.\eqno(25.18)$$

\medskip\noindent
{\bf Example 6.}\enspace
Specializing Example 5 further, let
$$S(w,z)=e^{U(w,z)+V(w,z)}\,\sum_{r=0}^Rw^rE_r(wz)\,,\eqno(25.19)$$
where $R$ is any nonnegative integer. Then the stopping configurations
$\nabla S$ represent the time when an evolving multigraph first
acquires excess $R+1$. Expression (25.18) becomes almost trivial because
$\nabla H_r$ is zero for all $r\neq R+1$; we have
$$\nabla
S(w,z)=w^{R+1}e^{U(w,z)}\,\vartheta_z^2e^{V(wz)}E_R(wz)\,.\eqno(25.20)$$

This family ${\cal S}$ has the property that $\Phi_n\nabla S=1$, by
(25.12), because a multigraph surely acquires excess $R+1$ at some
time $m\leq n+R+1$. We can write the identity $\Phi_n\nabla S=1$ more
explicitly, using our known formula for~$E_R$, and using $r$ in place
of~$R$:
$$\Phi_n\left(w^{r+1}e^{U(wz)/w}\,\vartheta_z^2\,
\sum_{d=0}^{2r}e_{rd}\,{T(wz)^{2r-d}\over
\bigl(1-T(wz)\bigr)^{3r-d+1/2}}\right)=1\,,\eqno(25.21)$$
for all $n\geq 1$ and $r\geq 0$. Moreover, we can write (25.20) in the
form
$$\nabla
S(w,z)=2w^{R+1}e^{U(w,z)+V(wz)}\bigl(R+1+(1-T)\vartheta_z\bigr)
E_{R+1}(wz)\,,$$
using (5.3). Setting $r=R+1$ and applying (20.9) gives us another way
to express (25.21),
$$\Phi_n\left(w^re^{U(wz)/w}\,\sum_{d=0}^{2r}(6r-2d)e_{rd}\,{T(wz)^{2r-d}\over
\bigl(1-T(wz)\bigr)^{3r-d+3/2}}\right)=1\,,\eqno(25.22)$$
for all $n\geq 1$ and $r\geq 1$.

For example, the case $r=1$ of (25.22) is
$$\Phi_n\left(we^U\left({5\over 4}\,{T^2\over (1-T)^{9/2}}+{1\over
2}\,{T\over (1-T)^{7/2}}\right)\right)=1\,.\eqno(25.23)$$
The operator $\Phi_n$ is defined in (25.3) to be a sum over~$m$, and
the $m$\/th term of (25.23) is 
$$\eqalignno{&{2^{m-1}(m-1)!\,n!\over n^{2m}}\;[w^mz^n]\,
we^{U(w,z)}f\bigl(T(wz)\bigr)\cr
\noalign{\smallskip}
&\qquad\qquad ={1\over 2m(n-m+1)}\,{2^mm!\,n!\over
n^{2m}(n-m)!}\;[z^n]\,U(z)^{n-m+1}f\bigl(T(z)\bigr)\cr 
\noalign{\smallskip}
&\qquad\qquad ={1\over 4m(n-m+1)}\;{2^mm!\,n!\over
n^{2m}(n-m)!}\;[z^n]\,U(z)^{n-m}g\bigl(T(z)\bigr)\,,&(25.24)\cr}$$
where $f(T)={5\over 4}T^2\!/(1-T)^{9/2}+\half T/(1-T)^{7/2}$ and
$g(T) =(2-T)Tf(T)$. We can write
$$g(T)={5/4\over (1-T)^{9/2}}-{2\over (1-T)^{7/2}}-{1/2\over
(1-T)^{5/2}}+{2\over (1-T)^{3/2}}-{3/4\over (1-T)^{1/2}}\,,$$
so we can evaluate (25.24) by summing five applications of formula
(10.1). 
The value is negligibly small unless $m$ is $\half n+O(n^{2/3})$,
hence the factor $4m(n-m+1)$ can be assumed to equal $n^2+O(n^{5/3})$.
The five terms of~$g$ yield values of order $n^{4/3}$, $n$, $n^{2/3}$,
$n^{1/3}$, and~1 respectively, according to (10.1); thus the leading
term ${5\over 4}/(1-T)^{9/2}$ must be responsible for the major
contribution to (25.23), and the $m$\/th term of the sum when
$m=\half n+\half \mu n^{2/3}$ will be
$$\textstyle{{5\over 4}n^{-2/3}\sqrt{2\pi}\,A({9\over 2}
\,,\mu)+O(n^{-1})\,.}$$ 
Summing over $m$ yields 1. Therefore it must be true that
$$\int_{-\infty}^{\infty}A\left({9\over
2}\,,\mu\right)\,d\mu={8/5\over \sqrt{2\pi}}\,.$$
This integral formula is not at all obvious from the definition of
$A(y,\mu)$ in (10.2), and it would be interesting to find a direct
proof. 

The argument we have just given can be extended to arbitrary~$r$,
starting with (25.22), and it implies the following remarkable result:
$$\int_{-\infty}^{\infty}A(3r+{\textstyle{3\over 2}}\,,\mu)\,d\mu
={1\over 3re_r\sqrt{2\pi}}\,,\qquad\hbox{integer }r\geq
1\,.\eqno(25.25)$$
By (8.17) we can also write
$$\int_{-\infty}^{\infty}A(3r+{\textstyle{3\over 2}}\,,\mu)\,d\mu
={1\over 3}\,\left({2\over 3}\right)^r\,{\Gamma(r)\,\sqrt{2\pi}\over
\Gamma(r+{5\over 6})\,\Gamma(r+{1\over 6})}\,,\quad\hbox{integer }r\geq
1\,.\eqno(25.26)$$

We have just proved that, if $M_{r,n}=\half n+{1\over
2}U_{r,n}n^{2/3}$ is the number of edges when the excess first
reaches~$r$, then 
$$\Pr(M_{r,n}=m)\;=\;6re_r\sqrt{2\pi}\,A\bigl(3r+{\textstyle{3\over
2}}\,,\mu)n^{-2/3}+O(n^{-1})\,;\eqno(25.27)$$
hence $U_{r,n}\ra U_r$ in distribution, where $U_r$ has the density
function 
$$f_r(\mu)=3re_r\sqrt{2\pi}\,A\bigl(3r+{\textstyle{3\over
2}}\,,\mu\bigr)\,,\qquad-\infty <\mu<\infty\,.\eqno(25.28)$$
Combining this formula with (13.17), we have
$$\eqalign{\sqrt{2\pi}\,e_rA(3r+{\textstyle\half }\,,\mu)=
\lim_{n\rightarrow\infty}\,\Pr({\cal E}_r)
&=\lim_{n\rightarrow\infty}\,\Pr(M_{r,n}\leq m<M_{r+1,n})\cr
\noalign{\smallskip}
&=\int_{-\infty}^{\mu}\bigl(f_r(u)-f_{r+1}(u)\bigr)\,du\,,\cr}$$
whence
$$\sqrt{2\pi}\,e_r\,A'(3r+{\textstyle{1\over
2}}\,,\mu)=f_r(\mu)-f_{r+1}(\mu)\,.\eqno(25.29)$$
In fact, (25.29) can be derived also by setting $y=3r+\half $ in
the formula
$$A'(y,\mu)=(y-{\textstyle\half })A(y+1,\mu)-{\textstyle{1\over
2}}\,y(y+2)A(y+4,\mu)\,,\eqno(25.30)$$
which is a consequence of (10.22) and (10.23).

\bigbreak\noindent
{\bf 26. Continuous excess.}\enspace
Let $I(y)$ be the integral in (25.25) when the parameter~$r$ is not
necessarily an integer:
$$I(y)=\int_{-\infty}^{\infty}A(y,\mu)\,d\mu\,.\eqno(26.1)$$
It is natural to conjecture that formula (25.26) holds for $y$ in
general:
$$I(y)={2^{y/3}\,\sqrt{\pi}\,\,\Gamma(y/3-1/2)\over
3^{y/3+1/2}\,\Gamma(y/3+1/3)\,\Gamma(y/3-1/3)}\,,\qquad
y>{3\over 2}\,.\eqno(26.2)$$
$\bigl($The condition $y>{3\over 2}$ is necessary and sufficient for
convergence of the integral, because of (10.3) and (10.4).$\bigr)$
And indeed, this conjecture is true.

\proclaim
Theorem 14. The integral (26.1) has the closed form (26.2).

\proof
Let $I_0(y)$ be the right-hand side of (26.2); we wish to show that
$I(y)=I_0(y)$. Clearly
$$I_0(y+3)={2y-3\over (y+1)(y-1)}\,I_0(y)\,,\qquad y>{3\over
2}\,.\eqno(26.3)$$ 
Since $\int_{-\infty}^{\infty}A'(y,\mu)\,d\mu=0$ for $y>\half $,
by (10.3) and (10.4), we can integrate (25.30) and replace~$y$ by
$y-1$ to get the same recurrence for $I(y)$:
$$I(y+3)={2y-3\over (y+1)(y-1)}\,I(y)\,,\qquad y>{3\over
2}\,.\eqno(26.4)$$
Therefore $I(y)/I_0(y)$ is a periodic function, and we need only prove
asymptotic equivalence $I(y)\sim I_0(y)$ as $y\ra\infty$ in order to
verify strict equality $I(y)=I_0(y)$ for all $y>{3\over 2}$.

The duplication and triplication formulas for the Gamma function
provide us with an alternate expression for $I_0(y)$:
$$I_0(3y)=\left({9\over
2}\right)^{y-1}\,{\Gamma(2y-1)\over\Gamma(3y-1)}\sim
{1\over\sqrt{6}}\,\left({2e\over 3y}\right)^y\,.\eqno(26.5)$$
To show that $I(y)$ has the same asymptotic behavior, we break the
integral into two parts,
$$I(y)=\int_{-\infty}^0A(y,\mu)\,d\mu+\int_0^{\infty}A(y,\mu)\,d\mu
=I_{-}(y)+I_{+}(y)\,.\eqno(26.6)$$
By definition (10.2) we have
$$I_{+}(y)=
{1\over 3^{(y+1)/3}}\,\int_0^{\infty}\,\sum_{k\geq 0}\,
{e^{-\mu^3\!/6}\bigl(\half 3^{2/3}\mu\bigr)^k\,d\mu\over
k!\,\Gamma\bigl((y+1-2k)/3\bigr)}\,;\eqno(26.7)$$
we will show that the asymptotic value of $I_{+}(y)$ can be obtained
by interchanging summation and integration, then estimating the
resulting sum.

Let $a_k$ be the $k$th term after integration,
$$a_k=\int_0^{\infty}\,{e^{-\mu^3\!/6}\bigl({1\over
2}3^{2/3}\mu\bigr)^k\,d\mu\over k!\ \Gamma\bigl((y+1-2k)/3\bigr)}
={2^{(1-2k)/3}3^{k-2/3}\,\Gamma\bigl((k+1)/3\bigr)\over k!\
\Gamma\bigl((y+1-2k)/3\bigr)}\,.\eqno(26.8)$$
If $a_k=0$ then $a_{k+3}=0$; otherwise we have
$${a_{k+3}\over a_k}={(2k+5-y)(2k+2-y)\over
4(k+2)(k+3)}\,,\eqno(26.9)$$
which is greater than 1 when $k<{1\over 4}y-{17\over 8}-{5\over
4}(y+{3\over 2})^{-1}$, less than~1 when $k$ exceeds that value, and
nonnegative except for one or two values of~$k$ near $\half y$. So
the largest terms $a_k$ occur when $k$ is near~${1\over 4}y$.
If $y>5$ and $k>y/2$, we have
$$\left\vert{a_{k+3}\over a_k}\right\vert\;\leq\;\left(1-{y-5\over
2k}\right)^2\leq \bigl(e^{3/k}\bigr)^{(5-y)/3}\leq\left({k+3\over
k}\right)^{(5-y)/3}\,,$$ 
and it follows that $a_k=O(k^{(5-y)/3})$ as $k\to\infty$. Therefore $\sum\vert
a_k\vert$ exists, and the interchange of summation and integration is
justified, at least for large~$y$.

Let $k={1\over 4}y+x$, where $\vert x\vert\leq y^{1/2+\epsilon}$.
Then Stirling's formula tells us that
$$\ln a_k={y+3\over 3}\,\ln 2+{y-2\over 3}\,\ln 3-{2y+3\over 6}\,\ln y
+{y\over 3}-\ln\,\sqrt{\pi}-{8x^2\over
3y}+O(y^{3\epsilon-1/2})\,.\eqno(26.10)$$ 
If $0<\epsilon<{1\over 6}$, this implies that the sum of all terms for
$\vert x\vert >y^{1/2+\epsilon}$ is superpolynomially small in relation to
the sum of terms for $\vert x\vert\leq y^{1/2+\epsilon}$; hence
$$\sum_{k=0}^{\infty}a_k\;\sim\;{2^{(y+3)/3}3^{(y-2)/3}e^{y/3}\over
y^{(2y+3)/6}\;\sqrt{\pi}\,}\;\int_{-\infty}^{\infty}e^{-8x^2\!/(3y)}\,dx$$
and we have
$$I_{+}(y)={1\over
3^{(y+1)/3}}\,\sum_{k=0}^{\infty}a_k\;\sim\;{1\over\sqrt{6}}\,
\left({2e\over y}\right)^{y/3}\;\sim\;I_0(y)\,.\eqno(26.11)$$

The proof of (26.2) will therefore be complete if we can show that
$I_{-}(y)/I_{+}(y)\ra 0$ as $y\ra\infty$. For this we can use (10.9)
to show that
$$A(y,-\alpha)\leq{1\over
2\pi}\,\alpha^{1/2-y}\int_{-\infty}^{\infty}\,
\left\vert1+{it\over\alpha^{3/2}}\right\vert^{1-y}
e^{-t^2\!/2}\,dt\leq{\alpha^{1/2-y}\over\sqrt{2\pi}}\,;$$
therefore the first portion of $I_{-}(y)$ is quite small,
$$\int_{-\infty}^{-y^{1/3}}A(y,\mu)\,d\mu\leq{1\over
\sqrt{2\pi}}\,\int_{y^{1/3}}^{\infty}
\alpha^{1/2-y}\,d\alpha=O(y^{-1/2-y/3})\,.$$
On the other hand when $-y^{1/3}\leq\mu\leq 0$ we can integrate (10.7)
from $y^{1/3}-i\infty$ to $y^{1/3}+i\infty$, obtaining
$$\eqalign{A(y,\mu)
&\leq{1\over 2\pi}\,\int_{-\infty}^{\infty}\,\vert
y^{1/3}+it\vert^{1-y}\,
e^{\Re\,K(\mu,y^{1/3}+it)}\,dt\cr
\noalign{\smallskip}
&\leq{1\over
2\pi}\,y^{(1-y)/3}\,e^{K(\mu,y^{1/3})}\,
\int_{-\infty}^{\infty}\,\exp\left(-\left(y^{1/3}+{\mu\over
2}\right)t^2\right)\,dt\cr
\noalign{\smallskip}
&={y^{(1-y)/3}\exp\bigl(y/3+\mu(3y^{2/3}-\mu^2)/6\bigr)\over
\sqrt{2\pi\,(2y^{1/3}+\mu)}}\leq{y^{1/6}\over\sqrt{2\pi}}\,\left({e\over
y}\right)^{y/3}\,;\cr}$$
hence
$$\int_{-y^{1/3}}^0A(y,\mu)\,d\mu=O\left(y^{1/2}\left({e\over
y}\right)^{y/3}\right)$$ 
and $I_{-}(y)\ll I_{+}(y)$ as desired.\quad\pfbox

\medskip
Theorem 14 sheds further light on the results of [\FKP], where the
first cycle of a random multigraph was shown to have average length
asymptotic to $\sqrt{\pi/2}\,I(2)\,n^{1/6}$. According to a lengthy
numerical calculation sketched there, this coefficient was determined
to be 2.0337, correct to four decimal places. Sure enough, equation
(26.2) now confirms that the exact value is
$${\pi^{1/2}\,\Gamma(1/3)\over 2^{1/6}\,3^{2/3}}=
2.03369\;20140\,63898\,89186\,17247\,01028\,49830\,16693{-}\,.\eqno(26.12)$$

Section 7 of [\FKP] also proves implicitly that, if the random
variables~$L$ and~$S$ are respectively the length of the first cycle and
the size of the component containing that cycle, we have
$$\eqalignno{{\rm E}_n\,L^k\;
&\sim\;\sqrt{\pi\over 2}\,k!\,I(k+1)\,n^{k/3-1/6}\,;&(26.13)\cr
\noalign{\smallskip}
{\rm E}_n\,S^k\;
&\sim\;2^{k-1/2}\,\Gamma\bigl(k+{\textstyle{1\over
2}}\bigr)\,I(2k+1)\,n^{2k/3-1/6}&(26.14)\cr}$$
In particular, the variance of $L$ is asymptotically $\sqrt{2\pi n}$;
the asymptotic mean and variance of $S$ are $\sqrt{\pi n/2}$ and
$Kn^{7/6}$, where $K$ is the constant in (26.12). For graphs instead
of multigraphs, these coefficients should all be multiplied
by~$e^{3/4}$.

Notice that $I(3)=1$. Hence the function $A(3,\mu)$, which is
expressible in terms of Airy series or Bessel functions $\bigl($see
(10.32)$\bigr)$, defines a probability density.

Let $V_y$ be a random variable with density function $A(y,\mu)/I(y)$,
when $y>{3\over 2}$. Then, by (10.22),
$${\rm E}\,V_y=\int_{-\infty}^{\infty}\,{\mu\,A(y,\mu)\,d\mu\over
I(y)}={y\,I(y+2)-I(y-1)\over I(y)}={(y-3)I(y-1)\over
(y-2)I(y)}\,,\eqno(26.15)$$
if $y>{5\over 2}$. In particular, the variable $U_r$ of (25.28), which
is $V_{3r+3/2}$, has the mean value
$${(3r-3/2)I(3r+1/2)\over (3r-1/2)I(3r+3/2)}=\left({3\over
4}\right)^{1/3}\,{\Gamma(2r-2/3)\over\Gamma(2r-1)}\,.\eqno(26.16)$$
This is the limit as $n\ra\infty$ of E$\,U_{r,n}$, which represents
the mean waiting time for a graph or multigraph to reach excess~$r$.
The values are 0.8113, 1.2621, 1.5191, 1.7104, 1.8666, 2.0002, 2.1181,
2.2241, 2.3209, 2.4102 when $1\leq r\leq 10$.

Similarly, (10.23) implies that
$${\rm E}\,V_y^2={I(y-2)\over I(y)}={(y-2)\over
6^{1/3}}\;{\Gamma\bigl((2y-7)/3\bigr)\over \Gamma(2y/3-2)}\,,\quad
y>{\textstyle{7\over 2}}\,.\eqno(26.17)$$
Hence E$\,V_y=(y/2)^{1/3}\bigl(1-{7\over 6}y^{-1}+O(y^{-2})\bigr)$, 
E$\,V_y^2=(y/2)^{2/3}\bigl(1-{2\over 3}y^{-1}+O(y^{-2})\bigr)$, and we
have
$$\hbox{Var}\,V_y={5\over
2^{2/3}\,3}\,y^{-1/3}+O(y^{-4/3})\,.\eqno(26.18)$$

Let us now set $\mu=(y/2)^{1/3}+\sigma z$, where
$$\sigma^2={5\over 2^{2/3}\,3}\,y^{-1/3}\,.\eqno(26.19)$$
An argument similar to the derivation of (26.11) proves that
$${A(y,\mu)\over
I(y)}\;\sim\;{1\over\sqrt{2\pi\sigma^2}}\,e^{-z^2\!/2}\,,
\qquad z=O(1)\,,\quad y\ra\infty\,.\eqno(26.20)$$
Therefore $(\hbox{Var}\,V_y)^{-1/2}(V_y-{\rm E}\,V_y)$ approaches the
normal distribution $N(0,1)$ as $y\ra\infty$. In particular, this
establishes a kind of asymptotic normality of~$U_{r,n}$
(and~$M_{r,n}$), if we first let $n\ra\infty$ and then $r\ra\infty$.

\bigbreak\noindent
{\bf 27. Proof of the top-line conjecture.}\enspace
We are almost ready to settle the conjecture that was made in
section~18, but first we should carry out the promised refinement of
our estimates (23.5) and (23.9) for the sizes of the acyclic and
unicyclic parts of a random multigraph.

The first step is to consider the quantity (23.3), when $m={1\over
2}n(1+\mu)$ and $k=\kappa n$. If $k\geq m$ or $k\geq n$, expression
(23.3) is zero; otherwise $0\leq\kappa <\min\bigl({1+\mu\over
2},1\bigr)$, and Stirling's approximation yields
$$e^k\,\sqrt{m{-}k\over m}\,\sqrt{n{-}k\over
n}\,{2^k\,m^{\underline{k}}\,n^{\underline{k}}\over
(n{-}k)^{2k}}\left(n{-}k\over n\right)^{2m}
=\exp\!\left(nf(\kappa,\mu)+O\left(1\over m{-}k\right)+
O\left(1\over n{-}k\right)\right),\eqno(27.1)$$
where
$$f(\kappa,\mu)={1+\mu\over 2}\,\ln(1+\mu)-{(1+\mu-2k)\over
2}\,\ln(1+\mu-2\kappa)+(\mu-\kappa)\,\ln(1-\kappa)-\kappa\,.\eqno(27.2)$$
Notice that
$$\eqalign{{\partial f(\kappa,\mu)\over
\partial\kappa}&=\ln(1+\mu-2\kappa)-\ln(1-\kappa)
+{\kappa-\mu\over 1-\kappa}\,,\cr
\noalign{\smallskip}
{\partial^2\!f(\kappa,\mu)\over \partial\kappa^2}&={(1-\mu)(\mu-\kappa)\over
(1+\mu-2\kappa)(1-\kappa)^2}\,,\cr}$$
so both first and second derivatives vanish when $\kappa=\mu$. The
first derivative is $\leq 0$ when $\kappa=0$; if $0<\mu<1$ it
increases to zero when $\kappa=\mu$, then becomes negative; if
$\mu\leq 0$ or $\mu\geq 1$ it decreases steadily. Thus $f(\kappa,\mu)$
is a decreasing function of~$\kappa$, as claimed in section~23.

We also have
$${\partial f(\kappa,\mu)\over \partial\mu}={\ln(1+\mu)\over
2}-{\ln(1+\mu-2\kappa)\over 2}+\ln(1-\kappa)\,;$$
this derivative decreases steadily, passing through zero when
$\mu=\kappa/(2-\kappa)$. Therefore we have
$$f(\kappa,\mu)\leq f\left(\kappa,\,{\kappa\over 2-\kappa}\right)=
-\left({\kappa^3\over 24}+{\kappa^4\over 24}+{11\kappa^5\over
320}+\cdots +{1-2j/2^j\over
j(j-1)}\,\kappa^j+\cdots\,\right)\,,\eqno(27.3)$$
for all $\mu>2\kappa-1$. In particular, we can conclude that terms like
(23.1) and (23.8) are superpolynomially small for all $k\geq
n^{2/3+\epsilon}$, since they are
$O\bigl(\exp(-n^{3\epsilon}/24)\bigr)$ when $k=n^{2/3+\epsilon}$.

Our next goal is to estimate the sum of (23.8) for $k\geq 1$ when
$\mu\geq n^{-1/3}$. This sum $V(m,n)$ is the expected number of
vertices in unicyclic components after $m$~steps of the multigraph
process.
The formulas above allow us to write
$$\eqalignno{V(m,n)
&=\sum_{k\leq n^{2/3+\epsilon}}\,\half \;{k^k Q(k)\over k!\;e^k}\,
\sqrt{{m\over m-k}}\,\sqrt{{n\over n-k}}\,e^{nf(k/n,\mu)+O(n^{-1})}
+O(e^{-n^\epsilon})\cr
\noalign{\smallskip}
&=\sum_{k\leq n^{2/3+\epsilon}}\,\half \;{k^kQ(k)\over k!\;e^k}\,\exp
\left(k\bigl(\ln(1+\mu)-\mu\bigr)+{\mu k^2\over 2n}-{k^3\over 6n^2}\right.\cr
&\hskip10em\left.+O
\left({\mu^2k^2\over n}+{k^4\over n^3}+{k\over
n}\right)\right)\,+\,O(e^{-n^\epsilon})\,.&(27.4)\cr}$$

Let $\mu=\alpha n^{-1/3}$, so that $\alpha$ is the quantity we called
$\mu$ in sections 10--20 above. We will assume that $\alpha\geq 1$,
and also that $\alpha\leq cn^{1/3}$ (hence $\mu\leq c$), where $c$ is
a sufficiently small constant. The terms of $V(m,n)$ are negligible
for $k\geq n^{2/3+\epsilon}$, regardless of the value of~$\mu$; and
when $n^{-1/3}\leq\mu\leq c$ we can in fact ignore all terms for
$k>\alpha^{\epsilon}\mu^{-2}$. The reason is that
$$\eqalign{k\bigl(\ln(1+\mu)-\mu\bigr)+{\mu\,k^2\over 2n}-{k^3\over
6n^2}
&={-\mu^2k\over 2}\left({1\over 4}+{1\over 3}\left({3\over
2}-{k\over\mu n}\right)^2\right)+O(k\mu^3)\cr
\noalign{\smallskip}
&\leq{-\mu^2k\over 8}\,\bigl(1+O(\mu)\bigr)\leq{-\mu^2k\over
100}\cr}$$
if we choose $c$ small enough. The sum of
$O\bigl(e^{-\mu^2k/100}\bigr)$ for $\alpha^{\epsilon}/\mu^2<k<\infty$
is then $O\bigl(\mu^{-2}e^{-\alpha^{\epsilon}/100}\bigr)$, which is
dominated by the error bounds we will encounter below.

When $k\leq \alpha^{\epsilon}\mu^{-2}$, we have $\mu
k^2\!/2n\leq\alpha^{2\epsilon-3}\!/2\leq 1/2$. Therefore we are
justified in moving terms out of the exponent in (27.4):
$$\interdisplaylinepenalty=10000
\eqalignno{\kern-2em V(m,n)&=\sum_{k\geq 1}\,{1\over
2}{k^kQ(k)\over k!}\;
\bigl((1+\mu)e^{-(1+\mu)}\bigr)^k\left(1+{\mu k^2\over 2n}-{k^3\over
6n^2}\right.\cr
\noalign{\smallskip}
&\hskip10em\left.+
O\left({\mu k^2\over 2n}-{k^3\over 6n^2}\right)^2
+O\left({\mu^2k^2\over n}+{k^4\over n^3}+{k\over n}\right)\right)\cr
\noalign{\smallskip}
&=\sum_{k\geq 1}\,\half {k^k Q(k)\over k!}\,
\bigl((1-\sigma)e^{-(1-\sigma)}\bigr)^k
\left(1+{\mu k^2\over 2n}+O(\alpha^{4\epsilon-6})
+O({\alpha^{2\epsilon-2}\over n^{1/3}})\right).&(27.5)\cr}$$
Here $\sigma$ is the ``shadow'' of $\mu$ as in (23.4) and (23.10), and
the error bounds are computed under the assumption
$k\leq\alpha^{\epsilon}\!/\mu^2$. The trick of (23.5) and (23.9) now applies,
using (23.7), and we have
$$\eqalignno{\kern-2em V(m,n)
&={1\over
2}\,\left(1+O(\alpha^{4\epsilon-6})+O({\alpha^{2\epsilon-2}\over n^{1/3}})
+{\mu\over
2n}\vartheta^2\right){T\bigl((1-\sigma)e^{-(1-\sigma)}\bigr)\over
\bigl(1-T\bigl((1-\sigma)e^{-(1-\sigma)}\bigr)\bigr)^2}\cr
\noalign{\smallskip}
&={1\over
2}\,\left((1+O(\alpha^{4\epsilon-6})+O({\alpha^{2\epsilon-2}\over
n^{1/3}})\bigr)
{1{-}\sigma\over\sigma^2}+{\mu\over 2n}\!
\left({8{-}17\sigma{+}11\sigma^2{-}2\sigma^3\over\sigma^6}\right)\right)
\,.&(27.6)\cr}$$
If we had expanded the summand further, we would have obtained still
more accuracy; therefore we are allowed to set $\epsilon=0$ in (27.6).
The term $O(\alpha^{-2}n^{-1/3})$ dominates $O(\alpha^{-6})$ when
$\alpha\geq n^{1/12}$; it comes from both $O(\mu^2k^2\!/n)$ and
$O(k/n)$ in (27.5).

We are assuming that $\mu$ is small, hence
$\sigma=\mu\bigl(1+O(\mu)\bigr)$. Thus (27.6) can be simplified to
$$V(m,n)=\left({1\over
2\alpha^2}+{2\over\alpha^5}+O\left({1\over\alpha^8}\right) 
\right)\,n^{2/3}\bigl(1+O(\mu)\bigr)\,,$$
and with an extension of the same approach we obtain an asymptotic
expansion that begins
$$V(m,n)=\left({1\over
2\alpha^2}+{2\over\alpha^5}+{20\over\alpha^8}+{320\over\alpha^{11}} 
+{7040\over\alpha^{14}}
+O\left({1\over\alpha^{17}}\right)
\right)n^{2/3}
\bigl(1+O(\mu)\bigr)\,.\eqno(27.7)$$
This expansion is readily computed if we note that
$$\vartheta^k\,{T(z)\over\bigl(1-T(z)\bigr)^2}={2^k\,k!\over
\bigl(1-T(z)\bigr)^{2k+2}}+\cdots\,,\eqno(27.8)$$
where the remaining terms $a_{k1}/\bigl(1-T(z)\bigr)^{2k+1}
+a_{k2}/\bigl(1-T(z)\bigr)^{2k}+\cdots$ are negligible
when we replace $T(z)$ by $1-\mu-O(\mu^2)$. The asymptotic series in
(27.7) is obtained also from the integral
$${\textstyle{1\over 4}}\int_0^{\infty}e^{-\alpha^2t/2+\alpha
t^2\!/2-t^3\!/6} \,dt={\textstyle{1\over 4}}\int_0^{\infty}
e^{(\alpha-t)^3\!/6-\alpha^3\!/6}\,dt\,,\eqno(27.9)$$
because we can expand $e^{\alpha t^2\!/2-t^3\!/6}$ into powers of $t$
and use the formula
$$\int_0^{\infty}e^{-\alpha^2t/2}t^k\,dt={2^{k+1}\,k!\over
\alpha^{2k+2}}\,,\eqno(27.10)$$
which matches (27.8).
The coefficients of (27.7) follow a simple pattern; for example,
$7040=22\cdot 16\cdot 10\cdot 4\,/\,2$. Thus we are led to
conjecture the asymptotic series
$$\int_0^{\infty}e^{(\alpha-t)^3\!/6-\alpha^3\!/6}\,dt\;\sim\;
2F\bigl({\textstyle{2\over 3}},1;;6/\alpha^3\bigr)/\alpha^2\quad{\rm
as}\quad \alpha\ra\infty\,;\eqno(27.11)$$
the right-hand side here is a formal power series that diverges for
all finite~$\alpha$.
And indeed, this conjecture is true, as we will see momentarily.

A similar calculation allows us to estimate $U(m,n)$, the number of
vertices in trees. The analog of (27.5) is
$$U(m,n)={n\over 1+\mu}\,\sum_{k\geq 1}{k^{k-1}\over
k!}\,\bigl((1-\sigma)e^{-(1-\sigma)}\bigr)^k
\left(1+{\mu k^2\over
2n}+O(k\alpha^{3\epsilon-4}n^{-2/3})+O(k\alpha^\epsilon n^{-1})\right)
;\eqno(27.5^{\prime})$$
we leave a factor of $k$ in the $O$ terms because it will lead to a
better final estimate. Then the analogs of (27.6)--(27.10) are
$$\eqalignno{U(m,n)
&={n\over
1+\mu}\!\left(1-\sigma+O(\alpha^{3\epsilon-5}n^{-1/3})+O(\alpha^{\epsilon-1}n^{-2/3})
+{\mu\over
2n}\!\left({1{-}\sigma\over\sigma^3}\right)\right);&(27.6^{\prime})\cr
\noalign{\medskip}
U(m,n)
&=n+\left(-2\alpha+{1\over 2\alpha^2}+{11\over 8\alpha^5}+{175\over
16\alpha^8}+{19005\over 128\alpha^{11}}+{735735\over
256\alpha^{14}}\right.\cr
\noalign{\smallskip}
&\hskip6em\null+O\left(\left.{1\over\alpha^{17}}\right)
\right)n^{2/3}\bigl(1+O(\mu)\bigr)\,;&(27.7^{\prime})\cr
\noalign{\smallskip}
\vartheta^kT(z)&={2^{k-1}\,\Gamma(k-1/2)\over\Gamma(1/2)}\;
{1\over\bigl(1-T(z)\bigr)^{2k-1}}+\cdots\;,\quad k\geq
1\,;&(27.8^{\prime})\cr}$$
\vskip-8pt
$$\eqalignno{
{1\over\sqrt{2\pi}}\int_0^{\infty}{e^{-\alpha^2t/2}(e^{\alpha
t^2\!/2-t^3\!/6}{-}1)\,dt\over t^{3/2}}
&={1\over\sqrt{2\pi}}\int_0^{\infty}
{e^{(\alpha-t)^3\!/6-\alpha^3\!/6}{-}1\over
t^{3/2}}\,dt+\alpha\,;&(27.9^{\prime})\cr
\noalign{\smallskip}
{1\over\sqrt{2\pi}}\,\int_0^{\infty}\,e^{-\alpha^2t/2}t^{k-3/2}\,dt
&={2^{k-1}\,\Gamma(k-1/2)\over\sqrt{\pi}\;\alpha^{2k-1}}\,,\quad k\geq
1\,. &(27.10^{\prime})\cr}$$
The asymptotic series (27.7) and (27.7$'$) for $\alpha\ra\infty$ blend
perfectly with the results obtained in [\LPW] when $\alpha$ is any
constant (positive, negative, or zero):
$$\eqalignno{V(m,n)
&={{1\over
4}}\left(\int_0^{\infty}e^{(\alpha-t)^3\!/6-\alpha^3\!/6}\,dt\right)
\,n^{2/3}+O(n^{1/3})\,;&(27.12)\cr
\noalign{\smallskip}
U(m,n)
&=n+\left(-\alpha+{1\over\sqrt{2\pi}}\,\int_0^{\infty}\,
{e^{(\alpha-t)^3\!/6-\alpha^3\!/6}-1\over t^{3/2}}\;dt\right)\,
n^{2/3}+O(n^{1/3})\,.&(27.12^{\prime})\cr}$$
These integrals are entire functions of $\alpha$,
$$\eqalignno{
\int_0^{\infty}\,e^{(\alpha-t)^3\!/6-\alpha^3\!/6}\,dt
&={6^{1/3}\,\Gamma(1/3)\over 3}\,e^{-\alpha^3\!/6}+\alpha\,F
(1;\,{\textstyle{4\over 3}}\,;-\alpha^3\!/6)\,;&(27.13)\cr
\noalign{\medskip}
\int_0^{\infty}\,{e^{(\alpha-t)^3\!/6-\alpha^3\!/6}-1\over
t^{3/2}}\;dt
&=-e^{-\alpha^3\!/6}\bigl(\bigl(6^{5/6}\,\Gamma({\textstyle{5\over
6}})/3\bigr)\,F({\textstyle\half },\,{\textstyle{5\over 6}};
\,{\textstyle{1\over 3}},\,{\textstyle{2\over 3}};
\,\alpha^3\!/6)\cr
\noalign{\smallskip}
&\qquad\null-\bigl(6^{1/2}\,\Gamma({\textstyle{1\over
2}})/6\bigr)\,\alpha\, 
F({\textstyle{5\over 6}},\,{\textstyle{7\over 6}};
\,{\textstyle{2\over 3}},\,{\textstyle{4\over 3}};\,\alpha^3\!/6)\cr
\noalign{\medskip}
&\qquad\null+\bigl(6^{1/6}\,\Gamma({\textstyle{1\over 6}})/8\bigr)\,
\alpha^2\,F({\textstyle{7\over 6}},\,{\textstyle{3\over 2}};
\,{\textstyle{4\over 3}},\,{\textstyle{5\over 3}};
\,\alpha^3\!/6)\bigr)\,.&(27.13^{\prime})\cr}$$
Equation (27.13) is proved by observing that if
$g(\alpha)=\int_0^{\infty} e^{(\alpha-t)^3\!/6}\,dt$, then
$g'(\alpha)=\int_0^{\infty} \,{(\alpha-t)^2\over
2}\,e^{(\alpha-t)^3\!/6}\,dt =e^{\alpha^3\!/6}$.
It implies (27.11) by well-known properties of confluent
hypergeometric series. Equation (27.13$'$) is proved by setting
$h(\alpha)=\int_0^{\infty}(e^{(\alpha-t)^3\!/6}-e^{\alpha^3\!/6})\,t^{-3/2}\,dt
=-\int_0^{\infty}(\alpha-t)^2\,e^{(\alpha-t)^3\!/6}t^{-1/2}\,dt$ and
proving that $h'''(\alpha)=\half \,\alpha^2h''(\alpha)+{5\over
2}\,\alpha h'(\alpha)+{15\over 8}\,h(\alpha)$, hence
$[\alpha^{k+3}]\,h(\alpha) =[\alpha^k]\,h(\alpha)(k+{3\over
2})(k+{5\over 2})/\bigl(2(k+1)(k+2)(k+3)\bigr)$. Recall that we
enumerated $n-U(m,n)-V(m,n)$, the expected number of vertices in
complex components, using a complementary approach in (15.13), by
summing over the excess~$r$. 

\proclaim
Lemma 9.
Let $V_{mn}$ be the number of vertices in unicyclic components of a
random multigraph with $m$~edges and $n$~vertices. If $m={1\over
2}n(1+\mu)$ and $\mu\geq n^{-1/3}$, the expected value
of~$V^{l}_{mn}$ is $O(\mu^{-2l})$, for every fixed integer $l\geq
1$.

\proof
Equation (27.7) proves this for $l=1$ and $n^{-1/3}\leq\mu\leq c$,
where $c$ is some positive constant. A~similar argument applies for
arbitrary~$l$, because the generating function
$\vartheta^{l}e^{V(z)}$ is $e^{V(z)}/\bigl(1-T(z)\bigr)^{2l}$
times a polynomial in $T(z)$; this means we are summing terms like
(27.5), but with $Q(k)$ replaced by a semipolynomial in~$k$ of degree
$l-\half $. (See the proof of Theorem~3 in section~8.) The
analog of (27.6) will then be $O(\sigma^{-2l})$, which is
$O(\mu^{-2l})$ if $\mu\leq c$. Incidentally, for this range of $\mu$
we will have
$${\rm E}\,V_{mn}^{l}\;=\;{\Gamma(l+1/4)\over\Gamma(1/4)}\;
\left({2\over\sigma^2}\right)^{l}\bigl(1+O(\mu)+
O(\mu^{-3}n^{-1})\bigr)\,.\eqno(27.14)$$

If $c\le\mu\le n^\epsilon$, with $\epsilon<{1\over4}$, let
$0<\delta<1-\ln(1+c)/c$. Then each term in the analog of (27.4) with
$k\le n^{3/4}$ is $O\bigl(k^{l-1}\exp\bigl(k(\ln(1+\mu)-\mu)+
O(\mu^2k^2\!/n)\bigr)\bigr)=O\bigl(k^{l-1}\exp(-k\delta\mu)\bigr)
=O(k^{l-1}e^{-\delta\mu-\delta ck})$. Hence E$\,V^l_{mn}=O(e^{-\delta
\mu})$.

Finally, if $\mu\geq n^{\epsilon}$ the value of E$\,V_{mn}^{l}$ is
superpolynomially small, for it is a sum of $n$~terms each of which is
bounded by a polynomial in~$m$ and~$n$ times $(1-1/n)^{2m}$, which is
$O(\mu^de^{-\mu})$ for some finite degree~$d$.\quad\pfbox

\proclaim
Corollary.
The probability that a random multigraph never acquires a new complex
component after it has gained $m=\half (n+\alpha n^{2/3})>{1\over
2}n$ edges is $1-O(\alpha^{-3})$.

\proof
We may assume that $\alpha\geq 1$. A new complex component must be
bicyclic. A~multigraph gains a new bicyclic component if and only if
the endpoints of a new edge both fall in unicyclic components. The
probability that this occurs at time $m=\half (n+\mu n)$ is
E$\,V_{mn}^2\!/n^2=O(\mu^{-4}n^{-2})$, by the lemma. Summing for
$m\geq\half (n+\alpha n^{2/3})$ gives $O(\alpha^{-3})$ as an upper
bound on the probability that at least one new bicyclic component
appears after time $\half (n+\alpha n^{2/3})$.\quad\pfbox

\proclaim Theorem 15. The probability that an evolving graph or multigraph on
$n$~vertices never has more than one complex component throughout its
evolution approaches ${5\pi\over 18}$ as $n\rightarrow\infty$.

\proof
Let $\epsilon>0$ be fixed. By the corollary just proved, there exists
a number~$\alpha$, independent of~$n$, such that the probability of a
random multigraph obtaining a new complex component after time $m={1\over
2}(n+\alpha n^{2/3})$ is less than~$\epsilon$.

By section 14 and the corollary of section 13, there is a number~$R$,
independent of~$n$, such that the probability of having excess $>R$ at
this time~$m$ is less than~$\epsilon$. So the probability that a
random multigraph leaves the top line after excess~$R$ is $<2\epsilon$.
(Either it reaches excess~$R$ before time~$m$, or it leaves the top
line after time~$m$.)

But the probability that a random multigraph leaves the top line
before excess~$R$ is $1-{5\pi\over 18}+O(R^{-1})+O(n^{-1/3})$, by
(18.2). We may choose $R$ sufficiently large that this $O(R^{-1})$ is
less than~$\epsilon$; then we may choose $n$ sufficiently large that
the $O(n^{-1/3})$ is less than~$\epsilon$. The probability that a
random multigraph leaves the top line for such $n$ is therefore
between
$1-{5\pi\over 18}-2\epsilon$ and $1-{5\pi\over
18}+4\epsilon$.

For graphs, we note that an evolving graph may be constructed from an evolving
multigraph by ignoring all new edges that would be loops or parallel to
an existing edge. Since this reduction preserves or decreases both the
excess and the number of complex components, it follows that if the graph
leaves the top line after excess~$R$, then the multigraph does too.
Hence this event likewise has probability $<2\epsilon$, and the proof
is completed as for multigraphs.\quad\pfbox

\proclaim
Theorem 16.
Given any set $S$ of infinite paths in Figure~1, the probability that
the evolution of a random multigraph follows a path in~$S$ converges
as $n\ra\infty$ to the corresponding probability for the Markov chain
with the transition probabilities given in Theorem~9.
Similarly, if the evolution of a random graph, which stops at excess
${n\choose2}-n$ when the complete graph is reached, is continued along
the top line to an infinite path in Figure~1, then the probability that
this path lies in~$S$ converges to the same limit.

\proof
Given $\epsilon >0$, let $R$ be as in the preceding proof so that a random
graph or multigraph leaves the top line after excess~$R$ with probability
$<2\epsilon$. We can also choose $R$ large enough that
$c_R>e_R(1-\epsilon)$, by (8.7). Since $c_R/e_R$ is the sum of all
Markov transition probabilities for paths that intersect the top line
at excess~$R$, if we cut Figure~1 at excess~$R$, the Markov
probabilities for paths in~$S$ that do not have this property must sum
to less than~$\epsilon$. When $R$ is large enough, the sum of Markov
probabilities for all paths that diverge from the top line after
excess~$R$ is likewise less than~$\epsilon$, because it
is $O\bigl(\sum_R^\infty r^{-2}\bigr)=O(R^{-1})$.

Let $P_n(S)$ be the probability that the evolution of a random graph or
multigraph on $n$ vertices follows a path in~$S$, and let $P_\infty(S)$
denote the corresponding Markov probability.
If $S_R$ is the subset of~$S$ having all paths on the top line
when the excess is $\geq R$, then $0\le P_n(S)-P_n(S_R)<2\epsilon$ for
all $n\le\infty$. Similarly, if $S'_R$ is the set of all paths that
follow a path in $S_R$ up to excess~$R$, but afterwards are arbitrary,
then $0\le P_n(S'_R)-P_n(S_R)<2\epsilon$, for $n\le\infty$.
Finally, by Theorem~10, $\bigl\vert P_n(S'_R)-P_\infty(S'_R)\bigr\vert
<\epsilon$ if $n$ is large enough, and we have
$\bigl\vert P_n(S)-P_\infty(S)\bigr\vert<5\epsilon$.\quad\pfbox

\medskip
Theorem 16 says that the evolutionary path, regarded as a random
element of the set of
all paths in Figure~1, converges in distribution to the Markov
process. There are uncountably many paths, but the theorem needs no
measurability restriction since the distributions for finite~$n$ and
for the limit are concentrated on the countable set of paths that
eventually follow the top line.
Note that we cannot strengthen the statement for random graphs
to deduce the limiting probability that the evolution follows a path
in~$S$ until it stops at excess~${n\choose2}-n$; for example, if
$S$ is the set of all paths that do {\it not\/} eventually follow
the top line, the Markov probability $P_\infty(S)$ is zero, while
$P_n(S)=1$ for all finite~$n$.

\proclaim
Corollary.
The probability that an evolving graph or multigraph never has more than $l$
complex components converges to a limit~$P_{@l}$.\quad\pfbox

Closed form expressions for $P_{@l}$ might not exist when $l\ge2$, but
the values can be estimated from below using the following related
probabilities:

\proclaim
Corollary.
The probability that an evolving graph or multigraph acquires exactly
$l\ge1$ new complex components during the evolution converges to
$$p'_l=\Pr\biggl(\sum_{r=0}^\infty I_r=l\biggr)
      =\Pr\biggl(\sum_{r=1}^\infty I_r=l-1\biggr)\,,\eqno(27.15)$$
where $I_0$, $I_1$, $I_2$, $I_3$, \dots\ are independent Bernoulli distributed
random variables with $\Pr({I_r=1})=1-\Pr({I_r=0})=5/(6r+1)(6r+5)$.

\vskip-3pt
\noindent In other words, the number of new complex components converges
in distribution to $\sum_{r=0}^\infty I_r$.

\medskip
\proof
Let $I_r=1$ if the Markov process acquires a new bicyclic component
when the excess goes from $r$ to $r+1$, and $I_r=0$ otherwise;
in particular $I_0=1$ always. By Theorem~9, $\Pr(I_r=1)=5/(6r+1)(6r+5)$
independently of the previous history, and thus the variables are
independent.\quad\pfbox

\medskip
The probabilities $p'_l$ have a surprisingly simple generating function:
We have
$$\eqalignno{p'_l&=[z^l]\,\prod_{r=0}^\infty\left(1+(z-1){5\over(6r+1)
(6r+5)}\right)\cr
\noalign{\medskip}
&=[z^l]\,\prod_{r=0}^\infty{(r+\half +{1\over6}\sqrt{@9-5z}\,)
(r+\half -{1\over6}\sqrt{@9-5z}\,)\over(r+{1\over6})(r+{5\over6})}\cr
\noalign{\medskip}
&=[z^l]\,{\Gamma\bigl({1\over6}\bigr)\,\Gamma\bigl({5\over6}\bigr)\over
\Gamma\bigl(\half +{1\over6}\sqrt{@9-5z}\,\bigr)\,
\Gamma\bigl(\half -{1\over6}\sqrt{@9-5z}\,\bigr)}\cr
&=[z^l]\,\cos\left({\pi\over6}\sqrt{@9-5z}\,\right)\bigg/
\cos{\pi\over3}\,.&(27.16)\cr}$$
Computing the coefficients of the Taylor series for
$\cos\bigl({\pi\over6}\sqrt{@9-5z}
\bigr)$, we find that the numbers $p'_l$ are rational polynomials
in $\pi$:
$$\openup2\jot
\eqalign{p'_1&={5\pi\over18}\approx0.87266\,;\cr
p'_2&={50\pi\over6^4}\approx0.12120\,;\cr
p'_3&={500\pi\over6^6}\left(1-{\pi^2\over12}\right)\approx0.00598\,;\cr
p'_4&={6250\pi\over6^8}\left(1-{\pi^2\over10}\right)\approx0.00015\,.\cr}$$
Let $P'_l=\sum_{j=1}^l p'_j$; numerically we have
$P_2>P'_2\approx0.99387$, $P_3>P'_3\approx0.99985$,
$P_4>P'_4\approx0.999998$.

The number of new complex components is also studied in~[\Jan], where
further results are given. The methods of~[\Jan] do not, however, seem
to yield the sharp results obtainable with generating functions.

\bigbreak\noindent
{\bf 28. Empirical data.}\enspace
Computer simulations of random multigraphs tend to confirm the theoretical
results derived above, although there are a few surprises apparently
due to the slow convergence of some asymptotic formulas. In this
section we will discuss some of the statistics computed during 1000
trials of the multigraph process on 20,000 vertices, so that readers can
obtain a feel for the way in which random multigraphs actually evolve
in practice. The data was divided into two groups of 500 runs each,
and both groups exhibited essentially the same behavior; therefore the
full set of 1000 runs is being treated as a unit~here.

When a statistic is given in the form `$x\pm y$' below, $x$ is the
sample mean and $y$ is the sample standard deviation divided by $\sqrt{1000}$.
The sample standard deviation has been computed by taking the square root
of an unbiased estimate of the variance. The ``time'' of an event is the
number of edges present when that event occurred.

The first cycle was formed at time $6769\pm96$; this agrees reasonably
well with the asymptotic formula $n/3$ found in [\FKP, Corollary~3].
The size of the first unicyclic component was $188\pm14$. 
According to (26.14), the mean should be approximately $\sqrt{\pi
n/2}\approx 177$.

The length of the first cycle was $3.9\pm0.1$; in fact, the histogram was
$$
\vcenter{\halign{\hfil#\ =&&\quad\hfil$#$\hfil\cr
length&1&2&3&4&5&6&7&\ge8\cr
actual&321&132&89&88&78&86&60&146\cr
theoretical&333&133&76&51&37&28&23&318\cr}}$$
The distribution has infinite mean, approximately $2.03 n^{1/6}+O(n^{3/22})$,
and its standard deviation is of order~$n^{1/4}$ by (26.13), so the
length of the first cycle should not be expected to be a robust statistic.
However, the marked deviation in the histogram for cycle lengths $\ge4$
was unexpected. Apparently $n$ must become quite large before the
asymptotic probability of first cycle length~$k$ will assert itself.

Several people have suggested in conversation that the ``last cycle''
ought to have the same statistical characteristics as the first. The
last cycle is the last unicyclic component that is present during a multigraph's
evolution: After it is absorbed into a component of higher complexity,
no further unicycles exist, and no further unicycles
are formed. (If two cycles disappear simultaneously
when the edge $\langle x,y\rangle$ is added, we say that the cycle
containing~$y$ was the last to go.) The manner in which the giant component
swallows other structures is rather like the initial stages of evolution
but in reverse: First the unicycles tend to go, then the larger trees,
and finally only isolated vertices are left (see Bollob\'as [\Biii, sections
VI.3 and VII.1]). A strong formulation of this symmetry principle was
proved by \L uczak~[\Li]; the phenomenon can be explained by the symmetry
between $T(z)$ and $2-T(z)$ in $U(z)$. However, the length of the last
cycle has a distinctly different distribution from the length
of the first cycle (see [\JL]). In these computer runs it had
the following histogram:
$$
\vcenter{\halign{\hfil#\ =&&\quad \hfil$#$\hfil\cr
length&1&2&3&4&5&6&7&\ge8\cr
observed&423&144&107&79&63&62&40&82\cr}}$$
with mean $3.1\pm0.1$.

The total number of unicyclic components formed during
the entire evolution was
$$
\vcenter{\halign{\hfil#\ =&&\quad \hfil$#$\hfil\cr
number&1&2&3&4&5&6&7&\ge8\cr
observed&53&148&221&219&178&98&44&39\cr}}$$
with mean $4.0\pm0.1$.

The excess of the multigraph changed from 0 to 1 at time $10331\pm13$. The
number of unicyclic components present was about 2.7 just before this event,
and about 1.5 just after. As soon as the excess became positive it began a steady rise:
$$\vcenter{\halign{\hfil#&&\quad\hfil#&${}\pm#$\hfil\cr
&&\omit&\multispan2\quad\hfil unicyclic size\hfil&
\multispan2\quad\hfil unicyclic size\hfil&
\multispan2\quad\hfil complex size\hfil&
\multispan2\quad\hfil complex size\hfil\cr
excess&\multispan2\quad\hfil time\hfil&
\multispan2\quad\hfil just before\hfil&
\multispan2\quad\hfil just after\hfil&
\multispan2\quad\hfil just before\hfil&
\multispan2\quad\hfil just after\hfil\cr
\noalign{\vskip2pt}
1&10331&13&1606&22&163&9&&\omit0\hfil&1442&21\cr
2&10501&10&265&14&132&7&1779&22&1912&20\cr
3&10603&8&168&9&111&7&2166&19&2222&19\cr
4&10675&8&132&8&90&5&2433&18&2475&17\cr
5&10738&8&105&6&85&5&2659&17&2680&17\cr
6&10789&7&95&6&76&5&2825&17&2844&16\cr
7&10835&7&83&5&69&4&2980&16&2994&16\cr
8&10880&7&77&5&66&4&3126&16&3137&16\cr
9&10920&7&72&5&62&4&3253&15&3263&15\cr
10&10955&7&66&4&58&4&3371&15&3379&15\cr}}$$
The value of $n^{2/3}$ is approximately
737 when $n=20000$, so each additional edge increases
the parameter $\mu$ of Lemma~3 by approximately $0.0027$. The value of~$\mu$
when $m=10955$ is approximately 2.59; then ${2\over3}\mu^3+1+{5\over24}\mu^{-3}
+{15\over16}\mu^{-6}\approx12.6$, so the excess is not quite keeping up
with the expected value in Theorem~6. 
Similarly, formula (26.16) predicts that the excess will reach~1 when
$m\approx 10299$, and~10 when $m\approx 10888$; random multigraphs for
finite~$n$ seem to become complex a bit ``late.''
It is interesting to note that the
observed standard deviations kept decreasing as the excess increased,
while the discrepancy from (26.16) kept increasing.

\vfill\eject % temporary patch...

The random 
multigraphs followed paths in Figure~1 with the frequencies shown in
Figure~3. 
When the excess changed from 9 to~10, the transition was from a
single $C_9$ 
to $C_{10}$ in 977~cases, from $C_9$ to $(C_1,C_9)$ in 2~cases,
from $(C_1,C_8)$ to $C_{10}$ in 8~cases, and from $(C_1,C_8)$ to
$(C_1,C_9)$ 
in the remaining 13~cases. Altogether 897 of the 1000 random
multigraphs remained
 on the top line of Figure~1 throughout their evolution.

%\vfill\eject

\unitlength=1cm
\centerline{\beginpicture(13,10)(0,-5)
\put(0,0){\bmath\makebox(0,0){$[0]$}}
\put(3,0){\bmath\makebox(0,0){$[1]$}}
\put(6,1){\bmath\makebox(0,0){$[0,1]$}}
\put(6,-1){\bmath\makebox(0,0){$[2]$}}
\put(9,2){\bmath\makebox(0,0){$[0,0,1]$}}
\put(9,0){\bmath\makebox(0,0){$[1,1]$}}
\put(9,-2){\bmath\makebox(0,0){$[3]$}}
\put(12,4){\bmath\makebox(0,0){$[0,0,0,1]$}}
\put(12,2){\bmath\makebox(0,0){$[1,0,1]$}}
\put(12,0){\bmath\makebox(0,0){$[0,2]$}}
\put(12,-2){\bmath\makebox(0,0){$[2,1]$}}
\put(12,-4){\bmath\makebox(0,0){$[4]$}}
\def\\#1//{\makebox(0,0){#1}}%
\put(0.4,0){\line(1,0){.5}}
\put(1.4,0){\\1000//}
\put(1.9,0){\vector(1,0){.8}}
\put(3.3,.1){\line(3,1){.85}}
\put(4.55,.53){\\962//}
\put(4.95,.7){\vector(3,1){.55}}
\put(3.3,-.1){\line(3,-1){.9}}
\put(4.5,-.6){\\38//}
\put(4.8,-.6){\vector(3,-1){.8}}
\put(6.6,1.2){\line(3,1){.5}}
\put(7.5,1.55){\\948//}
\put(7.85,1.667){\vector(3,1){.45}}
\put(6.35,.7){\line(3,-1){.25}}
\put(6.9,0.55){\\14//}
\put(7.2,.467){\vector(3,-1){1.1}}
\put(6.3,-.95){\line(3,4){.4}}
\put(6.85,-.2){\\14//}
\put(7.,-.02){\vector(3,4){1.15}}
\put(6.4,-1.016){\line(3,1){.5}}
\put(7.3,-.75){\\22//}
\put(7.6,-.616){\vector(3,1){.9}}
\put(6.4,-1.233){\line(3,-1){.4}}
\put(7.2,-1.5){\\2//}
\put(7.5,-1.6){\vector(3,-1){1.1}}
\put(9.6,2.4){\line(3,2){.4}}
\put(10.33,2.9){\\953//}
\put(10.6,3.033){\vector(3,2){.8}}
\put(9.65,2){\line(1,0){.2}}
\put(10.1,2){\\9//}
\put(10.4,2){\vector(1,0){.9}}
\put(9.6,.6){\line(3,4){.4}}
\put(10.2,1.4){\\9//}
\put(10.35,1.6){\vector(3,4){1.45}}
\put(9.6,.3){\line(3,2){.5}}
\put(10.475,0.85){\\21//}
\put(10.8,1.1){\vector(3,2){.7}}
\put(9.475,0){\line(1,0){.2}}
\put(10.0,0.07){\\6//}
\put(10.3,0){\vector(1,0){1.1}}
\put(9.06,-.24){\line(3,-2){.34}}
\put(9.7,-.65){\\0//}
\put(9.883,-.783){\vector(3,-2){1.5}}
\put(9.3,-1.8){\line(3,4){1.}}
\put(10.45,-0.3){\\2//}
\put(10.575,-.1){\vector(3,4){1.3}}
\put(9.3,-2){\line(1,0){.4}}
\put(10.05,-2){\\0//}
\put(10.4,-2){\vector(1,0){1.0}}
\put(9.15,-2.3){\line(3,-2){.55}}
\put(10,-2.9){\\0//}
\put(10.2,-3.0){\vector(3,-2){1.4}}
\put(12.9,4.3){\vector(3,2){.5}} % corrected
\put(12.9,3.8){\vector(3,-1){.5}}
\put(12.7,2.3){\vector(1,1){.6}} % corrected
\put(12.7,2.0){\vector(3,1){.6}} % corrected
\put(12.7,1.9){\vector(3,-2){.6}}
\put(12.65,1.8){\vector(3,-4){.65}}
\put(12.6,.15){\vector(3,1){.7}}
\put(12.6,0){\vector(1,0){.7}}
\put(12.6,-.15){\vector(3,-1){.7}}
\put(12.5,-1.84){\vector(1,1){.8}}
\put(12.5,-1.93){\vector(2,1){.8}}
\put(12.5,-2.0){\vector(1,0){.8}}
\put(12.5,-2.1){\vector(2,-1){.8}}
\put(12.5,-2.2){\vector(1,-1){.8}}
\put(12.3,-3.9){\vector(2,1){1.0}}
\put(12.3,-4.0){\vector(1,0){1.0}}
\put(12.3,-4.2){\vector(3,-2){1.0}} % corrected
\endpicture
}
\bigskip
{\narrower\narrower\smallskip\noindent
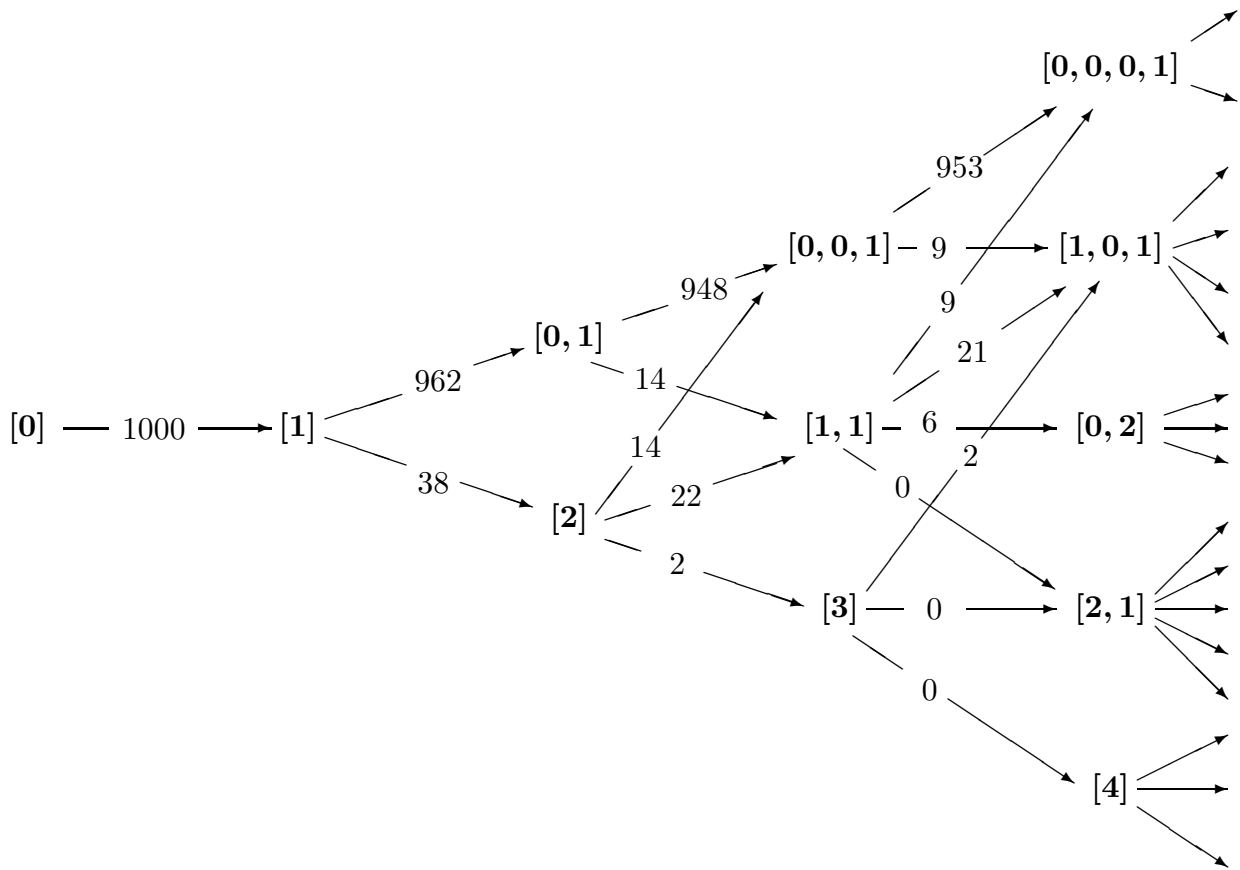
{\bf Figure 3.}\enspace The number of times the paths in Figure~1 were
actually traced,
 when 1000 random multigraphs on 20000 vertices were generated in
experimental tests.
\bigskip\bigskip}

\goodbreak
There comes a time
 when the giant component first succeeds in annihilating
everything except isolated vertices, after which it remains the only
component with edges. In these runs that time was $58352\pm224$.
The number of isolated vertices still remaining was then $71\pm1$.

The multigraph finally became connected at time $105294\pm404$. The
expected time for an evolving multigraph to have no isolated vertices
is $\half nH_n=\half n\ln n+{1\over
 2}\gamma n+{1\over4}+O(n^{-1})$,
which is approximately 104807 when $n=20000$.

\bigbreak\noindent
{\bf 29. Open problems.}\enspace
The topics discussed 
in this paper raise a host of interesting questions, and the
answers to those questions will no doubt bring additional striking
patterns to light.

But the reader may have noticed that this paper is already rather long.
Therefore it seems
 wise to stop at this point, with the hope that researchers
all over the world
 will enjoy exploring the tantalizing questions that remain.

For example, it would be interesting to
find a basis for as many linear combinations of terms $w^r T^a\!/(1-T)^b$ as
possible such that
$$\Phi_n\,w^r\,e^U\,T^a\!/(1-T)^b$$
has a known value, as in (25.22). We can find many linear combinations of
such functions for which $\Phi_n$ gives~0, because $\Phi_n\nabla S$ is
usually 0 or~1. Notice that
$${T^a\over(1-T)^{b+1}}={T^a\over(1-T)^b}+
{T^{a+1}\over(1-T)^{b+1}}\,;\eqno(29.1)$$
hence terms of excess $r+1$ can be expressed as combinations of terms of
excess~$r$. Conversely, we can go from excess~$r$ to excess~$r+1$, because
$${T^a\over(1-T)^b}={T^a\over(1-T)^{b+1}}-{T^{a+1}\over(1-T)^{b+2}}+
{T^{a+2}\over(1-T)^{b+3}}-\cdots\eqno(29.2)$$
is an infinite series that always ``converges'' under application of~$\Phi_n$;
all terms after a certain point are multiples of $T^{n+1}$, so they do
not change the coefficient of~$z^n$.

The stopping configuration machinery suggests many further problems of
interest. For example, we should be able to deduce more about the nature
of a random multigraph when its deficiency first exceeds a given number~$d$.

The discussion in section 23
 characterizes the stochastic behavior of $r$ and~$d$
when $\mu=o(1)$; what happens 
thereafter? Relations (23.12) and (23.13) may well continue to
describe the approximate mean values of~$r$ and~$d$
 as $\mu\to\infty$. The shadow point~$\sigma$ defined in~(23.2) will
approach~0, but it remains an analytic function of~$\mu$, and $1-\sigma$
remains a saddle point of the contour integral for $[z^n]\,U^{n-m+r}
T^{2r-d}/(1-T)^{3r-d+1/2}$.

The analytic function $T(z)$ has an interesting Riemann surface: There is
a quadratic singularity at
 $z=e^{-1}$, and if we travel around that point we get
to a second sheet in which there is a logarithmic singularity at $z=0$.
Winding around that logarithmic singularity takes us to infinitely many other
sheets having no finite singularities besides~0. It may be possible to
work out a theory under which contour integrals of importance in the study
of random graphs could be evaluated by paths that pass through the point
$1+\mu$, which lies on the ``wrong side''
of the quadratic singularity of~$T(z)$; $1+\mu$ turns out
to be a saddle point for several important generating functions.

Identity (8.15)--(8.16) suggests that the generating functions for random
multigraphs might have interesting continued fraction forms. Such
expressions could well be of special importance, because they often converge
when power series do not.

The fact that the recurrence for the coefficients $e_{rd}$ can be ``solved''
to yield (7.3)--(7.5) should prove to be a good challenge for computer
systems that are now being constructed to solve recurrence relations
automatically. The similar recurrence for the coefficients
$e'_{rd}$, discussed in (7.24) and (7.25), will probably be an even greater
challenge; at least, no simple derivation of (7.21) from (7.26) 
is known. 

The solution to the recurrence for $e_{rd}$ in section~7 relies on the
introduction of a ``half excess'' stage, in which the polynomials must be
evaluated at integers plus~$1\over2$ although the recurrence in which
they are used involves integers only. In section~20 we found, similarly,
that it was fruitful to break the process of adding an edge into
stages in which ``half-edges'' were added. Perhaps the theory of
fractional differentiation will be of value in future investigations.
However, the operators $D^{1/2}$ and $\vartheta^{1/2}$ do not seem
to transform the basic functions $T^a\!/(1-T)^b$ very nicely.

Is there an equation (27.11$^{\prime}$) analogous to (27.11)? There
must be a reason why the coefficients of (27.7$^{\prime}$) tend to
have small prime factors.

We have seen numerous examples in which the multigraph process leads to
formulas that are mathematically cleaner than the analogous formulas
for the graph process. This suggests that an analogous theory be introduced
in place of the alternative ``${\bf G}_{n,p}$'' model of random graphs:
Instead of saying that each edge is present with probability~$p$,
the multiplicity of each edge should be allowed to have a Poisson
distribution with mean~$p$. Readers are encouraged to experiment
with such an approach.

Convergence to limiting distributions often appears to be monotonic.
For example, the probability that an evolving multigraph on
$n$~vertices stays on the top line appears to be strictly decreasing
as $n$ increases. How could this be proved?

Our proof of the top-line probability in Theorem~15 was independent of
the difficult analyses in Lemma~7 and Theorem~13 about the behavior of
random multigraphs with more than $\half (n+n^{2/3+\epsilon})$
edges; moreover, it did not use the stopping-configuration machinery of
sections 24--26, although that theory was in fact motivated by
attempts to prove Theorem~15 in a sharper form via generating
functions. The top-line phenomenon may perhaps be understood more
deeply if we use a generating-function-based approach, and the
following ideas may therefore prove to be useful.
Let $S(w,z)$ be the bgf for all multigraphs that never leave the top
line of Figure~1, where each multigraph is weighted by the probability
of having a purely top-line history as discussed in section~17. The
discussion of sections~19 and~20 shows that
$$S(w,z)=e^{U(w,z)+V(w,z)}H(w,z)\,,\eqno(29.3)$$
where $H(w,z)$ satisfies a differential equation almost like the
equation (5.1) that defines $E(w,z)$:
$${\textstyle{1\over
w}}\,(\vartheta_w-T\vartheta_z)H={\textstyle\half }
\,e^{-V}\vartheta_z^2e^VH-{\textstyle{1\over
2}}\,e^{-V}(\vartheta_z^2e^V) (H-1)\,.\eqno(29.4)$$
The subtracted term $\half e^{-V}(\vartheta_z^2e^V)(H-1)$ accounts
for the forbidden case that a new edge marked by~$\vartheta_z^2$ lies
entirely in the unicyclic part generated by~$e^V$; a~second complex
component arises if and only if this happens. The correction applies
to $H-1$, not~$H$, because the very first complex component does not
violate the top-line condition.

Expressing $H(w,z)$ in the form (25.16), we have $H_1=E_1$, but $H_2$
is smaller than~$E_2$:
$$H_2={5\over 16}\;{T^4\over (1-T)^6}+{25\over 48}\;{T^3\over(1-T)^5}+
{11\over 48}\;{T^2\over (1-T)^4}+{1\over 48}\;{T\over (1-T)^3}\,.$$
In general we can write
$$H_r=\sum h_{rd}\,{T^{2r-d}\over (1-T)^{3r-d}}\eqno(29.5)$$
for appropriate coefficients $h_{rd}$. The special case $\mu=\nu=0$ of
(20.7) tells us that
$$\vartheta^2e^V=\vartheta^2(1-T)^{-1/2}={\textstyle{1\over
2}}\,T(1-T)^{-7/2} +{\textstyle{5\over
4}}\,T^2(1-T)^{-9/2}\,;\eqno(29.6)$$
therefore we can compute the coefficients $h_{rd}$ by making a slight
change to the rule for computing~$e_{rd}$ that is expressed in
(20.11): Subtract~5 from the numerator of the first coefficient term
in (20.11), and subtract~1 from the numerator of the second
coefficient. The first coefficient now simplifies to
$${(6r-2d+5)(6r-2d+1)-5\over 8(3r-d+3)}={3r-d\over 2}\,.$$
In particular, when $d=0$ we have $h_{(r+1)0}={3\over 2}\,rh_{r0}$;
hence $h_{r0}$ is the number we called $k_r$ in (24.3).

Equation (25.17) now gives us a useful expression for the stopping
configurations,
$$\eqalignno{\nabla S&=e^{U(w,z)}\sum_{r\geq
2}w^r(\vartheta_z^2e^V)H_{r-1}(wz)\cr
\noalign{\smallskip}
&=e^{U(w,z)}\sum_{r\geq 2}w^r\left(\half \;{T(wz)\over
\bigl(1-T(wz)\bigr)^{7/2}} +{5\over 4}\;{T(wz)^2\over
\bigl(1-T(wz)\bigr)^{9/2}}\right)H_{r-1}(wz)\,.&(29.7)\cr}$$
The probability that an evolving multigraph on $n$ vertices leaves the
top line of Figure~1 is $\Phi_n\nabla S$.

For fixed $r$ we can evaluate the contribution made to $\Phi_n\nabla
S$ by the $r$\/th term of (29.7), to within $O(n^{-1/3})$, because the
leading coefficient $h_{(r-1)0}$ controls the asymptotic behavior.
Indeed, we know from (25.22) and the subsequent discussion that
$$\Phi_n\left(w^r\,{e^{U(w,z)}T(wz)^{2r}\over\bigl(1-T(wz)\bigr)^{3r+3/2}}\right)
={1\over 6re_r}+O(n^{-1/3})\eqno(29.8)$$
for all fixed $r$. Therefore when $\Phi_n$ is applied to the $r$\/th
term of (29.7) we get
$$\Phi_ne^U\mskip-1muw^r\left(\half \,{T\over (1-T)^{7/2}}+{5\over
4}\,{T^2\over (1-T)^{9/2}}\right)H_{r-1}={5k_{r-1}\over
24re_r}+O(n^{-1/3})\,.\eqno(29.9)$$ 
When $r=2$, the limit is ${5\over 77}$; when $r>2$, (7.1) and (24.3) imply that
$${5k_{r-1}\over 24re_r}=\left({5k_{r-2}\over 24(r-1)e_{r-1}}\right)
\left({36(r-1)(r-2)\over (6r-1)(6r-5)}\right)\,.$$
It follows by induction that
$${5k_{r-1}\over 24re_r}={5\over
36(r-1)r}\,\prod_{k=1}^{r-1}\,{k(k+1)\over (k+{1\over 6})(k+{5\over
6})}=\prod_{k=1}^{r-2}\,{k(k+1)\over(k+{1\over 6})(k+{5\over 6})}
-\prod_{k=1}^{r-1}\,{k(k+1)\over (k+{1\over 6})(k+{5\over 6})}\,.$$
So the sum over $r$ is a telescoping series,
$$\sum_{r\geq 2}\,{5k_{r-1}\over 24re_r}=1-\prod_{k=1}^{\infty}\,
{k(k+1)\over (k+{1\over 6})(k+{5\over 6})}=1-{5\pi\over
18}\,.\eqno(29.10)$$
In other words, convergence to the top-line probability depends entirely on the
sum over~$r$ of the error term in (29.9).

The number of challenging and potentially fruitful questions that remain
unanswered seems to be almost endless. But we shall close this list
of research problems by stating what seems to be the single most important
related area ripe for investigation at the present time. Wright~[\Wsii]
gave a procedure for computing the number of {\it strongly connected
labeled digraphs\/} of excess~$r$, analogous to his formulas for connected
labeled undirected graphs.
Random directed multigraphs are of great importance in computer applications,
and it is shocking that so little attention has been given to
their study so far. Karp [\Kar] carried Wright's investigations further
and discovered a beautiful theorem:
A~random digraph with $n(1+\mu)$ directed arcs almost surely
has a giant strong component of size $\sim {\mit\Theta}(\mu)^2n$,
when ${\mit\Theta}(\mu)$ is the factor such that
an undirected graph with $\half n(1+\mu)$ edges almost surely
has a giant component of size $\sim{\mit\Theta}(\mu)n$.
(The function ${\mit\Theta}(\mu)$ is $(\mu+\sigma)/(1+\mu)$,
according to (23.11). Karp's investigation was based on ${\bf D}_{n,p}$,
in which every directed arc is present with probability~$p$, but a
similar result surely holds for other models of random digraphs.)
A complete analysis of the random {\it directed\/} multigraph process
is clearly called for, preferably based on generating functions
so that extensive quantitative information can be derived without
difficulty.

Here is a sketch of how such an investigation might begin.
The {\it directed multigraph process\/} consists of adding directed
arcs $x\to y$ repeatedly to an initially empty multiset of arcs
on the vertices $\{1,2,\ldots,n\}$, where $x$ and $y$ are independently
and uniformly distributed between 1 and~$n$. The {\it compensation
factor} $\kappa(M)$ of a multidigraph~$M$ with $m_{xy}$ arcs from
$x$ to~$y$ is $1\big/\prod_{x=1}^n\prod_{y=1}^n m_{xy}!\,$; we can use it to
compute bivariate generating functions as in~(2.1). The bgf for all
possible multidigraphs is $\sum_{n\ge0}e^{n^2w}z^n\!/n!=G(2w,z)$.

Let $\cal A$ be the family of all multidigraphs such that all vertices
are reachable from vertex~1 via a directed path, and let $A(w,z)$ be
the corresponding bgf.  There is a nice relation between $A(w,z)$ and
the bgf $C(w,z)$ for connected undirected multigraphs, (2.10): If
$A(w,z)=\sum_{n\ge1}a_n(w)z^n\!/n!$, we have
$$\sum_{n\ge1}a_n(w)\,e^{-n^2w/2}\,{z^n\over n!}\;=\;C(w,z)\,.\eqno(29.11)$$
This can be proved by replacing $z$ by $ze^{-w/2}$ and noting that
$C(w,ze^{-w/2})$ is the bgf for connected multigraphs without self-loops,
and by showing that all members of~$\cal A$ are obtainable from such
connected multigraphs~$M$ by the following reversible construction:
Define a linear ordering $\prec$ on the vertices $\{1,2,\ldots,n\}$ by
saying that $x\prec y$ if $d(x)<d(y)$ or $d(x)=d(y)$ and $x<y$, where
$d(x)$ is the distance from 1 to~$x$ in~$M$. Then define a multidigraph
$D\in\cal A$ by arcs $x\to y$ whenever $x\rbar y$ in~$M$ and $x\prec y$;
include arbitrary additional arcs $x\to y$ for all pairs of vertices
with $x\succeq y$. The construction is reversible because $d(x)$ is
easily seen to be the distance from 1 to~$x$ in~$D$, regardless of
the choice of additional arcs. The additional arcs correspond to a
multiplicative factor $e^{{n+1\choose2}w}=e^{n^2w/2}(e^{w/2})^n$ in
an $n$-vertex multigraph, with one factor $e^w$ for each of the
$n+1\choose2$ vertex pairs $x\succeq y$.

Let $\cal S$ be the family of all strongly connected multidigraphs,
and let $S(w,z)=s_1(w)z+s_2(w)z^2\!/2!+s_3(w)z^3\!/3!+\cdots$ be the
corresponding bgf. A nontrivial identity discovered by Wright~[\Wsi]
implies that we can calculate the coefficients $s_n(w)$ by using the formula
$$\sum_{n\ge1}s_n(w)\,e^{-n^2w/2}\,{z^{n-1}\over(n-1)!}\,{G(w,z)\over
G(w,ze^{-nw})}\;=\;C'(w,z)\,,\eqno(29.12)$$
where the prime in $C'(w,z)$ denotes differentiation with respect
to~$z$. Notice that our generating function~$G(w,z)$ satisfies
$$\eqalignno{&G'(w,z)=e^{w/2}G(w,ze^w)\,,\quad
G''(w,z)=e^{2w}G(w,ze^{2w})\,,\cr
&\hskip6em\ldots\,,\quad G^{(n)}(w,z)=e^{n^2w/2}G(w,ze^{nw})\,,\quad\ldots\,;
&(29.13)\cr}$$
thus the denominator $G(w,ze^{-nw})$ in (29.12) is essentially an
$n$-fold integral of $G(w,z)$.

Wright [\Wsii] proved that the number of strongly connected digraphs with
$n+r$ arcs on $n$ vertices, disallowing self-loops and multiple arcs,
is $n!$ times a polynomial in~$n$ of degree $3r-1$, when $n>r>0$.
His proof can be adapted to multidigraphs, and everything becomes
much simpler, just as formula (9.4) for multigraphs is simpler than
formula (9.20) for graphs. The analogs of (2.11) and (3.4) are
$$S(w,z)=w^{-1}S_{-1}(wz)+S_0(wz)+wS_1(wz)+w^2S_2(wz)+\cdots\,,\eqno(29.14)$$
where
$$\eqalignno{S_{-1}(z)&=z\,,&(29.15)\cr
S_0(z)&=-\ln(1-z)\,,&(29.16)\cr}$$
and $S_r(z)$ for $r\ge1$
can easily be shown to be $(1-z)^{-3r}$ times a polynomial
in~$z$ of degree $<3r$. For example, the multidigraphs enumerated by
$wS_1(wz)$ all arise by inserting (``uncancelling'') vertices in the
arcs of the reduced multidigraphs
$$\unitlength=10pt
\vcenter
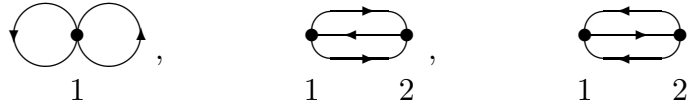
{\hbox{\beginpicture(4,3)(0,-1)
\put(2,1){\disk{.4}}
\put(2,-.75){\makebox(0,0){1}}
\put(1,1){\circle2}
\put(3,1){\circle2}
\put(0,.75){\vector(0,-1){0}}
\put(4,1.25){\vector(0,1){0}}
\endpicture}}
\,\;,\qquad\qquad
\vcenter{\hbox{\beginpicture(4,3)(0,-1)
\put(.5,1){\disk{.4}}
\put(.5,-.75){\makebox(0,0){1}}
\put(3.5,1){\disk{.4}}
\put(3.5,-.75){\makebox(0,0){2}}
\put(2,1){\oval(3,1.5)}
\put(.5,1){\line(1,0)3}
\put(2.5,1.75){\vector(1,0){0}}
\put(1.5,1){\vector(-1,0){0}}
\put(2.5,0.25){\vector(1,0){0}}
\endpicture}}
\,,\qquad\qquad
\vcenter{\hbox{\beginpicture(4,3)(0,-1)
\put(.5,1){\disk{.4}}
\put(.5,-.75){\makebox(0,0){1}}
\put(3.5,1){\disk{.4}}
\put(3.5,-.75){\makebox(0,0){2}}
\put(2,1){\oval(3,1.5)}
\put(.5,1){\line(1,0)3}
\put(1.5,1.75){\vector(-1,0){0}}
\put(2.5,1){\vector(1,0){0}}
\put(1.5,0.25){\vector(-1,0){0}}
\endpicture}}
$$
whose generating functions are respectively $\half w^2z$,
${1\over4}w^3z^2$, ${1\over4}w^3z^2$. The operation of uncancelling
corresponds to replacing $w$ by $w/(1-wz)$, as in Lemma~1; so
$wS_1(wz)=\half w^2z/(1-wz)^2+\half w^3z^2\!/(1-wz)^3=
\half w^2z/(1-wz)^3$, and $S_1(z)=\half z/(1-z)^3$.

In fact, the numerator of $S_r(z)$ turns out to have a surprisingly
small degree. Computer calculations indicate that we can write
$$S_r(z)={s_{r0}z^{2r-1}\over(1-z)^{3r}}+
         {s_{r1}z^{2r-2}\over(1-z)^{3r-1}}+\cdots+
         {s_{r(2r-2)}z\over(1-z)^{r+2}}\,,\eqno(29.17)$$
a formula analogous to (8.4), at least when $r\le5$. The coefficients are
$$\advance\baselineskip5pt
\centerline{\vbox{\halign{\hfil$#$&&\ \ \hfil$#$\hfil\cr
d=&0&1&2&3&4&5&6&7&8\cr
s_{1d}=
&{ 1\over 2}\cr
s_{2d}=
&{ 17\over 8}
&{ 13\over 8}
&{ 1\over 6}\cr
s_{3d}=
&{ 275\over 12}
&{ 427\over 12}
&{ 391\over 24}
&{ 13\over 6}
&{ 1\over 24}\cr
s_{4d}=
&{ 26141\over  64}
&{ 61231\over  64}
&{ 51299\over  64}
&{ 18473\over  64}
&{ 6047\over 144}
&{ 263\over 144}
&{  1\over 120}\cr
s_{5d}=
&{ 1630711\over   160}
&{ 1276481\over   40}
&{ 3125933\over   80}
&{ 2840093\over   120}
&{ 3546283\over   480}
&{ 6743\over  6}
&{ 25307\over  360}
&{ 43\over 36}
&{  1\over 720}\cr}}}$$
No reason why $S_r(z)$ should have the simple form (29.17) is apparent;
this phenomenon cries out for explanation, if it is indeed true for
all $r>0$, and the explanation will probably lead to new theorems
of interest. It can be shown that this conjecture is equivalent
to the assertion that the sum of $(-1)^\nu \kappa/\nu!$, over all
labelled, reduced, strongly connected multidigraphs of excess~$r$,
is zero; or in other words, if we choose a labelled, reduced, strongly
connected multidigraph of excess~$r$ at random, with probabilities
weighted in the natural way by the compensation factor~$\kappa$, then
the probability is $\half$ that
there will be an even number of vertices.

Is there a simple recurrence governing the leading coefficients
$s_{10}$, $s_{20}$, $s_{30}$, \dots, perhaps analogous to the relation
we observed for ordinary connected components in (8.5)?

\bigbreak\noindent
{\bf Acknowledgments.}\enspace
 The authors wish to thank Prof.~Richard Askey for
helpful correspondence relating to this research.

\bigbreak\noindent
{\bf Appendix.}\enspace
Here is a list of corrections to the related paper [\FKP].

\halign{\qquad#\hfil\cr
Page 175, line 10: $(1+t)^N$ should be $(1+t)^{-N}$\cr
Page 175, line 11: (3.5) should be (3.6)\cr
Page 182, (4.21): $\sqrt{3t}$ should be $\sqrt{3}\,t$\cr
Page 183, line 18: $\half \,\sqrt{3t}$ should be ${i\over
2}\,\sqrt{3t}$\cr 
Page 183, line 24: (4.27) should be (4.25)\cr
Page 184, (5.6): $l=1$ should be $l-1$\cr
Page 185, line 17: $l=2$ should be $l=3$\cr
Page 189, lines 4 and 9: $\half l(l-1)$ should be $\half l(l+1)$\cr
Page 192, (7.13): ${31\over 45}$ should be ${1\over 45}$;
$2+3\hat{p}_3$ should be $\hat{p}_3$\cr
Page 194, line 15: `than $\Re h(\lambda)-\lambda 
    -(1-\half \lambda)(\ln(1-\half \lambda)-\ln(1+\half \lambda))$\cr
\hskip10em $<\Re h(\lambda)-{1\over3}\lambda^2$ when'\cr
Page 205, line 7: delete `number of'\cr
Page 207, (11.9): delete commas in denominator\cr
Page 209, first line of (A.6): $ixt-it^3\!/3$ should be
$ixt+it^3\!/3$\cr
Page 213, the argument for enveloping series is incomplete\cr
Page 215, (11.12) and (11.14): delete commas in denominators\cr}

\vfill\eject

\centerline{\bf Bibliography}

\smallskip
\bib [\Bag]\enspace
G. N. Bagaev, ``Slucha\u\i nye grafy so stepen'\t\i u sv\t\i aznosti~2,''
{\sl Diskretny\u\i\ Analiz\/ \bf 22} (1973), 3--14.

\smallskip
\bib [\BD]\enspace
G. N. Bagaev and E. F. Dmitriev, ``Perechislenie sv\t\i aznykh otmechennykh
dvudol'nykh grafov,'' {\sl Doklady Akademi\t\i a Nauk BSSR\/ \bf28} (1984),
1061--1063.

\smallskip
\bib [\BCM]\enspace
Edward A. Bender, E. Rodney Canfield, and Brendan D. McKay, ``The asymptotic
number of labeled connected graphs with a given number of vertices and
edges,'' {\sl Random Structures and Algorithms\/ \bf1} (1990), 127--169.

\smallskip
\bib [\Bi]\enspace
B\'ela Bollob\'as, ``A probabilistic proof of an asymptotic formula
for the number of labelled regular graphs,'' {\sl European Journal
of Combinatorics\/ \bf 1} (1980), 311--316.

\smallskip
\bib [\Bii]\enspace
B\'ela Bollob\'as, ``The evolution of random graphs,''
{\sl Transactions of the American Mathematical Society\/ \bf 286}
(1984), 257--274.

\smallskip
\bib [\Biii]\enspace
B\'ela Bollob\'as, {\sl Random Graphs\/} (London: Academic Press, 1985).

\smallskip
\bib [\BF]\enspace
B. Bollob\'as and A. Frieze, ``On matchings and Hamiltonian cycles in
random graphs,'' in {\sl Random Graphs~'83}, edited by Micha\l\
Karo\'nski and Andrzej Ruci\'nski, {\sl Annals of Discrete
Mathematics\/ \bf 28} (1985), 23--46.

\smallskip
\bib [\Bor]\enspace
C. W. Borchardt, ``Ueber eine der Interpolation entsprechende Darstellung
der Elim\-inations-Resultante,''
{\sl Journal f\"ur die reine und angewandte Mathematik\/ \bf 57} (1860),
111--121.

\smallskip
\bib [\Bri]\enspace
V. E. Britikov, ``O strukture slucha\u\i nogo grafa vblizi kritichesko\u\i\
tochki,'' {\sl Diskretna\t\i a Matematika\/ \bf1},3 (1989), 121--128.
English translation, ``On the random graph structure near the critical
point,'' {\sl Discrete Mathematics and Applications\/ \bf1},3 (1991),
301--309.

\smallskip
\bib [\Cay]\enspace
A. Cayley, ``A theorem on trees,'' {\sl Quarterly Journal of Pure and
Applied Mathematics\/ \bf 23} (1889), 376--378. Reprinted in his {\sl
Mathematical Papers\/ \bf 13}, 26--28.

\smallskip
\bib [\Eis]\enspace
Gotthold Eisenstein, ``Entwicklung von $\alpha^{\alpha^{\alpha^{\dddots}}}$,''
{\sl Journal f\"ur die reine und angewandte Mathematik\/ \bf 28} 
(1844), 49--52.

\smallskip
\bib [\ERo]\enspace
P. Erd{\H o}s and A. R\'enyi, ``On random graphs~I,'' {\sl
Publicationes Mathematicae (Debrecen) \bf 6} (1959), 
290--297. Reprinted in {\sl
Paul Erd{\H o}s: The Art of Counting\/} (MIT Press, 1973), 561--568;
and in {\sl Selected Papers of Alfr\'ed R\'enyi\/} (Akad\'emiai Kiad\'o,
1976), 308--315.

\smallskip
\bib [\ER]\enspace
P. Erd{\H o}s and A. R\'enyi, ``On the evolution of random graphs,''
{\sl A Magyar Tudom\'anyos Akad\'emia Matematikai Kutat\'o
Int\'ezet\'enek K\"ozlem\'enyei\/
\bf 5} (1960), 17--61.
Reprinted in {\sl
Paul Erd{\H o}s: The Art of Counting\/} (MIT Press, 1973), 574--618;
and in {\sl Selected Papers of Alfr\'ed R\'enyi\/} (Akad\'emiai Kiad\'o,
1976), 482--525.

\smallskip
\bib [\FKP]\enspace
Philippe Flajolet, Donald E. Knuth, and Boris Pittel,
``The first cycles in an evolving graph,'' {\sl Discrete Mathematics\/
\bf 75} (1989), 167--215.

\smallskip
\bib [\FKG]\enspace
C. M. Fortuin, P. W. Kasteleyn, and J. Ginibre, ``Correlation
inequalities on some partially ordered sets,'' {\sl Communications in
Mathematical Physics\/ \bf 22} (1971), 89--103.

\smallskip
\bib [\GJ]\enspace
I. P. Goulden and D. M. Jackson, {\sl Combinatorial Enumeration\/}
(New York: Wiley, 1983).

\smallskip
\bib [\CM]\enspace
Ronald L. Graham, Donald E. Knuth, and Oren Patashnik, {\sl Concrete
Mathematics\/} (Reading, Massachusetts: Addison-Wesley, 1989).

\smallskip
\bib [\Hen]\enspace
Peter Henrici, {\sl Applied and Computational Complex Analysis}, volume~2,
Wiley, 1977.

\smallskip
\bib [\Jan]\enspace
Svante Janson, ``Multicyclic components in a random graph process,''
{\sl Random Structures and Algorithms\/ \bf4} (1993), 71--84.

\smallskip
\bib [\JL]\enspace
Svante Janson and Tomasz \L uczak, ``The size of the last cycle in the
random graph process,'' {\sl Abstracts of Papers Presented to the
American Mathematical Society\/ \bf13} (1992), 354, abstract 875-05-131.

\smallskip
\bib [\Kar]\enspace
Richard M. Karp, ``The transitive closure of a random digraph,'' {\sl Random
Structures and Algorithms\/ \bf1} (1990), 73--93.

\smallskip
\bib [\AOC]\enspace
Donald E. Knuth, ``An analysis of optimum caching,'' 
{\sl Journal of Algorithms\/ \bf 6} (1985), 181--199.

\smallskip
\bib [\CP]\enspace
Donald E. Knuth, ``Convolution polynomials,'' {\sl The Mathematica
Journal\/ \bf2},4 (Fall  1992), 67--78.

\smallskip
\bib [\KP]\enspace
Donald E. Knuth and Boris Pittel, ``A recurrence related to trees,''
{\sl Proceedings of the American Mathematical Society\/ \bf 105}
(1989), 335--349.

\smallskip
\bib [\Li]\enspace
Tomasz \L uczak, ``Component behavior near the critical point of the random
graph process,'' {\sl Random Structures and Algorithms\/ \bf1} (1990),
287--310.

\smallskip
\bib [\Lii]\enspace
Tomasz \L uczak, ``Cycles in a random graph near the critical point,''
{\sl Random Structures and Algorithms\/ \bf 2} (1991), 421--439.

\smallskip
\bib [\LW]\enspace
Tomasz \L uczak and John C. Wierman, ``The chromatic number of random
graphs at the double-jump threshold,'' {\sl Combinatorica\/ \bf 9}
(1989), 39--49.

\smallskip
\bib [\LPW]\enspace
Tomasz \L uczak, Boris Pittel, and John C. Wierman, ``The
structure of a random graph at the point of phase transition,''
to appear in {\sl Transactions of the American Mathematical Society}.

\smallskip
\bib [\Mit]\enspace
D. S. Mitrinovi\'c, {\sl Analytic Inequalities\/} (Springer-Verlag, 1970).

\smallskip
\bib [\Ram]\enspace
S. Ramanujan, ``Questions for solution, number 294,'' {\sl Journal of
the Indian Mathematical Society\/ \bf3} (1911), 128; {\bf4} (1912), 151--152.

\smallskip
\bib [\Ren]\enspace
Alfred R\'enyi, ``Some remarks on the theory of trees,''
{\sl A Magyar Tudom\'anyos Aka\-d\'emia Matematikai Kutat\'o
Int\'ezet\'enek K\"ozlem\'enyei\/  \bf 4} (1959), 73--85. Reprinted in
{\sl Selected Papers of Alfr\'ed R\'enyi\/ \bf 2}, 363--374.

\smallskip
\bib [\Rid]\enspace
Robert James Riddell, Jr., {\sl Contributions to the Theory of Condensation}.
Dissertation,
 University of Michigan, 1951. (The main results of this dissertation
were published as R. J. Riddell, Jr., and G. E. Uhlenbeck, ``On the theory of
the virial development of the equation of state of monoatomic gases,''
{\sl Journal of Chemical Physics\/ \bf 21} (1953), 2056--2064.)

\smallskip
\bib [\Sei]\enspace
G. Seitz, ``Une remarque aux in\'egalit\'es,'' 
{\sl Aktuarsk\'e V{\v e}dy\/ \bf 6} (1936/37), 167--171.

\smallskip
\bib [\Sla]\enspace
L. J. Slater, ``Expansions of generalized Whittaker functions,'' {\sl
Proceedings of the Cambridge Philosophical Society\/ \bf 50} (1954),
628--630.

\smallskip
\bib [\Si]\enspace
V. E. Stepanov, ``Neskol'ko teorem otnositel'no slucha\u{\i}nykh grafov,''
{\sl Vero\t\i atnostnye metody v diskretno\u\i\ matematike\/}
(Karel'ski\u\i\ filial Akademi\t\i a Nauk SSSR, Petrozavodsk, 1983),
90--92.

\smallskip
\bib [\Sii]\enspace
V. E. Stepanov, ``O nekotorykh osobennost\t\i akh stroeni\t\i a\
 slucha\kern.05em{\u\i}nogo grafa
vblizi kritichesko\kern.05em\u\i\ tochki,'' 
{\sl Teoriya Veroyatnostei i ee Primeneni\t\i a\/ 
\bf 32} (1988), 633--657.
English translation, ``On some features of the structure of a random graph
near a critical point,'' {\sl Theory of Probability and Its
Applications\/ \bf 32} (1988), 573--594.

\smallskip
\bib [\Syl]\enspace
J. J. Sylvester, ``On the change of systems of independent
variables,'' {\sl Quarterly Journal of Pure and Applied Mathematics\/
\bf 1} (1857), 42--56. Reprinted in his {\sl Mathematical Papers\/ \bf
2}, 65--85.

\smallskip
\bib [\Vob]\enspace
V. A. Vobly\u\i, ``O koeffitsientakh Ra{\u\i}ta i Stepanova-Ra{\u\i}ta,''
{\sl Matematicheskie Zametki\/} 42 (1987), 854--862. English translation,
V. A. Voblyi, ``Wright and Stepanov-Wright coefficients,'' {\sl Mathematical
Notes\/ \bf 42} (1987), 969--974.

\smallskip
\bib [\Wo]\enspace
E. M. Wright, ``A relationship between two sequences,'' {\sl
Proceedings of the London Mathematical Society
 \bf 17} (1967), 296--304, 547--552.

\smallskip
\bib [\Wsi]\enspace
E. M. Wright, ``The number of strong digraphs,'' {\sl Bulletin of the
London Mathematical Society\/ \bf3} (1971), 348--350.

\smallskip
\bib [\Wi]\enspace
E. M. Wright, ``The number of connected sparsely edged graphs,''
{\sl Journal of Graph Theory\/ \bf 1} (1977), 317--330.

\smallskip
\bib [\Wsii]\enspace
E. M. Wright, ``Formulae for the number of sparsely-edged strong
labelled digraphs,'' {\sl Quarterly Journal of Mathematics},
Oxford (2), {\bf 28} (1977), 363--368.

\smallskip
\bib [\Wii]\enspace
E. M. Wright, ``The number of connected sparsely edged graphs. II.
Smooth graphs and blocks,''
{\sl Journal of Graph Theory\/ \bf 2} (1978), 299--305.

\smallskip
\bib [\Wiii]\enspace
E. M. Wright, ``The number of connected sparsely edged graphs. III.
Asymptotic results,''
{\sl Journal of Graph Theory\/ \bf 4} (1980), 393--407.

\smallskip
\bib [\Wiv]\enspace
E. M. Wright, ``The number of connected sparsely edged graphs. IV.
Large nonseparable graphs,''
{\sl Journal of Graph Theory\/ \bf 7} (1983), 219--229.

\bye